\documentclass[letterpaper,11pt,twosides,onecolumn,final]{report}

\usepackage{multirow}

\usepackage[table]{xcolor}

\usepackage{stmaryrd}
\usepackage{amsfonts,amsmath,amssymb}

\usepackage{color}
\usepackage[table]{xcolor}

\usepackage{rotating}
\usepackage{mathabx}

\usepackage{cite}

\newcommand{\tabsize}{\fontsize{8}{9pt}\selectfont}
\usepackage{relsize}
\newcommand{\ltimesx}{\mathlarger{\ltimes}}

\usepackage[vlined,ruled,commentsnumbered]{algorithm2e}

\definecolor{aqua}{rgb}{0 , 0.4118, 0.5686}
\definecolor{violet}{rgb}{0.4706 ,   0.3137 ,   0.6667}
\definecolor{violetdark}{rgb}{0.3686  ,  0.2353   , 0.5216}
\definecolor{milkblue}{rgb}{0.2,0.2,0.7}
\definecolor{redaccent}{rgb}{0.3882  ,  0.1412  ,  0.1373}
\definecolor{bluedark}{rgb}{0.2,0.3,0.6}
\definecolor{QuartzPurple}{rgb}{0.302,	0.235,	0.573}
\definecolor{SapphireBlue}{rgb}{0.000,	0.388,	0.671}
\definecolor{DarkBlue}{rgb}{0.090,	0.212,	0.365}
\definecolor{BrightBlue}{rgb}{0.9137 ,   0.9686,    0.9922}

\colorlet{algshading_color}{BrightBlue}
\renewcommand\algcolor  {\color{black}}
\renewcommand\algcaptioncolor {\color{QuartzPurple}}

\usepackage{algorithmic}
\usepackage{booktabs}

\usepackage{url}

\usepackage{colortbl}

\definecolor{shadingcolor}{rgb} {.87,    0.92,    .98}
\newcommand{\shadingbox}[1]{   
    \fboxsep 0pt
    \colorbox{shadingcolor}{
        {\hskip -2pt #1}\hskip -2.5pt
    }
}

\newcommand{\minitab}[2][l]{\begin{tabular}{@{}#1}#2\end{tabular}}

\newcommand{\skron}{\,|\!\!\otimes\!\!|\,}   
\newcommand{\Cprod}{\,|\!\bullet\!|\,}   

\input{def1b.set}

\usepackage{mathpazo}

\usepackage[mathscr]{eucal}

\newcommand{\matn}[2][n]{\ensuremath{\mathbf{#2}^{(#1)}}}

\makeatletter
\newsavebox{\@brx}
\newcommand{\llangle}[1][]{\savebox{\@brx}{\(\m@th{#1\langle}\)}%
  \mathopen{\copy\@brx\kern-0.5\wd\@brx\usebox{\@brx}}}
\newcommand{\rrangle}[1][]{\savebox{\@brx}{\(\m@th{#1\rangle}\)}%
  \mathclose{\copy\@brx\kern-0.5\wd\@brx\usebox{\@brx}}}
\makeatother

\usepackage[titletoc,title]{appendix}
\usepackage[british]{babel}

\title{Low-Rank Tensor Networks  for  Dimensionality  Reduction and
Large-Scale Optimization  Problems:  Perspectives and Challenges PART 1
{\footnote{Copyright  A.Cichocki {\it et al.} Please make reference to: A. Cichocki, N. Lee, I. Oseledets, A.-H. Phan, Q. Zhao and D.P. Mandic (2016), ``Tensor Networks for Dimensionality Reduction and Large-scale Optimization: Part 1 Low-Rank Tensor Decompositions'', Foundations and Trends in Machine Learning: Vol. 9: No. 4-5, pp 249-429. }} \\
\vspace{0.9cm}
A. Cichocki\hspace{0.9cm}  N. Lee, \\ I.V. Oseledets,  \hspace{0.9cm} A-H. Phan, \\
 Q. Zhao, \hspace{0.6cm}  D. Mandic \vfill\newpage}

\author{Andrzej CICHOCKI\\
RIKEN Brain Science Institute (BSI), Japan and SKOLTECH, Russia \\
cia@brain.riken.jp\\
\\
Namgil LEE \\
RIKEN BSI,   namgil.lee@riken.jp\\
\\
Ivan OSELEDETS\\
Skolkovo Institute of Science and Technology (SKOLTECH), and\\
Institute of Numerical Mathematics of Russian Academy of Sciences,\\
 Russia\\
i.oseledets@skolkovotech.ru\\
\\
Anh-Huy PHAN \\
RIKEN BSI, phan@brain.riken.jp\\
\\
Qibin ZHAO \\
RIKEN BSI,  qbzhao@brain.riken.jp\\
\\
Danilo P. MANDIC\\
Imperial College, UK\\
d.mandic@imperial.ac.uk
}

\date{}

\emergencystretch=\maxdimen
\usepackage{hyperref}

\usepackage{atbegshi}

\begin{document}

\sloppy

\pagenumbering{gobble}
\begingroup
\makeatletter
\let\ps@empty\ps@plain 
\maketitle
\endgroup




\pagenumbering{arabic}
\begin{abstract}
\thispagestyle{plain}

Machine learning and data mining  algorithms are becoming increasingly important in analyzing  large volume,  multi-relational and multi--modal datasets, which are often conveniently represented as multiway arrays or tensors.  It is therefore timely and valuable for the multidisciplinary research community to review  tensor decompositions and tensor networks as emerging tools for large-scale data analysis and data mining. We provide  the mathematical and graphical representations and interpretation of tensor networks, with the main focus on the Tucker and Tensor Train (TT) decompositions and their extensions or generalizations.

To make the material self-contained, we also address the concept of tensorization which allows for the creation of very high-order tensors from lower-order  structured datasets represented by vectors or matrices. Then, in order to combat the curse of dimensionality and  possibly obtain  linear or even sub-linear complexity of storage and computation,  we address super-compression of tensor data through low-rank tensor networks. Finally, we demonstrate how such approximations can be used to solve a wide class of huge-scale linear/ multilinear dimensionality reduction and related optimization problems that are far from being tractable when using classical numerical methods.

The challenge for huge-scale optimization problems is therefore to develop methods which scale linearly or sub-linearly (i.e., logarithmic complexity) with the size of datasets, in order to benefit from the well-- understood optimization frameworks for smaller size problems.  However, most efficient optimization algorithms are  convex and do not scale well with data volume, while linearly scalable algorithms typically only apply to very specific scenarios. In this review, we address this problem through the concepts of low-rank tensor network approximations, distributed tensor networks, and the associated learning algorithms. We then elucidate how these concepts can be used to convert otherwise intractable huge-scale optimization problems into a set of much smaller linked and/or distributed sub-problems of affordable size and complexity. In doing so, we highlight the ability of tensor networks to account for the couplings between the multiple variables, and for multimodal, incomplete and noisy data.

The methods and approaches discussed in this work can be considered both as an alternative and a complement to emerging methods for huge-scale optimization, such as the random coordinate descent (RCD) scheme, subgradient methods, alternating direction method of multipliers (ADMM) methods, and proximal gradient descent methods. This is PART1 which consists of Sections 1-4.\\

\end{abstract}

\setcounter{page}{3}

\chapter{Introduction and Motivation}
\label{chap:intro}

\vspace{0cm}

This monograph aims to present a coherent account of ideas and methodologies related to tensor decompositions (TDs) and
tensor networks models (TNs).
Tensor decompositions (TDs) decompose principally data tensors  into  factor matrices, while tensor networks (TNs)
 decompose higher-order tensors into sparsely interconnected  small-scale low-order core tensors.
 These low-order core tensors are
 called ``components'', ``blocks'', ``factors'' or simply ``cores''. In this way, large-scale data can be approximately
 represented in highly compressed and distributed formats.

 In this monograph,  the TDs and TNs are treated  in a unified way, by considering TDs as simple tensor networks or sub-networks; the terms
   ``tensor decompositions'' and  ``tensor networks'' will  therefore be used interchangeably.
  Tensor networks can be thought of  as special graph structures which break down high-order tensors into a set of sparsely interconnected low-order core tensors, thus allowing for both enhanced interpretation and computational advantages.
  Such an approach is valuable in many application contexts which require the computation of eigenvalues and the
  corresponding  eigenvectors of extremely high-dimensional linear or nonlinear operators.  These operators typically describe
  the coupling between many degrees of freedom within real-world physical systems; such degrees of freedom are often only  weakly coupled.
Indeed, quantum physics provides  evidence that couplings between multiple data channels usually do not exist among all
the degrees of freedom but mostly locally, whereby ``relevant'' information, of relatively low-dimensionality, is embedded
into very large-dimensional measurements \cite{verstraete08MPS,Schollwock13,Orus2013,Murg_TTNS15}.

 Tensor networks offer a theoretical and computational framework for the analysis of computationally prohibitive large
 volumes of data, by  ``dissecting'' such data into the  ``relevant'' and ``irrelevant'' information, both of lower dimensionality.
In this way, tensor network representations  often allow for super-compression  of datasets as large as  $10^{50}$ entries,
down to the affordable levels of $10^7$ or even less entries
\cite{OseledetsTT09,QTT_Tucker,Kazeev_Toeplitz13,Kazeev14CME,KressnerEIG2014,vervliet2014breaking,Dolgov15CME,TPA2015,Multigrid_Amen2016}.

With the emergence of the big data paradigm, it is therefore both timely and important to provide
the multidisciplinary machine learning and data analytic communities
with a comprehensive overview of tensor networks, together with an example-rich guidance on their application
in several generic optimization problems
for huge-scale structured data. Our aim is also to unify the terminology, notation, and algorithms for tensor
decompositions and tensor networks which are being developed not only in machine learning, signal processing,
numerical analysis and scientific computing,
  but also in quantum physics/ chemistry for the representation of, e.g., quantum many-body systems.

 \section{Challenges in Big Data Processing}

\begin{figure}
(a)
\begin{center}
\includegraphics[width=7.8cm]{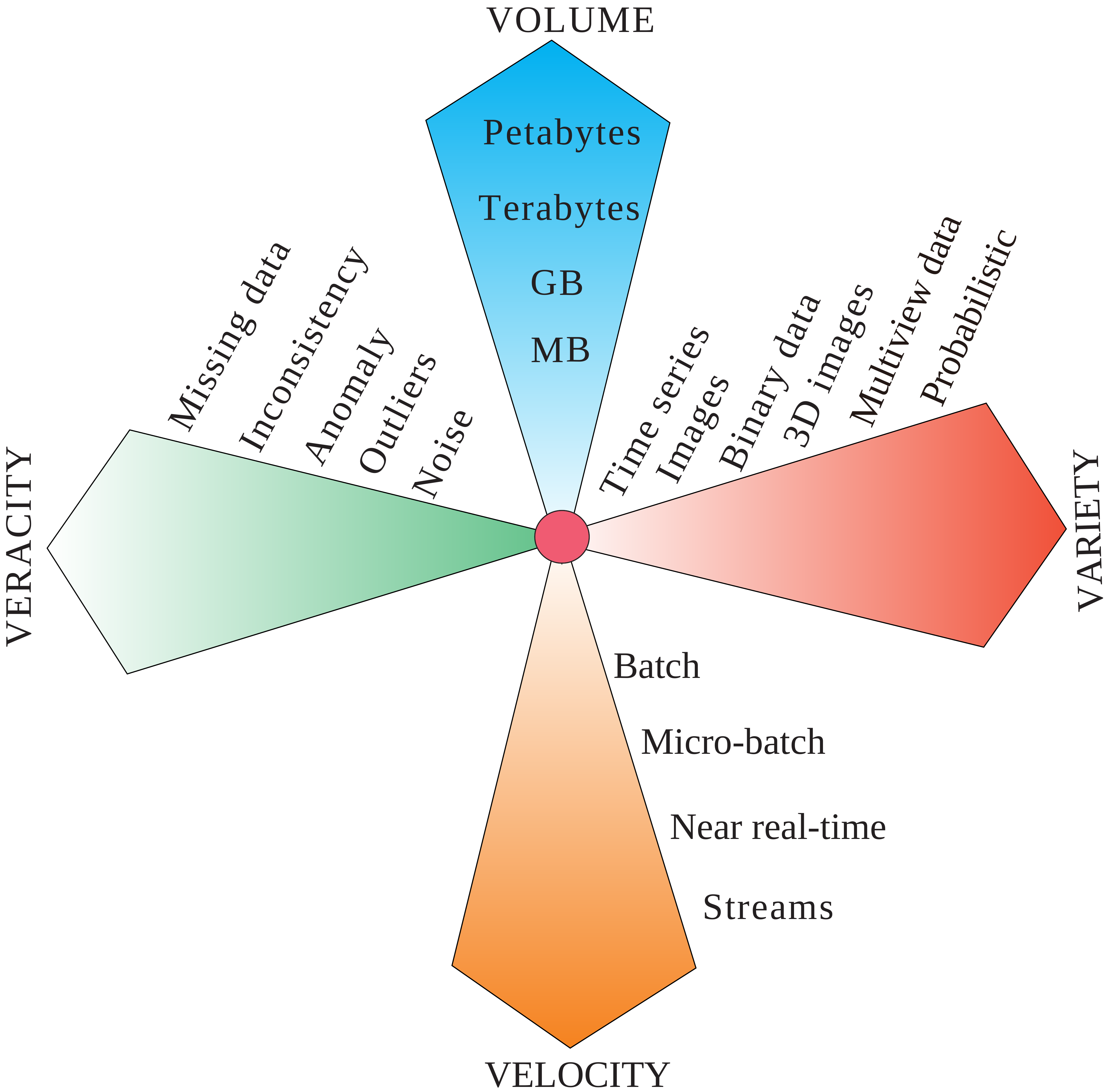}
\end{center}
(b)
\vspace{-0.2cm}
\begin{center}
\includegraphics[width=10.79cm]{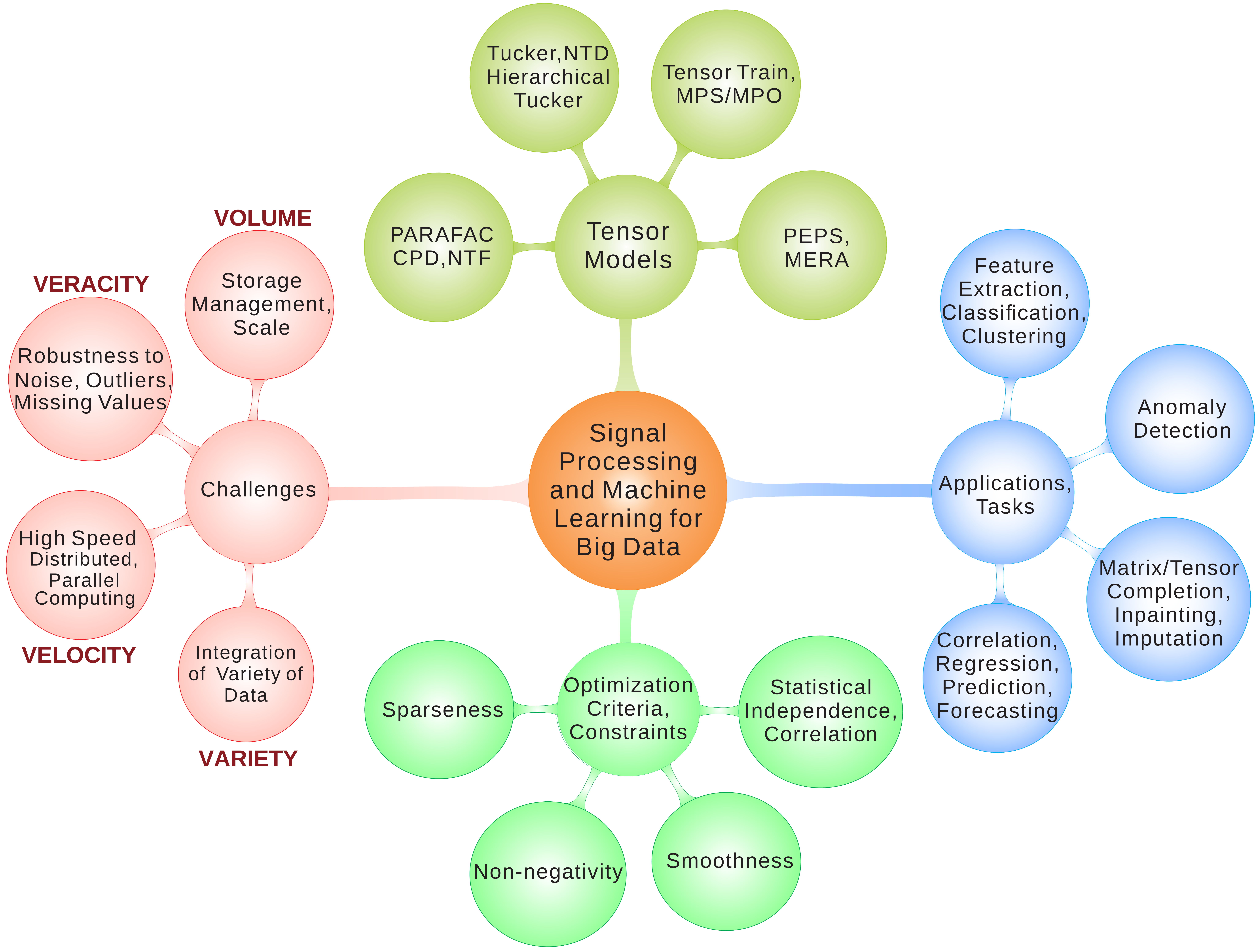}
\end{center}
\caption{{\small A framework for extremely large-scale data analysis. (a) The 4V challenges for big data.
  (b) A unified framework for the 4V challenges and the potential applications based on tensor decomposition approaches.}}
\label{Fig:BData}
\end{figure}

The  volume and structural complexity of modern datasets are becoming exceedingly  high, to the extent
which renders standard analysis methods and algorithms inadequate.
Apart from the huge Volume, the other  features which characterize big data include Veracity, Variety and Velocity
(see Figures~\ref{Fig:BData}(a) and~(b)). Each of the ``V features'' represents
a research challenge in its own right. For example, high volume implies the need for algorithms that are scalable; high
Velocity requires the processing
of big data streams in near real-time; high Veracity calls for robust and predictive algorithms for noisy, incomplete
and/or inconsistent data; high Variety
 demands the fusion of different data types, e.g.,  continuous,  discrete, binary, time series, images, video, text,
 probabilistic or multi-view.
 Some applications give rise to additional ``V challenges'', such as  Visualization,  Variability and Value.
The Value feature is particularly interesting and refers to  the  extraction of high quality and consistent information,
from which meaningful and interpretable results can be obtained.

Owing to the increasingly affordable recording
devices, extreme-scale volumes and variety of data are becoming ubiquitous across the science and engineering  disciplines. In the case of  multimedia (speech, video), remote sensing
and medical / biological data, the analysis also
requires a paradigm shift in order to efficiently process massive datasets within tolerable time (velocity).
Such massive datasets  may have billions of entries and are typically represented  in the form of huge block matrices
and/or tensors. 
This has spurred a  renewed interest in the development of  matrix / tensor algorithms that are suitable for very large-scale datasets.
We show that tensor networks provide a natural sparse and distributed representation for
big data, and address both established and  emerging methodologies for tensor-based representations and optimization.
Our particular focus is on  low-rank tensor network representations, which allow for
huge data tensors to be approximated (compressed) by
 interconnected  low-order core tensors.

\section{Tensor Notations and Graphical Representations}

\begin{figure}[t!]
\centering
\includegraphics[width=4.9cm]{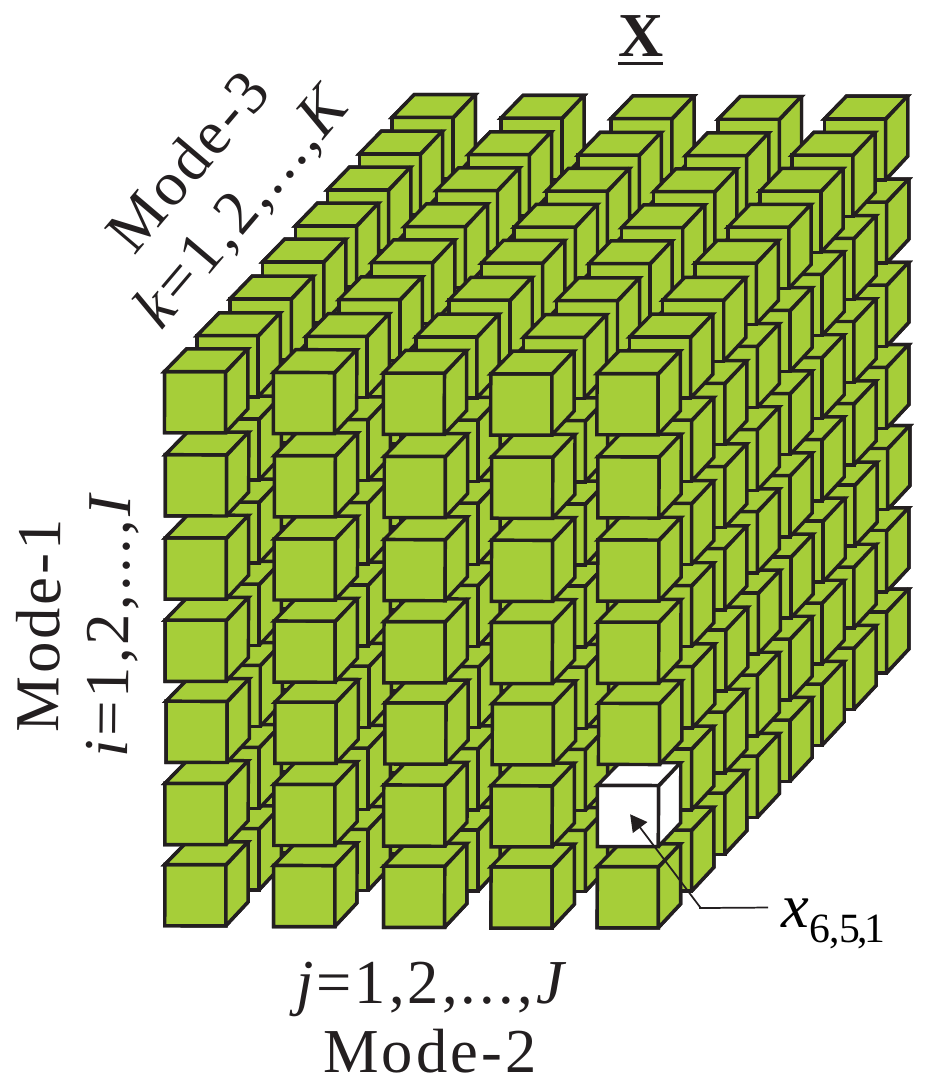}\\
\vspace{0.3cm}
\includegraphics[width=8.6cm]{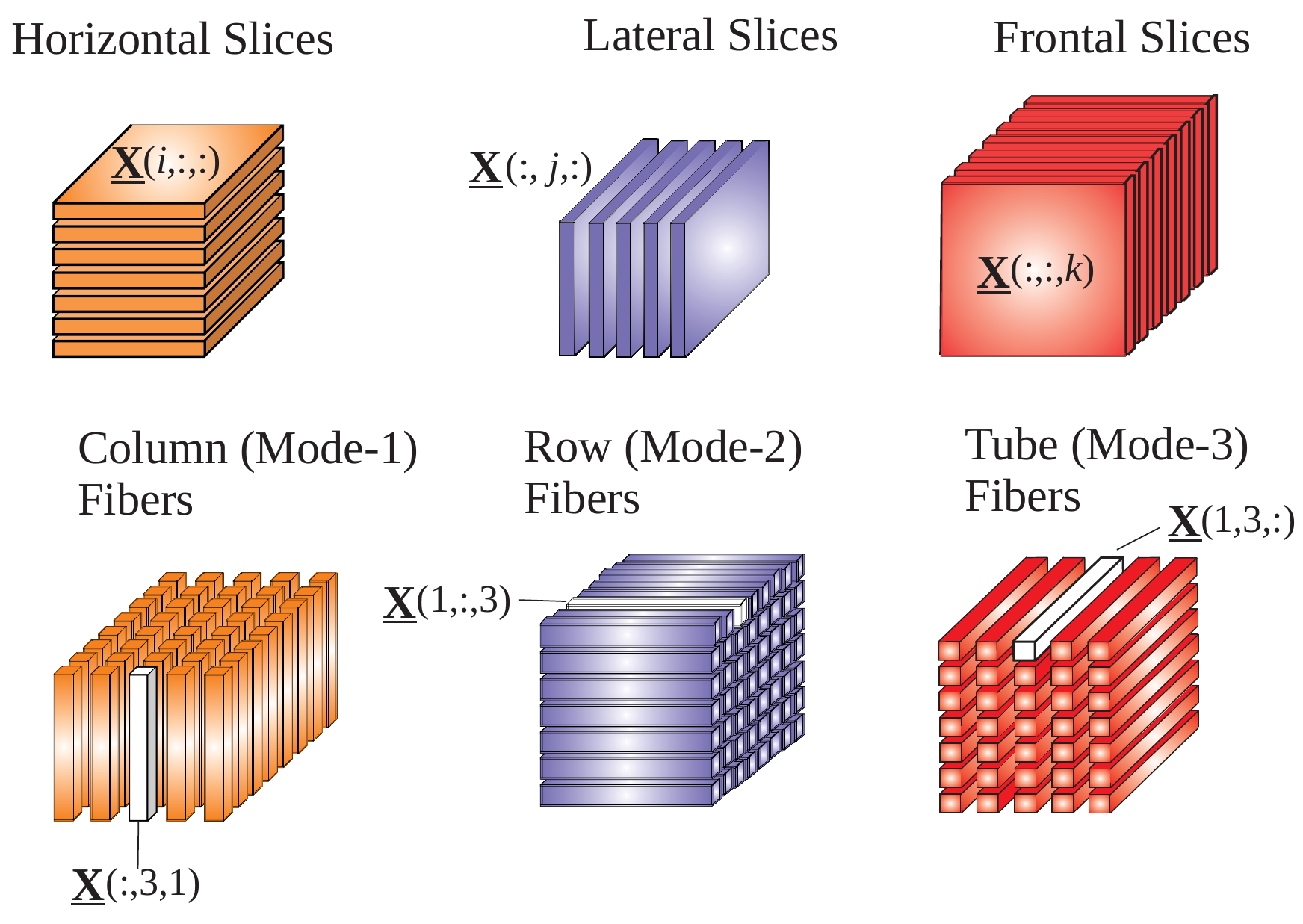}
\caption{A 3rd-order tensor $\underline \bX \in \Real^{I \times J \times K}$, with
entries $x_{i, j, k}=\underline \bX(i,j,k)$, and its subtensors: slices (middle) and  fibers (bottom).
All  fibers are treated as column vectors.}
\label{Fig:fibers}
\end{figure}
\begin{figure}[t]
\centering
\includegraphics[width=6.6cm]{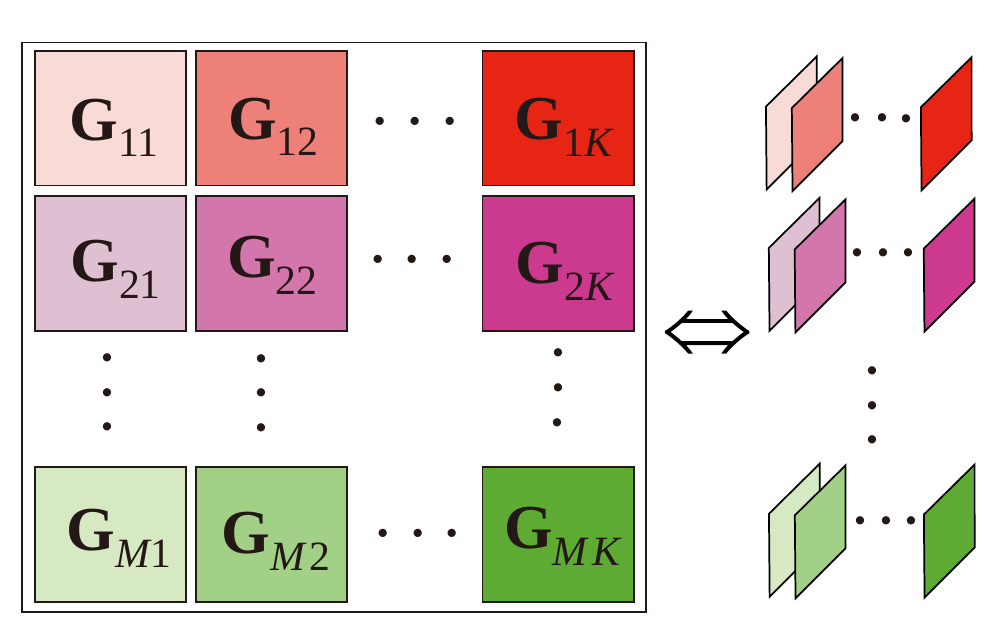}
\caption{A block matrix  and its representation as a 4th-order tensor, created by
reshaping (or a projection) of blocks  in the rows into lateral slices of 3rd-order tensors.}
\label{Fig:Tens_block_matr}
\end{figure}
\begin{figure}[ht!]
\centering
\includegraphics[width=11.5cm]{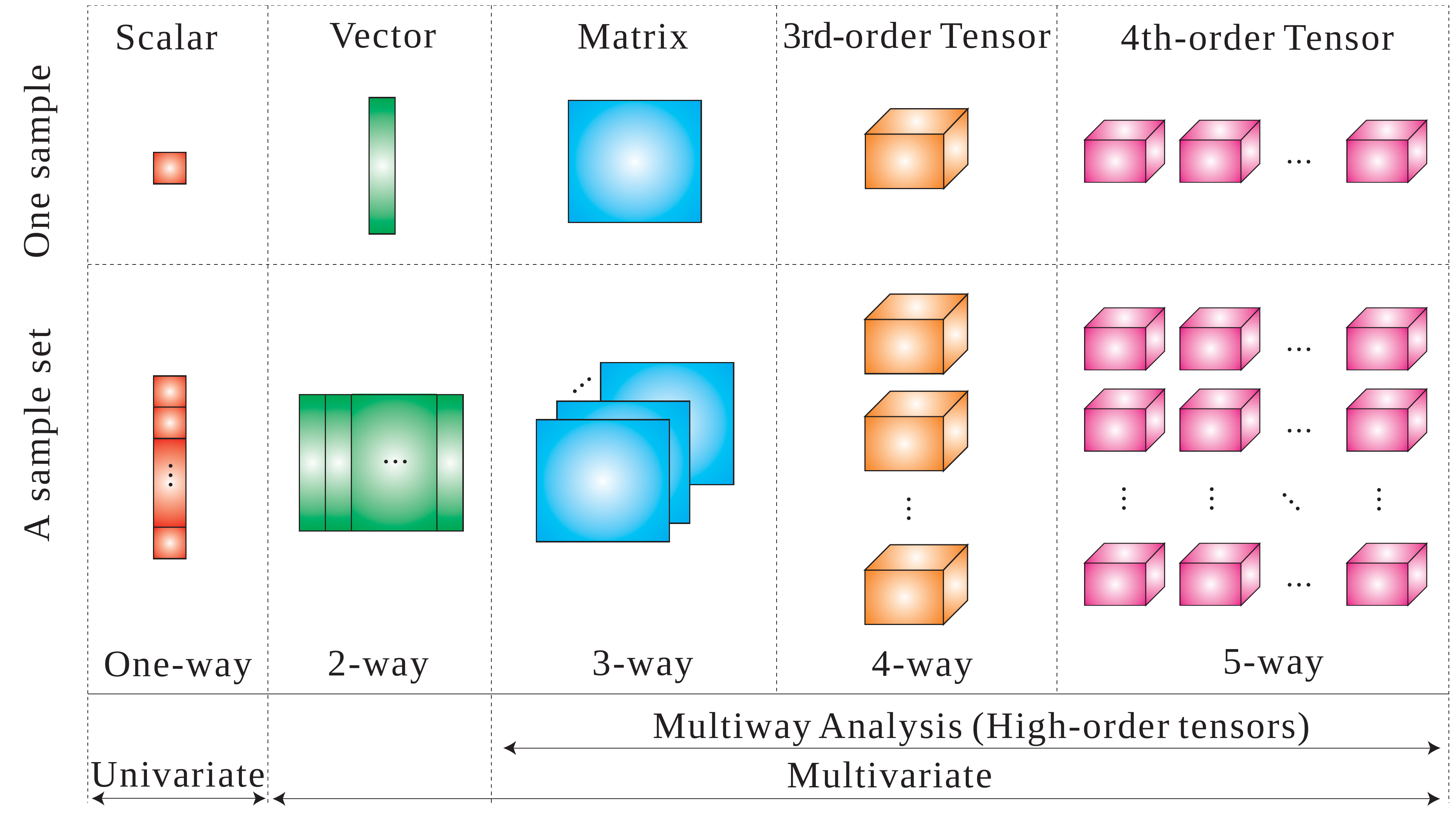}
\caption{Graphical representation of multiway array (tensor) data  of increasing structural complexity and ``Volume''
(see \cite{OlivieriMW2008} for more detail).}
\label{Fig:Multiway}
\end{figure}

\begin{figure}[t!]
(a)\vspace{-0.1cm}
\begin{center}
\includegraphics[width=9.2cm,trim = 0 0.2cm 0 0, clip  = true]{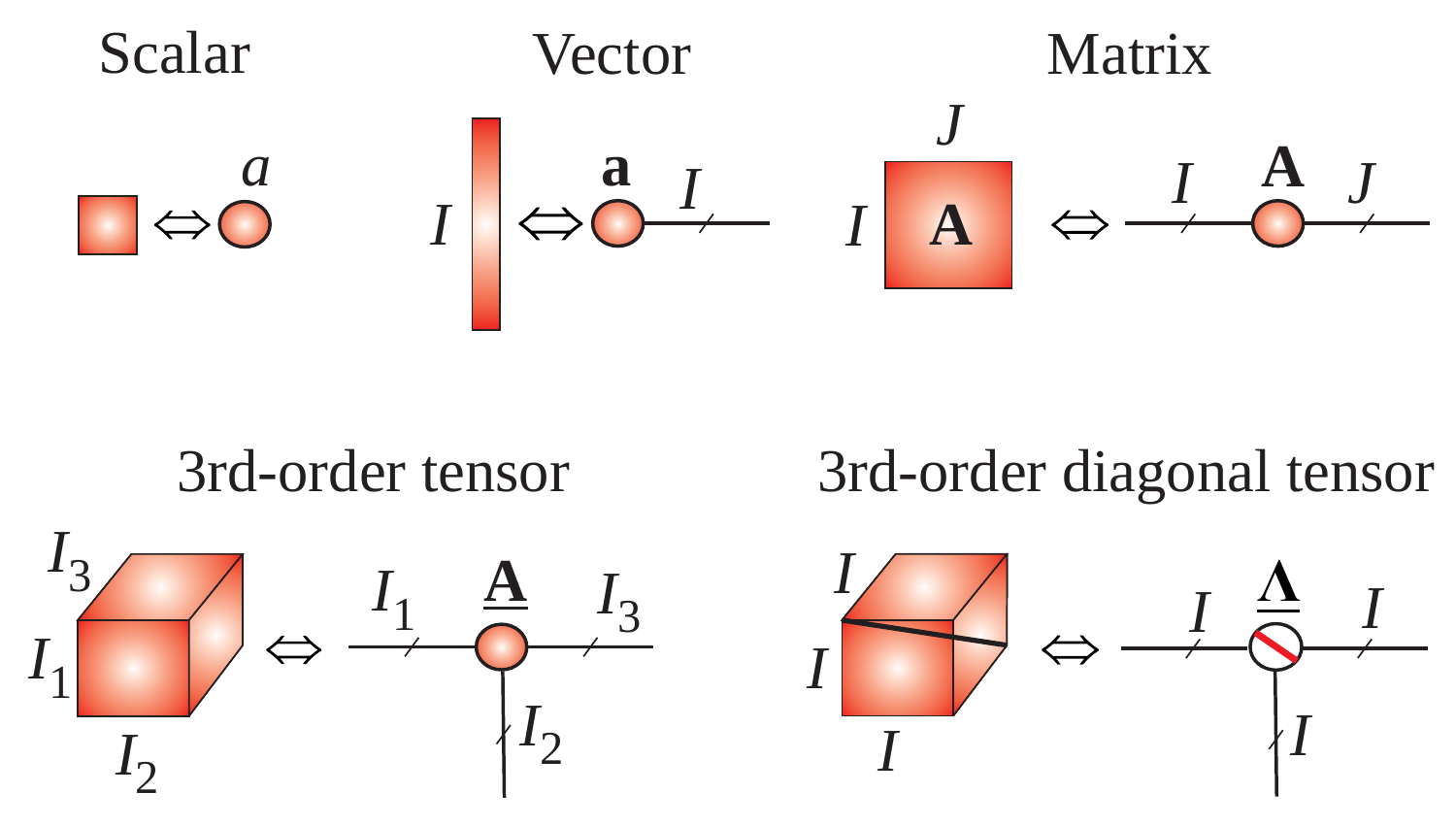}
\end{center}
(b)
\vspace{-0.1cm}
\begin{center}
\includegraphics[width=6.2cm, trim = 0 0.5cm 0 0, clip  = true]{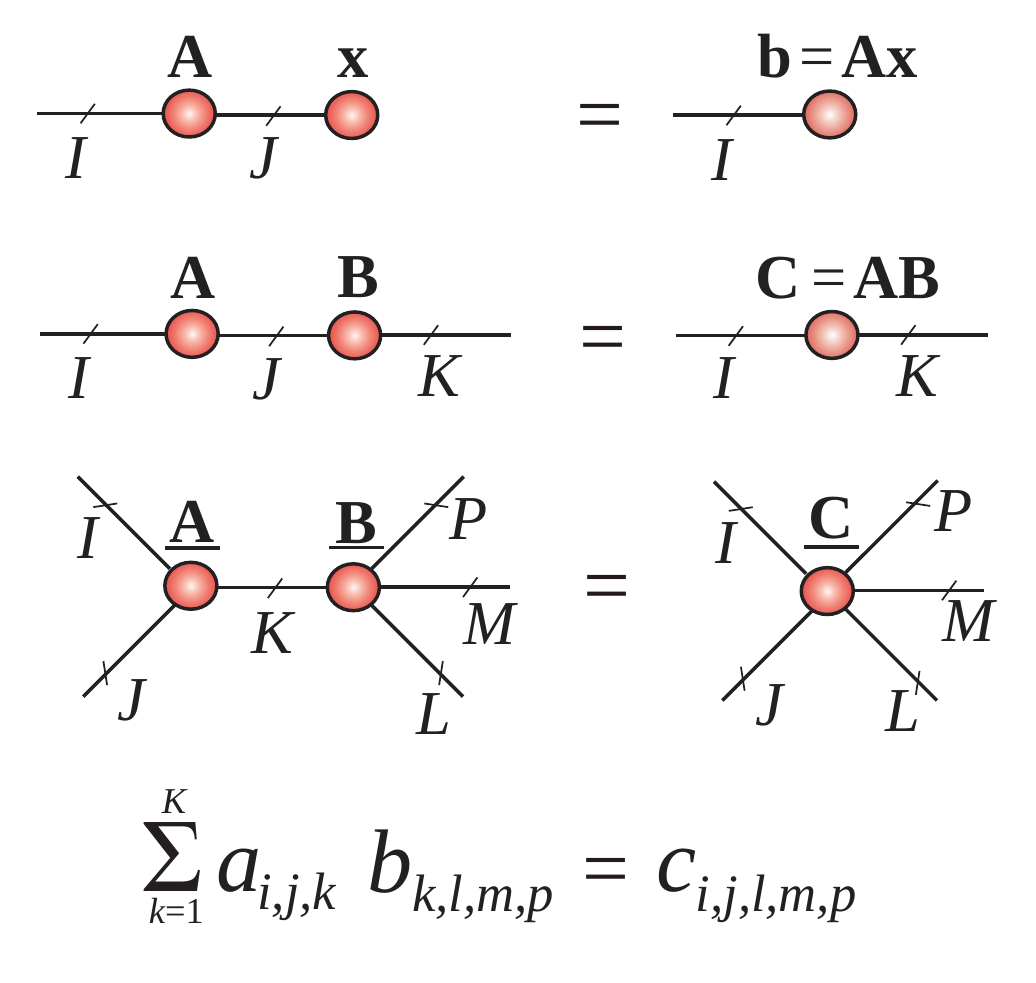}
\end{center}
\vspace{-0.1cm}
\caption{Graphical representation of tensor manipulations. (a) Basic building blocks for tensor network diagrams. (b)  Tensor network diagrams for matrix-vector  multiplication (top),  matrix by matrix multiplication (middle) and  contraction of two tensors (bottom). The order of reading of indices is anti-clockwise, from the left  position.
}
\label{Fig:symbols}
\end{figure}

\begin{figure}[t!]
\centering
4th-order tensor \\
\includegraphics[width=6.2cm]{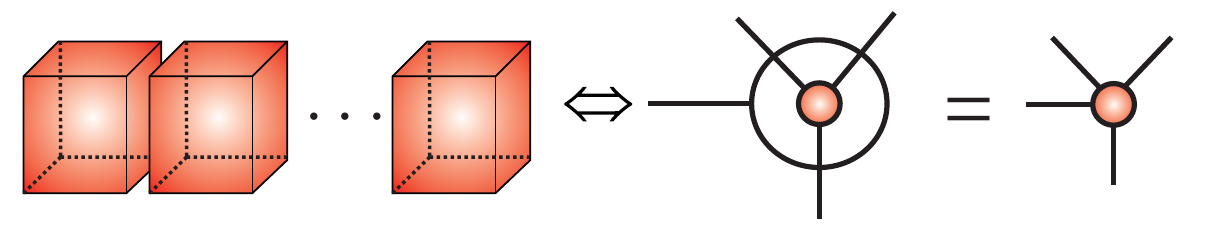}\\
\vspace{0.2cm}
5th-order tensors\\
\includegraphics[width=5.9cm]{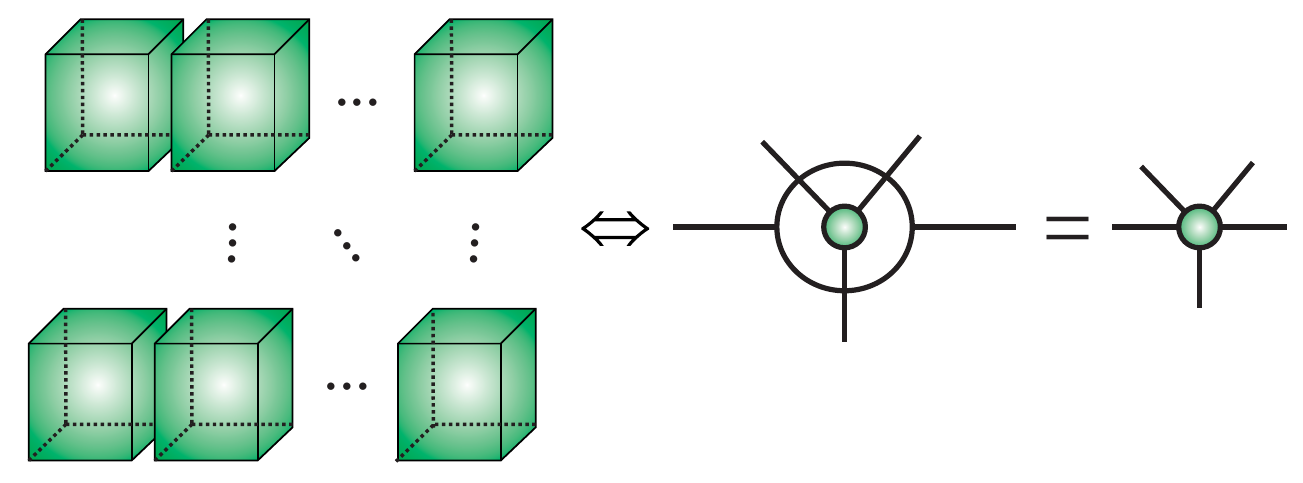}
\hspace{0.1cm}
\includegraphics[width=5.2cm]{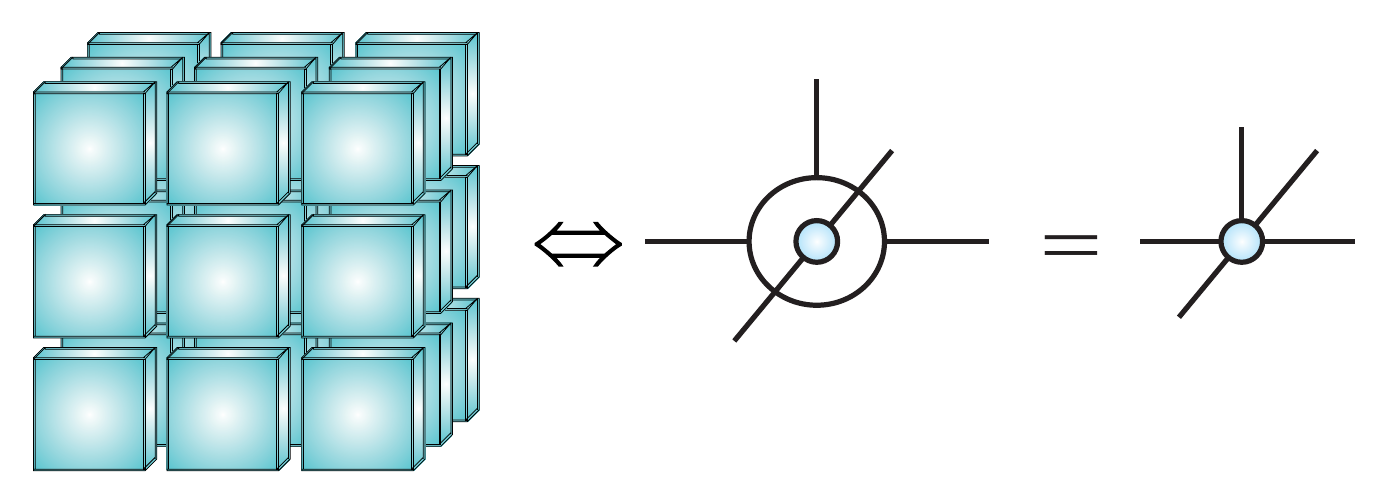}\\
\vspace{0.3cm}
6th-order tensor\\
\vspace{0.3cm}
\includegraphics[width=6.0cm]{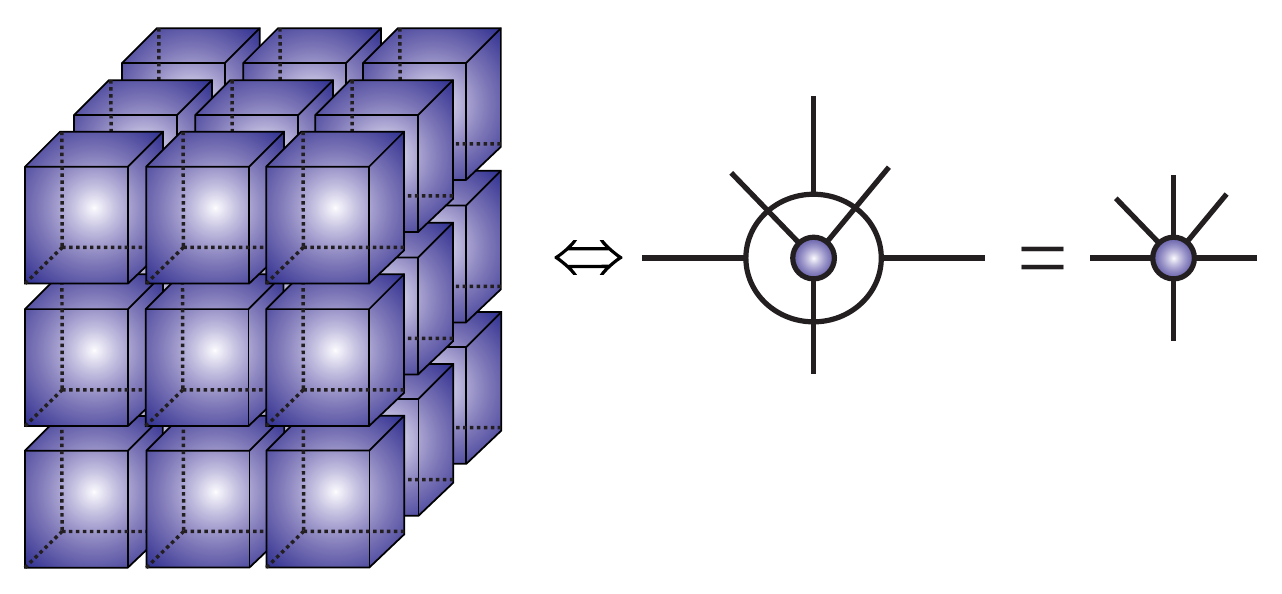}
\caption{Graphical representations and symbols for higher-order block tensors. Each block represents either a 3rd-order tensor or a 2nd-order tensor.
The outer circle indicates a global structure of the block tensor (e.g. a vector, a matrix, a 3rd-order  block tensor), while the inner circle
 reflects the structure of each element within the block tensor. For example, in the top diagram a vector of 3rd order tensors is represented by an outer circle with one edge (a vector) which surrounds an inner circle with three edges (a 3rd order tensor), so that the whole structure designates a 4th-order tensor.}
\label{Fig:symbols2}
\end{figure}

\begin{figure}[t!]
(a)
\vspace{-0.1cm}
\begin{center}
\includegraphics[width=11.0cm]{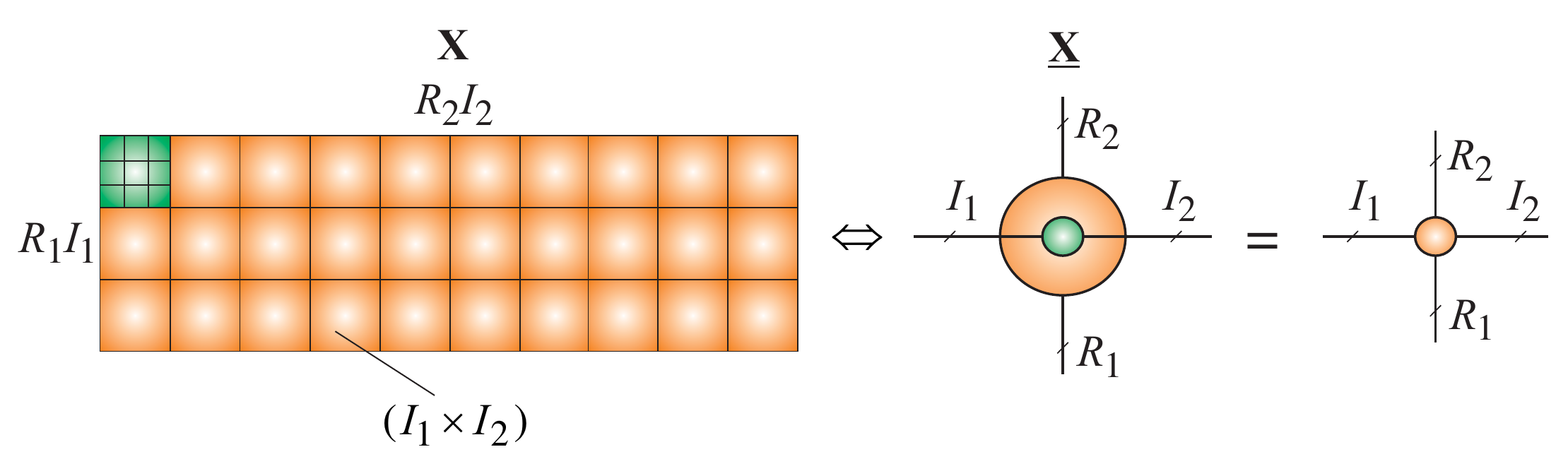}
\end{center}
(b)
\vspace{-0.1cm}
\begin{center}
\includegraphics[width=9.5cm]{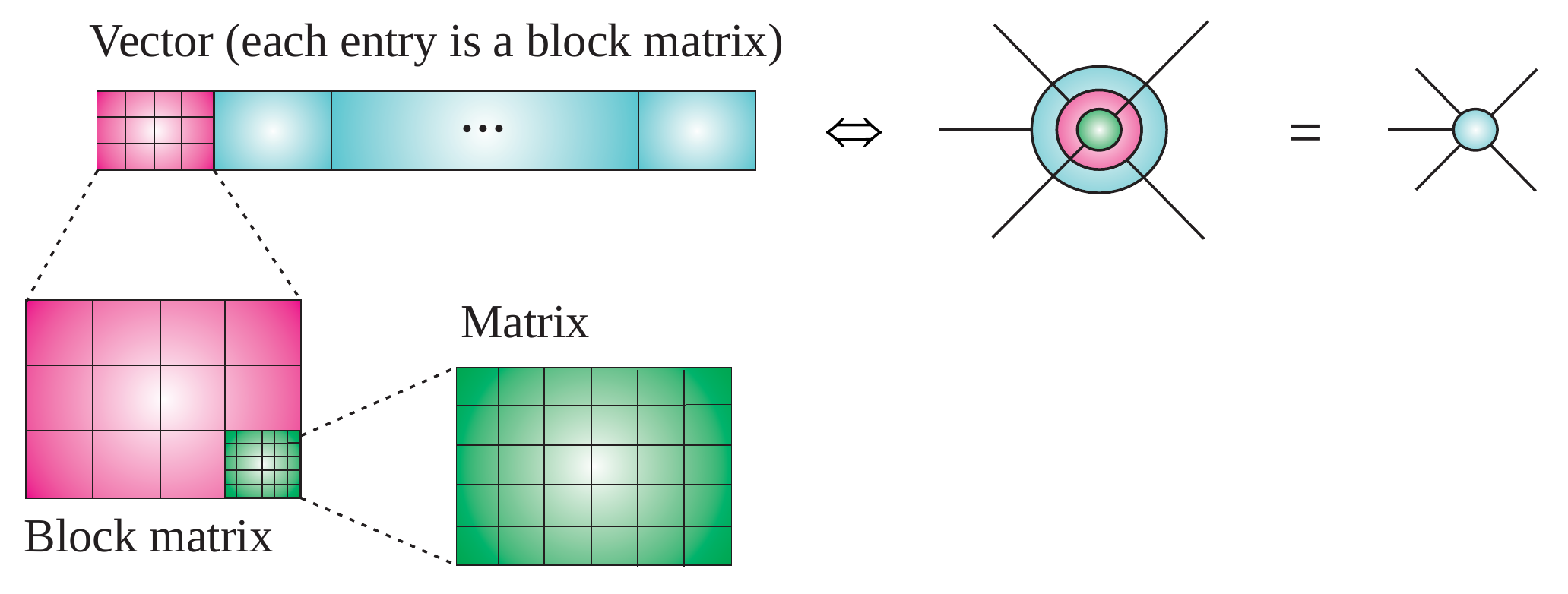}
\end{center}
(c)
\vspace{-0.1cm}
\begin{center}
\hspace{-0.2cm} \includegraphics[width=9.0cm]{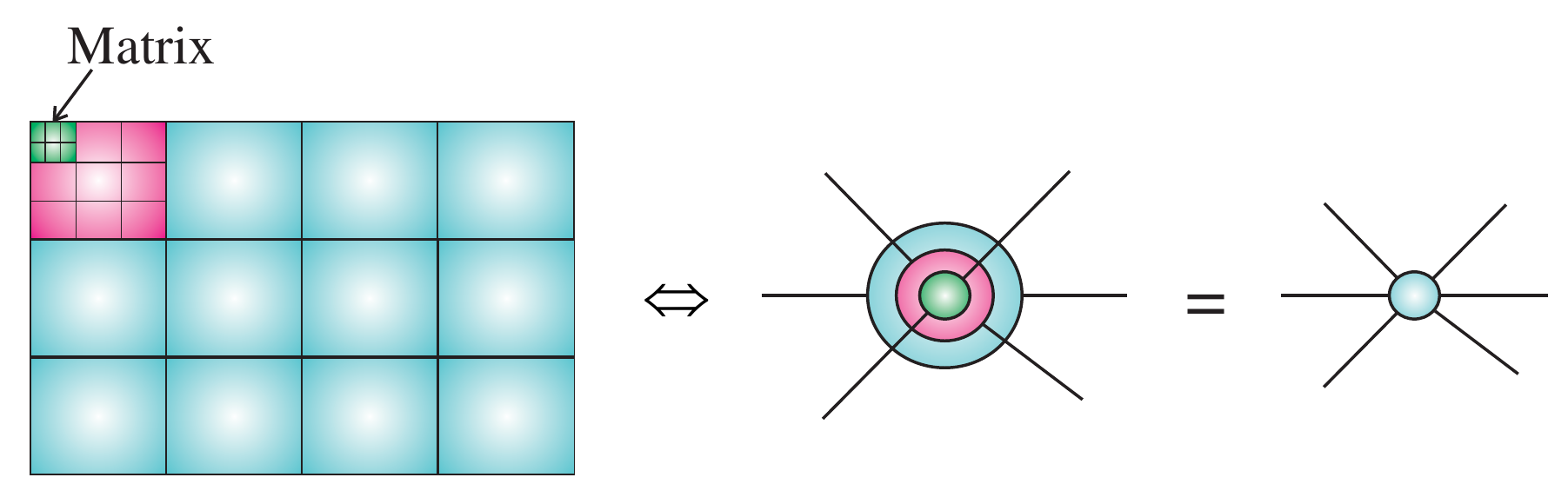}
\end{center}
\caption{Hierarchical  matrix structures and their symbolic representation as tensors.
(a) A 4th-order tensor representation for a block matrix $\bX \in \Real^{R_1 I_1 \times R_2 I_2}$ (a matrix of matrices),  which comprises block matrices $\bX_{r_1,r_2} \in \Real^{I_1 \times I_2}$. (b) A 5th-order tensor. (c) A 6th-order tensor.}
\label{Fig:Tens_block5_8}
\end{figure}

 Tensors are multi-dimensional generalizations of matrices. A matrix (2nd-order tensor) has two modes,
 rows and columns, while an $N$th-order tensor has $N$ modes (see Figures \ref{Fig:fibers}--\ref{Fig:Tens_block5_8});
 for example, a 3rd-order tensor (with three-modes) looks like a cube (see Figure  \ref{Fig:fibers}).
 Subtensors are formed when a subset of tensor indices is fixed.
 Of particular interest are {\it fibers} which are vectors obtained by fixing every tensor index but one, and
{\it matrix slices} which are two-dimensional sections (matrices) of a tensor,
obtained by fixing all the tensor indices but two. It should  be noted that block matrices can also be represented by tensors,
 as illustrated in Figure \ref{Fig:Tens_block_matr} for 4th-order tensors.

We adopt the notation whereby tensors (for $N\geq 3$) are denoted by bold underlined capital letters,
 e.g., $\underline \bX  \in \Real^{I_{1} \times I_{2} \times \cdots \times I_{N}}$. For simplicity, we assume that
 all tensors are real-valued, but it is, of course, possible to define tensors as complex-valued or over arbitrary fields.
Matrices  are denoted by boldface capital letters, e.g., $\bX \in \Real^{I \times J}$, and vectors (1st-order tensors) by boldface lower
case letters, e.g., $\bx \in \Real^J$.  For example, the columns of the matrix $\bA=[\ba_1,\ba_2, \ldots,\ba_R]  \in \Real^{I \times R}$ are the vectors denoted by $\ba_r \in \Real^I$, while the elements of a matrix (scalars) are denoted by lowercase letters, e.g., $a_{ir}=\bA(i,r)$ (see Table \ref{table_notation1}).

A specific entry of an $N$th-order tensor $\underline \bX  \in \Real^{I_{1} \times I_{2} \times \cdots \times I_{N}}$ is  denoted by
$x_{i_1,i_2,\ldots,i_N} = \underline \bX(i_1,i_2, \ldots, i_N) \in \Real$.
The order of a tensor is the number of its ``modes'',  ``ways'' or ``dimensions'', which can include
 space, time, frequency,  trials, classes, and dictionaries.
The term {\emph ``size''} 
stands for the number of values that an index can take in a particular mode.
For example, the tensor $\underline \bX  \in \Real^{I_{1} \times I_{2} \times \cdots \times I_{N}}$ is of order $N$ and size $I_n$ in all modes-$n$ $(n=1,2,\ldots,N)$.
Lower-case letters  e.g, $i, j$  are used for the subscripts in
running indices and capital letters  $I,J$ denote the upper bound of an index, i.e., $i=1,2,\ldots,I$ and $j=1,2,\ldots, J$.
For a positive integer $n$,  the shorthand notation $<n>$  denotes the set of indices $\{1,2,\ldots,n\}$.

\minrowclearance 2ex
\begin{table}
\vspace{-2.5cm}
\caption{Basic matrix/tensor notation and symbols.}
\centering
 {\tabsize \normalsize \shadingbox{
    \begin{tabular*}{\linewidth}[t]{@{\extracolsep{\fill}}@{\hspace{1ex}}l@{\hspace{1ex}}l} \hline
{$\underline \bX \in \Real^{I_1 \times I_2 \times \cdots \times I_N}$} & {$N$th-order tensor of size $I_1\times  I_2 \times \cdots \times I_N$} \\
%
%
{$x_{i_1,i_2,\ldots,i_N}=\underline \bX(i_1,i_2,\ldots,i_N)$} &
{$(i_1,i_2,\ldots,i_N)$th entry of $\underline \bX$} \\
{$x, \; \bx, \; \bX $}& {scalar, vector and matrix} \\ 
{$\underline \bG,  \; \underline \bS, \; \underline \bG^{(n)}, \; \underline \bX^{(n)}$}   & {core tensors} \\  [-2ex]
{$\underline {\mbi {\Lambda}} \in \Real^{R \times R \times \cdots \times R}$}   & \minitab[p{.50\linewidth}]{\hspace{-0.6em} $N$th-order diagonal core tensor with nonzero entries $\lambda_r$ on the main diagonal} \\  [-2ex]
 ${\bA}^{\text{T}}$,   ${\bA}^{-1}$, ${\bA}^{\dag}$ & \minitab[p{.5\linewidth}]{transpose, inverse and Moore--Penrose \hbox{pseudo-inverse} of a matrix $\bA$} \\[-2ex]
{$\bA = [\ba_1,\ba_2,\ldots,\ba_R] \in \Real^{I \times R} $}& \minitab[p{.50\linewidth}]{matrix with $R$ column vectors $\ba_r \in \Real^I$, with entries $a_{ir}$} \\
{$\bA, \, \bB, \,\bC, \; \bA^{(n)}, \bB^{(n)}, \; \bU^{(n)} $}& {component (factor) matrices} \\[1ex]
$\bX_{(n)} \in \Real^{I_n \times I_1 \cdots  I_{n-1} I_{n+1}  \cdots  I_N}$ & {mode-$n$ matricization of $\underline \bX \in \Real^{I_1 \times \cdots \times I_N}$}\\
$\bX_{<n>} \in \Real^{I_1 I_2 \cdots I_n \times  I_{n+1} \cdots I_N}$ &  {mode-($1,\ldots,n$) matricization of $\underline \bX \in \Real^{I_1 \times \cdots \times I_N}$}\\[-1ex]
$\underline \bX({:,i_2,i_3,\ldots,i_N}) \in \Real^{I_1}$ & \minitab[p{.54\linewidth}]{mode-1 fiber of a tensor $\underline \bX$ obtained by fixing all indices but one (a vector)}\\[-3ex]
$\underline \bX({:,:,i_3,\ldots,i_N}) \in \Real^{I_1 \times I_2}$ & \minitab[p{.50\linewidth}]{slice (matrix) of a tensor $\underline \bX$ obtained by fixing all indices
but two}\\[-2ex]
$\underline \bX({:,:,:,i_4,\ldots,i_N})$ & \minitab[p{.50\linewidth}]{subtensor of $\underline \bX$, obtained by fixing several indices}\\[-1ex]
$R, \; (R_1, \ldots, R_N)$& tensor rank $R$ and multilinear rank \\[-2ex]
\minitab[p{.28\linewidth}]{$\circ\;, \; \odot\;, \; \otimes $ \\$ \otimes_{L}\:, \; \skron$}  & \minitab[p{.58\linewidth}] {outer, Khatri--Rao,   Kronecker products \\
Left Kronecker, strong Kronecker products}\\
$\bx=\vtr{\underline \bX}$& vectorization of $\underline \bX$  \\
$\tr( \bullet) $& trace of a square matrix \\
$\diag(\bullet)$& diagonal matrix \\
\hline
    \end{tabular*}
    }}
\label{table_notation1}
\vspace{30pt}
\end{table}
\minrowclearance 0ex

\begin{table}
\setlength{\tabcolsep}{2pt}
\renewcommand{\arraystretch}{1.5}
\centering
\caption{Terminology used for tensor networks across the machine learning / scientific computing  and quantum physics / chemistry communities.}
{
\shadingbox{
\begin{tabular}{p{4.5cm}|p{7cm}}
\hline
Machine Learning &  Quantum Physics \\
\hline
$N$th-order tensor &  rank-$N$ tensor \\

high/low-order tensor &  tensor of high/low dimension\\
ranks of TNs & bond dimensions of TNs\\
unfolding, matricization & grouping of indices \\
tensorization & splitting of indices\\
core   & site  \\
variables & open (physical) indices \\

ALS Algorithm & one-site DMRG or DMRG1\\
MALS Algorithm & two-site DMRG or DMRG2 \\
column vector $\bx \in \Real^{I \times 1}$ & ket $| \Psi \rangle$ \\

row vector $\bx^{\text{T}} \in \Real^{1 \times I}$ & bra  $\langle  \Psi |$ \\

inner product $\langle\bx,\bx \rangle =\bx^{\text{T}} \bx$ &  $\langle  \Psi |  \Psi\rangle $ \\

Tensor Train (TT) & Matrix Product State (MPS) (with Open Boundary Conditions (OBC))\\

Tensor Chain (TC) & MPS with Periodic Boundary Conditions (PBC)\\

Matrix TT & Matrix Product Operators (with OBC)\\

Hierarchical Tucker (HT) & Tree Tensor Network State  (TTNS) with rank-3 tensors\\
\hline
\end{tabular}
}}
\label{Table:ML_QC}
\vspace{30pt}
\end{table}

Notations and terminology used for tensors and tensor networks differ across the scientific communities (see Table~\ref{Table:ML_QC}); to this end we employ a unifying notation particularly suitable for machine learning and signal processing research, which is summarized in Table~\ref{table_notation1}.

Even with the above notation conventions, a precise description of tensors and tensor operations is often tedious and cumbersome, given the multitude of indices involved. To this end, in this monograph, we grossly simplify the description of tensors and their mathematical operations through diagrammatic representations borrowed from physics and quantum chemistry (see \cite{Orus2013} and references therein).
In this way, tensors are represented graphically by nodes of any geometrical shapes (e.g., circles, squares, dots), while each outgoing line (``edge'', ``leg'',``arm'') from a node represents the indices of a specific mode (see Figure~\ref{Fig:symbols}(a)).
In our adopted notation, each scalar (zero-order tensor), vector (first-order tensor), matrix
(2nd-order tensor), 3rd-order tensor or higher-order tensor is represented by a circle (or rectangular), while the order of a
tensor is determined by the number of lines (edges) connected to it.
According to this notation, an $N$th-order tensor $\underline \bX \in \Real^{I_1 \times \cdots \times I_N}$ is represented by a circle (or any shape)  with
$N$ branches each of size $I_n, \; n=1,2, \ldots, N$ (see Section~\ref{chap:basic_operations}).
An interconnection between two circles designates a contraction of tensors, which is a summation of products over a common index
(see Figure~\ref{Fig:symbols}(b) and Section~\ref{chap:basic_operations}).

Block tensors, 
 where each entry (e.g., of  a matrix or a vector) is an individual subtensor,   can  be represented in a similar graphical form, as illustrated in Figure \ref{Fig:symbols2}.
Hierarchical (multilevel block) matrices are also naturally represented by tensors and vice versa, as illustrated in Figure  \ref{Fig:Tens_block5_8} for 4th-, 5th- and 6th-order tensors.
All mathematical operations on tensors can be therefore equally performed on block matrices.

 In this monograph, we make extensive use of tensor network diagrams as an intuitive and visual way to efficiently represent tensor decompositions. Such graphical notations
 are of great help in studying and implementing sophisticated tensor operations. We highlight the significant advantages of such diagrammatic notations in the description of tensor manipulations, and show that most tensor operations can be visualized through changes in the architecture of a tensor network diagram.


%


\section{Curse of Dimensionality and Generalized Separation of Variables for Multivariate Functions}
\markright{\thesection.\quad Curse of Dimensionality}

\subsection{Curse of Dimensionality}

The term {\it curse of dimensionality} was coined by \cite{bellman1961adaptive} to indicate that the number of samples needed to estimate an
arbitrary function with a given level of accuracy grows exponentially with the number
of variables, that is, with the dimensionality of the function.
In a general context of machine learning and the underlying optimization problems, the ``curse of dimensionality'' may also refer to an exponentially increasing number of parameters required to describe the data/system or an extremely large number of degrees of freedom. The term ``curse of dimensionality'', in the context of tensors, refers to the phenomenon whereby the number of elements, $I^N$,
of an $N$th-order tensor of size $(I \times I \times \cdots \times I)$  grows exponentially with the tensor order, $N$.
Tensor volume can therefore easily become prohibitively big for multiway arrays for which the number of dimensions (``ways'' or ``modes'') is very high, thus requiring enormous  computational and memory resources to process such data.
The understanding and handling of the inherent dependencies among the excessive degrees of freedom create both  difficult to solve problems and fascinating new opportunities, but comes at a price of  reduced accuracy, owing to  the necessity to involve various approximations.

We show that the curse of dimensionality can be alleviated or even fully dealt with through tensor network representations; these naturally cater for the  excessive volume, veracity and variety of data (see Figure \ref{Fig:BData}) and are supported by efficient tensor decomposition algorithms which involve  relatively simple mathematical operations.  
Another desirable aspect of tensor networks is their relatively small-scale and low-order \emph{core tensors}, which act as ``building blocks'' of tensor networks. These core tensors are relatively easy to handle and visualize, and enable  super-compression of the  raw, incomplete, and noisy huge-scale datasets.
This also suggests a solution to a more general quest for  new technologies for processing of  exceedingly large datasets within affordable computation times. 
%

 To address the curse of dimensionality, this work mostly focuses on
approximative low-rank representations of tensors, the so-called low-rank tensor approximations (LRTA) or low-rank tensor network decompositions.

\section{Separation of Variables and Tensor Formats}

A tensor is said to be in a \emph{full format} when it is represented as an original (raw) multidimensional array \cite{Klus-EVD}, however,  distributed storage and processing of high-order tensors in their full format is infeasible due to the curse of dimensionality. The \emph{sparse  format} is a variant of the full tensor format which stores only the nonzero entries of a tensor, and is used extensively in software tools such as the Tensor Toolbox \cite{KoldaTensorToolbox} and in the sparse grid approach \cite{garcke2001data,bungartz2004sparse,Hackbush2012}.

 As already mentioned, the problem of  huge dimensionality can  be alleviated  through various distributed and compressed tensor network formats,  achieved by low-rank tensor network approximations.
  The underpinning idea is that by employing tensor networks formats, both computational costs and storage requirements  may be dramatically reduced through distributed storage and computing resources. It is important to note that, except for very special data structures, a tensor cannot be compressed without incurring some compression error, since a low-rank tensor representation is only an approximation of the original tensor.

The concept of compression of multidimensional large-scale data by tensor network decompositions can be intuitively explained as follows.
 Consider  the approximation of an $N$-variate function
$f(\bx)= f(x_1,x_2,\ldots,x_N)$  by  a finite sum of products of individual functions,  each depending on only one or a very few variables \cite{bebendorf2011adaptive,Dolgovth,cho2016numerical,Trefethen}.
In the simplest scenario, the function $f(\bx)$ can be (approximately) represented in the following separable form
\be
f(x_1,x_2,\ldots,x_N) \cong f^{(1)}(x_1) f^{(2)} (x_2) \cdots f^{(N)}(x_N).
\label{simp-separability}
\ee
In practice, when an $N$-variate function $f(\bx) $ is discretized into an $N$th-order array, or a tensor, the approximation in (\ref{simp-separability}) then corresponds to the representation by  rank-1 tensors, also called  elementary tensors
(see Section~\ref{chap:basic_operations}). Observe that with $I_n, \;n=1,2,\ldots,N$ denoting the size of each mode
 and $I=\max_n \{I_n\}$, the memory requirement to store such a full tensor is $\prod_{n=1}^N I_n \leq I^N$, which grows exponentially with $N$.
 On the other hand, the separable representation in (\ref{simp-separability}) is completely defined by its factors,
 $f^{(n)}(x_n),  \;(n=1,2, \ldots,N$), and requires only $\sum_{n=1}^N I_n \ll I^N$ storage units.
 If $x_1,x_2, \ldots, x_N$
are statistically independent random  variables, their joint probability density function is equal to the product of marginal probabilities, $f(\bx) = f^{(1)}(x_1) f^{(2)} (x_2) \cdots f^{(N)}(x_N)$,  in an exact analogy to outer products of elementary tensors.
Unfortunately, the form of separability in (\ref{simp-separability}) is rather rare  in practice.

The concept of tensor networks rests upon  generalized (full or  partial) separability of the variables
of a high dimensional function.  This  can be achieved in different tensor formats, including:
\begin{itemize}

\item The Canonical Polyadic (CP) format (see Section~\ref{sect:CPD}), where
\begin{equation}
f(x_1,x_2,\ldots,x_N)\cong \sum_{r=1}^R f_{r}^{(1)}(x_1) \, f_{r}^{(2)} (x_2) \cdots f_{r}^{(N)}(x_N),
\label{eq:CPfun}
\end{equation}
in an exact analogy to (\ref{simp-separability}). In a discretized form, the above CP format  can be written as an $N$th-order  tensor
\be
\underline \bF \cong \sum_{r=1}^R {\boldf}_{\;r}^{(1)} \circ {\boldf}_{\;r}^{(2)} \circ \cdots \circ  {\boldf}_{\;r}^{(N)}
\in \Real^{I_1 \times I_2 \times \cdots \times I_N},
\label{eq:CPfun2}
\ee
where ${\boldf}_{\;r}^{(n)} \in \Real^{I_n}$  denotes a discretized version  of the univariate function $f_r^{(n)}(x_n)$, symbol $\circ$  denotes the outer product, and $R$ is the tensor rank.
\item The Tucker format, given  by
\begin{equation}
f(x_1,\ldots,x_N) \cong \sum_{r_1=1}^{R_1} \cdots \sum_{r_N=1}^{R_N} g_{r_1,\ldots,r_N} \; f_{r_1}^{(1)}(x_1)  \cdots f_{r_N}^{(N)}(x_N), \label{eq:HT_function} 
\end{equation}
and its distributed  tensor network variants  (see Section~\ref{sect:Tucker}),

\item The Tensor Train (TT) format (see Section~\ref{sect:TT}), in the form
\be
f(x_1,x_2,\ldots,x_N) &\cong& \sum_{r_1=1}^{R_1} \sum_{r_2=1}^{R_2} \cdots  \sum_{r_{N-1}=1}^{R_{N-1}} f_{r_1}^{(1)}(x_1) \, f_{r_1 \, r_2}^{(2)} (x_2)  \cdots \nonumber  \\
&\cdots& f_{r_{N-2} \, r_{N-1}}^{(N-2)} (x_{N-1})\; f_{r_{N-1}}^{(N)}(x_N),
 \label{eq:TT_def}
\ee
with the equivalent compact matrix representation
\be
f(x_1,x_2,\ldots,x_N) \cong \bF^{(1)} (x_1) \, \bF^{(2)} (x_2) \cdots \bF^{(N)} (x_N), 
\ee
where $\bF^{(n)}(x_n) \in \Real^{R_{n-1} \times R_n}$, with $R_0=R_N=1$.

\item The Hierarchical Tucker (HT) format  (also known as the Hierarchical Tensor format)
can be expressed via a hierarchy of nested separations in the following way. Consider nested nonempty disjoint subsets $u$, $v$, and $t  = u \cup v  \subset \{1,2,\ldots,N\}$, then for  some $1\leq N_0<N$,
with $u_0 = \{1,\ldots,N_0\}$ and $v_0 = \{N_0+1,\ldots,N\}$, the HT format can be expressed as
	\begin{eqnarray}
	f(x_1,\ldots,x_N)
	&\cong& \sum_{r_{u_0}=1}^{R_{u_0}} \sum_{r_{v_0}=1}^{R_{v_0}}
	g^{(12\cdots N)}_{r_{u_0},r_{v_0}} \,
	f^{({u_0})}_{r_{u_0}}( \bx_{u_0} ) \, f^{({v_0})}_{r_{v_0}}( \bx_{v_0} ) , \notag
	\\
	f^{(t)}_{r_t} ( \bx_t )
	&\cong&  \sum_{r_u=1}^{R_u} \sum_{r_v=1}^{R_v}
	g^{(t)}_{r_u,r_v,r_t} \,
	f^{(u)}_{r_u}( \bx_u ) \, f^{(v)}_{r_v}( \bx_v ),  \notag
	\end{eqnarray}
where $\bx_t = \{x_i :  i \in t\}$.
See Section~\ref{sect:HT} for more detail.

{\bf Example.}  In a particular case for $N$=4,
the HT format can be expressed by
\begin{eqnarray}
	f(x_1,x_2,x_3,x_4) &\cong&
		\sum_{r_{12}=1}^{R_{12}} \sum_{r_{34}=1}^{R_{34}}
		g^{(1234)}_{r_{12},r_{34}} \,
		f^{(12)}_{r_{12}} (x_1, x_2) \,
		f^{(34)}_{r_{34}} (x_3, x_4) , \notag \\
	f^{(12)}_{r_{12}} (x_1, x_2) &\cong&
		\sum_{r_{1}=1}^{R_{1}} \sum_{r_{2}=1}^{R_{2}}
		g^{(12)}_{r_{1},r_{2},r_{12}} \,
		f^{(1)}_{r_{1}} (x_1) \,
		f^{(2)}_{r_{2}} (x_2) , \notag \\
	f^{(34)}_{r_{34}} (x_3, x_4) &\cong&
		\sum_{r_{3}=1}^{R_{3}} \sum_{r_{4}=1}^{R_{4}}
		g^{(34)}_{r_{3},r_{4},r_{34}}\,
		f^{(3)}_{r_{3}} (x_3) \,
		f^{(4)}_{r_{4}} (x_4) .\notag
	\end{eqnarray}
%

The Tree Tensor Network States (TTNS) format,  which is an extension of the HT format,  can be obtained by generalizing the two subsets, $u,v$, into a larger number of disjoint subsets $u_1,\ldots,u_m$, $m\geq 2$. In other words, the TTNS can be obtained by more flexible separations of variables through products of larger numbers of functions at each hierarchical level (see Section~\ref{sect:HT} for graphical illustrations and more detail).
\end{itemize}
All the above approximations adopt the form of ``sum-of-products'' of single-dimensional functions, a procedure which plays a key role in all tensor factorizations and  decompositions.

 Indeed, in many applications based on multivariate functions, very good approximations are obtained with  a surprisingly small number of factors; this number corresponds to the
 tensor rank, $R$, or tensor network ranks, $\{R_1,R_2,\ldots, R_N\}$ (if the representations are exact and minimal).
 However, for some specific  cases  this approach  may fail to obtain  sufficiently good low-rank TN approximations.
%
The concept of generalized separability has already been explored in
 numerical methods for high-dimensional density function equations \cite{TPA2015,Trefethen,cho2016numerical} and within a  variety of  huge-scale optimization problems (see Part 2 of this monograph).

To illustrate how tensor decompositions address excessive volumes of data, if all computations are performed on a CP tensor format in  (\ref{eq:CPfun2}) and not on the raw $N$th-order data tensor itself, then instead of the original, \emph{exponentially growing}, data dimensionality of  $I^N$, the number of parameters in a CP representation reduces to
$N I R$, which \emph{scales linearly} in the tensor order $N$ and size $I$ (see Table \ref{table_complexity}). For example, the discretization of a $5$-variate function over $100$ sample points
on each axis would yield the difficulty to manage $100^5 = 10,000,000,000$ sample points, while a rank-$2$ CP representation would require only $5 \times 2 \times 100 = 1000$ sample points.

Although the CP format in (\ref{eq:CPfun}) effectively bypasses the curse of dimensionality, the CP approximation may involve numerical problems for very high-order tensors, which in addition to the intrinsic uncloseness of the CP format (i.e., difficulty to arrive at a canonical format), the corresponding algorithms for CP decompositions are often ill-posed \cite{deSilva-Lim08}.
 As a remedy, greedy approaches may be considered which, for enhanced  stability,  perform consecutive rank-1 corrections \cite{Lim2010}.
On the other hand, many efficient and stable algorithms exist for the more  flexible Tucker format in (\ref{eq:HT_function}),  however, this format is 
 not practical for tensor orders $N>5$  because the number of entries of both the original data tensor and the  core tensor
(expressed in (\ref{eq:HT_function}) 
by elements $g_{r_1,r_2,\ldots,r_N} $)
scales exponentially in the tensor order $N$ (curse of dimensionality).

In contrast to CP decomposition algorithms, TT tensor network formats in (\ref{eq:TT_def}) exhibit both very good numerical properties
and the ability to control the error of approximation, so that 
a desired accuracy of approximation is obtained relatively easily.
%
The main advantage of the TT format over the CP decomposition is the  ability to provide stable quasi-optimal rank reduction, 
achieved through, for example, truncated singular value decompositions (tSVD) or adaptive cross-approximation
\cite{oseledets2010tt,bebendorf2011adaptive,Khoromskij16efficient}.
This makes the TT format one of the most stable and simple approaches to separate latent variables in a sophisticated way,
 while the associated TT decomposition algorithms  provide full control over  low-rank TN approximations{\footnote{Although similar
approaches have been known in quantum physics for a long time, their rigorous mathematical analysis is still a work in progress (see \cite{OseledetsTT11,Orus2013} and references therein).}}.
%
%
In this monograph, we  therefore make extensive use of the TT format for low-rank TN approximations and employ the TT toolbox software for efficient implementations \cite{oseledets2012tt}.
The TT format will also serve as a basic prototype for high-order tensor representations, while we also consider the Hierarchical Tucker (HT) and the Tree Tensor Network States (TTNS) formats (having more general tree-like structures) whenever advantageous in applications.

Furthermore, we address in depth the concept of 
tensorization of structured vectors and matrices  to convert a wide class
of huge-scale optimization problems into much smaller-scale interconnected optimization sub-problems which can be solved by existing optimization methods (see Part 2 of this monograph).

  The tensor network  optimization framework is  therefore performed through the  two main steps:

 \begin{itemize}

\item  Tensorization of data vectors and matrices into a high-order tensor, followed by a  distributed approximate representation of a cost function in a specific
low-rank tensor network format.

\item Execution of all computations and analysis in tensor  network formats (i.e., using only core tensors) that scale linearly, or even sub-linearly (quantized tensor networks),  in  the tensor order $N$. This yields both  the reduced computational complexity  and distributed
memory requirements.


 \end{itemize}


 \section{Advantages of  Multiway Analysis  via Tensor Networks}

In this monograph, we focus on  two main challenges in huge-scale data analysis which are  addressed by  tensor networks: (i) an approximate representation of a specific  cost (objective) function by a tensor network while maintaining the desired accuracy of approximation, and (ii) the  extraction of  physically meaningful latent variables from data in
a  sufficiently accurate and computationally affordable way.
The benefits of multiway (tensor) analysis methods  for  large-scale datasets then include:
\begin{itemize}

\item  Ability to perform all mathematical operations in tractable tensor network formats; 

\item  Simultaneous and flexible distributed representations of both the  structurally rich data and complex optimization tasks;

\item Efficient compressed formats of large multidimensional data achieved  via tensorization and low-rank tensor decompositions  into low-order factor matrices and/or  core tensors;


 \item Ability to operate with noisy and  missing  data  by virtue of numerical stability and robustness
 to noise of low-rank tensor / matrix approximation algorithms;

\item  A flexible  framework which naturally  incorporates various diversities and constraints,  thus seamlessly extending the standard, flat view, Component Analysis (2-way CA) methods to  multiway  component analysis;

\item Possibility to analyze linked (coupled) blocks of large-scale matrices and tensors in order to separate  common / correlated from  independent / uncorrelated
components  in the observed  raw data;

\item Graphical representations of tensor networks allow us to express mathematical operations on tensors (e.g., tensor contractions and reshaping)  in a simple and intuitive way, and without the explicit use of complex  mathematical expressions.
\end{itemize}
In that sense, this monograph both reviews current research in this area and complements optimisation methods,
such as the Alternating Direction Method of Multipliers (ADMM) \cite{BoydFTML2011}.

Tensor decompositions (TDs) have been already adopted in  widely diverse disciplines, including 
 psychometrics, chemometrics, biometric, quantum physics / information, quantum chemistry, signal and image processing, machine learning,  and brain science
\cite{Smilde,tao2007general,Kroonenberg,Kolda08,Hackbush2012,Favier-deAlmeida14,NMF-book,Cich-Lath}.
This is largely due to their advantages  in the analysis of data that exhibit not only large volume but also very high variety (see Figure \ref{Fig:BData}), as in   the case in bio- and neuroinformatics and in computational neuroscience, where various forms of  data collection
 include sparse tabular structures and graphs or hyper-graphs.

Moreover, tensor networks have the ability to efficiently parameterize, through structured compact representations, very general high-dimensional spaces
 which arise in modern   applications
 \cite{KressnerSV2013,Cichocki2014optim,zhang-Osel15Anova,Zorin-TT,Litsarev16Integral,Khoromskij16efficient,Benner16reduced}.
%
Tensor networks also naturally account for intrinsic multidimensional and distributed patterns present in data, and thus provide the opportunity to develop very sophisticated models for capturing multiple interactions  and couplings in data -- these are more physically insightful and interpretable  than   standard pair-wise interactions.

\section{Scope and Objectives}

Review and tutorial papers \cite{Kolda08,MSLSurvey2011,Grasedyck-rev,Cich-Lath,Favier-deAlmeida16,Sidiro-Lath2016,papalexakis2016tensors,Bachmayr2016}  and books \cite{Smilde,Kroonenberg,NMF-book,Hackbush2012} dealing with TDs and TNs already exist, however, they typically focus on standard models, with no explicit links to very large-scale data processing topics or connections to a wide class of optimization problems.
The aim of this  monograph is therefore to extend beyond the standard  Tucker and CP tensor decompositions, and to demonstrate the perspective of TNs in extremely large-scale data analytics,
together with their role as a mathematical backbone in the discovery of hidden structures in prohibitively  large-scale data. Indeed, we show that TN models provide a framework for the analysis of linked (coupled) blocks of tensors with millions and even billions of non-zero entries.

We also demonstrate that TNs provide natural extensions of  2-way (matrix) Component Analysis (2-way CA) methods to multi-way component analysis (MWCA), which deals with the extraction of desired components from multidimensional  and multimodal data.
This paradigm shift requires  new models and associated algorithms capable of identifying core relations among the different tensor modes, while guaranteeing  linear / sub-linear scaling with the size of datasets{\footnote{Usually, we  assume that  huge-scale problems operate on at least
 $10^7$ parameters.}}.

Furthermore, we review tensor decompositions and the  associated algorithms for very large-scale linear / multilinear dimensionality reduction problems.
The related  optimization problems  often involve structured matrices and vectors with  over a billion entries
 (see \cite{Grasedyck-rev,Dolgovth,Garreis2016} and references therein).
In particular, we focus on Symmetric Eigenvalue Decomposition (EVD/PCA) and Generalized Eigenvalue Decomposition (GEVD) \cite{dolgovEIG2013,KressnerEIG2014,KressnerEVD2016},
 SVD  \cite{Lee-SIMAX-SVD},
solutions of overdetermined and undetermined  systems of linear algebraic equations  \cite{Oseledets-Dolgov-lin-syst12,DolgovAMEN2014},
 the Moore--Penrose pseudo-inverse of structured matrices \cite{Lee-TTPI}, and Lasso problems \cite{Lee-Lasso}.
 Tensor networks for extremely large-scale multi-block (multi-view) data are also discussed, especially TN models for orthogonal Canonical Correlation
  Analysis (CCA) and related Partial Least Squares (PLS) problems. 
  For convenience, all these problems are reformulated as constrained optimization problems which are then, by virtue  of low-rank tensor
  networks reduced to manageable lower-scale optimization sub-problems. The enhanced tractability and scalability is achieved through tensor network contractions and other tensor network transformations.

  The methods and approaches discussed in this work can be considered a both  an alternative and  complementary
 to other emerging methods for huge-scale
optimization problems like random coordinate  descent (RCD) scheme \cite{Nesterov2012,PCoDM2015}, sub-gradient methods
\cite{Nesterov2014subgradient}, alternating direction method of multipliers (ADMM) \cite{BoydFTML2011},  and proximal gradient descent methods \cite{Proximal2014} (see also \cite{Cevher2014,Hong2016} and references therein).

This monograph  
systematically introduces  TN models and the associated algorithms for TNs/TDs and illustrates  many potential applications of TDs/TNS.
 The dimensionality reduction and optimization frameworks (see Part 2 of this monograph) are considered in detail, and we also illustrate the use of TNs in  other challenging problems for huge-scale datasets which can be  solved using the tensor network approach,
including  anomaly detection, tensor completion, compressed sensing, clustering,  and classification.

\chapter{Tensor Operations and Tensor Network Diagrams}
\chaptermark{Tensor Operations and Tensor Network Architectures}
\label{chap:basic_operations}

\vspace{0cm}

Tensor operations benefit from the power  of multilinear algebra which is structurally much richer than linear algebra, and even some basic properties, such as the rank, have a more complex meaning. We next introduce the background on fundamental mathematical operations in multilinear algebra, a prerequisite for the understanding of higher-order tensor decompositions. A unified account of both the  definitions and properties of  tensor network operations is provided, including
the outer,  multi-linear, Kronecker, and  Khatri--Rao  products.
 For clarity,  graphical illustrations are provided,  together with an example rich guidance  for tensor network operations
 and their
 properties.  To avoid any confusion that may arise given the numerous options on  tensor reshaping, both mathematical operations and their properties are expressed directly in their native multilinear contexts,   supported by graphical visualizations.

\minrowclearance 2ex
\begin{table}
\vspace{-2.0cm}
\caption{Basic tensor/matrix operations.} \centering
  {\shadingbox{
    \begin{tabular*}{1.00\linewidth}[t]{@{\extracolsep{\fill}}@{\hspace{2ex}}ll} \hline
\\[-3em]
$\underline \bC = \underline \bA \times_n \bB$ & \minitab[p{.6\linewidth}]{\\[-3.2ex] Mode-$n$ product of a tensor $\underline \bA \in \Real^{I_1 \times I_2 \times \cdots \times I_N}$ and a matrix $\bB \in \Real^{J \times I_n}$ yields a tensor
$\underline \bC \in  \Real^{I_1 \times \cdots \times I_{n-1} \times J \times I_{n+1} \times \cdots \times I_N}$,
with entries $c_{\,i_1, \ldots, i_{n-1}, \, j, \, i_{n+1}, \ldots, i_N} = \sum_{i_n=1}^{I_n} \; a_{i_1, \ldots, i_n, \ldots, i_N} \, b_{j, \, i_n}$} \\ [-2ex]
%
\\[-1em]
$\underline \bC = \llbracket \underline \bG; \bB^{(1)},  \ldots, \bB^{(N)}\rrbracket$ &\minitab[p{.61\linewidth}]
{Multilinear (Tucker) product of a core tensor, $\underline \bG$, and factor matrices $\bB^{(n)}$, which gives \\ $\underline \bC = \underline \bG \times_1 \bB^{(1)} \times_2 \bB^{(2)} \cdots \times_N \bB^{(N)}$}\\  [-2ex]
\\[-1em]
$\underline \bC = \underline \bA \,\bar{\times}_n \; \bb$ &\minitab[p{.58\linewidth}]
	{Mode-$n$  product of a tensor $\underline \bA \in \Real^{I_1\times\cdots\times I_N}$ and
	vector $\bb \in \Real^{I_n}$ yields a tensor $\underline \bC \in
   \Real^{I_1\times\cdots\times I_{n-1}\times I_{n+1}\times\cdots \times I_N}$,  with entries
   $c_{\,i_1,\ldots,i_{n-1},i_{n+1},\ldots,i_N}
   =\sum_{i_n=1}^{I_n} a_{i_1, \ldots, i_{n-1}, i_n, i_{n+1},\ldots,i_N} \; b_{i_n}$} \\  [-2ex]
$\underline \bC = \underline \bA \times_N^1 \underline \bB=\underline \bA \times^1 \underline \bB$ & \minitab[p{.6\linewidth}] {Mode-$(N,1)$ contracted product of
		tensors $\underline \bA \in\Real^{I_1\times I_2 \times \cdots\times I_N}$ and
	$\underline \bB \in \Real^{J_1 \times J_2 \times \cdots \times J_M}$,
	   with $I_N=J_1$, yields a tensor $\underline \bC \in \Real^{I_1\times\cdots\times I_{N-1}\times J_2\times\cdots \times J_M}$
		with entries  $c_{i_1,\ldots,i_{N-1},j_2,\ldots,j_M} =
		\sum_{i_N=1}^{I_N} a_{i_1,\ldots,i_N} \; b_{i_N,j_2,\ldots,j_M}$}\\ [-3ex]
$\underline \bC = \underline \bA \circ \underline \bB$ & \minitab[p{.61\linewidth}]{\\[-3.2ex] Outer product of tensors $\underline \bA \in \Real^{I_1 \times I_2 \times \cdots \times I_N}$ and $\underline \bB \in \Real^{J_1 \times J_2 \times \cdots \times J_M}$ yields an $(N+M)$th-order tensor $\underline \bC$, with entries $c_{\; i_1, \ldots, i_N, \,j_1, \ldots, j_M} = a_{i_1, \ldots, i_N} \; b_{j_1, \ldots, j_M}$} \\  [-2ex]
{$\underline \bX = \ba \circ \bb \circ \bc \in \Real^{I \times J \times K}$}&\minitab[p{.58\linewidth}]{Outer product of vectors $\ba,  \bb$ and  $\bc$ forms a rank-1 tensor, $\underline \bX$, with entries $x_{ijk} = a_i \; b_j \; c_k$} \\  [-2ex]
{$\underline \bC = \underline \bA \otimes_L \underline \bB$}&\minitab[p{.6\linewidth}]{(Left) Kronecker product of tensors $\underline \bA \in
\Real^{I_1 \times I_2 \times \cdots \times I_N}$ and $\underline \bB \in \Real^{J_1 \times J_2 \times \cdots \times J_N}$ yields a tensor
$ \underline \bC \in \Real^{I_1 J_1 \times \cdots \times I_N J_N}$, with entries
$c_{\;\overline{i_1 j_1}, \ldots, \overline{i_N j_N}} = a_{i_1, \ldots,i_N} \: b_{j_1, \ldots,j_N}$} 
\\[-2ex]
{$\bC = \bA\odot_L \bB$}&\minitab[p{.6 \linewidth}]{(Left) Khatri--Rao product of matrices $\bA=[\ba_1, \ldots, \ba_J] \in \Real^{I\times J}$ and $\bB = [\bb_1, \ldots, \bb_J] \in \Real^{K\times J}$ yields a matrix $\bC \in \Real^{I K \times J}$, with columns $\bc_j = \ba_j \otimes_L \bb_j \in \Real^{IK}$}\\[-1ex]
\\[-6pt] \hline
    \end{tabular*}
    }}
\label{table_notation2}
\end{table}
\minrowclearance 0ex

\begin{figure}[t]
\centering
\includegraphics[width=7.6cm]{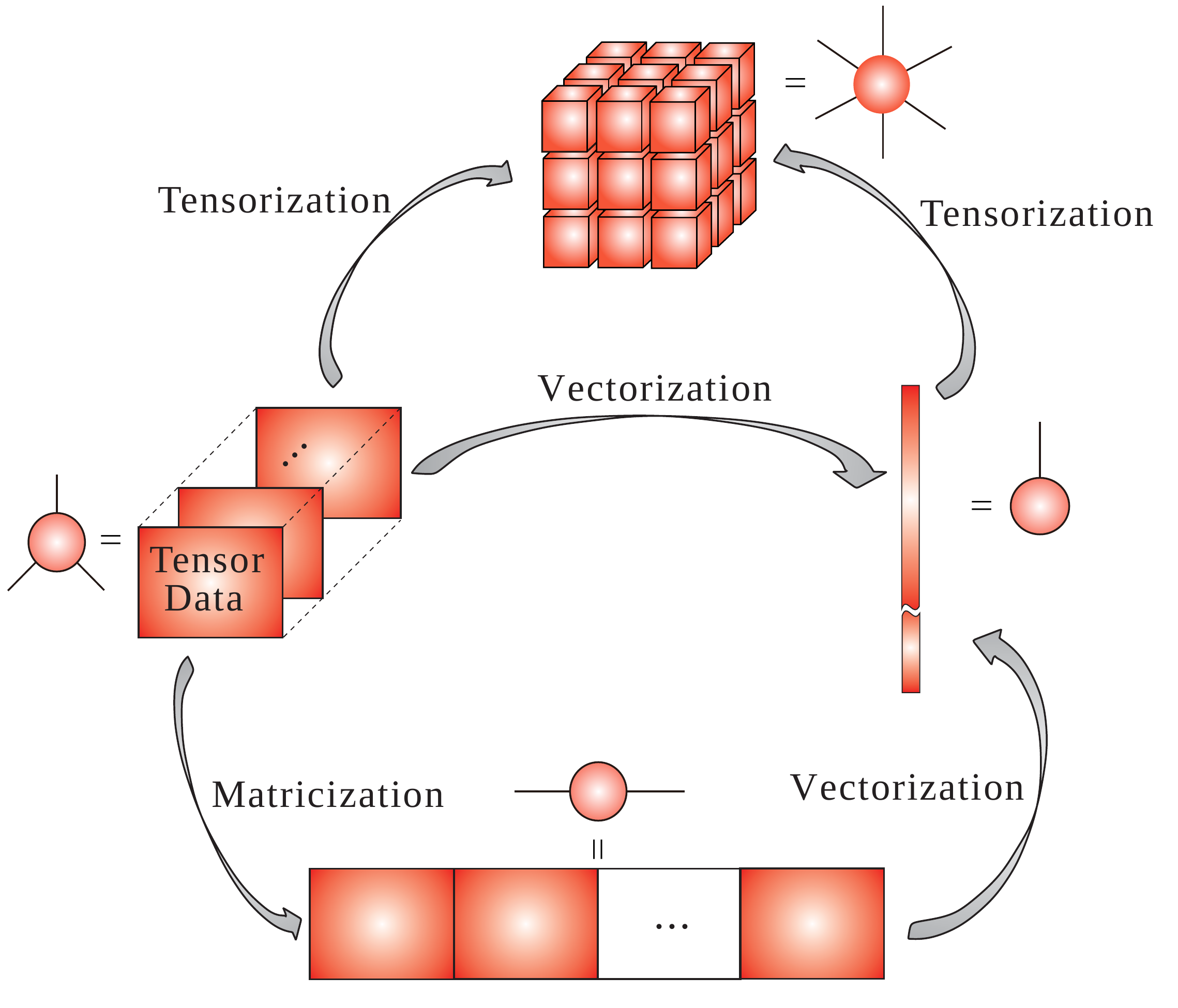}
\caption{Tensor  reshaping operations: Matricization, vectorization and tensorization.  Matricization refers to converting a tensor into a matrix, vectorization to converting a tensor or a matrix into a vector, while tensorization refers to converting   a vector, a matrix   or a low-order tensor
into a higher-order tensor.}
\label{Fig:Multiway2}
\end{figure}

\section{Basic Multilinear  Operations}

\begin{figure}
(a)
\vspace{-0.01cm}
\begin{center}
\includegraphics[width=10.35cm]{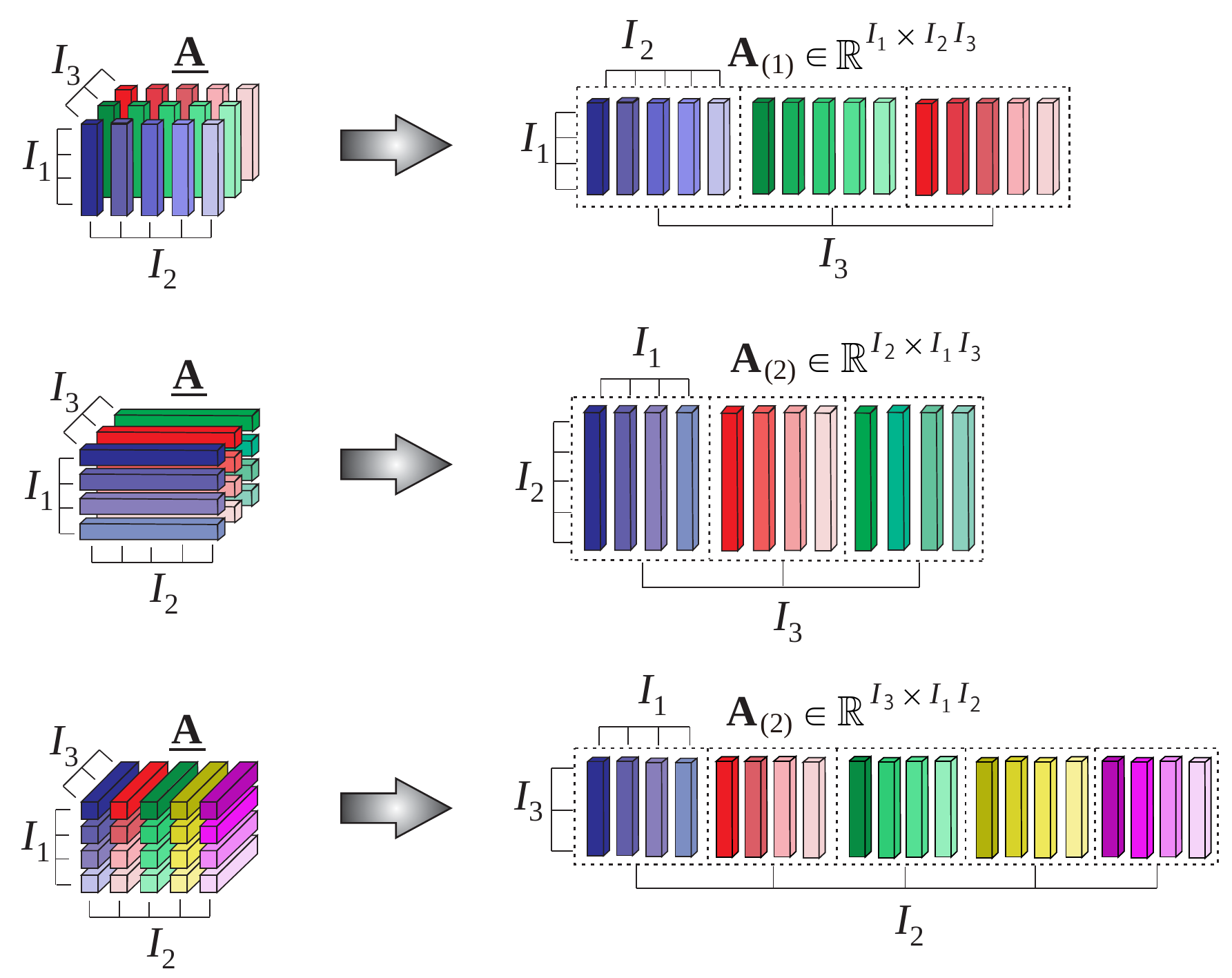}
\end{center}
(b)
\vspace{-0.5cm}
\begin{center}
\includegraphics[width=7.74cm]{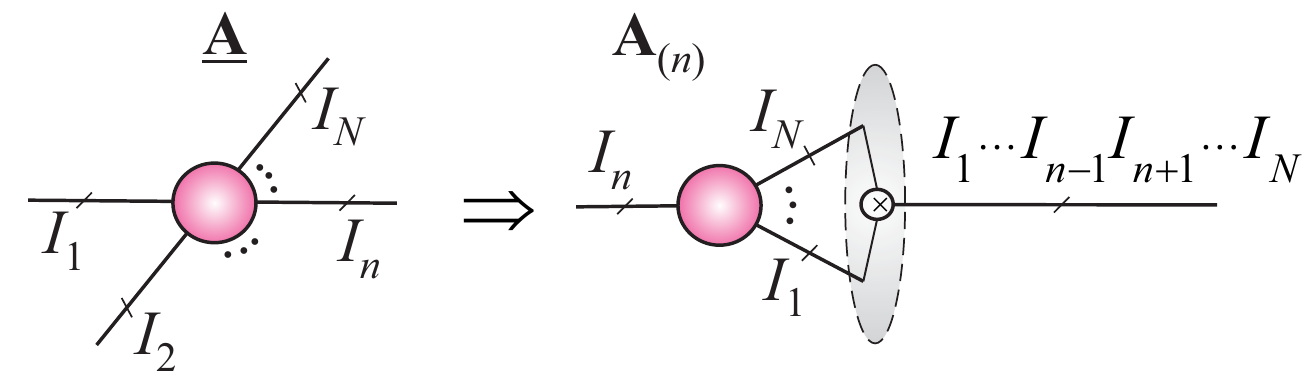}
\end{center}
(c)
\vspace{-0.5cm}
\begin{center}
\includegraphics[width=6.84cm]{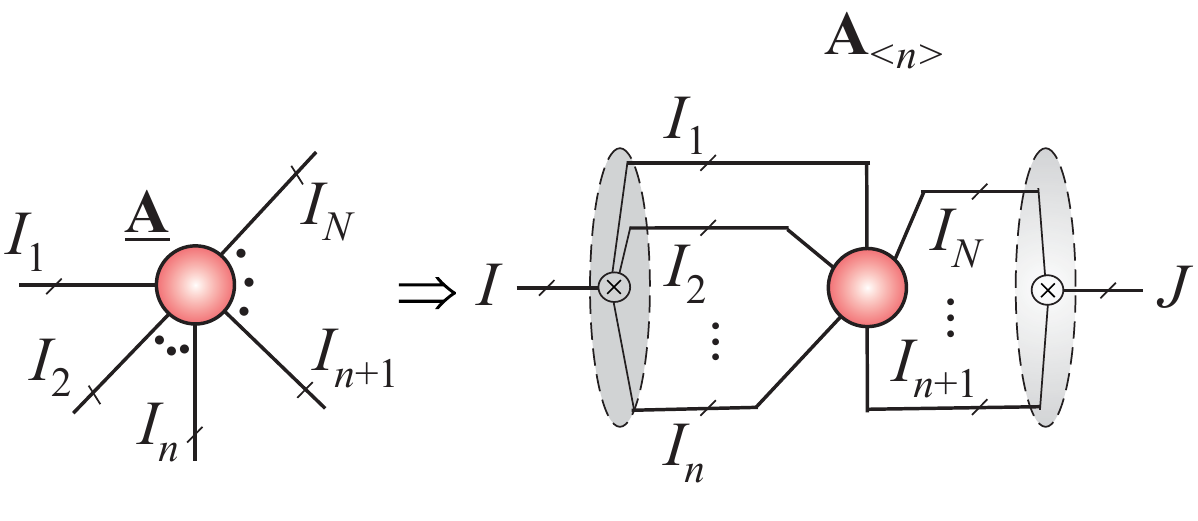}
\end{center}
\caption{{\small Matricization (flattening, unfolding) used in tensor reshaping. (a) Mode-1, mode-2, and mode-3 matricizations of a 3rd-order tensor, from the top to the bottom panel. (b) Tensor network diagram for the  mode-$n$ matricization of an $N$th-order tensor,  $\underline \bA \in \Real^{I_1  \times  I_2 \times \cdots \times I_N}$, into a short and wide matrix,   $\bA_{(n)} \in \Real^{I_n  \;\times \; I_1\cdots I_{n-1} I_{n+1} \cdots I_N}$. (c)  Mode-$\{1,2,\ldots,n\}$th (canonical) matricization of an $N$th-order tensor, $\underline \bA$, into a matrix $\bA_{<n>} = \bA_{(\overline{i_1 \cdots i_n}\; ; \; \overline{i_{n+1} \cdots  i_N})} \in \Real^{I_1 I_2 \cdots I_n \; \times  \;I_{n+1} \cdots I_N}$.}}
\label{Fig:unfolding}
\end{figure}

\begin{figure}[t]
\centering
\includegraphics[width=10.8cm,trim = 0 0.25cm 0 0.4cm, clip  = true]{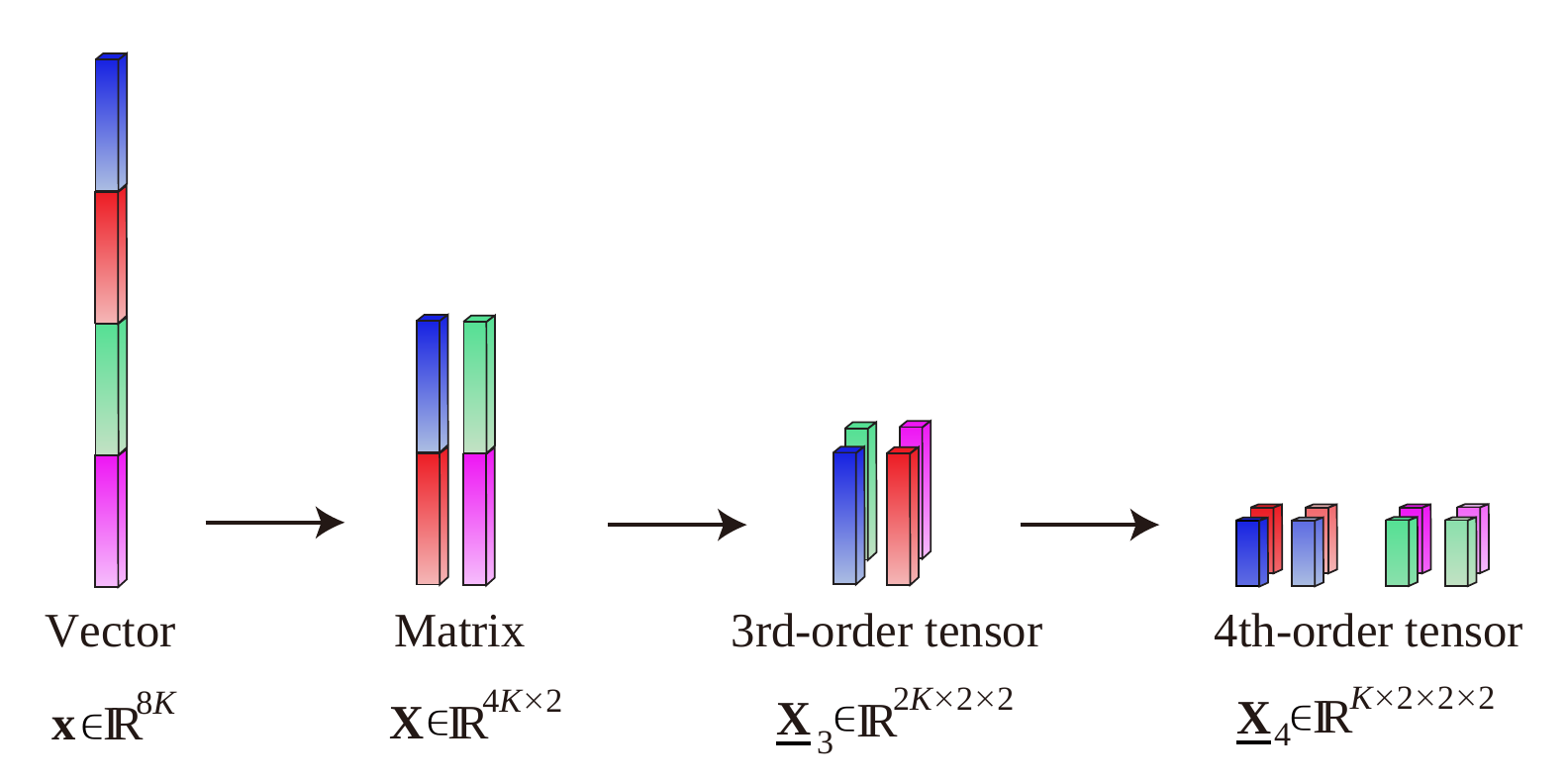}
\caption{Tensorization of a vector into a matrix, 3rd-order tensor and 4th-order tensor.}
\label{Fig:Tensorization1}
\end{figure}

\begin{figure}[t]
\centering
\includegraphics[width=9.99cm,trim = 0 0.2cm 0 .2cm, clip = true]{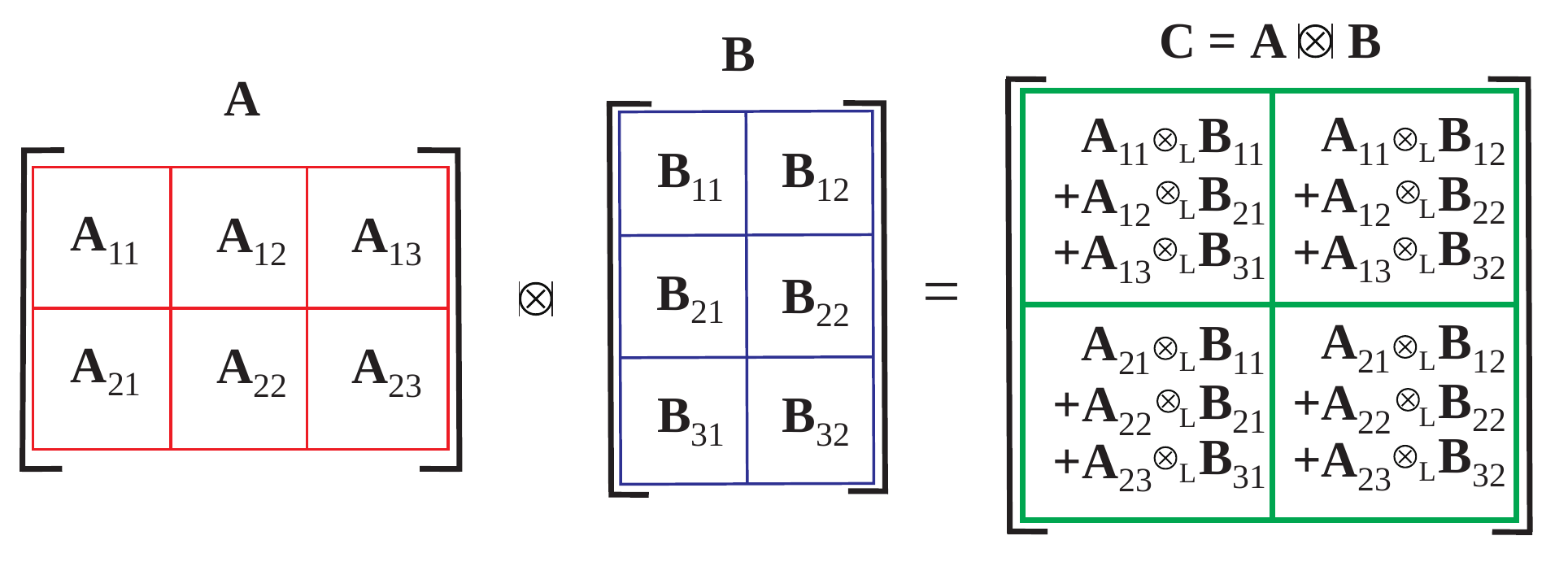}
\caption{Illustration of the strong Kronecker product of two block matrices, $\bA =[\bA_{r_1,r_2}] \in \Real^{R_1 I_1 \times R_2 J_1}$ and $\bB = [\bB_{r_2,r_3}] \in \Real^{R_2 I_2 \times R_3 J_2}$, which is defined as a block matrix $\bC =\bA \skron \bB \in \Real^{R_1 I_1 I_2 \times R_3 J_1 J_2}$, with the blocks $\bC_{r_1,r_3} =\sum_{r_2=1}^{R_2} \bA_{r_1,r_2} \otimes_L \bB_{r_2,r_3} \in \Real^{ I_1 I_2 \times J_1 J_2}$, for $r_1=1,\ldots, R_1,\; r_2=1,\ldots,R_2$ and $r_3=1,\ldots,R_3$.
}
\label{Fig:strongKron}
\end{figure}

\begin{figure}[t]
\begin{center}
 \includegraphics[width=7.9cm,trim = 0 0.0cm 0 0.3cm, clip = true]{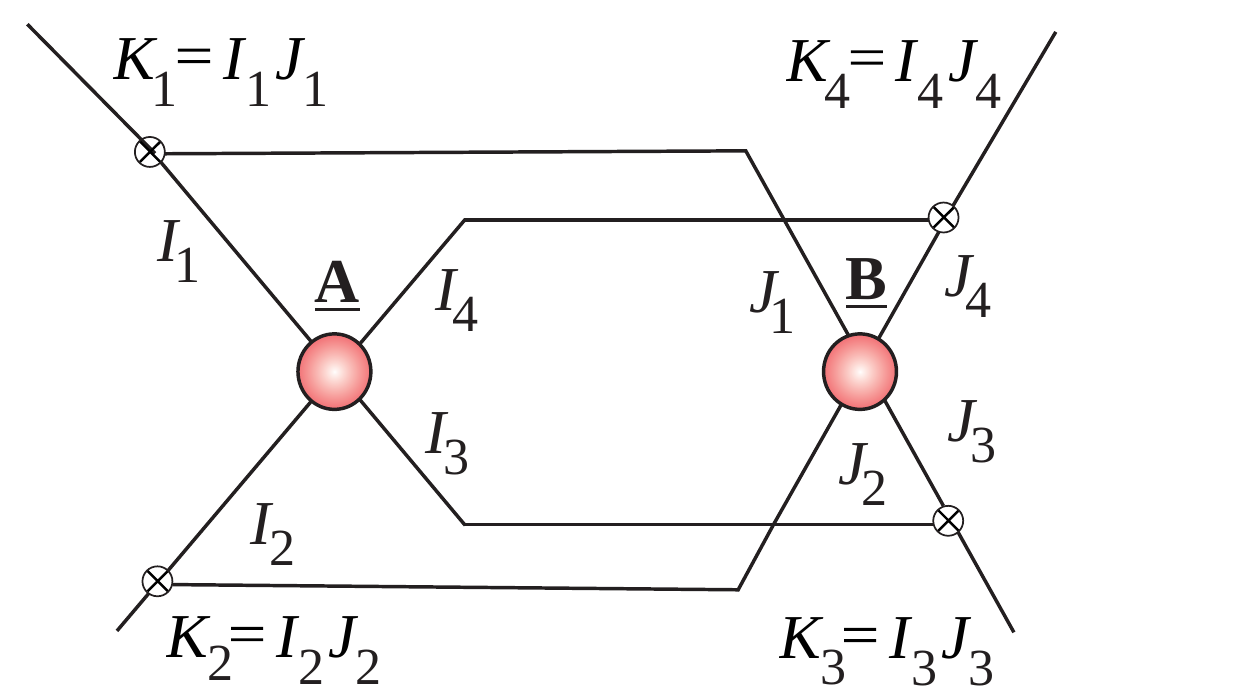}
 \end{center}
\caption{The left Kronecker product of two 4th-order tensors, $\underline \bA$ and $\underline \bB$, yields a 4th-order tensor,   ${\underline \bC} = {\underline \bA} \otimes_L {\underline \bB} \in \Real^{I_1 J_1 \times \cdots \times I_4 J_4}$, with entries $c_{k_1,k_2,k_3,k_4}  = a_{i_1, \ldots,i_4} \: b_{j_1, \ldots,j_4}$, where $k_n= \overline{i_n j_n} $ ($n=1,2,3,4$). Note that the order of tensor $\underline \bC$ is the same as the order of $\underline \bA$ and $\underline \bB$, but the size in every mode within $\underline \bC$ is a product of the respective sizes of $\underline \bA$ and $\underline \bB$.}
\label{Fig:KP}
\end{figure}

\begin{figure}
\begin{center}
 \includegraphics[width=11.7cm]{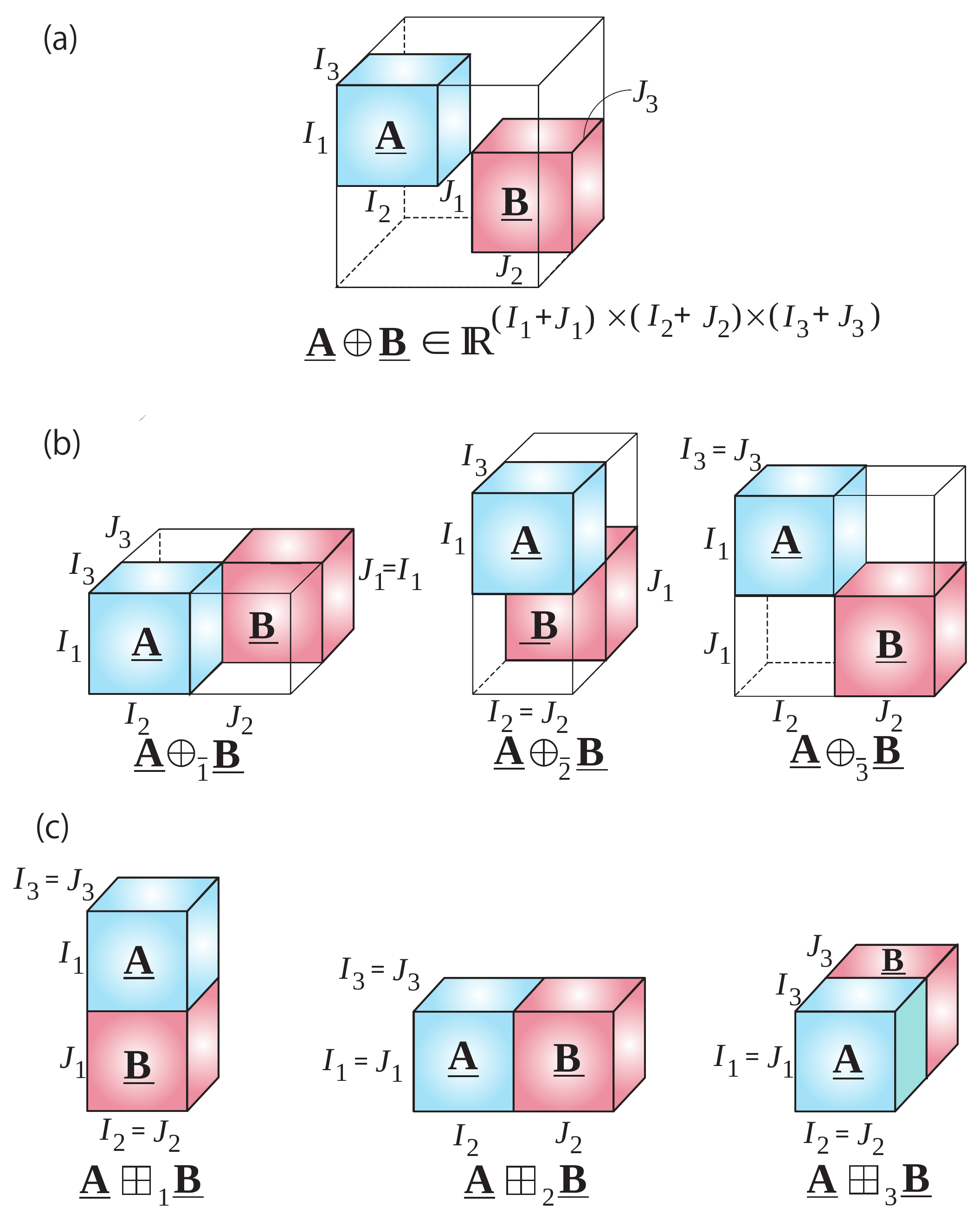}
 \end{center}
\caption{Illustration of the direct sum, partial direct sum and concatenation operators  of two 3rd-order tensors.
(a) Direct sum. (b) Partial (mode-1, mode-2, and mode-3) direct sum. (c) Concatenations along mode-1,2,3.}
\label{Fig:direct_sum}
\vspace{30pt}
\end{figure}

\begin{figure}
\begin{center}
 \includegraphics[width=9.9cm]{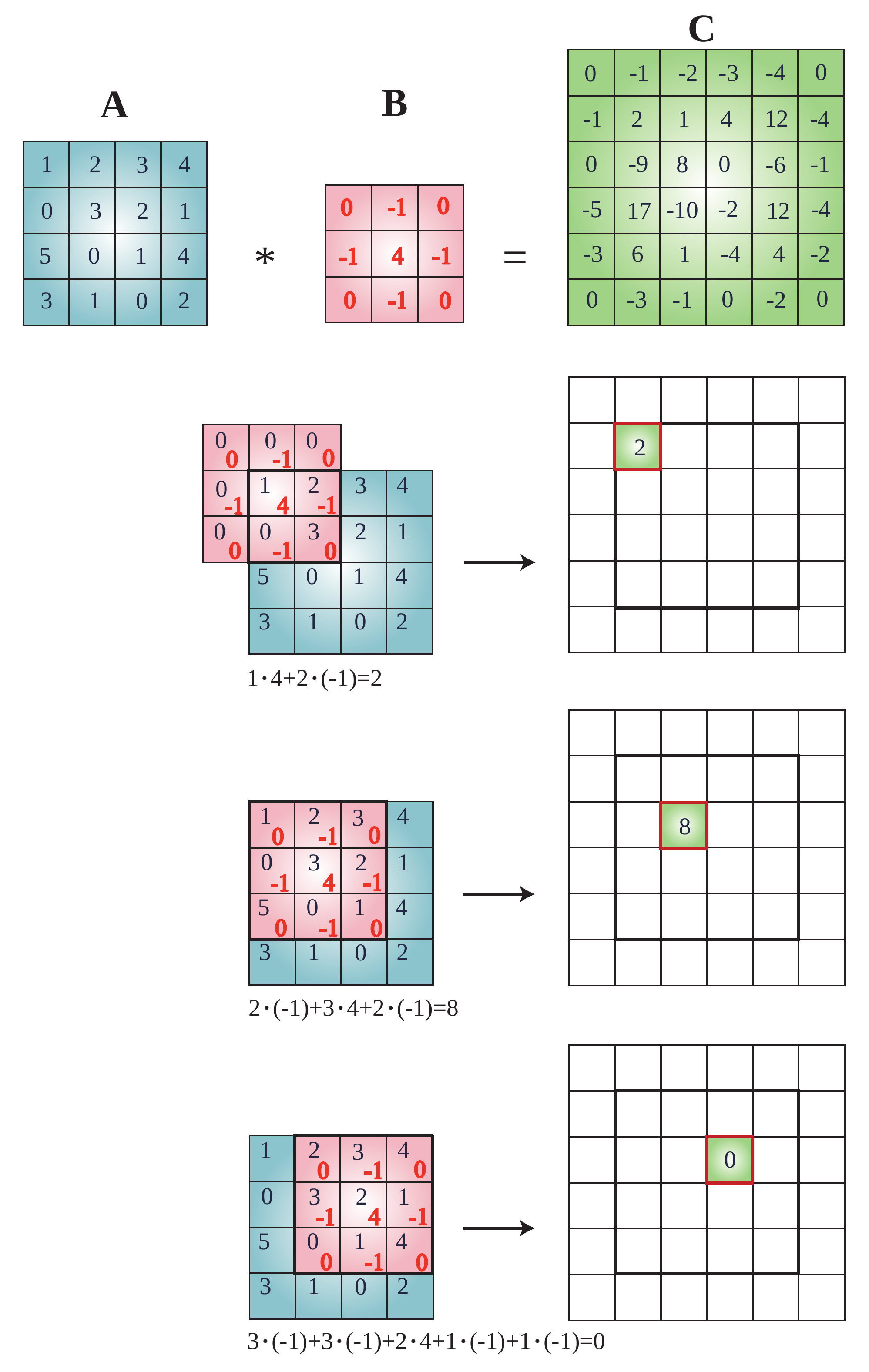}
 \end{center}
\caption{Illustration of the 2D convolution operator, performed through a sliding window
 operation  along both the horizontal and vertical index.}
\label{Fig:2Dconv}
\end{figure}

\begin{figure}[t]
\centering
 \includegraphics[width=11.9cm,trim = 0 0.2cm 0 0.4cm , clip = true]{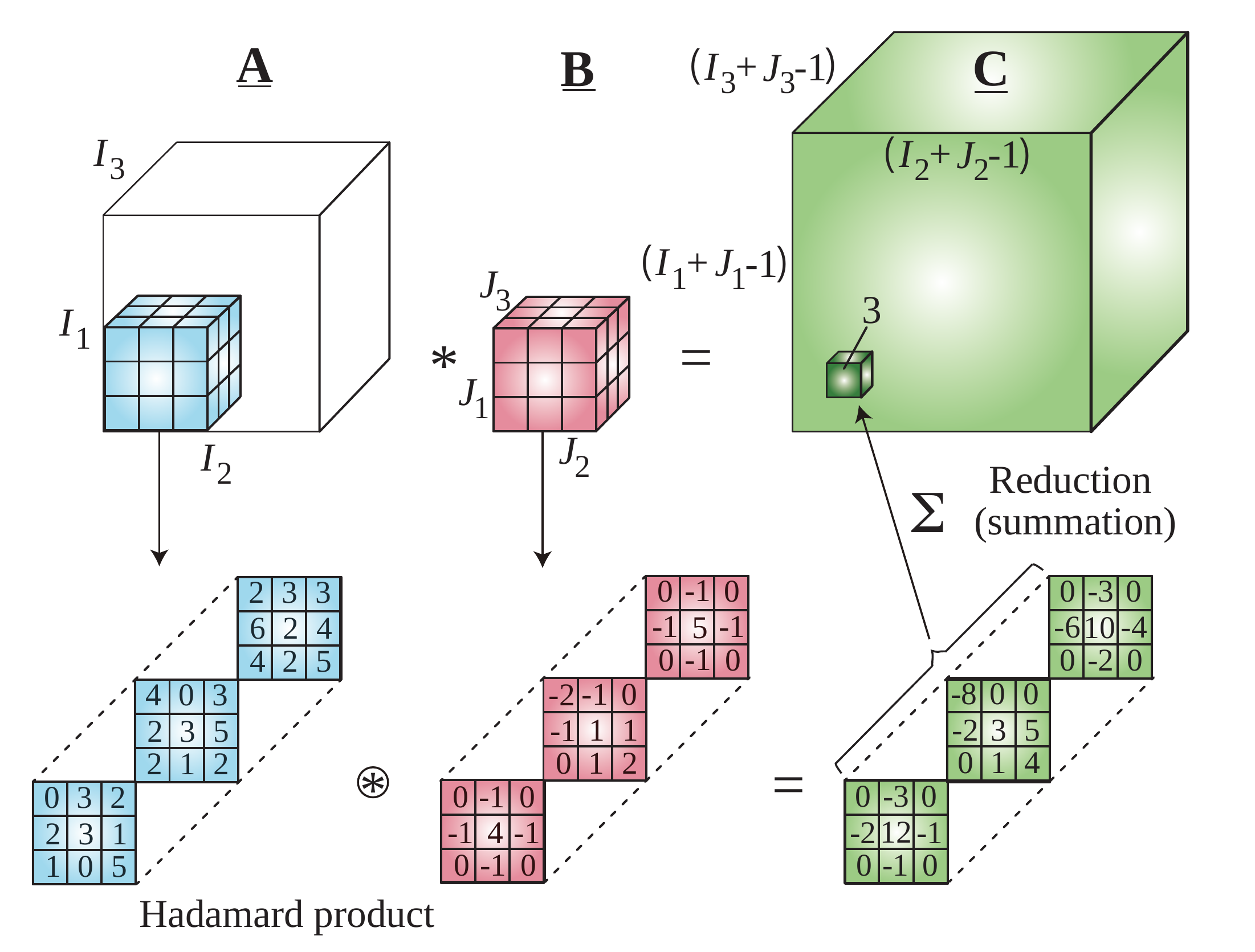}
\caption{Illustration of the 3D convolution operator, performed through a sliding window operation along all three indices.}
\label{Fig:3Dconv}
\end{figure}

\begin{figure}[t]
 \centering
 \includegraphics[width=8.1cm, trim = 0 0.5cm 0 0.1cm, clip = true]{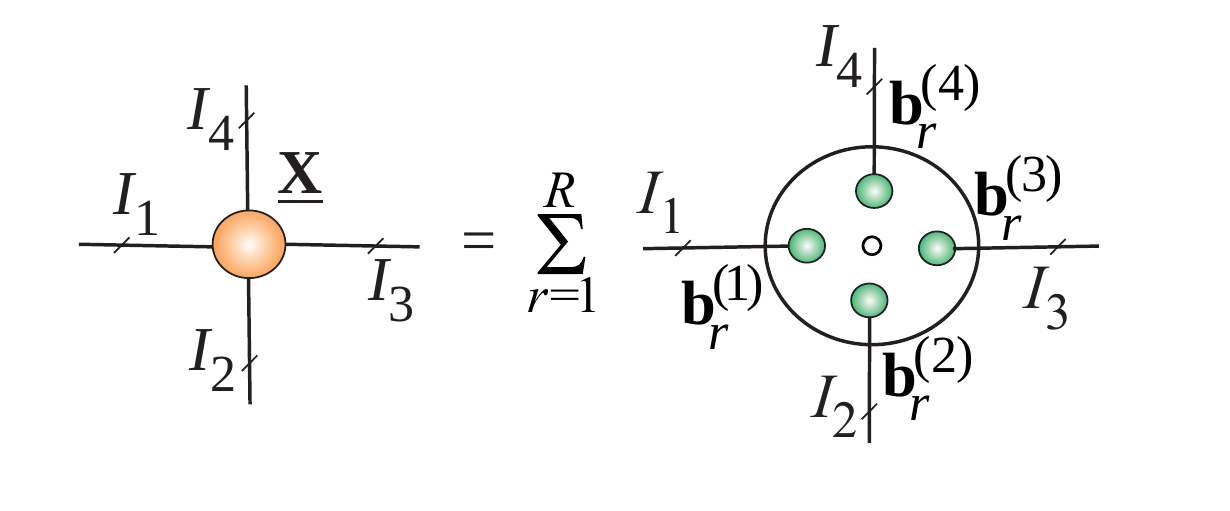}
\caption{The CP decomposition for a 4th-order tensor $\underline \bX$ of rank $R$. Observe that the rank-1 subtensors  are formed through the outer products of the vectors $\bb_r^{(1)}, \ldots,\bb_r^{(4)} $, $r = 1, \ldots, R$.}
\label{Fig:CPD4}
\end{figure}

\begin{figure}[t!]
(a)
\vspace{-0.1cm}
\begin{center}
\includegraphics[width=10.99cm, trim  = 0 0.5cm 0 0 , clip = true]{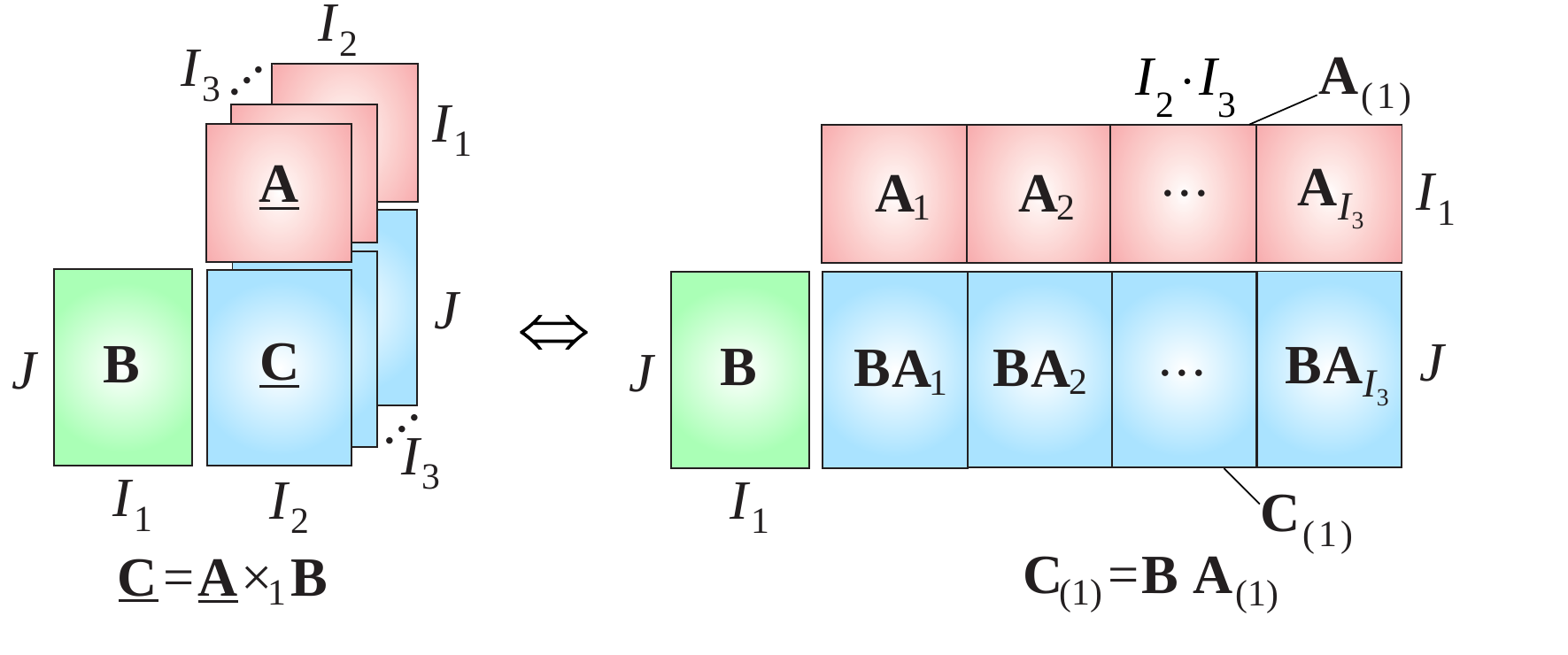}
\end{center}
(b)
\vspace{-0.1cm}
\begin{center}
\includegraphics[width=10.99cm]{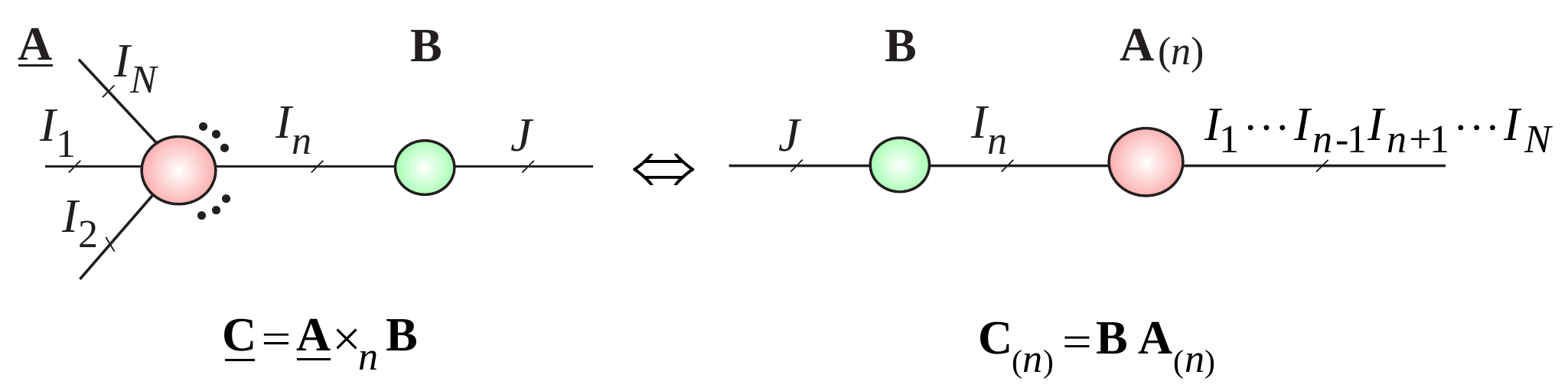}
\end{center}
\caption{Illustration of the multilinear mode-$n$ product, also known as   the TTM (Tensor-Times-Matrix) product, performed in the tensor format (left) and the matrix format (right). (a) Mode-$1$ product of a 3rd-order tensor,  $\underline \bA \in \Real^{I_1 \times I_2 \times I_3}$, and a factor (component) matrix, $\bB \in \Real^{J \times I_1}$,  yields a tensor $\underline \bC = \underline \bA \times_1 \bB \in \Real^{J \times I_2 \times I_3}$. This is equivalent to a simple matrix multiplication formula, $\bC_{(1)}=\bB \bA_{(1)}$. (b) Graphical representation of a mode-$n$ product of an $N$th-order tensor, $\underline \bA \in \Real^{I_1 \times I_2 \times \cdots \times I_N}$,  and a factor matrix, $\bB \in \Real^{J \times I_n}$.}
\label{Fig:AxnB}
\end{figure}

\begin{figure}
(a)
\begin{center}
\includegraphics[width=7.36cm]{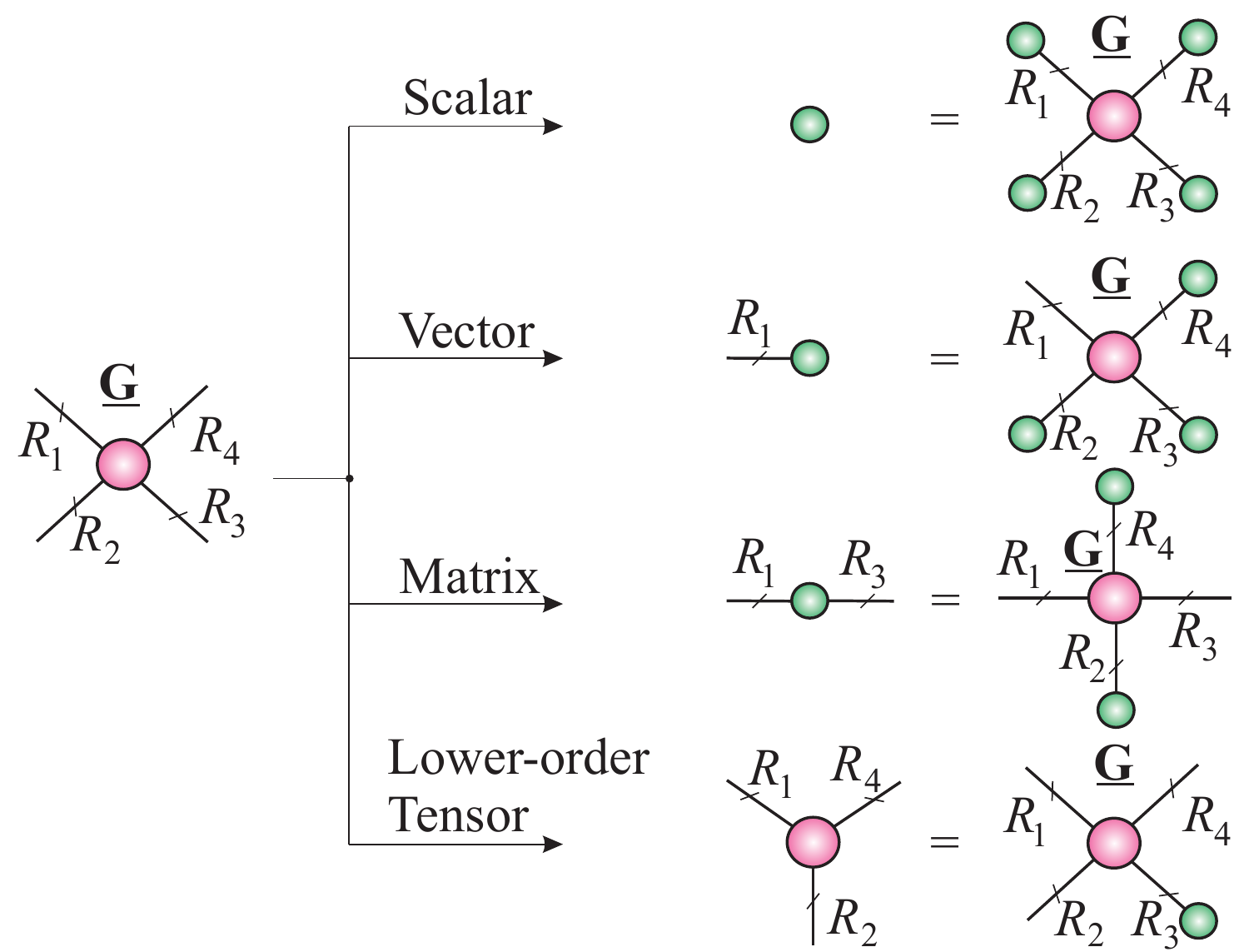}
\end{center}
(b) \hspace{5.9cm} (c)
\begin{center}
\includegraphics[width=9.2cm]{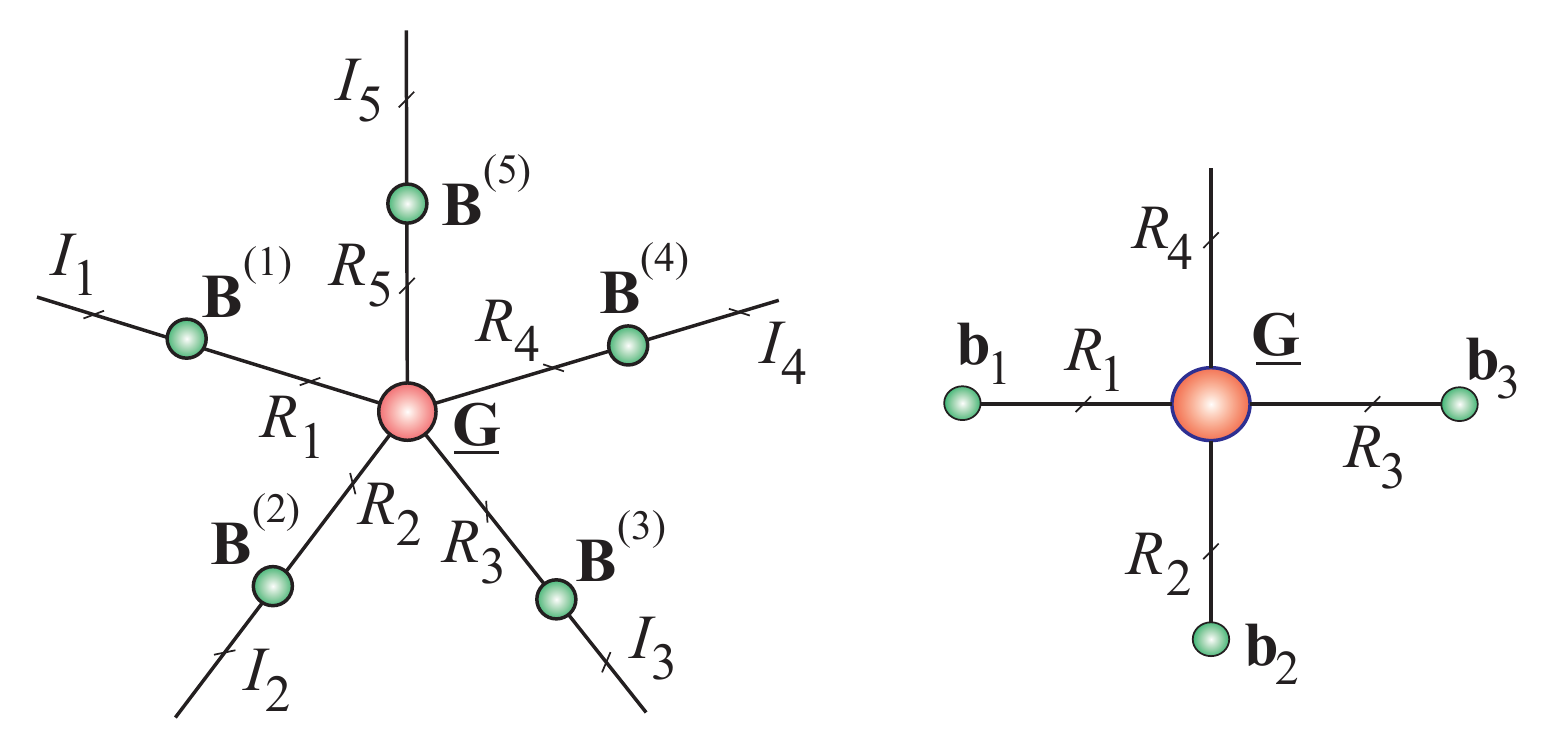}
\end{center}
\caption{Multilinear tensor products in a compact tensor network notation. (a) Transforming and/or compressing a 4th-order tensor, $\underline \bG \in \Real^{R_1 \times R_2 \times R_3 \times R_4}$, into a scalar, vector, matrix and 3rd-order tensor, by  multilinear products of the tensor and vectors. Note that a mode-$n$ multiplication of a tensor by a matrix does not change the order of a tensor, while a multiplication of a tensor by a vector reduces its order by one. For example, a multilinear product of a 4th-order tensor and four vectors (top diagram) yields a scalar. (b) Multilinear product of a tensor, $\underline \bG \in \Real^{R_1 \times R_2 \times \cdots \times R_5}$, and five factor (component) matrices, $\bB^{(n)} \in \Real^{I_n \times R_n}$ ($n=1,2,\ldots,5$), yields the tensor $\underline \bC = \underline \bG \times_1 \bB^{(1)}  \times_2 \bB^{(2)} \times_3\bB^{(3)}  \times_4 \bB^{(4)}   \times_5 \bB^{(5)} \in \Real^{I_1 \times I_2 \times \cdots \times I_5}$. This corresponds to the Tucker format.  (c) Multilinear product of a 4th-order tensor, $\underline \bG \in \Real^{R_1 \times R_2 \times R_3 \times R_4}$,  and three vectors, $\bb_n \in \Real^{R_n}$ $(n=1,2,3)$, yields  the vector $\bc= \underline \bG \bar \times_1 \bb_1  \bar \times_2 \bb_2 \bar \times_3  \bb_3 \in \Real^{R_4} $.}
\label{Fig:TO1}
\end{figure}

\begin{figure}[t]
 \centering
\includegraphics[width=10.9cm]{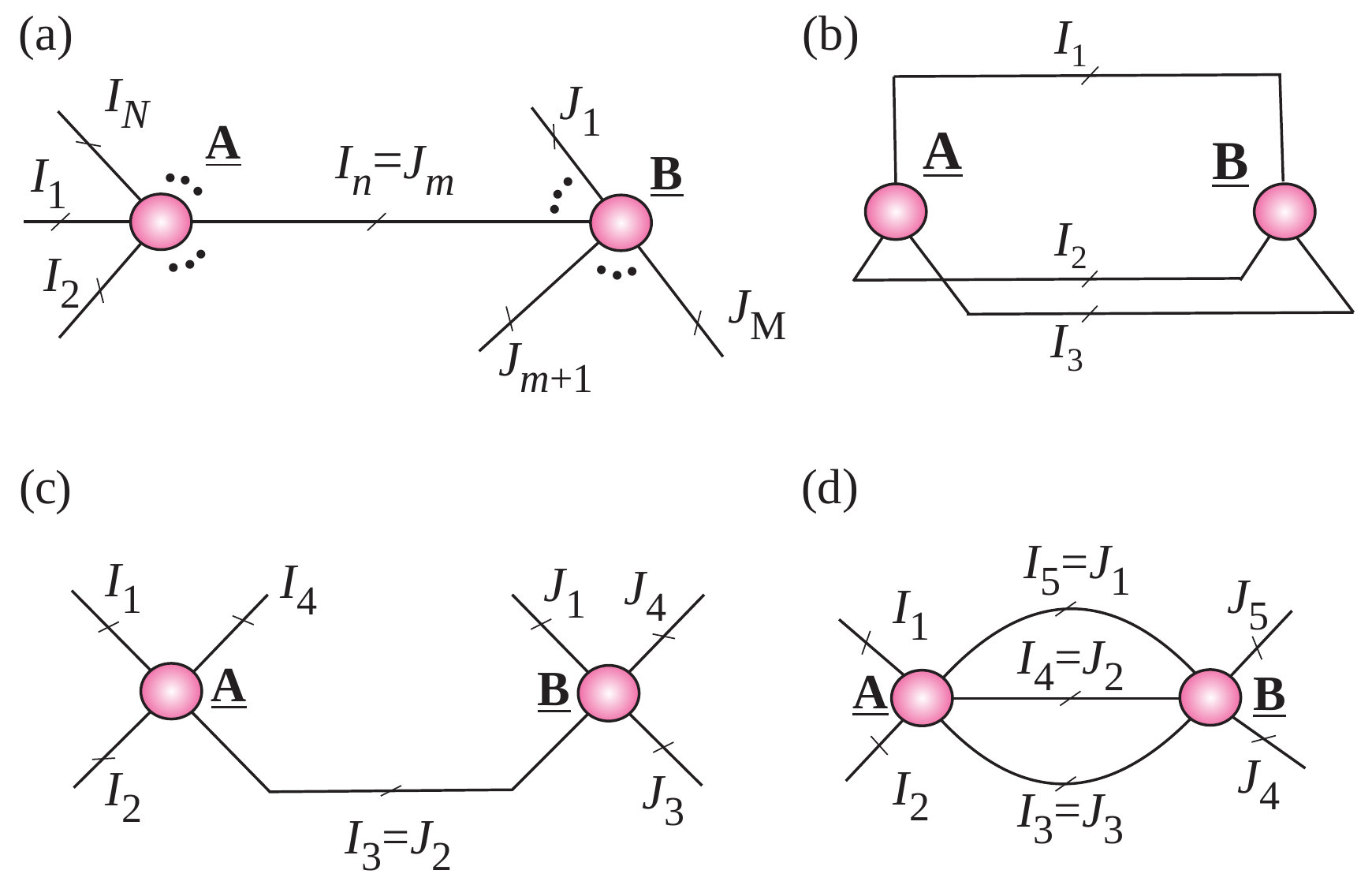}
\caption{Examples of  contractions of two tensors. (a) Multilinear product of two tensors is denoted by  $ \underline \bA \; \times_n^m \; \underline \bB$. (b) Inner product of two 3rd-order tensors yields a scalar $c=\langle\underline \bA, \underline \bB\rangle=\underline \bA \; \times_{1,2,3}^{1,2,3} \; \underline \bB  =\underline \bA \; \bar \times \; \underline \bB=\sum_{i_1,i_2,i_3} \; a_{i_1,i_2,i_3} \; b_{i_1,i_2,i_3}$.  (c) Tensor contraction of two 4th-order tensors, along mode-3 in $\underline \bA$ and mode-2 in $\underline \bB$, yields a 6th-order tensor, $\underline \bC = \underline \bA \; \times_{3}^{2} \; \underline \bB \in \Real^{I_1 \times I_2 \times I_4 \times J_1 \times J_3 \times J_4}$, with entries $c_{i_1,i_2,i_4,j_1,j_3,j_4}= \sum_{i_3} \; a_{i_1,i_2,i_3,i_4}  \; b_{j_1,i_3,j_3,j_4}$. (d) Tensor contraction of two 5th-order tensors along the modes $3,4,5$ in $\underline \bA$ and $1,2,3$ in $\underline \bB$ yields a 4th-order tensor, $\underline \bC = \underline \bA \; \times_{5, 4, 3}^{1, 2, 3} \; \underline \bB \in \Real^{I_1 \times I_2 \times J_4 \times J_5}$.}
\label{Fig:TO2}
\end{figure}

\begin{figure}[t]
\centering
 \includegraphics[width=5.99cm]{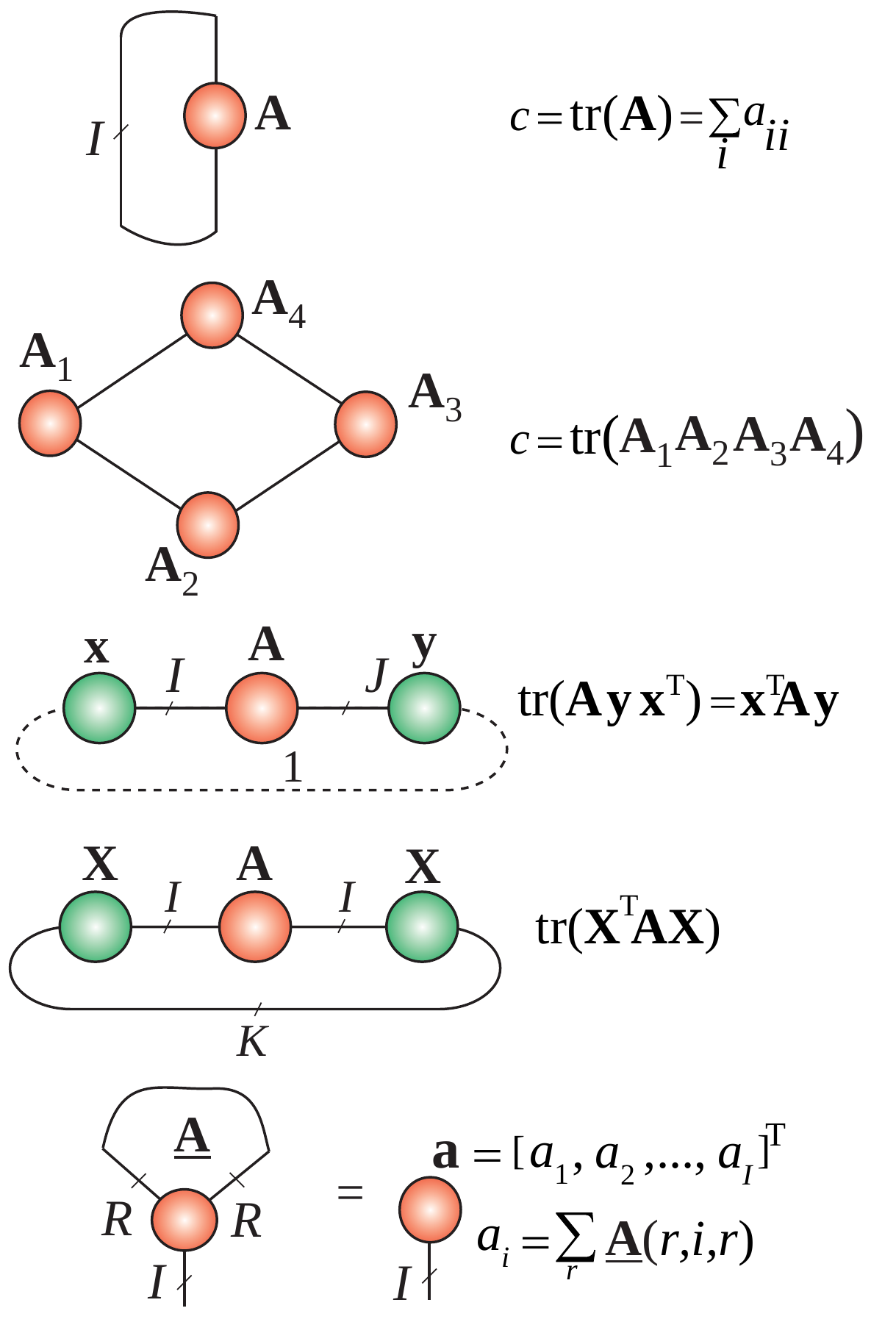}
\caption{Tensor network notation  for the traces of matrices (panels 1-4 from the top), and a (partial) tensor trace (tensor self-contraction) of a 3rd-order tensor (bottom panel). Note that graphical representations of the trace of matrices intuitively explain the permutation property of trace operator, e.g., $\tr(\bA_1 \bA_2 \bA_3 \bA_4) = \tr(\bA_3 \bA_4 \bA_1 \bA_2)$.}
\label{Fig:trace}
\end{figure}

The following symbols are used for most common tensor multiplications: $\otimes$ for the Kronecker product, $\odot$ for the Khatri--Rao product,  $\*$ for the Hadamard
(componentwise) product,  $\circ$ for the outer product  and  $ \times_n$ for the mode-$n$  product. Basic tensor  operations are summarized in Table \ref{table_notation2},
 and illustrated in Figures \ref{Fig:Multiway2}--\ref{Fig:trace}.
 We refer to \cite{Kolda08,NMF-book,Lee-TTfund1} for more detail regarding the basic
  notations and tensor operations.
For convenience, general operations, such as $\mbox{vec}(\cdot)$ or $\diag(\cdot)$, are defined similarly to the MATLAB syntax. \\

\noindent {\bf Multi--indices:} By a multi-index $i = \overline{i_1 i_2 \cdots  i_N}$ we refer to an index which takes all possible
combinations of values of indices, $i_1, i_2, \ldots, i_N$,  for $i_n = 1,2,\ldots, I_n$, $\;n=1,2,\ldots,N$ and in a specific order.
Multi--indices can be defined  using two  different conventions \cite{DolgovAMEN2014}:
\begin{enumerate}
\item Little-endian convention (reverse  lexicographic ordering)
\begin{equation}
\overline{i_1 i_2 \cdots i_N } = i_1 + (i_{2} - 1) I_1+ (i_{3} - 1) I_1 I_2 +
 \cdots   + (i_N - 1) I_1 \cdots  I_{N-1}. \notag
\end{equation}
\item Big-endian (colexicographic ordering)
\be
\overline{i_1 i_2 \cdots i_N } &=& i_N + (i_{N-1} - 1)I_N
+  (i_{N-2} - 1)I_N I_{N-1}+  \notag \\
&&   \cdots  + (i_1 - 1)I_2 \cdots  I_N. \notag
\ee
\end{enumerate}
The little-endian convention is used, for example, in Fortran and
MATLAB, while the big-endian convention is used in C language. Given the complex and non-commutative nature of tensors,  the basic definitions, such as the  matricization, vectorization and  the Kronecker  product,  should be consistent with the chosen convention{\footnote{
 Note that using the colexicographic ordering, the vectorization of an outer product of two vectors, $\ba$ and $\bb$, yields their Kronecker product, that is,   $\mbox{vec}(\ba \circ  \bb) = \ba \otimes \bb$,
 while using the  reverse lexicographic ordering, for the same operation, we need to use the Left Kronecker product,  $\mbox{vec}(\ba \circ  \bb) = \bb \otimes \ba=\ba \otimes_L \bb $.}}.
In this monograph, unless otherwise stated, we will use  little-endian notation.\\

\noindent {\bf Matricization.} The matricization operator, also known as the unfolding or flattening, reorders the elements of a tensor into a
matrix  (see Figure \ref{Fig:unfolding}). Such a matrix is re-indexed according to the choice of multi-index described  above, and the following two fundamental matricizations are used extensively.\\

\noindent {\bf The mode-$n$ matricization.} For a fixed index $n \in \{1,2,\ldots,N\}$, the mode-$n$ matricization of an $N$th-order tensor, $\underline \bX \in \Real^{I_1\times \cdots \times I_N}$, is defined as the (``short'' and ``wide'') matrix 
	\begin{equation} \label{expr_mode_matricization}
	\bX_{(n)}
	\in \Real^{I_n\times I_1 I_2 \cdots I_{n-1} I_{n+1} \cdots I_N},
	\end{equation}
with $I_n$ rows and $I_1 I_2 \cdots I_{n-1} I_{n+1} \cdots I_N$ columns, the entries of which are
	$$
	(\bX_{(n)})_{i_n, \overline{i_1 \cdots i_{n-1} i_{n+1} \cdots i_N}}
	= x_{i_1,i_2,\ldots,i_N}.
	$$
Note that the columns of a mode-$n$ matricization, $\bX_{(n)}$, of a tensor $\underline \bX$ are the mode-$n$ fibers of $\underline \bX$.\\

\noindent {\bf The mode-$\{n\}$ canonical matricization.} For a fixed index $n \in \{1,2,\ldots,N\}$, the mode-$(1,2,\ldots, n)$ matricization, or  simply mode-$n$ canonical matricization, of a tensor $\underline \bX \in \Real^{I_1\times \cdots \times I_N}$ is defined as the matrix 
	\begin{equation}
	\bX_{<n>}
	\in \Real^{I_{1} I_{2} \cdots I_{n} \times
			   I_{n+1} \cdots I_{N}},
	\end{equation}
with  $I_{1} I_{2} \cdots I_{n}$ rows and  $I_{n+1} \cdots I_{N}$ columns, and the entries
	$$
	(\bX_{<n>})_{\overline{i_1 i_2 \cdots i_n}, \, \overline{i_{n+1}  \cdots i_N}} = x_{i_1,i_2,\ldots,i_N}.
	$$
The matricization operator in the  MATLAB notation (reverse lexicographic) is given by
\begin{equation}
\bX_{<n>} =\mbox{reshape}\left( \underline \bX, I_1I_2 \cdots I_n,   I_{n+1} \cdots I_N \right).
\end{equation}
 As special cases we immediately have (see Figure \ref{Fig:unfolding})
\begin{equation}
\bX_{<1>} =\bX_{(1)}, \quad \bX_{<N-1>} =\bX^{\text{T}}_{(N)}, \quad \bX_{<N>} =\text{vec} (\bX).
\end{equation}

The tensorization of a vector or a matrix can be considered as a reverse process  to the vectorization or matricization
(see  Figures ~\ref{Fig:Multiway2} and~\ref{Fig:Tensorization1}).\\

\noindent {\bf Kronecker, strong Kronecker, and  Khatri--Rao products of matrices and  tensors.}
For an $I \times J$ matrix $\bA$ and a $K \times L$ matrix $\bB$, the standard (Right) Kronecker product, $\bA \otimes \bB$, and the Left Kronecker product, $\bA \otimes_L \bB$, are the following $IK \times JL$ matrices
\begin{gather}
\bA \otimes \bB  =  \left[ \begin{matrix}
   a_{1,1} \bB & \cdots & a_{1,J} \bB \\
    \vdots & \ddots & \vdots \\
   a_{I,1}\bB & \cdots & a_{I,J} \bB
    \end{matrix}  \right],
    \quad
    \bA \otimes_L \bB
    =  \left[ \begin{matrix}
    \bA b_{1,1} &\cdots&  \bA b_{1,L}\\
    \vdots&\ddots&\vdots\\
   \bA b_{K,1} & \cdots &  \bA b_{K,L}
    \end{matrix}  \right].\notag  
  \end{gather}
  Observe that $\bA \otimes_L \bB = \bB \otimes \bA$, so that the Left Kronecker product will be used in most cases in this monograph as it is consistent with the little-endian notation.

\noindent Using Left Kronecker product,  the strong Kronecker product of two block matrices, $\bA \in \Real^{R_{1} I \times R_2 J}$ and
$\bB  \in \Real^{R_2 K \times R_{3} L}$, given by
\be
{\bA} = \begin{bmatrix}\bA_{1,1}&\cdots&  \bA_{1,R_2}\\
				\vdots & \ddots & \vdots \\
				\bA_{R_{1},1}&\cdots&  \bA_{R_{1},R_2}
		\end{bmatrix}, \notag
\qquad
{\bB} = \begin{bmatrix} \ \bB_{1,1}&\cdots& \bB_{1,R_{3}}\\
				\vdots & \ddots & \vdots \\
				 \bB_{R_{2},1}&\cdots& \bB_{R_{2},R_{3}}
		\end{bmatrix}, \notag
 \ee
 can be defined as a block matrix (see Figure \ref{Fig:strongKron} for a graphical illustration)
\be
\bC =  {\bA} \; \skron \;   {\bB} \in \Real^{R_1 I K \times R_3 J L},
\ee
with blocks $\bC_{r_{1},r_{3}}=\sum_{r_2=1}^{R_2}  \bA_{r_{1},r_2} \otimes_L   \bB_{r_2,r_{3}} \in \Real^{I K \times J L}$,    where $\bA_{r_{1},r_2} \in \Real^{I \times J}$ and  $\bB_{r_{2},r_{3}} \in \Real^{K \times L}$ are the blocks of matrices within ${\bA}$ and $ {\bB}$, respectively  \cite{Seberry94,Kazeev_Toeplitz13,Kazeev2013LRT}.
Note that the strong Kronecker product is  similar to the standard block matrix multiplication, but performed using  Kronecker products of the blocks instead of the standard matrix-matrix products. The above definitions of Kronecker products  can be naturally extended to tensors \cite{Phan2012-Kron} (see Table~\ref{table_notation2}), as shown below.\\

\noindent{\bf The Kronecker product of tensors.} The (Left) Kronecker product of two $N$th-order tensors, $\underline \bA \in \Real^{I_1 \times I_2 \times \cdots \times I_N}$ and $ \underline  \bB \in \Real^{J_1 \times J_2 \times \cdots \times J_N}$, yields a tensor $\underline  \bC = \underline \bA \otimes_L \underline \bB \in \Real^{I_1 J_1 \times \cdots \times I_N J_N}$ of the same order but enlarged in size, with entries $c_{\overline{i_1  j_1},\ldots, \overline{i_N  j_N}} = a_{i_1, \ldots, i_N} \: b_{j_1, \ldots,j_N}$  as illustrated in Figure \ref{Fig:KP}.\\

\noindent{\bf The mode-$n$  Khatri--Rao product of tensors.} The Mode-$n$ Khatri--Rao product of two $N$th-order tensors, $\underline \bA \in \Real^{I_1 \times I_2 \times \cdots \times I_n \times  \cdots \times I_N}$ and $ \underline  \bB \in \Real^{J_1 \times J_2 \times \cdots \times J_n \times \cdots \times J_N}$, for which $I_n=J_n$,  yields a tensor $\underline  \bC = \underline \bA \odot_{\;n} \underline \bB \in \Real^{I_1 J_1 \times \cdots \times I_{n-1}J_{n-1} \times I_n \times I_{n+1}J_{n+1} \times \cdots \times I_N J_N}$, with subtensors $\underline \bC(:,\ldots:,i_n,:,\ldots,:) = \underline \bA(:,\ldots:,i_n,:,\ldots,:) \otimes \underline \bB(:,\ldots:,i_n,:,\ldots,:)$.\\

\noindent {\bf The mode-$2$ and mode-1  Khatri--Rao product of matrices.} The above definition simplifies to the standard Khatri--Rao (mode-2) product of two matrices, $\bA= [\ba_1, \ba_2, \ldots, \ba_R] \in \Real^{I \times R}$ and
$\bB =[\bb_1, \bb_2, \ldots, \bb_R] \in \Real^{J \times R}$,
or in other words a ``column-wise Kronecker product''. Therefore,   the standard Right and Left Khatri--Rao products for matrices are  respectively given by{\footnote{For simplicity, the mode $2$ subindex is usually neglected, i.e., $\bA \odot_2 \bB=\bA \odot \bB$.}}
\be
\bA \odot \bB &=& [\ba_1 \otimes \bb_1, \ba_2 \otimes \bb_2, \ldots, \ba_R  \otimes \bb_R] \in \Real^{IJ \times R}, \\
\bA \odot_L \bB &=& [\ba_1 \otimes_L \bb_1, \ba_2 \otimes_L \bb_2, \ldots, \ba_R  \otimes_L \bb_R] \in \Real^{IJ \times R}.
\ee

Analogously, the mode-1 Khatri--Rao product of two matrices  $\bA\in \Real^{I \times R}$ and
$\bB\in \Real^{I \times Q}$, is defined as
\be
\bA \odot_1 \bB  = \left[\begin{matrix} \bA(1,:) \otimes \bB(1,:) \\ \vdots \\ \bA(I,:) \otimes \bB(I,:)  \end{matrix} \right] \in \Real^{I \times R Q}.
\ee

\noindent {\bf Direct sum of tensors.}
A direct sum of $N$th-order tensors $\underline \bA \in \Real^{I_1\times\cdots\times I_N}$ and
	$\underline \bB \in \Real^{J_1\times\cdots\times J_N}$ yields
	a tensor $\underline \bC =\underline \bA \oplus \underline \bB \in \Real^{(I_1+J_1) \times \cdots \times (I_N+J_N)}$, with entries
	$\underline \bC(k_1,\ldots,k_N) = \underline \bA(k_1,\ldots,k_N)$ if $1\leq k_n \leq I_n$, \, $\forall n$,
	$\underline \bC(k_1,\ldots,k_N) = \underline \bB (k_1-I_1,\ldots,k_N-I_N)$ if $I_n< k_n \leq I_n+J_n$, \, $\forall n$,
	and $\underline \bC(k_1,\ldots,k_N) = 0$, otherwise (see Figure~\ref{Fig:direct_sum}(a)). \\

\noindent {\bf Partial (mode-$n$) direct sum of tensors.}
A partial direct sum of  tensors
	   $\underline \bA \in \Real^{ I_1\times\cdots\times I_N}$ and
	$\underline \bB \in \Real^{J_1\times\cdots\times J_N}$, with $I_n=J_n$,
	yields a tensor $\underline \bC=\underline \bA \oplus_{\;\overline n} \underline \bB \in \Real^{(I_1+J_1)\times \cdots \times (I_{n-1}+J_{n-1})\times I_n \times (I_{n+1}+J_{n+1})\times \cdots \times (I_N+J_N)}$, where
	$\underline \bC(:,\ldots,:,i_n,:, \ldots, :) = \underline \bA(:,\ldots,:,i_n,:, \ldots, :)\oplus \underline \bB(:,\ldots,:,i_n,:, \ldots, :)$, as illustrated in Figure~\ref{Fig:direct_sum}(b). \\

\noindent {\bf Concatenation of $N$th-order  tensors.} A concatenation along mode-$n$ of  tensors
	   $\underline \bA \in \Real^{I_1\times\cdots\times I_N}$ and
	$\underline \bB \in \Real^{J_1\times\cdots\times J_N}$, for which $I_m=J_m$, $\;\forall m \neq n$
	yields a tensor $\underline \bC=\underline \bA \boxplus_{\; n} \underline \bB \in \Real^{I_1 \times \cdots \times I_{n-1} \times (I_n+J_n) \times I_{n+1}\times \cdots \times (I_N)}$, with subtensors
	$\underline \bC(i_1,\ldots,i_{n-1},:,i_{n+1}, \ldots, i_N) = \underline \bA(i_1,\ldots,i_{n-1},:,i_{n+1}, \ldots, i_N)\oplus \underline \bB(i_1,\ldots,i_{n-1},:,i_{n+1}, \ldots, i_N)$, as illustrated in Figure \ref{Fig:direct_sum}(c).  For a concatenation of two tensors of suitable dimensions along mode-$n$, we will use equivalent notations
$\underline \bC=\underline \bA \boxplus_{\; n} \underline \bB  =\underline \bA \frown_{\; n} \underline \bB$. \\


\noindent {\bf 3D Convolution.} For simplicity, consider two 3rd-order tensors  $\underline \bA \in \Real^{I_1\times I_2 \times I_3}$ and
	$\underline \bB \in \Real^{J_1\times J_2 \times J_3}$. Their 3D Convolution
	yields a tensor $\underline \bC=\underline \bA \ast  \underline \bB \in \Real^{(I_1 +J_1-1) \times (I_2+J_2-1) \times  (I_3 +J_3-1)}$,
with entries: \\
	$\underline \bC(k_1, k_2, k_{3}) = \sum_{j_1} \sum_{j_2} \sum_{j_3} \; \underline \bB(j_1, j_2, j_{3}) \; \underline \bA(k_1-j_1, k_2 -j_2, k_3-j_{3})$
as illustrated in Figure \ref{Fig:2Dconv} and Figure \ref{Fig:3Dconv}. \\

\noindent {\bf Partial (mode-$n$) Convolution.} For simplicity,  consider two 3rd-order tensors $\underline \bA \in \Real^{I_1\times I_2 \times I_3}$ and
	$\underline \bB \in \Real^{J_1\times J_2 \times J_3}$. Their mode-$2$ (partial) convolution
	yields a tensor $\underline \bC=\underline \bA \boxdot_2  \underline \bB \in \Real^{I_1 J_1 \times (I_2+J_2-1) \times  I_3 J_3}$,
 the subtensors (vectors) of which are
	$\underline \bC(k_1,:,k_{3}) = \underline \bA(i_1,:,i_{3}) \ast \underline \bB(j_1,:,j_{3}) \in \Real^{I_2+J_2-1}$, where $k_1=\overline{i_1 j_1}$,  and $k_3=\overline{i_3  j_3}$.\\

\noindent {\bf Outer product.}  The central operator in tensor analysis is the outer or tensor product, which for the
tensors $\underline \bA \in \Real^{I_1 \times \cdots \times I_N}$ and $\underline \bB \in \Real^{J_1 \times \cdots \times J_M}$ gives
the tensor $\underline \bC =\underline \bA \circ \underline \bB \in \Real^{I_1 \times \cdots \times I_N \times J_1 \times \cdots \times J_M }$ with entries
$c_{i_1,\ldots,i_N,j_1,\ldots,j_M} = a_{i_1,\ldots,i_N} \; b_{j_1,\ldots,j_M}$.

Note that for 1st-order tensors (vectors), the tensor product reduces to the standard  outer product of two nonzero vectors,
$\ba \in \Real^I$ and $\bb \in \Real^J$, which yields  a rank-1 matrix,
$\bX=\ba \circ \bb = \ba \bb^{\text{T}} \in \Real^{I \times J}$.
 The outer product of  three nonzero vectors, $\ba \in \Real^I, \; \bb \in \Real^J$ and  $\bc \in \Real^K$, gives a 3rd-order rank-1  tensor
(called  pure or elementary tensor),
$\underline \bX= \ba \circ \bb \circ \bc \in \Real^{I \times J  \times K}$,
with entries $x_{ijk} =a_i \; b_j \; c_k$.

\noindent{\bf Rank-1 tensor.} A tensor, $\underline \bX \in \Real^{I_1 \times I_2 \times \cdots \times I_N}$, is said to be of rank-1 if  it can be expressed exactly
as the outer product, $\underline \bX = \bb^{(1)} \circ \bb^{(2)} \circ \cdots  \circ \bb^{(N)}$ of nonzero vectors,  $\bb^{(n)} \in \Real^{I_n}$, with the tensor entries  given by $x_{i_1,i_2,\ldots,i_N} =
b^{(1)}_{i_1} b^{(2)}_{i_2} \cdots b^{(N)}_{i_N}$.\\

\noindent{\bf Kruskal tensor, CP decomposition.} For further discussion, it  is important to highlight  that any tensor can be expressed
as a finite sum of rank-1 tensors, in the form
\begin{equation}
\label{eq:CPD}
\underline \bX= \sum_{r=1}^R \bb_r^{(1)} \circ \bb_r^{(2)}  \circ \cdots \circ \bb_r^{(N)} = \sum_{r=1}^R \; \left(\overset{N} {\underset{n=1} {\circ}} \;\bb_r^{(n)}\right), \ \ \bb_r^{(n)} \in \Real^{I_n},
\end{equation}
 which is exactly the form of the Kruskal tensor, illustrated in   Figure \ref{Fig:CPD4}, also known under the names of CANDECOMP / PARAFAC,
  Canonical Polyadic Decomposition (CPD), or simply the CP decomposition in (\ref{eq:CPfun}). We will use the acronyms CP and CPD.\\

\noindent{\bf Tensor rank.} The tensor rank, also called the CP rank, is a natural extension of the matrix rank and is defined as a minimum  number, $R$, of rank-1 terms in an exact CP decomposition of the form in (\ref{eq:CPD}).

Although   the CP decomposition has already  found many practical applications,
its limiting theoretical property  is
 that the best rank-$R$ approximation of a given data tensor may not  exist (see  \cite{deSilva-Lim08} for more detail).
However, a rank-$R$ tensor can be approximated arbitrarily well by a  sequence of tensors for which the
CP ranks are strictly less than $R$.  For these reasons, the concept of  border
rank was proposed \cite{Bini1985BR}, which is defined as the minimum number of rank-1 tensors which
 provides the approximation of a given tensor with an arbitrary accuracy.\\

\noindent{\bf Symmetric tensor decomposition.} A symmetric tensor (sometimes called a super-symmetric tensor) is invariant to the permutations of its indices.
A  symmetric tensor of $N$th-order has equal sizes,  $I_n=I, \, \forall n$,  in all its
modes, and the same value of entries for every permutation of its indices. For example, for vectors $\bb^{(n)}= \bb \in \Real^{I}, \; \forall n$, the  rank-1 tensor, constructed by  $N$ outer products,
$\circ_{n=1}^N \bb^{(n)} =\bb \circ \bb \circ \cdots \circ \bb \in \Real^{I \times I \times \cdots \times I}$, is symmetric.
Moreover, every symmetric tensor can be expressed as a
linear combination of such symmetric rank-1 tensors through the so-called symmetric CP decomposition, given by
\begin{equation}
\underline \bX= \sum_{r=1}^R \lambda_r \bb_r \circ \bb_r \circ \cdots \circ \bb_r, \qquad \bb_r \in \Real^{I},
\label{eq:symCP}
\end{equation}
where $\lambda_r \in \Real$ are the scaling parameters for the unit length vectors $\bb_r$, while
the symmetric tensor rank  is the minimal number $R$ of rank-1 tensors that is necessary for its exact representation.\\
%

\noindent {\bf Multilinear products.} The  mode-$n$  (multilinear) product, also called  the  tensor-times-matrix product (TTM),
of a tensor,
$\underline \bA \in \Real^{I_{1} \times \cdots \times I_{N}}$, and a matrix,  $\bB \in \Real^{J \times I_n}$,  gives the tensor 
\begin{equation}
\underline \bC = \underline \bA \times_n \bB \in \Real^{I_1 \times \cdots \times I_{n-1} \times J \times I_{n+1} \times \cdots \times I_N},
\end{equation}
with entries
\begin{equation}
c_{i_1,i_2,\ldots,i_{n-1},j,i_{n+1},\ldots, i_{N}} =\sum_{i_n=1}^{I_n} a_{i_1,i_2,\ldots,i_N} \; b_{j,i_n}. \label{eq:TTM_output}
\end{equation}
From (\ref{eq:TTM_output}) and Figure \ref{Fig:AxnB}, the equivalent  matrix form is $\bC_{(n)} =\bB \bA_{(n)}$, which allows us  to employ established fast matrix-by-vector
and matrix-by-matrix multiplications
when dealing with very large-scale tensors. Efficient and optimized algorithms for  TTM are, however, still emerging \cite{TTM2015,Ballard15MbyMimproved,Ballard15MbyM}.\\

\noindent {\bf Full multilinear (Tucker) product.} A full multilinear product, also called the Tucker product,
 of an $N$th-order tensor, $\underline \bG \in \Real^{R_1 \times R_2 \times \cdots \times R_N}$, and a set of $N$ factor matrices,
  $\underline \bB^{(n)} \in \Real^{I_n \times R_n}$ for $n=1,2,\ldots, N$,
 performs the multiplications  in all the modes and  can be compactly written as
  (see Figure \ref{Fig:TO1}(b))
\be
\underline \bC &=& \underline \bG \times_1 \bB^{(1)}  \times_2 \bB^{(2)}  \cdots  \times_N \bB^{(N)}  \\
&=& \llbracket \underline \bG; \bB^{(1)}, \bB^{(2)}, \ldots, \bB^{(N)} \rrbracket \in \Real^{I_1 \times I_2 \times \cdots \times I_N}.\notag
\label{mprod}
\ee
Observe that this format corresponds to the Tucker decomposition \cite{tucker64extension,Tucker1966,Kolda08} (see Section~\ref{sect:Tucker}).\\

\noindent {\bf Multilinear product of a tensor and a vector (TTV).} In a similar way, the mode-$n$ multiplication of a tensor,
$\underline \bA \in \Real^{I_{1} \times \cdots \times I_{N}}$, and a vector,  $\bb \in \Real^{I_n}$ (tensor-times-vector, TTV)
  yields a tensor
\be
\underline \bC = \underline \bA \bar \times_n \bb \in \Real^{I_1 \times \cdots \times I_{n-1}  \times I_{n+1} \times \cdots \times I_N},
\ee
 with entries
 \begin{equation}
  c_{i_1,\ldots,i_{n-1},i_{n+1},\ldots, i_{N}} =\sum_{i_n=1}^{I_n} a_{i_1, \ldots,i_{n-1}, i_n, i_{n+1},\ldots,i_N}\; b_{i_n}.
\end{equation}
Note that the mode-$n$ multiplication of a tensor by a matrix does not change the tensor  order,
while the multiplication of a tensor by vectors reduces its order, with the mode $n$ removed (see Figure \ref{Fig:TO1}).

Multilinear products of   tensors by matrices or vectors  play a key role in  deterministic methods for the reshaping of tensors and dimensionality reduction, as well as in
probabilistic methods for randomization / sketching procedures and in  random projections of  tensors into matrices or vectors.
In other words, we can also perform reshaping  of a tensor through  random projections  that change its entries,  dimensionality or size of modes, and/or
 the tensor order. This is achieved by multiplying a tensor by random matrices or vectors,  transformations which preserve  its basic properties.
\cite{sun2006beyond,drineas2007randomized,MSLSurvey2011,li2012hashing,pham2013fast,SmolaNIPS15,kuleshov2015tensor,sorber2016exact}
 (see Section~\ref{sect:Sketching} for more detail).\\

\noindent {\bf Tensor contractions.}
Tensor contraction  is a fundamental and  the most important operation in tensor networks, and  can be considered a  higher-dimensional analogue of matrix multiplication,
inner product, and outer product. 

In a way similar to the  mode-$n$ multilinear product{\footnote{In the literature, sometimes the symbol $\times_n$ is replaced by $\bullet_n$.}},   the mode-$ (_n^m)$ product (tensor contraction) of two tensors,  $\underline \bA \in \Real^{I_1 \times I_2 \times \cdots \times I_N}$ and $\underline \bB \in \Real^{J_1 \times J_2 \times \cdots \times J_M}$,  with  common modes, $I_n=J_m$, yields   an $(N+M-2)$-order tensor, $ \underline \bC  \in \Real^{I_1 \times \cdots \times  I_{n-1} \times I_{n+1} \times  \cdots \times I_N \times J_{1} \times \cdots \times J_{m-1} \times J_{m+1} \times \cdots \times J_M}$, in the form (see Figure \ref{Fig:TO2}(a))
 \be \underline \bC =  \underline \bA \; {\times}_n^m \; \underline \bB,
 \ee
  for which the entries are computed as
   \be && c_{i_1,\, \ldots, \, i_{n-1}, \, i_{n+1}, \, \ldots, i_N, \, j_1, \, \ldots, \, j_{m-1}, \, j_{m+1}, \,\ldots, \, j_M} = \notag \\
   &&= \sum_{i_n=1}^{I_n} a_{i_1, \ldots, i_{n-1}, \, i_n, \, i_{n+1}, \, \ldots,\, i_N} \;  b_{j_1, \, \ldots, \,j_{m-1}, \, i_n, \, j_{m+1}, \, \ldots, \, j_M}.
   \ee
This operation is referred to as a {\it contraction of two tensors in single common mode}.

  Tensors can be contracted in several modes or even in all modes, as illustrated in Figure \ref{Fig:TO2}.
For convenience of presentation, the super- or sub-index, e.g., $m,n$, will be omitted in a few special cases. For example, the multilinear product of the tensors, $\underline \bA \in  \Real^{I_{1} \times I_{2} \times \cdots \times I_{N}}$ and $\underline \bB \in  \Real^{J_{1} \times J_{2} \times \cdots \times J_{M}}$,  with  common  modes, $I_N=J_1$, can be written   as
\begin{equation}
\underline \bC = \underline \bA \; \times_N^1 \; \underline \bB = \underline \bA \times^1 \underline \bB = \underline  \bA \bullet \underline \bB  \in   \Real^{I_{1} \times I_2 \times  \cdots \times I_{N-1} \times  J_2 \times \cdots \times J_{M}},
\end{equation}
for which the entries
\begin{equation}
 c_{i_1,i_2,\ldots,i_{N-1},j_2,j_3,\ldots,j_M}= \sum_{i_N=1}^{I_N} a_{i_1,i_2,\ldots,i_N} \; b_{i_N,j_2,\ldots,j_M}. \notag
\end{equation}

In this notation, the multiplications of matrices and vectors can be written as,  $\bA \times^1_2 \bB = \bA \times^1 \bB =\bA \bB$,
$\;\bA \times^2_2 \bB  =\bA \bB^{\text{T}}$,  $\;\bA \times^{1,2}_{1,2} \; \bB = \bA \bar \times \bB = \langle \bA, \bB\rangle$, and $\bA \times^1_2 \; \bx =\bA \times^1 \bx =\bA \bx$.

Note that tensor contractions are, in general not associative or commutative, since when contracting more than two tensors, the order has to be precisely specified (defined),
for example, $\underline \bA \times^b_a
(\underline \bB \times_c^d \underline \bC)$ for $b < c$.

It is also important to note that a  matrix-by-vector product, $\by=\bA \bx \in \Real^{I_1 \cdots I_N}$, with  $\bA \in \Real^{I_1 \cdots I_N \times J_1 \cdots J_N}$  and
$\bx \in \Real^{J_1 \cdots J_N}$, can be expressed in a tensorized
 form via the contraction operator as $\underline \bY =  \underline  \bA \bar \times \underline  \bX$, where the symbol $\bar \times$ denotes the contraction of all modes of the tensor $\underline \bX$ (see Section~\ref{sect:TT-oper}).

Unlike the  matrix-by-matrix multiplications for which several efficient parallel schemes have
been developed, (e.g. BLAS procedure) 
 the number of efficient algorithms for tensor contractions is rather limited. 
In practice,  due to the high computational complexity of tensor contractions, especially for tensor networks with loops,
 this operation is often  performed  approximately  \cite{ContracPEPS14,Di_Tens_Contrac2014,Pfeifer2014ncon,Uni10Kao2015}. \\

\noindent {\bf Tensor trace.}
Consider a tensor with partial self-contraction modes, where the outer (or open) indices represent physical modes of the tensor, while the inner indices indicate its contraction modes. The Tensor Trace operator performs the summation of all inner indices of the tensor \cite{gu2009tensor}.
For example, a tensor $\underline \bA$ of size $R \times I \times R$ has two inner indices, modes 1 and 3 of size $R$, and one open mode of size $I$.
Its tensor trace yields a vector of length $I$, given by
\begin{equation}
\ba = \mbox{Tr}(\underline \bA)  = \sum_{r}  \underline \bA(r,:,r) \,, \notag
\end{equation}
the elements of which are the traces of its lateral slices $\bA_i \in \Real^{R \times R}$  $\;(i=1,2,\ldots,I)$, that is, (see bottom of Figure  \ref{Fig:trace})
\be
\ba =  [\tr(\bA_1),  \ldots, \tr(\bA_i), \ldots, \tr(\bA_I)]^{\text{T}}.
\ee
A tensor can have more than one pair of inner indices, e.g., the tensor $\underline \bA$ of size $R \times I \times S \times S \times I \times R$ has two pairs of inner indices, modes 1 and 6, modes 3 and 4, and two open modes (2 and 5). The tensor trace of $\underline \bA$ therefore returns a matrix of size $I \times I$ defined as
\be
	\mbox{Tr}(\underline \bA) = \sum_{r} \sum_{s} \underline \bA(r,:,s,s,:,r) \,. \notag
\ee

A variant of Tensor Trace \cite{Lee-TTfund1} for the case of the partial tensor self-contraction considers a tensor $\underline \bA \in \Real^{R \times I_1 \times I_2 \times \cdots \times I_N \times R}$ and yields a reduced-order tensor  $\underline {\widetilde\bA}= \mbox{Tr} (\underline \bA)  \in \Real^{I_1 \times I_2 \times \cdots \times I_N }$,   with entries
\begin{equation}
\underline {\widetilde\bA}(i_1,i_2,\ldots,i_N) = \sum_{r=1}^R \underline \bA (r,i_1,i_2,\ldots,i_N,r),
\end{equation}
Conversions of tensors to scalars, vectors, matrices or  tensors with  reshaped modes and/or  reduced orders are illustrated in  Figures  \ref{Fig:TO1}-- \ref{Fig:trace}.

\section{Graphical Representation of Fundamental Tensor Networks}

 \begin{figure}[t]
 \centering
\includegraphics[width=11.75cm]{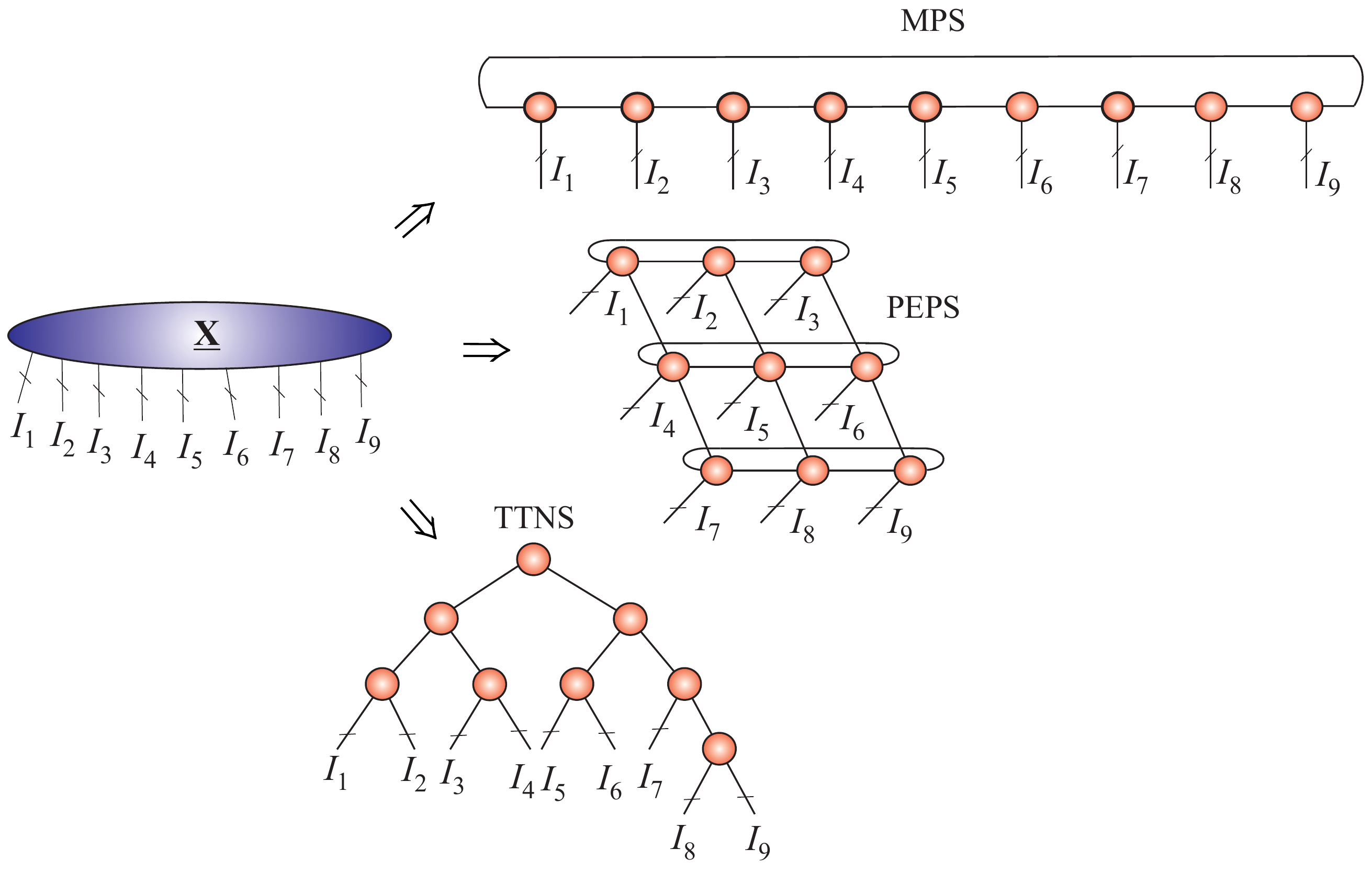}
\caption{Illustration of the decomposition of a 9th-order tensor, $\underline \bX \in \Real^{I_1 \times I_2 \times \cdots \times I_9}$, into different forms of tensor networks (TNs). In general, the objective is to decompose a very high-order tensor into sparsely (weakly) connected low-order and small size tensors, typically 3rd-order and 4th-order tensors called cores. Top: The Tensor Chain (TC) model, which is equivalent to the  Matrix Product State (MPS) with periodic boundary conditions (PBC). Middle:  The Projected Entangled-Pair States (PEPS), also with PBC. Bottom: The  Tree Tensor Network State (TTNS).}
\label{Fig:TN1}
\end{figure}

Tensor networks (TNs) represent a higher-order tensor  as a set of sparsely interconnected  lower-order tensors (see Figure  \ref{Fig:TN1}), and in this way provide computational and storage benefits.
The lines (branches, edges)  connecting core tensors correspond to the contracted modes while their  weights (or numbers of branches)  represent the rank of a tensor network{\footnote{Strictly speaking, the minimum set of internal indices $\{R_1,R_2,R_3,\ldots \}$ is called the  rank (bond dimensions) of a specific tensor network.}}, whereas the lines
which do not connect core tensors  correspond to the ``external''  physical variables (modes, indices) within the data tensor. In other words, the number of free (dangling) edges (with weights larger than one) determines the order of a data tensor  under consideration, while set of weights of internal branches represents  the TN rank.

\section{Hierarchical Tucker (HT) and Tree Tensor Network State  (TTNS) Models}
\label{sect:HT}

\begin{figure}[t]
\centering
  \includegraphics[width=10.8cm]{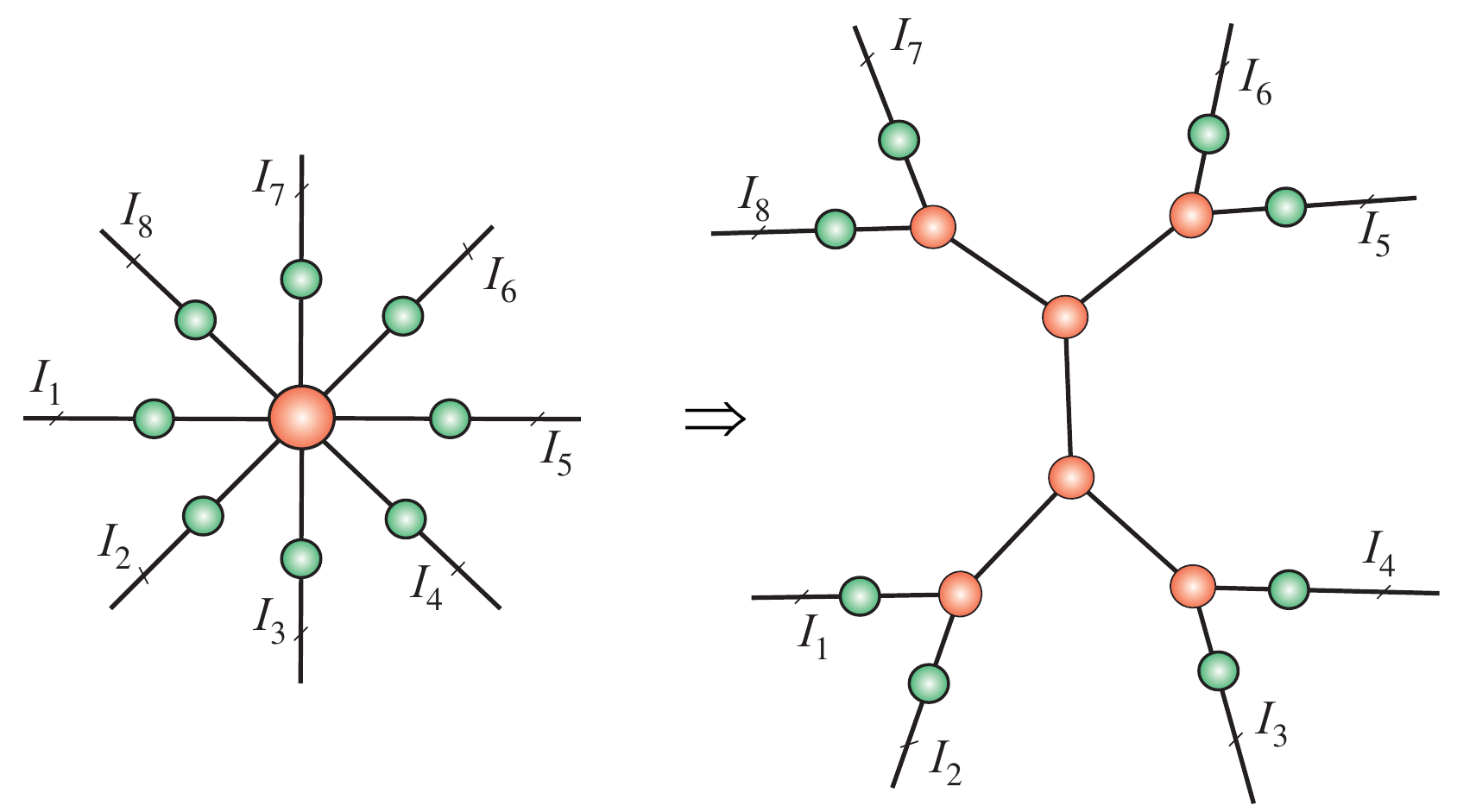}
\caption{The standard Tucker decomposition of an 8th-order
tensor into a core tensor (red circle) and eight factor matrices (green circles), and its transformation into an equivalent  Hierarchical Tucker (HT) model  using  interconnected smaller size  3rd-order core tensors and the same factor matrices.}
\label{Fig:HT8}
\end{figure}

\begin{figure}
\begin{center}
\includegraphics[width=9.54cm]{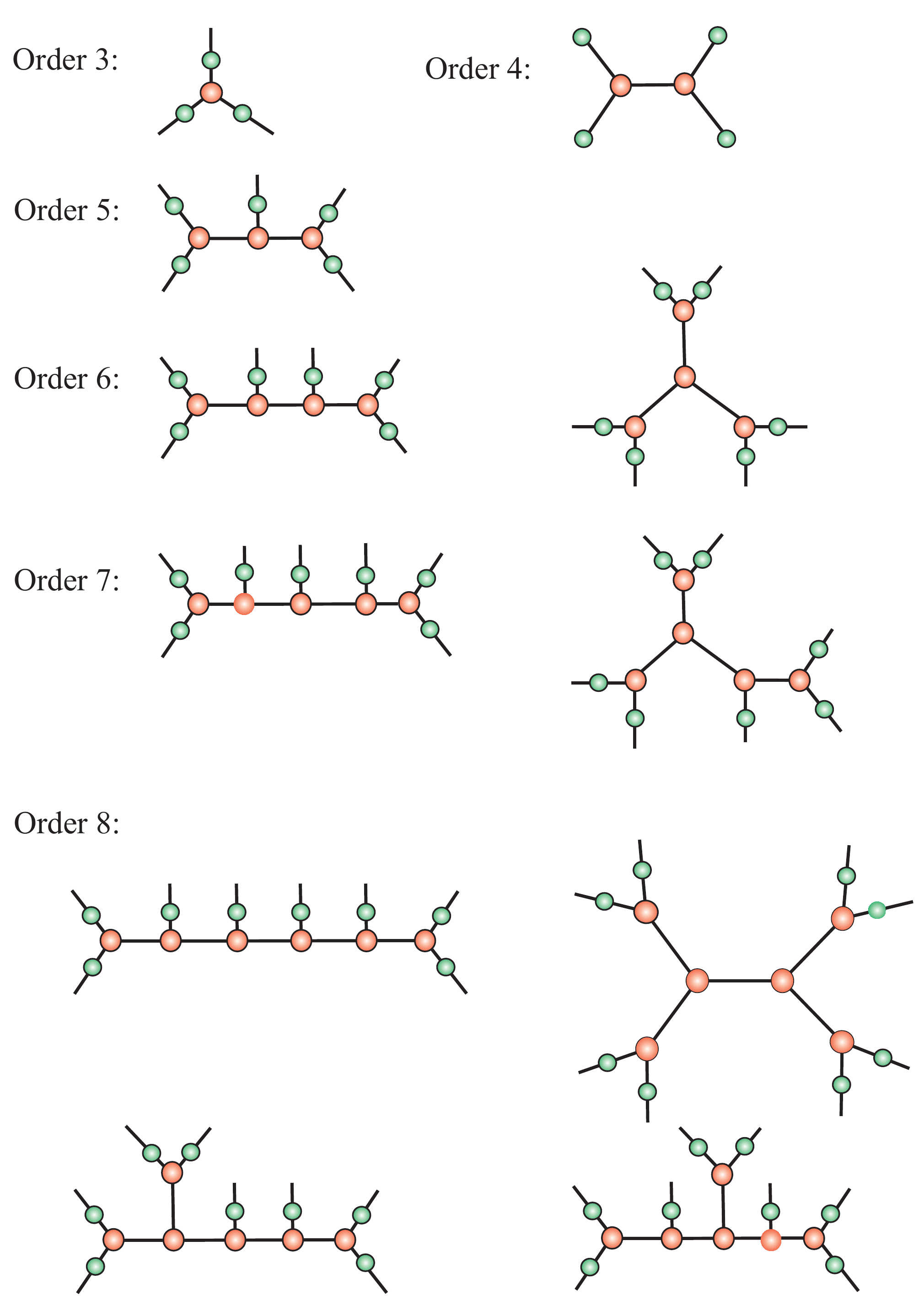}
 \end{center}
\caption{Examples of HT/TT models (formats) for distributed Tucker decompositions with 3rd-order cores, for different orders of data tensors. Green circles denote factor matrices (which can be absorbed by core tensors), while red circles indicate cores. Observe that the representations are not unique.}
\label{Fig:HT8b}
\end{figure}

\begin{figure}[t!]
\begin{center}
\includegraphics[width=9.4cm]{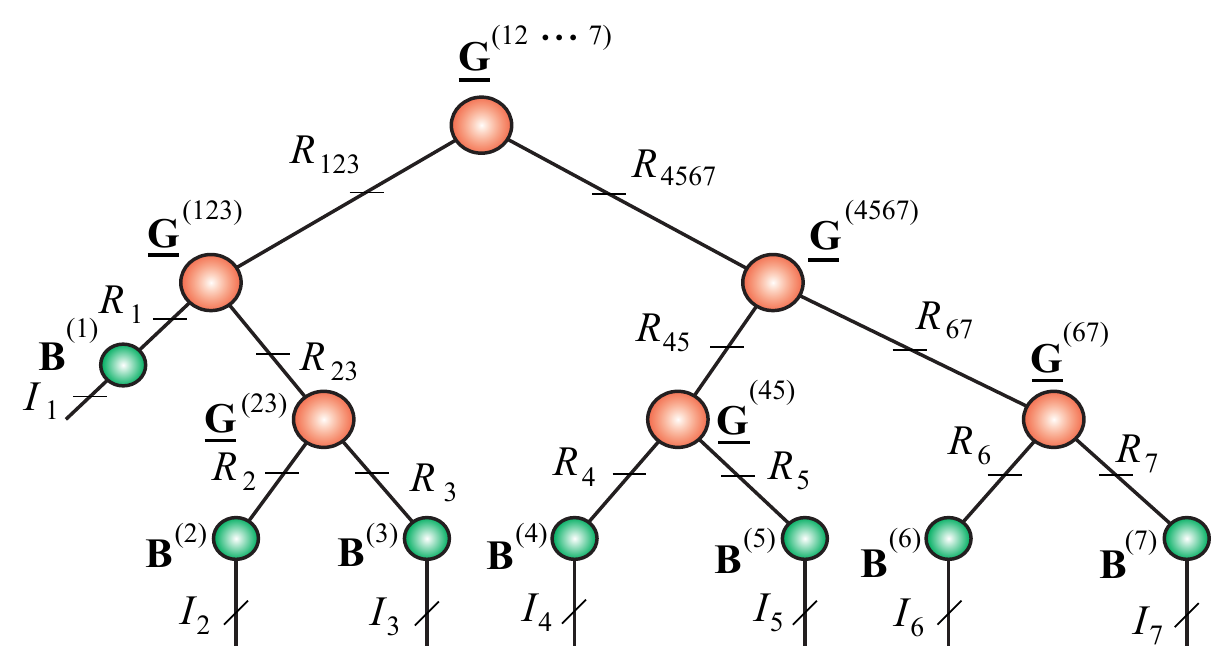}
\end{center}
\caption{Example illustrating the HT decomposition for a 7th-order data tensor.}
\label{Fig:HT7}
\end{figure}

Hierarchical Tucker (HT) decompositions (also called hierarchical tensor representation)
 have been introduced in   \cite{HackbuschHT09}
and also independently in \cite{hTucker1},
see also \cite{Hackbush2012,Lubich_Schneider13,Uschmajew_Vander2013,KressnerTobler14,Bachmayr2016}
and references therein{\footnote{The HT model was developed independently, from a different perspective,
in the  chemistry
 community under the name MultiLayer Multi-Configurational Time-Dependent Hartree method (ML-MCTDH) \cite{wang2003MCTDH}.
 Furthermore, the PARATREE model, developed independently for signal processing applications \cite{VisaSP-09},
 is quite similar to the HT model \cite{hTucker1}.}}.
Generally, the  HT decomposition requires splitting
the set of modes of a tensor in a hierarchical way,  which results in a binary tree containing  a subset
of modes at each branch (called a dimension tree); examples of binary trees
are given in Figures  \ref{Fig:HT8}, \ref{Fig:HT8b} and \ref{Fig:HT7}.
In tensor networks based on binary trees, all the cores are of order of three  or less. Observe that the HT model does not contain any cycles (loops), i.e., no edges connecting a node with itself. The splitting operation of the set of modes of the original data tensor by binary tree edges
 is performed through a suitable matricization. 

{\bf Choice of dimension tree.} The dimension tree within the HT format is chosen {\it a priori} and defines  the topology of the HT decomposition. Intuitively, the dimension
tree specifies which groups of  modes
are ``separated'' from other groups of modes, so that  a sequential HT decomposition can be performed
via a (truncated) SVD applied to a suitably  matricized tensor.
One of the simplest and most straightforward choices of a dimension tree
is the linear and unbalanced tree, which gives rise to the tensor-train (TT) decomposition, discussed in detail in    Section~\ref{sect:TT1} and
 Section~\ref{chap:TT} \cite{OseledetsTT11,OseledetsTT09}.

Using mathematical formalism, a dimension tree is a binary tree $T_N$, $N>1$, which satisfies that
\begin{enumerate}
\item[(i)] all nodes $t\in T_N$ are non-empty subsets of \{1, 2,\ldots, N\},
\item[(ii)] the set $t_{root} = \{1,2,\ldots,N\}$ is the root node of $T_N$, and
\item[(iii)] each non-leaf node has two children $u,v\in T_N$ such that $t$ is a disjoint union
$t = u \cup v$.
\end{enumerate}
The HT model is illustrated through the following Example.\\[1em]
{\bf Example.} Suppose that the dimension tree $T_7$ is given, which gives  the HT decomposition illustrated in Figure~\ref{Fig:HT7}.
The HT decomposition of a tensor $\underline \bX \in \mathbb{R}^{I_1\times\cdots\times I_7}$ with given set of integers $\{R_t\}_{t\in T_7}$
can be expressed in the tensor and vector / matrix forms as follows.
Let intermediate tensors
$\underline \bX^{(t)}$ with $t=\{n_1,\ldots, n_k\}\subset\{1,\ldots,7\}$ have the size $I_{n_1}\times I_{n_2}\times\cdots\times I_{n_k} \times R_t$. Let $\underline \bX^{(t)}_{r_t} \equiv \underline \bX^{(t)}(:,\ldots,:,r_t)$ denote the subtensor of $\underline \bX^{(t)}$ and
$\bX^{(t)} \equiv \bX^{(t)}_{<k>} \in\mathbb{R}^{I_{n_1}I_{n_2}\cdots I_{n_k} \times R_t}$ denote the corresponding unfolded matrix.
Let $\underline \bG^{(t)}\in\mathbb{R}^{R_u\times R_v\times R_t}$ be core tensors where $u$ and $v$ denote respectively the left and right children of $t$.

The HT model shown in Figure \ref{Fig:HT7} can be then described mathematically in the vector form as
\begin{align}
	&\text{vec}(\underline \bX) \cong ( \bX^{(123)} \otimes_L \bX^{(4567)} ) \; \text{vec}(\bG^{(12\cdots 7)}) , \notag \\
\notag \\
	&\bX^{(123)}
	\cong
		( \bB^{(1)} \otimes_L \bX^{(23)} ) \; \bG^{(123)}_{<2>} , \quad
	\bX^{(4567)}
	\cong
		( \bX^{(45)} \otimes_L \bX^{(67)} ) \; \bG^{(4567)}_{<2>} ,\notag \\
\notag \\
	&\bX^{(23)}
	\cong
		( \bB^{(2)} \otimes_L \bB^{(3)} ) \; \bG^{(23)}_{<2>} , \qquad
	\bX^{(45)}
	\cong
		( \bB^{(4)} \otimes_L \bB^{(5)} ) \; \bG^{(45)}_{<2>} , \notag \\
\notag\\
	& \bX^{(67)}
	\cong
		( \bB^{(6)} \otimes_L \bB^{(7)} ) \; \bG^{(67)}_{<2>} . \notag
	\end{align}

An equivalent, more explicit form,  using  tensor notations becomes
\begin{align}
	&\underline \bX \cong
		\sum_{r_{123}=1}^{R_{123}} \sum_{r_{4567}=1}^{R_{4567}}
		g^{(12\cdots 7)}_{r_{123},r_{4567}} \;
		{\underline \bX}^{(123)}_{r_{123}} \circ
		{\underline \bX}^{(4567)}_{r_{4567}} , \notag \\
\notag \\
	&{\underline \bX}^{(123)}_{r_{123}}
	\cong
		\sum_{r_{1}=1}^{R_{1}} \sum_{r_{23}=1}^{R_{23}}
		g^{(123)}_{r_{1},r_{23},r_{123}}
		{\bb}^{(1)}_{r_{1}} \circ
		{\bX}^{(23)}_{r_{23}} , \notag \displaybreak[3]\\
\notag \\
	&{\underline \bX}^{(4567)}_{r_{4567}}
	\cong
		\sum_{r_{45}=1}^{R_{45}} \sum_{r_{67}=1}^{R_{67}}
		g^{(4567)}_{r_{45},r_{67},r_{4567}}
		{\bX}^{(45)}_{r_{45}} \circ
		{\bX}^{(67)}_{r_{67}},  \notag \displaybreak[3]\\
\notag \\
	&{\bX}^{(23)}_{r_{23}}
	\cong
		\sum_{r_{2}=1}^{R_{2}} \sum_{r_{3}=1}^{R_{3}}
		g^{(23)}_{r_{2},r_{3},r_{23}} \;
		{\bb}^{(2)}_{r_{2}} \circ
		{\bb}^{(3)}_{r_{3}} , \notag \\
\notag \\
	&{\bX}^{(45)}_{r_{45}}
	\cong
		\sum_{r_{4}=1}^{R_{4}} \sum_{r_{5}=1}^{R_{5}}
		g^{(45)}_{r_{4},r_{5},r_{45}} \;
		{\bb}^{(4)}_{r_{4}} \circ
		{\bb}^{(5)}_{r_{5}} , \notag \\
\notag \\
	& {\bX}^{(67)}_{r_{67}}
	\cong
		\sum_{r_{6}=1}^{R_{6}} \sum_{r_{7}=1}^{R_{7}}
		g^{(67)}_{r_{6},r_{7},r_{67}} \;
		{\bb}^{(6)}_{r_{6}} \circ
		{\bb}^{(7)}_{r_{7}} . \notag
	\end{align}
The TT/HT decompositions lead  naturally  to  a distributed   Tucker decomposition, where a single core tensor is replaced by
interconnected cores of lower-order,
 resulting in  a distributed network  in  which  only some cores are connected  directly with factor matrices, as
 illustrated in Figure   \ref{Fig:HT8}.
 Figure  \ref{Fig:HT8b} illustrates exemplary HT/TT structures for data tensors of various orders \cite{Toblerthesis,KressnerTobler14}.
Note that for a 3rd-order tensor, there is only one HT  tensor network representation,
 while for a 5th-order we have 5, and for a 10th-order tensor there are 11 possible HT architectures.

\begin{figure}[t]
\begin{center}
\includegraphics[width=11.87cm]{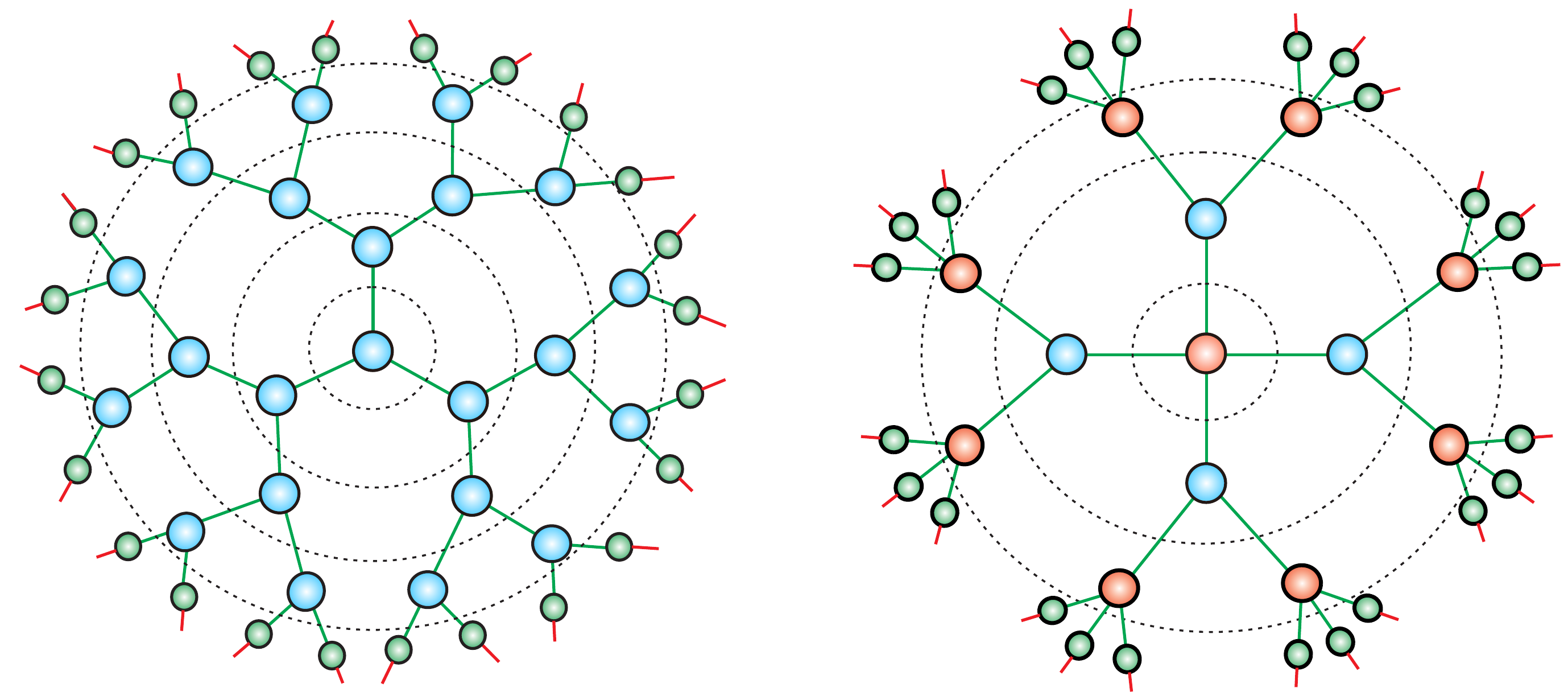}
\end{center}
\caption{The Tree Tensor  Network State (TTNS) with 3rd-order and 4th-order cores for the representation of 24th-order data tensors.
The TTNS can be considered both as a generalization of HT/TT format  and
 as a distributed model for the Tucker-$N$ decomposition (see Section~\ref{sect:Tucker}).}
\label{Fig:TTNS}
\end{figure}

A simple  approach to reduce the size  of a large-scale core tensor in the standard  Tucker decomposition (typically, for $N>5$)
would be to apply  the concept of distributed tensor networks (DTNs). The DTNs assume two kinds of cores (blocks):
(i)~the internal  cores (nodes) which are connected only to other cores and  have no free edges and (ii)~external cores which do
have free edges representing physical modes (indices) of a given data tensor (see also Section~\ref{sect:distrTN}).
Such distributed representations of   tensors are not unique.

The  tree tensor network state (TTNS) model, whereby all nodes are of 3rd-order or higher, can be considered as a generalization of
 the TT/HT decompositions, as illustrated  by two examples in Figure \ref{Fig:TTNS} \cite{DMRG2013}. 
 A more detailed mathematical description of the TTNS 
  is given in Section~\ref{sect:Tucker}.

\begin{figure}
(a)
\vspace{-0.15cm}
\begin{center}
 \includegraphics[width=11.76cm]{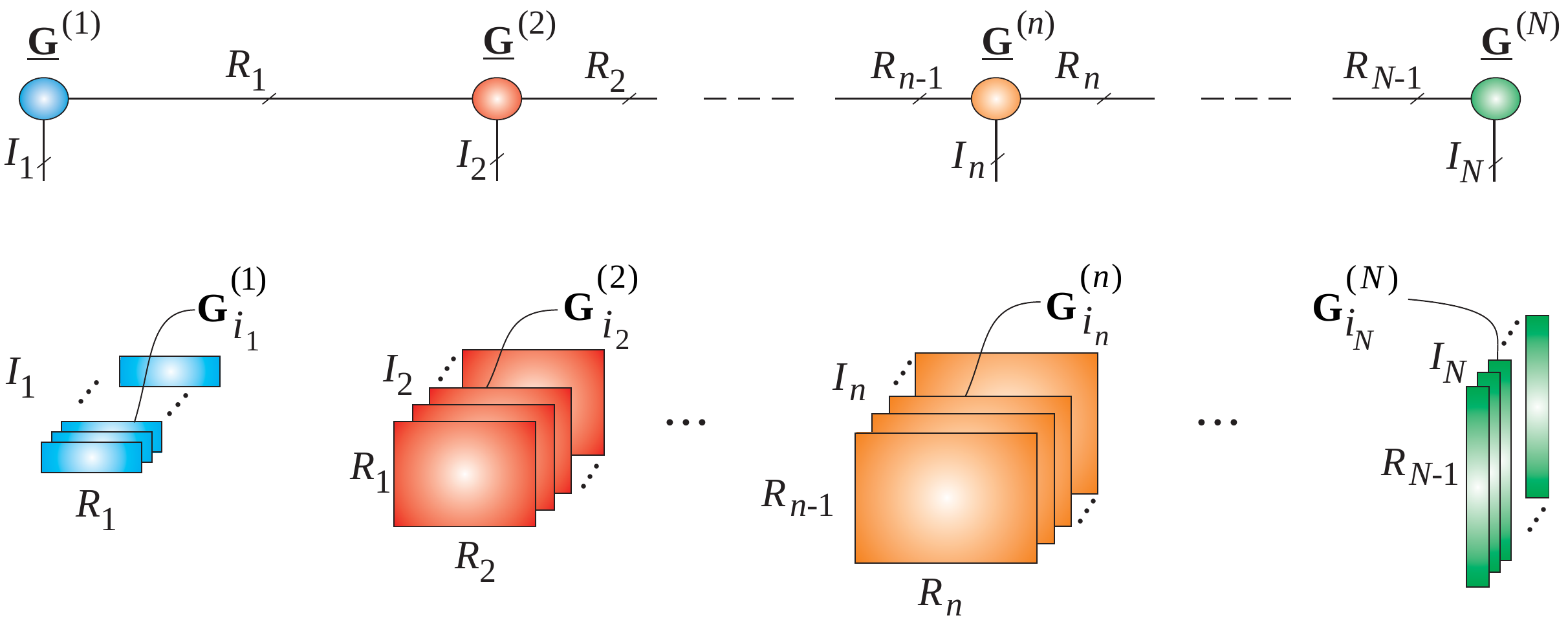}
\end{center}
(b)
\vspace{-0.1cm}
\begin{center}
 \includegraphics[width=11.27cm]{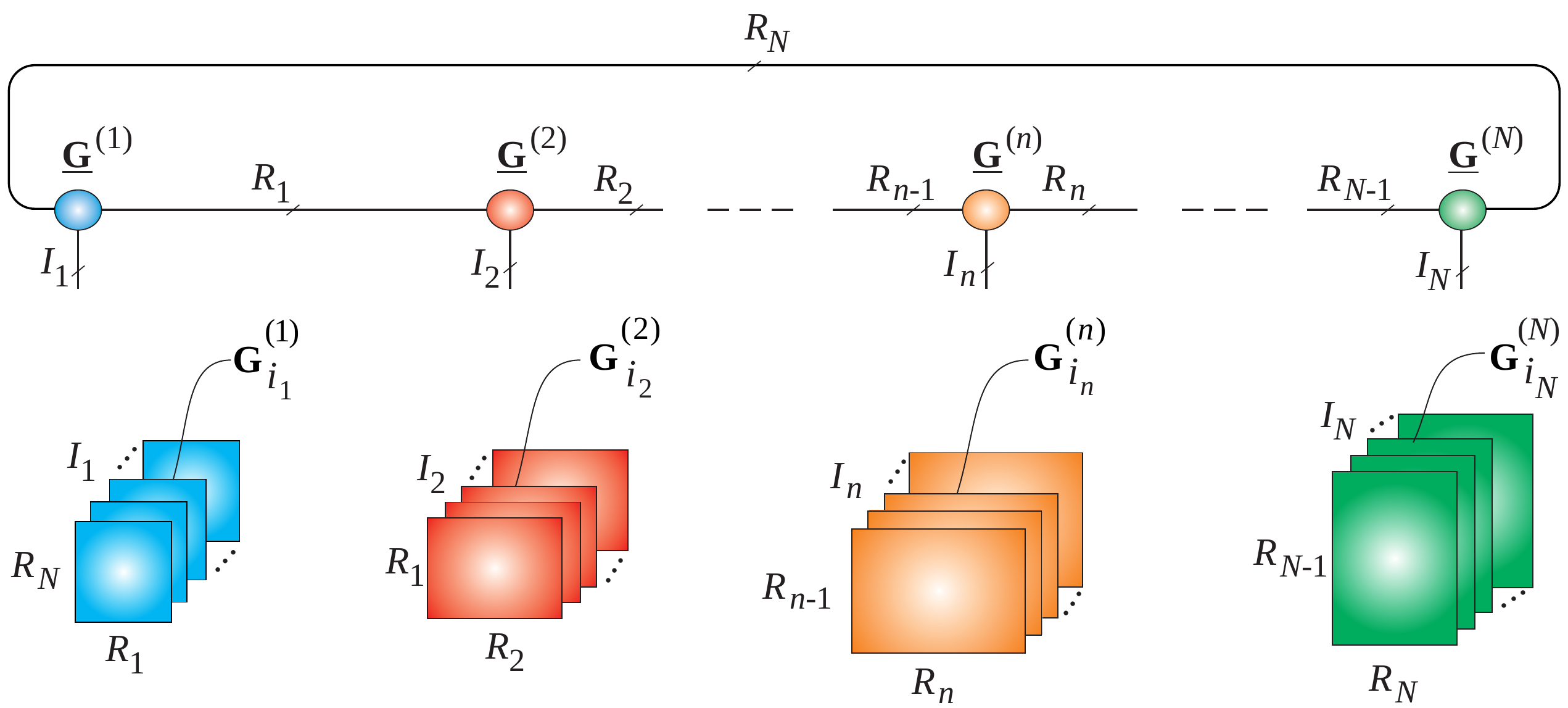}
\end{center}
\caption{{\small Concepts of the tensor train (TT) and tensor chain  (TC) decompositions (MPS with OBC and PBC, respectively) for an $N$th-order data tensor,
$\underline \bX \in \Real^{I_1 \times I_2  \times \cdots \times I_N}$. (a) Tensor Train (TT)
 can be mathematically described as $x_{i_1,i_2,\ldots,i_N} =
\bG^{(1)}_{i_1} \;  \bG^{(2)}_{i_2} \; \cdots \;  \bG^{(N)}_{i_N}$, where (bottom panel) the slice matrices of TT-cores
$\underline \bG^{(n)} \in \Real^{R_{n-1} \times I_n \times R_n}$ are defined as
$\bG^{(n)}_{i_n}= \underline \bG^{(n)}(:, i_n,:)\allowbreak \in \Real^{R_{n-1} \times R_n}$  with $R_0=R_N=1$.
(b) For the Tensor Chain (TC), the entries of a tensor are expressed as  $x_{i_1,i_2,\ldots,i_N} =
\tr \, (\bG^{(1)}_{i_1} \;  \bG^{(2)}_{i_2} \; \cdots \;  \bG^{(N)}_{i_N}) = \displaystyle \sum_{r_1=1}^{R_1} \sum_{r_2=1}^{R_2} \cdots \sum_{r_N=1}^{R_N}
 g^{(1)}_{\;r_N, \,i_1, \,r_1}  \; g^{(2)}_{\;r_{1}, \,i_2, \,r_2} \cdots g^{(N)}_{\;r_{N-1},\, i_N, \,r_N}$, where (bottom panel) the lateral slices of the TC-cores are defined as $\bG^{(n)}_{i_n}= \underline \bG^{(n)}(:, \, i_n,\, :) \in \Real^{R_{n-1} \times R_n}$  and $g^{(n)}_{\; r_{n-1}, \, i_n, \,r_n} = \underline \bG^{(n)} (r_{n-1}, \,i_n, \,r_n)$ for $n=1,2, \ldots,N$, with $R_0=R_N >1$. Notice that TC/MPS is effectively a TT with a single loop connecting the first and the last core, so that all TC-cores are of 3rd-order.}}
\label{Fig:TTouter}
\end{figure}

\section{Tensor Train (TT) Network}
\label{sect:TT1}

\begin{figure}[t]
\centering
\includegraphics[width=9.6cm]{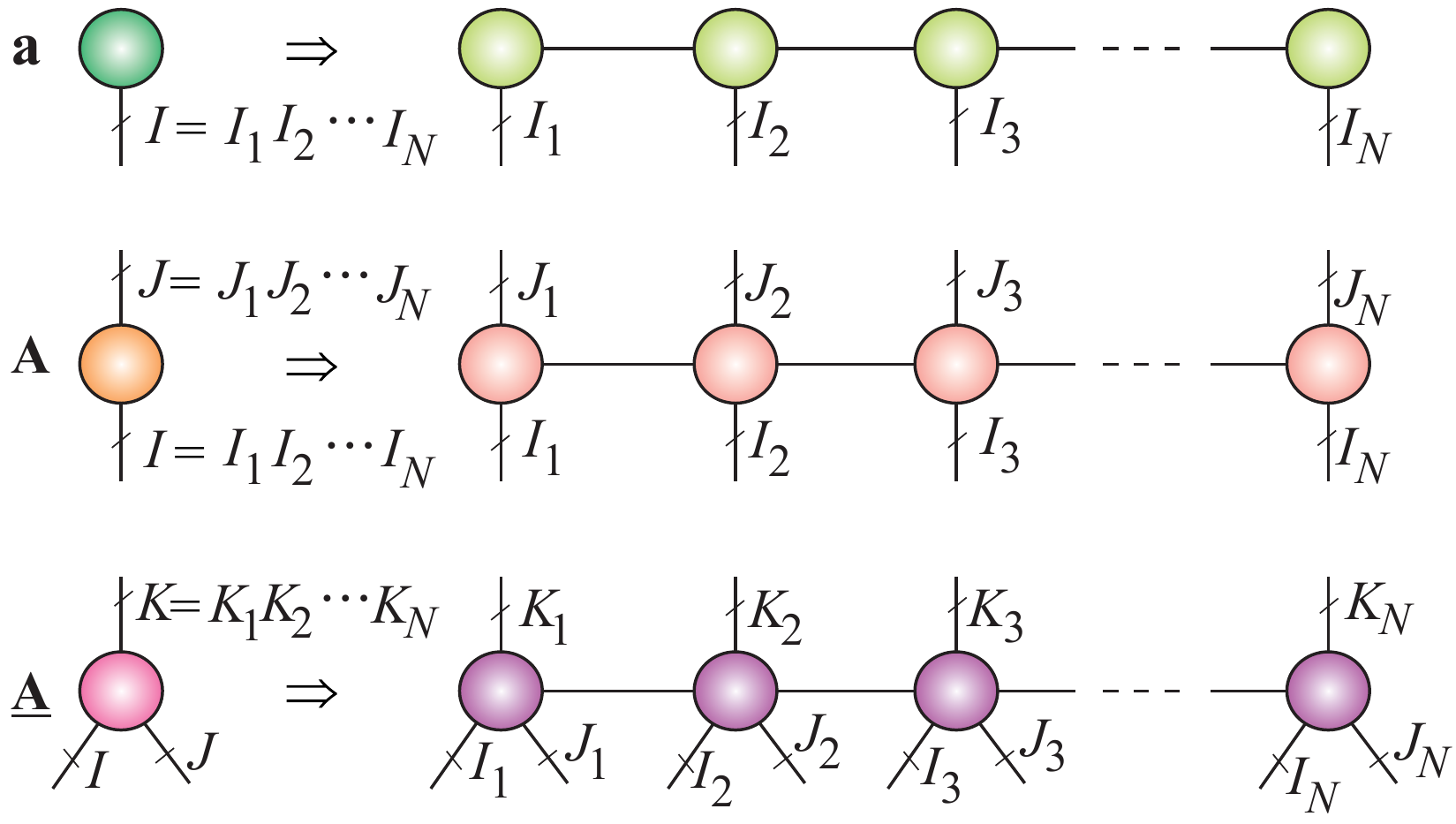}
\caption{Forms of tensor train decompositions for  a  vector, $\ba \in \Real^I$,
 matrix, $\bA \in \Real^{I \times J}$, and  3rd-order  tensor, $\underline \bA \in \Real^{I \times J \times K}$ (by applying a suitable tensorization).}
\label{Fig:VariousTT}
\end{figure}

The Tensor Train (TT) format  can be  interpreted as a special case of the HT format, where all nodes (TT-cores) of the underlying tensor network  are  connected in cascade (or train), i.e., they are aligned
while  factor matrices corresponding to the leaf modes are assumed to be identities and thus need not be stored.
The TT format was first proposed in  numerical analysis and scientific computing in  \cite{OseledetsTT11,OseledetsTT09}.
Figure \ref{Fig:TTouter} presents the concept of TT decomposition for an $N$th-order tensor, the entries of which can be  computed as a cascaded (multilayer) multiplication of appropriate matrices (slices of TT-cores).  The weights of internal edges (denoted by $\{R_1,R_2,\ldots, R_{N-1}\}$) represent the TT-rank. In this way,
the so aligned sequence of core tensors represents a ``tensor train'' where the role of ``buffers'' is played by TT-core connections. It is important to highlight that TT networks can be applied  not only for the approximation of tensorized vectors but also for scalar multivariate functions, matrices, and even large-scale low-order tensors, as illustrated in Figure \ref{Fig:VariousTT} (for more detail see Section~\ref{chap:TT}).

In the quantum physics community, the TT format is known as the Matrix Product
State (MPS) representation with the Open Boundary Conditions (OBC) and was introduced in 1987 as the ground state of the 1D AKLT model
\cite{Affleck_MPS1987}. It was subsequently extended  by many researchers{\footnote{In fact, the TT  was rediscovered several times under different names: MPS, valence bond states, and density matrix renormalization group (DMRG) \cite{White1993}. The DMRG usually refers  not only to a tensor network format but also the efficient computational algorithms (see also \cite{schollwock11-DMRG,Hubig15} and references therein). Also, in quantum physics the ALS algorithm is called the one-site DMRG, while the Modified ALS (MALS) is known as the two-site DMRG (for more detail, see  Part~2).}} (see \cite{White1993,Vidal03,MPS2007,verstraete08MPS,Schollwock13,Huckle2013,Orus2013} and references therein).\\

\noindent{\bf Advantages of TT formats.} An important advantage of the TT/MPS format over the HT format is its simpler practical implementation, as no binary tree needs to be determined (see Section~\ref{chap:TT}).
Another attractive property of the TT-decomposition is its simplicity when performing  basic mathematical  operations on tensors directly in the TT-format (that is, employing only core tensors).
 These include matrix-by-matrix and matrix-by-vector multiplications, tensor addition, and the entry-wise (Hadamard) product
of tensors. These operations produce tensors, also in the TT-format, which generally exhibit
 increased TT-ranks. A detailed description of basic  operations supported by the TT format is  given in Section~\ref{sect:TT-oper}.
Moreover, only TT-cores need to be stored and processed, which makes the number of parameters to scale linearly in the tensor  order, $N$, of a data tensor
and all mathematical operations are then performed only on the low-order and relatively small size core tensors.

 \begin{figure}[t]
\centering
\includegraphics[width=9.1cm, trim = 0 0.5cm 0 0, clip = true]{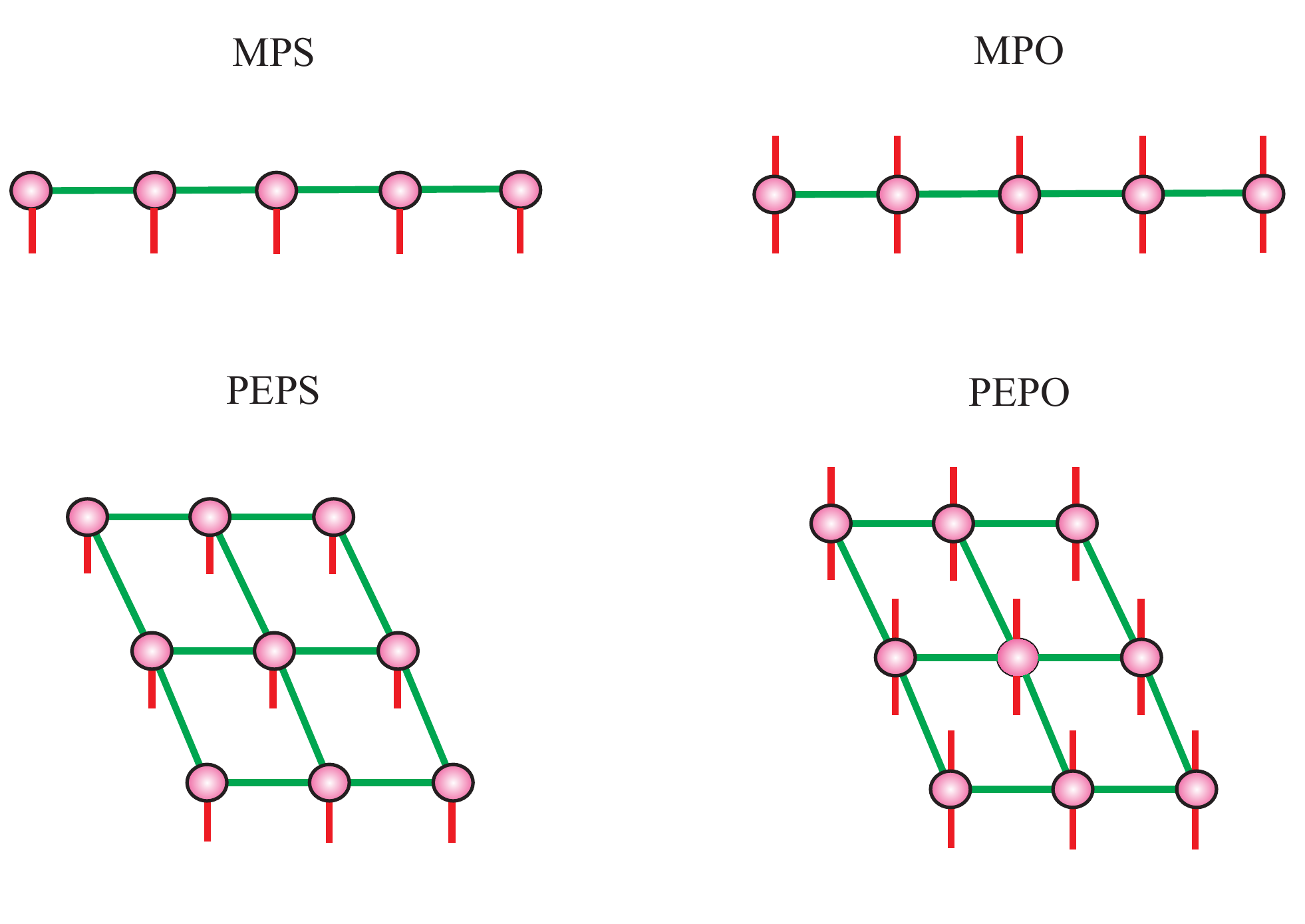}
\caption{Class of 1D and 2D tensor train networks with open boundary conditions (OBC):  the Matrix Product State (MPS) or (vector) Tensor Train (TT),
the Matrix Product Operator (MPO) or Matrix TT, the Projected Entangled-Pair States (PEPS) or Tensor Product State (TPS),
and the Projected Entangled-Pair Operators (PEPO).}
\label{Fig:PEPS}
\end{figure}

The TT rank is defined as an $(N-1)$-tuple of the form
\begin{equation}
\mbox{rank}_{\mbox{\small TT}} (\underline \bX)= \brr_{TT} = \{R_1, \ldots,R_{N-1}\}, \quad R_n=\mbox{rank}(\bX_{<n>}),
\end{equation}
 where $\bX_{<n>} \in \Real^{I_1\cdots I_n \times I_{n-1} \cdots I_N}$ is an $n$th canonical matricization of the tensor $\underline \bX$. Since the TT rank determines memory requirements of a tensor train, it has a strong impact on
the complexity, i.e., the suitability of tensor train representation for a given raw data tensor.
%

The number of data samples to be stored scales linearly in the tensor order, $N$, and the size, $I$, and  quadratically in the maximum TT rank bound, $R$, that is
\be
\sum_{n=1}^N R_{n-1} R_n I_n \sim {\cal{O}} (N R^2 I), \quad R := \max_{n}\{R_n\}, \quad I := \max_{n}\{I_n\}. \qquad
\ee
This is why it is crucially important to have low-rank TT approximations{\footnote{In the worst case scenario  the TT ranks can grow up to $I^{(N/2)}$ for an $N$th-order tensor.}}.
A drawback of the TT format is that the ranks of a tensor train decomposition
depend  on the ordering (permutation)  of the modes, which gives  different size of cores for different ordering.
To solve this challenging permutation problem, we can estimate mutual information between individual TT cores pairwise (see \cite{Barcza2011,ehlers2015entanglement}). The procedure can be arranged in the following three steps: (i) Perform a rough (approximate) TT decomposition with relative low TT-rank and calculate mutual information between all pairs of cores, (ii) order TT cores in such way that the mutual information matrix is close to a diagonal matrix, and finally, (iii) perform TT decomposition again using the so optimised order of TT cores (see also Part 2).

\section{Tensor Networks with Cycles: PEPS,  MERA and Honey-Comb Lattice (HCL)}
\label{sect:MERA}

An important issue in tensor networks is the rank-complexity trade-off in the design.
Namely, the main idea behind TNs is to dramatically reduce computational cost and provide distributed storage and computation  through low-rank TN approximation. However, the TT/HT ranks, $R_n$, of 3rd-order core tensors  sometimes increase rapidly with the order of a data tensor and/or increase of a desired  approximation accuracy,
for any choice of a tree of tensor network. 
The ranks can be often kept under control through hierarchical two-dimensional TT models  called the PEPS (Projected Entangled Pair States{\footnote{An ``entangled pair state'' is  a tensor that cannot be
represented as an elementary rank-1 tensor. The state is called ``projected'' because it is not a real physical state  but a projection onto some subspace. The term ``pair'' refers to the entanglement being  considered only for maximally entangled state pairs \cite{Orus2013,Handschuhth}.}})
 and  PEPO (Projected Entangled Pair Operators) tensor networks, which contain cycles, as shown in Figure \ref{Fig:PEPS}.
In the PEPS and PEPO, the ranks
 are kept considerably smaller  at a cost of  employing  5th- or even 6th-order core  tensors and the associated higher
 computational complexity with respect to the order  \cite{verstraete08MPS,Evenbly-Vidal09alg,schuch2010peps}.

\begin{figure}[t]
(a) \hspace{5.6cm} (b)
\vspace{-2em}
\begin{center}
\includegraphics[width=11.87cm]{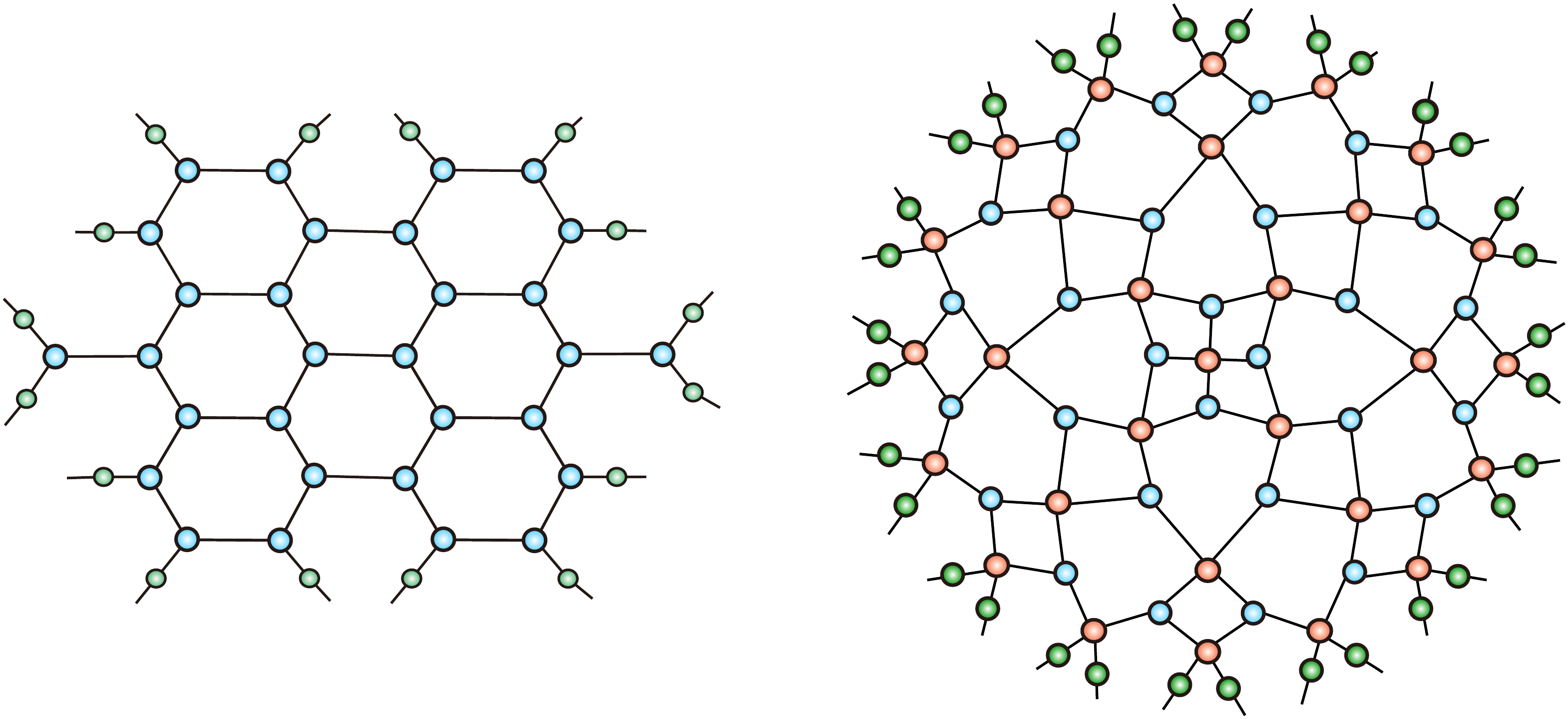}
\end{center}
\vspace{-.8em}
\caption{Examples of TN architectures  with loops. (a) Honey-Comb Lattice  (HCL) for a 16th-order tensor. (b) MERA for  a 32th-order tensor.}
\label{Fig:COMB}
\end{figure}

Even with the PEPS/PEPO architectures, for very high-order tensors, the ranks
 (internal size of cores) may increase rapidly with an increase in the desired
accuracy of approximation. For further control of  the ranks, alternative tensor networks can be  employed, such as: (1) the Honey-Comb Lattice (HCL) which uses 3rd-order cores, and (2)
  the Multi-scale Entanglement Renormalization Ansatz (MERA) which consist of both 3rd- and 4th-order core tensors (see Figure  \ref{Fig:COMB})
\cite{MERA08,Orus2013,matsueda2016MERA}.
 The ranks are often kept considerably small through special architectures  of such TNs, at the expense of  higher computational complexity with respect to tensor contractions due to many cycles.
%

Compared with the PEPS and PEPO formats, the  main advantage of the  MERA formats is that  the order and size of
each  core tensor in the internal tensor network structure is often much smaller, which dramatically reduces the number
of free parameters and provides more efficient distributed storage of huge-scale data tensors.
 Moreover, TNs with cycles, especially the MERA tensor network allow us to model more complex functions and interactions between variables.
%

\section{Concatenated (Distributed)  Representation of TT Networks}
\label{sect:distrTN}

\begin{figure}
\begin{center}
\includegraphics[width=11.21cm]{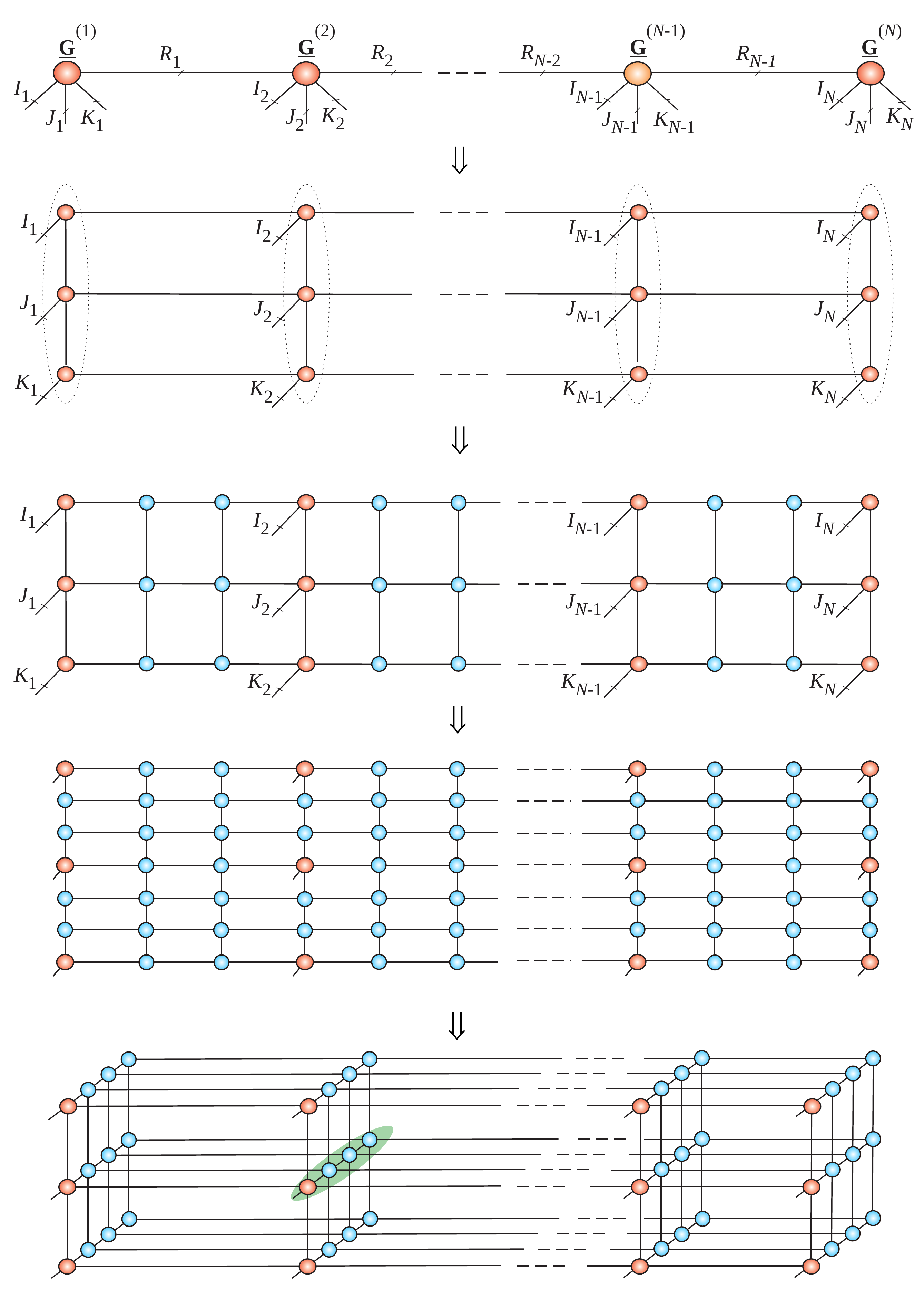}
\end{center}
\caption{Graphical  representation of  a large-scale data tensor via its TT model (top panel), the PEPS model of the TT (third panel), and  its transformation to a  distributed 2D (second from bottom panel) and 3D (bottom panel) tensor train networks.}
\label{Fig:TTPEPS}
\end{figure}

Complexity of algorithms for computation (contraction) on tensor networks typically scales polynomially with the rank, $R_n$, or size, $I_n$, of the core tensors, so that the computations quickly become intractable
 with the increase in $R_n$. A step towards  reducing storage and computational  requirements would be therefore to reduce the size (volume) of core tensors by increasing their number through  distributed tensor networks (DTNs),
 as illustrated   in Figure  \ref{Fig:COMB}.  The underpinning  idea  is that each  core tensor in an original TN is replaced by another TN (see Figure  \ref{Fig:TTPEPS} for TT networks),
resulting in a distributed  TN in which only some core tensors are associated with physical (natural) modes  of the original
data tensor \cite{hubener2010concatenated}.
A DTN  consists of two kinds of relatively  small-size cores (nodes),
internal nodes which have no free edges and external nodes which have free edges representing natural (physical) indices of a data tensor.

The  obvious advantage of DTNs  is that  the size of
each  core tensor in the internal tensor network structure is usually much smaller
than the size of the  initial core tensor; this allows for  a better management of distributed storage, and
 often in the reduction of the total number of network parameters  through distributed computing.
 However, compared to  initial tree structures, the contraction of the resulting distributed tensor network becomes much more difficult because of the loops in the architecture.

\section{Links between  TNs and Machine Learning Models}
\label{sect:analogy}

Table \ref{Table:Morton} summarizes the conceptual connections of  tensor networks with graphical and neural network models
in machine learning and statistics \cite{Morton12,Critchthesis,Morton14,Kazeev14CME,MRF-TT14,DNN_TF2016,Shashua-HT,Shashua_GTD2016,Mera_wavelets1}.
More research is needed to establish  deeper and more precise relationships.

\begin{table}[t]
\setlength{\tabcolsep}{2pt}
\renewcommand{\arraystretch}{1.5}
\centering
\caption{Links between tensor networks (TNs) and
graphical models used in Machine Learning (ML) and Statistics.
The corresponding categories are not exactly the same, but   have general analogies.}
{\shadingbox{
\begin{tabular}{p{3.5cm}|p{7.6cm}}
\hline
Tensor Networks & Neural Networks and Graphical Models in ML/Statistics \\
\hline
TT/MPS & Hidden Markov Models (HMM) \\
HT/TTNS &  Deep Learning Neural Networks, Gaussian Mixture Model (GMM) \\
PEPS & Markov Random Field (MRF), Conditional Random Field (CRF) \\
MERA & Wavelets, Deep Belief Networks (DBN) \\
ALS, DMRG/MALS Algorithms & Forward-Backward Algorithms, Block  Nonlinear Gauss-Seidel Methods \\
\hline
\end{tabular}
}}
\label{Table:Morton}
\end{table}

\section{Changing the Structure of Tensor Networks}
\label{sect:transform}

 An advantage of  the graphical (graph) representation of  tensor networks is that the graphs allow us  to perform  complex mathematical operations on core tensors
 in an intuitive and easy to understand way, without the need to resort to complicated mathematical expressions. Another important advantage
 is the ability to modify (optimize) the topology of a TN, while
 keeping the original physical modes intact.
The so optimized topologies yield simplified or more convenient graphical representations of a higher-order data tensor and facilitate practical applications  \cite{Zhao-Xie-RTNS10,hubener2010concatenated,Handschuhth}. In particular:
\begin{itemize}
  \item A change in topology to a HT/TT tree structure provides reduced computational complexity, through sequential contractions of core tensors and enhanced stability of the corresponding algorithms;
  \item Topology of TNs with cycles can be modified so as to completely eliminate the cycles or to reduce their number;
  \item  Even for vastly diverse original data tensors, topology modifications may produce identical or  similar TN structures  which make it easier to compare and jointly analyze block of interconnected data tensors. This provides opportunity to perform joint group (linked) analysis of tensors by decomposing them to TNs.
\end{itemize}

It is important to note that, due to the iterative way in which tensor contractions are performed, the computational requirements associated with tensor contractions are usually much smaller for tree-structured networks than for tensor networks containing many cycles.
Therefore, for stable computations, it is advantageous to transform a tensor network with cycles into a tree structure.\\

\begin{figure}[ht!]
(a)
\vspace{-0.1cm}
\begin{center}
\includegraphics[width=10.1cm]{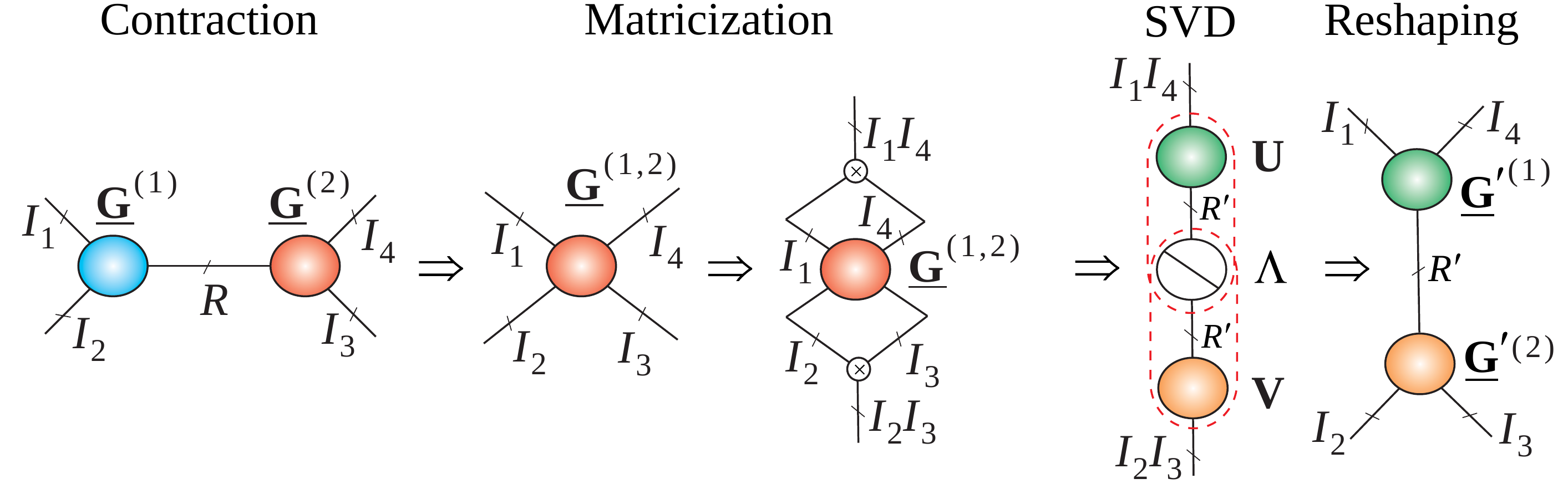}
\end{center}
\vspace{-0.1cm}
(b)
\vspace{-0.1cm}
\begin{center}
\includegraphics[width=10.5cm]{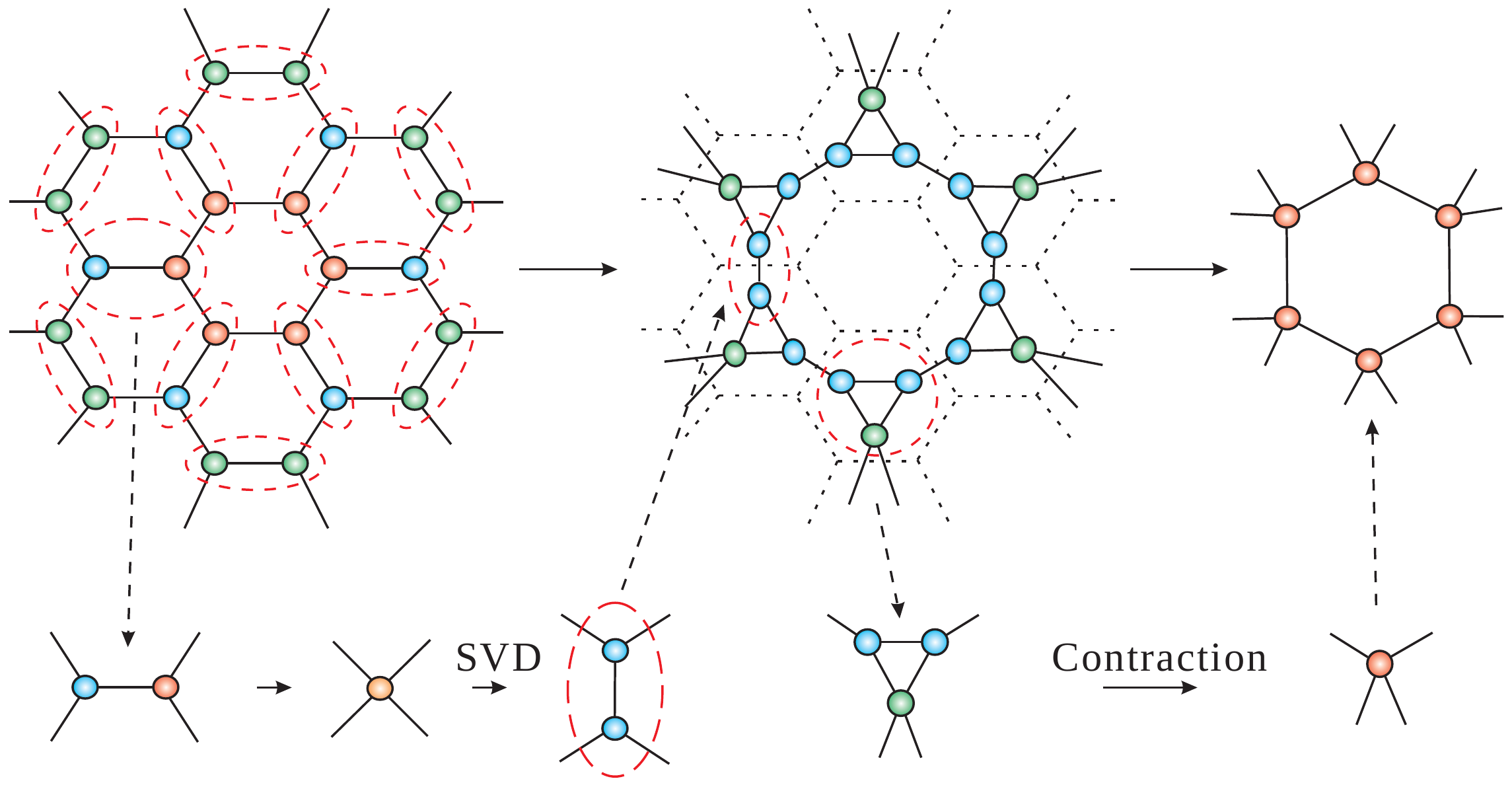}
\end{center}
\vspace{-0.1cm}
\caption{Illustration of  basic transformations on  a tensor network. (a) Contraction, matricization, matrix
factorization (SVD) and reshaping of matrices back into tensors. (b) Transformation of a Honey-Comb lattice into
a Tensor Chain (TC) via tensor contractions and the SVD.}
\label{Fig:TCTT}
\end{figure}

\noindent{\bf Tensor Network transformations.} In order to modify tensor network structures, we may perform sequential core contractions, followed by the unfolding of these contracted tensors into matrices,
 matrix factorizations (typically truncated SVD) and finally reshaping of such matrices back into  new core tensors,
as illustrated in Figures~\ref{Fig:TCTT}.

The example in Figure \ref{Fig:TCTT}(a) shows that,  in the
first step   a contraction of two core tensors, $\underline \bG^{(1)} \in \Real^{I_1 \times I_2 \times R}$ and $\underline \bG^{(2)}  \in \Real^{R \times I_3 \times I_4}$, is performed to give the tensor
\be
\underline \bG^{(1,2)} =\underline \bG^{(1)} \times^1 \; \underline \bG^{(2)} \in \Real^{I_1 \times I_2 \times I_3 \times I_4},
\ee
with entries  $g_{i_1,i_2,i_3,i_4}^{(1,2)} = \sum_{r=1}
^R g_{i_1,i_2,r}^{(1)} \; g_{r,i_3,i_4}^{(2)}$.
In the next step,   the tensor
$\underline \bG^{(1,2)}$ is transformed into a matrix via matricization, followed by a  low-rank matrix factorization  using the SVD, to give
\be
\bG^{(1,2)}_{\overline{i_1 i_4}, \,\overline{i_2 i_3}} \cong \bU  \bS \bV^{\text{T}} \in \Real^{I_1 I_4 \times I_2 I_3}.
\ee
In the final step, the factor matrices, $\bU \bS^{1/2} \in \Real^{I_1 I_4 \times R'}$ and $\bV \bS^{1/2} \in \Real^{R' \times I_2 I_3}$,
are reshaped  into new core tensors, $\underline \bG^{'(1)} \in \Real^{I_1 \times R' \times I_4 }$ and $\underline \bG^{'(2)} \in \Real^{R' \times I_2 \times I_3}$.

\begin{figure}[t]
\centering
\includegraphics[width=7.6cm]{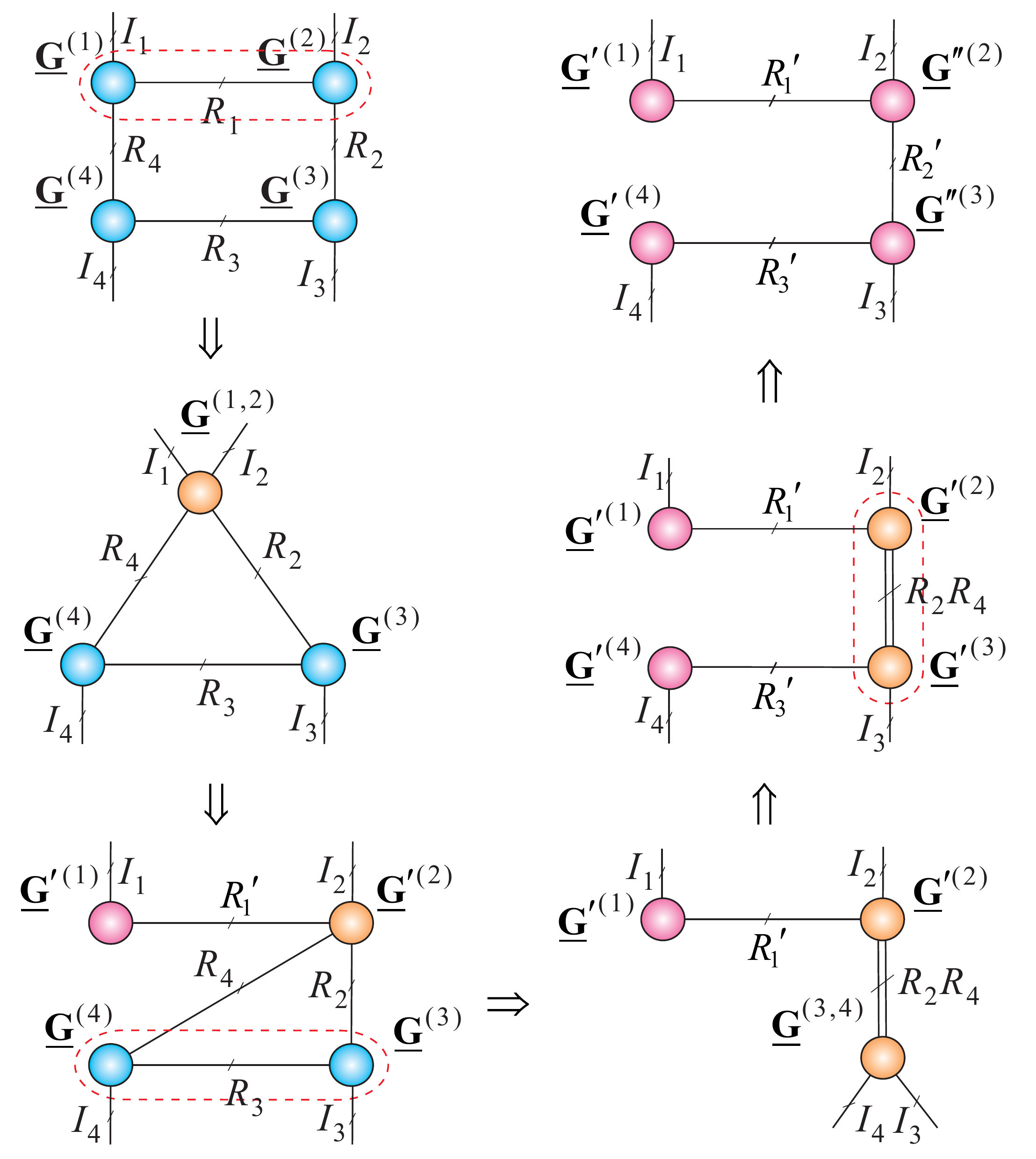}
\caption{Transformation of the closed-loop Tensor Chain (TC) into the open-loop Tensor Train (TT). This is achieved by suitable contractions, reshaping and decompositions of core tensors.}
\label{Fig:TCTT2}
\end{figure}

The above tensor transformation procedure is quite general, and is  applied  in Figure \ref{Fig:TCTT}(b)
to transform a Honey-Comb lattice  into a tensor chain (TC), 
while  Figure \ref{Fig:TCTT2} illustrates the conversion of a tensor chain (TC) into TT/MPS with OBC.

To convert a TC into TT/MPS, in the
first step, we perform  a contraction of two tensors, $\underline \bG^{(1)} \in \Real^{I_1 \times R_4 \times R_1}$ and $\underline \bG^{(2)}  \in \Real^{R_1 \times  R_2 \times I_2}$, as
\begin{equation}
\underline \bG^{(1,2)} = \underline \bG^{(1)} \times^1 \; \underline \bG^{(2)} \in \Real^{I_1 \times R_4 \times R_2 \times I_2},  \notag
\end{equation}
for which the entries  $g_{i_1,r_4,r_2,i_2}^{(1,2)} = \sum_{r_1=1}
^{R_1} g_{i_1,r_4,r_1}^{(1)} \; g_{r_1,r_2,i_2}^{(2)}$.
In the next step, the tensor
$\underline \bG^{(1,2)}$ is transformed into a matrix,  followed by a truncated SVD
\begin{equation}
\bG^{(1,2)}_{(1)} \cong \bU \bS \bV^{\text{T}} \in \Real^{I_1  \times R_4 R_2 I_2}. \notag
\end{equation}
Finally, the matrices, $\bU  \in \Real^{I_1 \times R'_1}$ and $\bV \bS \in \Real^{R'_1 \times R_4 R_2 I_2}$,
 are reshaped back into the core tensors, $\underline \bG^{'(1)}=\bU  \in \Real^{1 \times I_1 \times R'_1}$ and $\underline \bG^{'(2)} \in \Real^{R'_1 \times R_4 \times  R_2 \times I_2}$.
The procedure is repeated  all over again for different pairs of cores, as illustrated in Figure \ref{Fig:TCTT2}.

\section{Generalized Tensor Network Formats}
\label{sect:HOPTA}

The fundamental TNs considered so far assume that the links between the cores are expressed by tensor contractions.
In general, links between the core tensors (or tensor sub-networks) can also be expressed via  other mathematical linear/multilinear or nonlinear operators, such as the
outer (tensor) product, Kronecker product,  Hadamard product and convolution operator.
For example, the use of the outer product  leads to Block Term Decomposition
(BTD) \cite{Lath-BCM12,Lath-Nion-BCM3,LathauwerTBSS,Sorber2012b} and use
the Kronecker products yields  to the Kronecker Tensor Decomposition (KTD) \cite{Ragnarsson-PHD,Phan2012-Kron,Phan_BTDLxR}.
Block term decompositions (BTD) are closely related to constrained Tucker formats (with a sparse block Tucker core) and the Hierarchical Outer Product Tensor Approximation (HOPTA), which be employed for very high-order data tensors \cite{Cichocki2014optim}.

\begin{figure}[t]
\centering
\includegraphics[width=9.5cm]{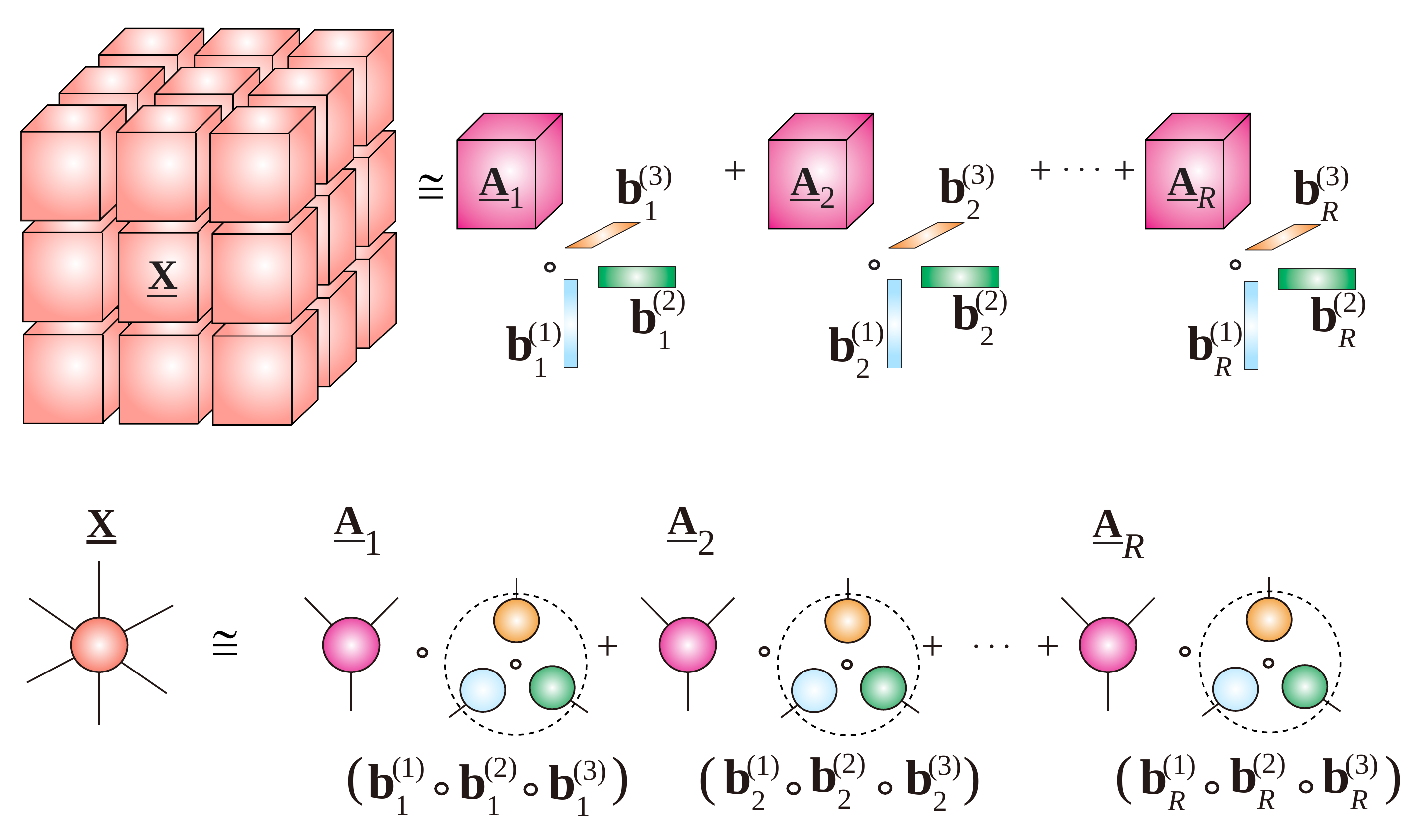}
\caption{Block term  decomposition (BTD)  of a 6th-order block tensor, to yield
$\underline \bX =\sum_{r=1}^R \underline \bA_r \circ \left(\bb_r^{(1)} \circ \bb_r^{(2)} \circ \bb_r^{(3)}\right)$ (top panel),
for more detail see \cite{Lath-BCM12,Sorber2012b}. BTD in the tensor network notation
(bottom panel). Therefore, the 6th-order tensor $\underline \bX$ is approximately represented as  a sum of $R$ terms, each of which is an outer product of a 3rd-order tensor, $\underline \bA_r$, and another a 3rd-order, rank-1 tensor, $\bb_r^{(1)} \circ \bb_r^{(2)} \circ \bb_r^{(3)}$ (in dashed circle), which itself is an outer product of three vectors.}
\label{Fig:BTDS}
\end{figure}
\begin{figure}
\begin{center}
\includegraphics[width=11.5cm]{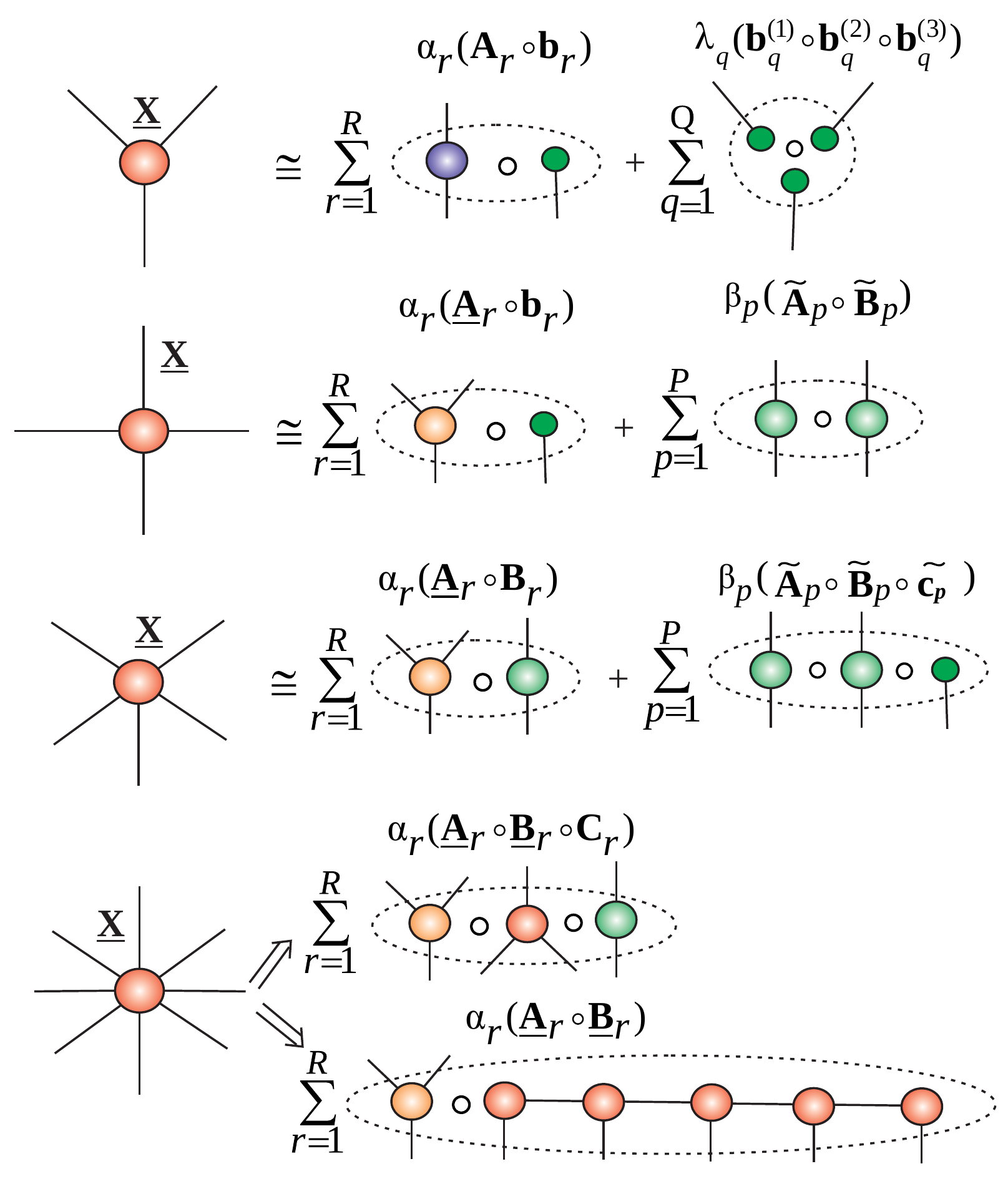}
\end{center}
\caption{Conceptual model of the HOPTA generalized tensor network, illustrated
for  data tensors of different orders. For simplicity, we use the standard outer (tensor) products, but  conceptually  nonlinear outer products (see Eq.~(\ref{eq:nouter})  and other tensor product operators (Kronecker, Hadamard) can also be employed.  Each component (core tensor), $\underline \bA_r$, $\underline \bB_r$  and/or $\underline \bC_r$, can be further hierarchically decomposed using suitable  outer products, so that the
 HOPTA models can be applied to very high-order tensors.}
\label{Fig:HOPTA}
\end{figure}

Figure~\ref{Fig:BTDS} illustrates  such a BTD model for a 6th-order tensor, where the links between the components
are expressed via outer products,  while  Figure  \ref{Fig:HOPTA} shows a  more flexible  Hierarchical Outer Product Tensor Approximation (HOPTA) model suitable for very high-order tensors.

Observe that the fundamental operator in the  HOPTA generalized tensor networks is outer (tensor) product,  which for two tensors $\underline \bA \in \Real^{I_1  \times \cdots \times I_N}$ and $\underline \bB \in \Real^{J_1 \times  \cdots \times J_M}$,  of arbitrary orders $N$ and $M$, is defined as  an $(N+M)$th-order tensor $\underline \bC =\underline \bA \circ \underline \bB \in \Real^{I_1 \times \cdots \times I_N \times J_1 \times \cdots \times J_M}$, with entries $c_{\; i_1, \ldots, i_N, \,j_1, \ldots, j_M} = a_{i_1, \ldots, i_N} \; b_{j_1, \ldots, j_M}$.
 This standard outer product of two tensors can be generalized to a nonlinear outer
product as follows
\be
\left(\underline \bA \circ_f \underline \bB\right)_{i_1,\ldots,i_N,j_1,\ldots,J_M} = f\left(a_{i_1,\ldots, i_N}, b_{j_1,\ldots,j_M}\right),
\label{eq:nouter}
\ee
where $f( \cdot, \cdot)$ is a  suitably designed nonlinear function with associative and commutative properties. In a similar way, we can define other \hbox{nonlinear} tensor products, for example, Hadamard, Kronecker  or Khatri--Rao products and employ them in generalized nonlinear tensor networks.
The advantage of the HOPTA model over other TN models is its flexibility and the ability to model more complex data structures by approximating  very high-order tensors through a relatively small number of low-order cores.

The BTD,  and KTD models can be expressed mathematically, for example, in simple nested (hierarchical) forms, given by
\be
\label{OPTD-model}
\textnormal{BTD}: \underline \bX &\cong& \sum_{r=1}^R (\underline \bA_r \circ \underline \bB_r),\\
\textnormal{KTD}: \underline {\tilde\bX} &\cong& \sum_{r=1}^R (\underline \bA_r \otimes \underline \bB_r),
\label{Phan-model2}
\ee
where, e.g., for BTD, each  factor tensor can be represented recursively as
$\underline \bA_r \cong  \sum_{r_1=1}^{R_1} (\underline \bA^{(1)}_{r_1} \circ \underline \bB^{(1)}_{r_1})$  or $\underline \bB_r \cong  \sum_{r_2=1}^{R_2} \underline \bA^{(2)}_{r_2} \circ \underline \bB^{(2)}_{r_2}$.

Note  that the $2N$th-order subtensors, $\underline \bA_r \circ \underline \bB_r$ and $\underline \bA_r \otimes \underline \bB_r$,  have the same elements,  just arranged differently. For example,  if $\underline \bX = \underline \bA \circ \underline \bB$ and $\underline \bX' = \underline \bA \otimes \underline \bB$, where $\underline \bA  \in \Real^{J_1 \times J_2 \times \cdots \times J_N}$ and  $\underline \bB  \in \Real^{K_1 \times K_2 \times \cdots \times K_N}$, then \\ $x_{j_1,j_2, \ldots,j_N,k_1,k_2, \ldots,k_N} = x'_{k_1+K_1(j_1-1), \ldots, k_N+ K_N(j_N-1)}$.

The definition of the tensor Kronecker product in the KTD model assumes that both core tensors, $\underline \bA_r$
and $\underline \bB_r$, have the same order. 
This is not a limitation, given that vectors and matrices can also be treated as tensors, e.g, a matrix of dimension $I \times J$
  as is also a 3rd-order tensor of dimension $I \times J \times 1$.
In fact, from the BTD/KTD models,  many existing and new TDs/TNs can be derived by
changing the structure and orders of factor tensors, $\underline \bA_r$ and $\underline \bB_r$.
 For example:
\begin{itemize}
\item If $\underline \bA_r$ are rank-1 tensors of size $I_1 \times I_2 \times \cdots \times I_N$, and
$\underline \bB_r$  are scalars, \, $\forall r$, then (\ref{Phan-model2}) represents the rank-$R$ CP decomposition;
\item If $\underline \bA_r$ are rank-$L_r$ tensors of size $I_1 \times I_2 \times \cdots \times I_R \times 1 \times \cdots \times 1$,  in the Kruskal (CP) format, and
$\underline \bB_r$ are rank-1  tensors of size $1 \times  \cdots \times 1 \times I_{R+1}\times \cdots \times I_N$,  \, $\forall r$,  then (\ref{Phan-model2}) expresses the rank-($L_r \circ 1$) BTD;
\item If $\underline \bA_r$ and $\underline \bB_r$ are expressed by KTDs, we arrive at the Nested Kronecker Tensor Decomposition (NKTD), a special case of which is the Tensor Train (TT) decomposition.
    Therefore, the  BTD model in (\ref{Phan-model2})
    can also be used for recursive TT-decompositions. 

\end{itemize}

The generalized tensor network approach caters for a large variety of  tensor decomposition models,  which
  may  find  applications in scientific computing, signal processing or deep learning
  (see, eg., \cite{LathauwerTBSS,CichockiSISA,Cichocki2014optim,Phan_ICASSP15,Shashua_GTD2016}).

In this monograph, we will  mostly focus on the more established  Tucker and  TT decompositions (and some of their extensions),
due to their conceptual simplicity, availability of stable and efficient algorithms for their computation and the  possibility to
naturally extend these models to more complex tensor networks.
In other words, the Tucker and  TT models are considered here as  simplest prototypes, which can then serve as building blocks
for more sophisticated tensor networks.

\chapter{Constrained Tensor Decompositions: From Two-way to Multiway Component Analysis}
\chaptermark{Multiway Component Analysis and Tensor Decompositions}
\label{chap:TDs}

\vspace{0cm}

The component analysis (CA) framework usually refers to the application of constrained matrix  factorization techniques to observed mixed signals
in order to extract  components with specific properties  and/or estimate the mixing matrix
\cite{Cichock_Amari2003,NMF-book,Comon-Jutten2010,Torre2012,HyvICA2013}.
%
In the machine learning practice, to aid the well-posedness and uniqueness of the problem,  component analysis methods exploit prior knowledge about the statistics and diversities of latent variables (hidden sources) within the data.
Here, by the diversities, we refer to different characteristics, features or morphology of  latent variables which allow us to extract the
desired components or features, for example, sparse or  statistically independent components.

\section{Constrained Low-Rank Matrix Factorizations}
\label{sect:2wayCA}

Two-way Component Analysis (2-way CA), in its simplest form,
can be formulated as a  constrained  matrix factorization of typically low-rank, in the form
\begin{equation}
 \label{2CA}
 \bX  = \bA \mbi \Lambda \bB^{\text{T}} + \bE = \sum_{r=1}^R \lambda_r \; \ba_r \circ \bb_r +\bE=\sum_{r=1}^R \lambda_r \, \ba_r \, \bb^{\text{T}}_r +\bE,
\end{equation}
where $\mbi \Lambda =\mbox{diag}(\lambda_1,\ldots,\lambda_R)$ is an optional   diagonal scaling   matrix. The potential constraints imposed on the factor matrices, $\bA$ and/or $\bB$, include orthogonality, sparsity, statistical independence, nonnegativity or smoothness.
In the bilinear 2-way CA  in (\ref{2CA}), $\bX \in \Real^{I \times J}$ is a known matrix of observed data, $\bE  \in \Real^{I \times J}$
represents residuals or noise,    $\bA = [ \ba_1, \ba_2, \ldots, \ba_R] \in \Real^{I \times R}$ is the unknown  mixing matrix with $R$ basis vectors $\ba_{r} \in  \Real^I $,
 and depending on application,  $\bB=[\bb_1,\bb_2,$ $\ldots,\bb_R]$ $ \in \Real^{J  \times R}$,  is the matrix of unknown components, factors, latent variables, or hidden sources, represented by vectors $\bb_r \in \Real^{J}$ (see Figure \ref{Fig:2CA_CPD}).

It should be noted that 2-way CA has an inherent symmetry. Indeed, Eq.~(\ref{2CA}) could also be written as $\bX^{\text{T}} \approx \bB \bA^{\text{T}}$,  thus interchanging the roles of sources and mixing process.

Algorithmic approaches to 2-way (matrix) component analysis are well established, and include  Principal Component Analysis (PCA), Robust PCA (RPCA), Independent Component Analysis (ICA), Nonnegative Matrix Factorization (NMF),  Sparse Component Analysis (SCA) and Smooth Component Analysis (SmCA)
\cite{BachICA2003,NMF-book,bruckstein2009sparse,Comon-Jutten2010,Hyv3way2015,Yokota-SmCA}. 
These techniques have become standard tools in blind source separation (BSS), feature extraction, and  classification paradigms. The columns of the matrix $\bB$, which represent different latent components, are then determined by specific chosen constraints and should be, for example, (i)~as statistically mutually independent as possible for ICA; (ii)~as sparse as possible for SCA; (iii)~as smooth as possible for SmCA;
(iv)~take only nonnegative values for NMF. 

Singular value decomposition (SVD) of the data matrix $\bX \in \Real^{I \times J}$ is a special, very important, case of the factorization in Eq.~(\ref{2CA}),
and is given by
\begin{equation}
\bX= \bU \bS \bV^{\text{T}}  =\sum_{r=1}^R \sigma_r \; \bu_r \circ \bv_r =\sum_{r=1}^R\sigma_r \; \bu_r \bv_r^{\text{T}},
\label{SVD1}
\end{equation}
where $\bU \in \Real^{I \times R} $ and $\bV  \in \Real^{J \times R}$ are column-wise orthogonal matrices and $\bS \in \Real^{R \times R}$ is a diagonal matrix  containing only nonnegative singular values $\sigma_r$ in a monotonically non-increasing order.

According to the well known Eckart--Young theorem, the truncated SVD provides the optimal, in the least-squares (LS) sense, low-rank matrix approximation{\footnote{\cite{Mirsky} has generalized this optimality to arbitrary unitarily invariant norms.}}. The SVD, therefore, forms the backbone of low-rank matrix approximations (and consequently low-rank tensor approximations).

Another virtue of component analysis comes from  the ability 
to perform simultaneous matrix factorizations
\be
\bX_k \approx \bA_k \bB_k^{\text{T}}, \qquad (k=1,2,\ldots,K),
\label{GBSS1}
\ee
 on several data matrices, $\bX_k$, which represent linked datasets,  subject to various constraints imposed on linked (interrelated) component (factor) matrices.
In the case of orthogonality or statistical  independence constraints, the problem in (\ref{GBSS1}) can be related to models of group PCA/ICA through suitable pre-processing, dimensionality reduction and post-processing procedures \cite{esposito2005independent,LinkedICA,Cichocki-SICE,HyvGPCA2014,Zhou-PIEEE}. 
The terms ``group component analysis'' and ``joint multi-block data analysis'' are used interchangeably to refer to methods which aim to identify links (correlations, similarities) between  hidden components in  data.
In other words, \emph{the objective of group component analysis is to analyze the correlation, variability, and consistency of  the latent components across  multi-block datasets}.
The field of 2-way CA is  maturing and has generated efficient algorithms for 2-way component analysis, especially for sparse/functional PCA/SVD, ICA, NMF and SCA \cite{BachICA2003,Cichock_Amari2003,Comon-Jutten2010,Zhou-Cichocki-SP,HyvICA2013}.

The rapidly emerging field of tensor decompositions is the next important step which naturally generalizes 2-way CA/BSS models and algorithms. Tensors, by virtue of multilinear algebra, offer enhanced flexibility in CA, in the sense that not all components need to be  statistically independent, and can be instead  smooth, sparse, and/or non-negative  (e.g., spectral components). Furthermore, additional constraints can be used to reflect  physical properties and/or diversities of spatial distributions, spectral and temporal patterns. We proceed to show how constrained matrix factorizations or 2-way CA models can be extended to multilinear models using  tensor decompositions, such as the Canonical Polyadic (CP)   and the Tucker decompositions, as illustrated in Figures
 \ref{Fig:CPD1}, \ref{Fig:2CA_CPD} and \ref{Fig:Tucker}.

\section{The CP Format}
\label{sect:CPD}

\begin{figure}
(a) {\small {Standard block diagram for CP decomposition of a 3rd-order tensor}}
\vspace{0.3cm}
\begin{center}
\includegraphics[width=8.64cm]{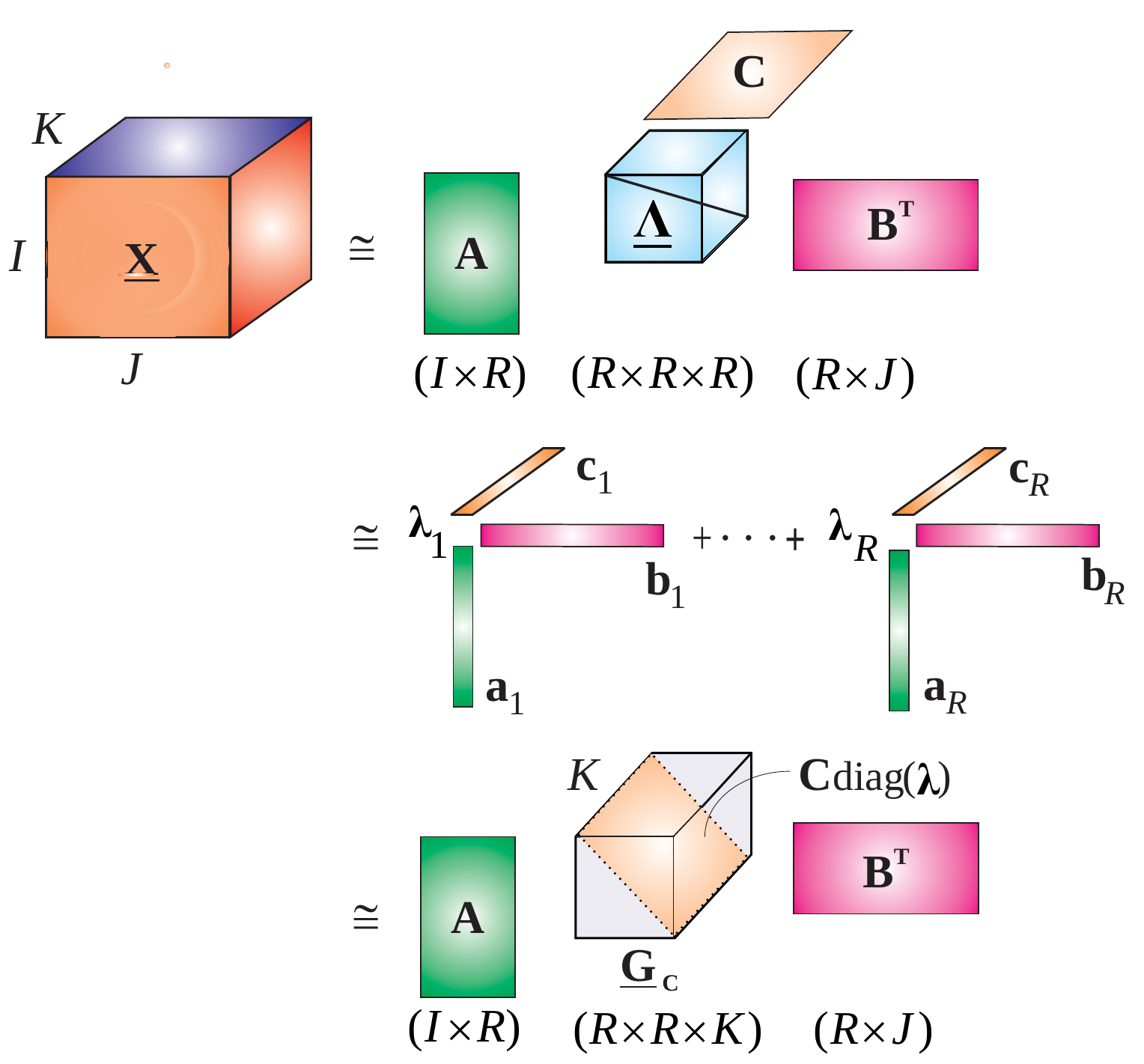}
 \end{center}
\vspace{0.2cm}
(b)  {\small {CP decomposition for a 4th-order tensor in the tensor network notation}}
\begin{center}
\includegraphics[width=6.84cm]{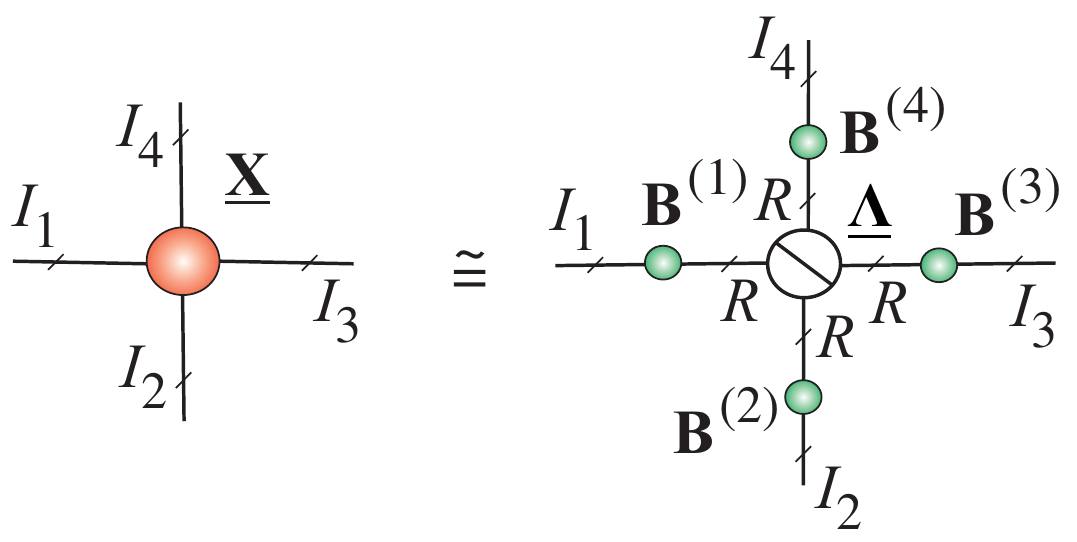}
\end{center}
\caption{Representations of the CP decomposition. The objective of the CP decomposition  is to estimate the factor matrices $\bB^{(n)} \in \Real^{I_n \times R}$ and  scaling coefficients $\{\lambda_1,\lambda_1,\ldots,\lambda_R\}$. (a) The CP decomposition of a 3rd-order tensor in the form, $\underline \bX \cong \underline {\mbi \Lambda} \times_1 \bA \times_2 \bB \times_3 \bC = \sum_{r=1}^R \lambda_r \; \ba_r \circ \bb_r \circ \bc_r = \underline \bG_c \times_1 \bA \times_2 \bB$, with $\underline \bG_c = \underline {\mbi \Lambda} \times_3 \bC$.  (b) The CP decomposition for a 4th-order tensor in the form $\underline \bX \cong  \underline {\mbi \Lambda} \times_1 \bB^{(1)} \times_2 \bB^{(2)} \times_3 \bB^{(3)} \times_4 \bB^{(4)} =\sum_{r=1}^R \lambda_r \; \bb^{(1)}_r \circ \bb^{(2)}_r \circ \bb^{(3)}_r \circ \bb^{(4)}_r$.} 
\label{Fig:CPD1}
\end{figure}

\begin{figure}[ht]
\centering
\includegraphics[width=11.75cm]{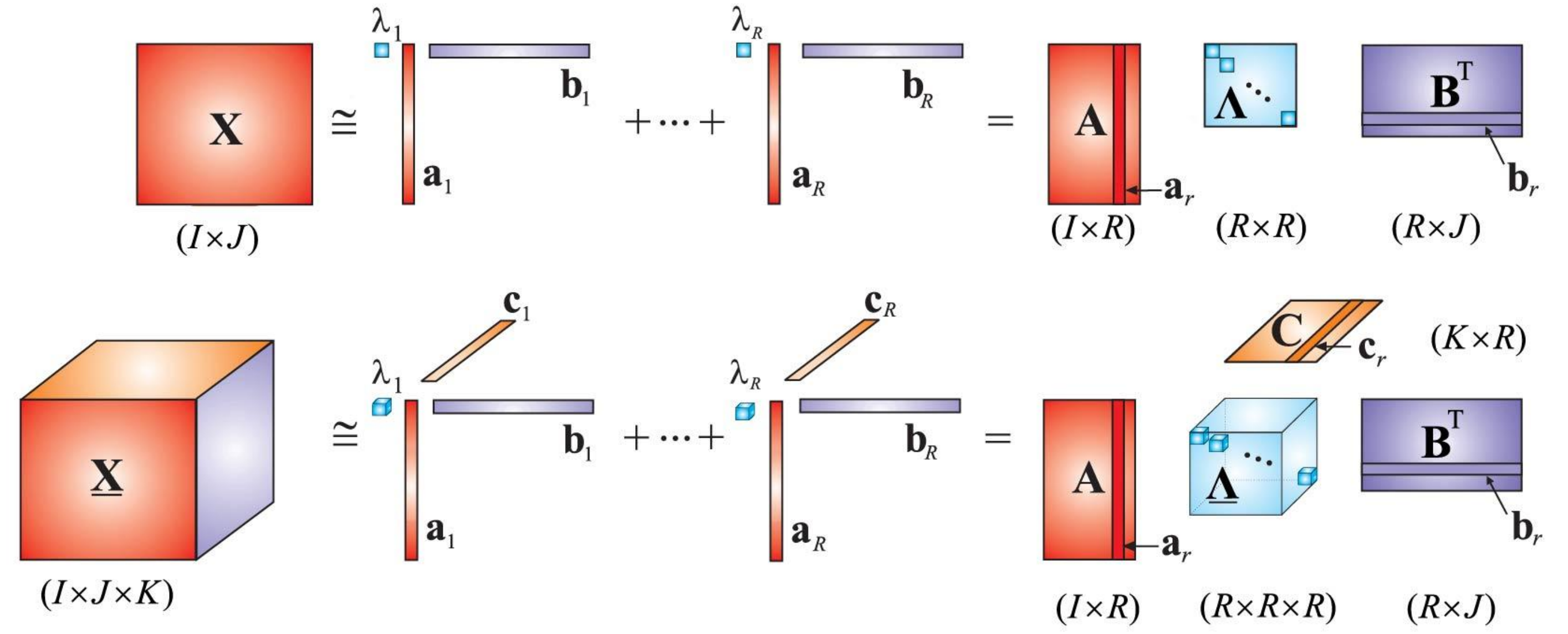}
\caption{Analogy between a low-rank matrix factorization, $\bX \cong \bA \mbi \Lambda \bB^{\text{T}} =\sum_{r=1}^R \lambda_r \; \ba_r \circ \bb_r$ (top), and a simple low-rank tensor factorization (CP decomposition),  $\underline \bX \cong \underline {\mbi \Lambda} \times_1 \bA  \times_2 \bB \times_3 \bC =\sum_{r=1}^R \lambda_r \; \ba_r \circ \bb_r \circ \bc_r$ (bottom).}
\label{Fig:2CA_CPD}
\end{figure}

The CP decomposition (also called the CANDECOMP, PARAFAC, or Canonical Polyadic decomposition) decomposes an $N$th-order tensor, $\underline \bX \in \Real^{I_1 \times I_2 \times  \cdots \times I_N}$, into a linear combination of terms, $\bb^{(1)}_r \circ  \bb^{(2)}_r  \circ \cdots \circ \bb^{(N)}_r$, which are rank-1 tensors,
and is given by \cite{Hitchcock1927,Harshman,PARAFAC1970Carroll}
\begin{equation}
\label{CPDN}
\begin{aligned}
 \underline \bX & \cong  \sum_{r=1}^R  \lambda_r \; \bb^{(1)}_r \circ  \bb^{(2)}_r  \circ \cdots \circ \bb^{(N)}_r \\
 &= \underline {\mbi \Lambda} \times_1 \bB^{(1)} \times_2 \bB^{(2)} \cdots
\times_N \bB^{(N)}\\
 &= \llbracket  \underline {\mbi \Lambda}; \matn[1]{B},  \matn[2]{B}, \ldots, \matn[N]{B}\rrbracket,
  \end{aligned}
\end{equation}
where $\lambda_r$  are   non-zero entries of the diagonal core tensor $\underline {\mbi {\Lambda}} \in \Real^{R \times R \times \cdots \times R}$  and
  $\bB^{(n)} = [\bb_1^{(n)},\bb_2^{(n)},\ldots,\bb_R^{(n)}] \in \Real^{I_n \times R}$ are factor matrices (see Figure \ref{Fig:CPD1}  and Figure \ref{Fig:2CA_CPD}).

Via the Khatri--Rao products (see Table \ref{table_notation2}), the CP decomposition can be equivalently expressed in a matrix/vector form as
\be
\label{CPD-KR2}
\bX_{(n)} &\cong&  \bB^{(n)} \mbi \Lambda (\bB^{(N)} \odot \cdots \odot \bB^{(n+1)} \odot \bB^{(n-1)} \odot \cdots \odot\bB^{(1)})^{\text{T}} \\
&=& \bB^{(n)} \mbi \Lambda (\bB^{(1)} \odot_L \cdots \odot_L \bB^{(n-1)} \odot_L \bB^{(n+1)} \odot_L \cdots \odot_L \bB^{(N)})^{\text{T}}  \notag
\ee
and
\be
\text{vec}(\underline \bX) &\cong& [\bB^{(N)} \odot \bB^{(N-1)} \odot \cdots \odot \bB^{(1)}] \; \mbi \lambda \\
 &\cong & [\bB^{(1)} \odot_L \bB^{(2)} \odot_L \cdots \odot_L \bB^{(N)}] \; \mbi \lambda, \notag
\ee
where  $\mbi \lambda =[\lambda_1,\lambda_2, \ldots,\lambda_R]^{\text{T}}$ and $\mbi \Lambda = \diag(\lambda_1, \ldots, \lambda_R)$ is a diagonal matrix.

The rank of a tensor $\underline \bX$ is  defined as the smallest $R$ for which the CP decomposition in (\ref{CPDN}) holds  exactly.\\

\noindent {\bf Algorithms to compute CP decomposition.} In real world applications, the signals of interest are corrupted by noise, so that the CP decomposition is rarely exact and has to be  estimated by  minimizing a suitable  cost function. Such cost functions are typically of the Least-Squares (LS) type, in the form of the Frobenius norm
\begin{equation}
J_2(\matn[1]{B},  \matn[2]{B}, \ldots, \matn[N]{B})=\|\underline \bX - \llbracket  \underline {\mbi \Lambda}; \matn[1]{B},  \matn[2]{B}, \ldots, \matn[N]{B}\rrbracket \|^2_F,
\end{equation}
 or Least Absolute Error (LAE) criteria \cite{Vorobyov}
\begin{equation}
J_1(\matn[1]{B},  \matn[2]{B}, \ldots, \matn[N]{B})=\|\underline \bX - \llbracket  \underline {\mbi \Lambda}; \matn[1]{B},  \matn[2]{B}, \ldots, \matn[N]{B}\rrbracket \|_1.
\end{equation}

The Alternating Least Squares (ALS) based algorithms minimize the cost function iteratively by individually optimizing each component  (factor matrix, $\bB^{(n)}$)), while keeping the other component matrices fixed  \cite{Harshman,Kolda08}.

To illustrate the ALS principle, assume that the diagonal matrix $\mbi \Lambda$ has been absorbed into one of the component matrices; then, by taking advantage of the Khatri--Rao structure in Eq.~(\ref{CPD-KR2}), the component matrices, $\bB^{(n)}$, can be updated sequentially as
\be
\bB^{(n)} \leftarrow \bX_{(n)} \left( \bigodot_{k\neq n} \bB^{(k)} \right) \left( \bigoasterisk_{k\neq n} (\bB^{(k)\;\text{T}} \bB^{(k)}) \right)^{\dagger}.
\label{ALS-CP}
\ee
 The main challenge (or bottleneck) in implementing ALS and Gradient Decent (GD) techniques for CP decomposition  lies  therefore in multiplying a matricized tensor and Khatri--Rao product (of factor matrices)  \cite{phan2013fast,Choi2014dfacto}  and in the computation of the pseudo-inverse of  $(R \times R)$ matrices (for the basic ALS see Algorithm \ref{alg:ALSCP}). \\

\begin{algorithm}[t]
{\small
\caption{\textbf{Basic ALS for the CP decomposition of a 3rd-order tensor}}
\label{alg:ALSCP}
\begin{algorithmic}[1]
\REQUIRE Data tensor $\underline \bX \in \Real^{I \times J \times K}$ and rank $R$
\ENSURE Factor matrices $\bA \in \Real^{I \times R}$, $\bB \in \Real^{J \times R}$, $\bC \in \Real^{K\times R}$,
 and scaling \\vector $\mbi \lambda \in \Real^R$
\STATE Initialize $\bA, \bB, \bC$
\WHILE {not converged or iteration limit is not reached}
    \STATE $\bA \leftarrow \bX_{(1)} (\bC \odot \bB) (\bC^{\text{T}} \bC \* \bB^{\text{T}} \bB  )^{\dagger}$
    \STATE Normalize column vectors of $\bA$ to unit length (by computing  the \\ norm of each column
    vector and dividing each element of a \\ vector by  its norm)
    \STATE $\bB \leftarrow \bX_{(2)} (\bC \odot \bA) (\bC^{\text{T}} \bC \* \bA^{\text{T}} \bA)^{\dagger}$
    \STATE Normalize column vectors of $\bB$ to unit length
    \STATE $\bC \leftarrow \bX_{(3)} (\bB \odot \bA) (\bB^{\text{T}} \bB \*
    \bC^{\text{T}} \bC)^{\dagger}$
    \STATE Normalize column vectors of $\bC$ to unit length, \\ store the  norms in vector $\mbi \lambda$
\ENDWHILE
\RETURN  $\bA, \bB, \bC$ and $\mbi \lambda$.
\end{algorithmic}
}
\end{algorithm}

 The ALS approach is attractive for its simplicity, and  often provides satisfactory performance for well defined problems with high SNRs and well separated and non-collinear components.
For ill-conditioned problems, advanced algorithms are required which typically exploit the rank-1 structure of the terms within CP decomposition to perform efficient computation and storage of the Jacobian and Hessian of the cost function \cite{Phan2012-Hess,Sorber2012b,phan2015tensor}.
Implementation of parallel ALS algorithm over distributed memory  for very large-scale tensors was  proposed in \cite{Choi2014dfacto,karlsson2015parallel}.
\\

\noindent{\bf Multiple random projections, tensor sketching and Giga-Tensor.}
Most of the existing algorithms for the computation of CP decomposition are based on the ALS or GD
approaches, however, these can be too computationally expensive for huge tensors.
Indeed, algorithms for tensor decompositions  have  generally not yet reached
the level of maturity and efficiency of  low-rank matrix factorization (LRMF)
methods. In order to employ efficient LRMF algorithms to tensors, we need to either: (i) reshape the tensor at hand into a set of matrices using
traditional   matricizations, (ii) employ reduced randomized unfolding matrices, or (iii) perform suitable  random  multiple projections of a data tensor onto lower-dimensional subspaces. The principles of the approaches (i) and (ii) are self-evident, while the approach (iii) employs a multilinear product of an $N$th-order tensor and $(N-2)$
random vectors, which are either chosen uniformly from a unit sphere or assumed to be i.i.d. Gaussian vectors
\cite{kuleshov2015tensor}.

For example, for a 3rd-order tensor, $\underline \bX \in \Real^{I_1 \times I_2 \times I_3}$,
 we can use the set of random projections, $\bX_{\bar 3} = \underline \bX \; \bar \times_3 \; \mbi \omega_3 \in \Real^{I_1 \times I_2}$,
$\bX_{\bar 2} = \underline \bX  \; \bar \times_2  \; \mbi \omega_2 \in \Real^{I_1 \times I_3}$ and $\bX_{\bar 1} = \underline \bX \; \bar \times_1 \; \mbi \omega_1 \in \Real^{I_2 \times I_3}$,
 where the vectors $\mbi \omega_n \in \Real^{I_n}$, $\, n=1,2,3$, are suitably chosen  random vectors.
%
%
Note that  random projections   in  such a case are non-typical -- instead of using projections for dimensionality
reduction, they are used to reduce a tensor of any order to matrices and  consequently transform  the CP decomposition problem
 to constrained matrix factorizations problem, which can be solved
  via simultaneous (joint) matrix diagonalization \cite{deLathauwer-JAD,chabriel2014JAD}.
It was shown that even a small number of random projections, such as ${\cal{O}}(\log R)$  is sufficient to preserve the spectral information in a tensor.
This mitigates the problem of the dependence
on the eigen-gap{\footnote{In linear algebra, the eigen-gap of a linear operator is the difference between two successive eigenvalues, where the eigenvalues are sorted in an ascending order.}} that plagued earlier tensor-to-matrix reductions.
Although a uniform random sampling may  experience problems for tensors with spiky elements, it often
outperforms the  standard CP-ALS  decomposition algorithms. 

Alternative  algorithms for the CP decomposition of huge-scale tensors include tensor sketching --
a  random mapping technique,  which exploits  kernel methods  and regression
\cite{pham2013fast,SmolaNIPS15}, and the class of distributed algorithms such as the DFacTo  \cite{Choi2014dfacto} and the GigaTensor which is based on Hadoop / MapReduce paradigm \cite{Giga-tensor}.\\

\noindent{\bf Constraints.} Under rather mild conditions, the CP decomposition is generally unique by itself  \cite{Kruskal1977,Sidiropoulos2000}.  It does not require additional constraints on the factor matrices to achieve uniqueness, which makes it a powerful and useful tool for tensor factroization. Of course, if the components in one or more modes are known to possess some properties,  e.g.,  they are known to be nonnegative, orthogonal, statistically independent or sparse, such prior knowledge may  be incorporated into the algorithms to compute CPD and at the same time  relax uniqueness conditions. More importantly, such constraints may enhance the accuracy and stability of CP decomposition algorithms and also facilitate better physical interpretability of the extracted components \cite{Nikos04,dhillon2009fast,sorensen,kim2014robust,Zhou2012-SPL,Sidiro2015parallel}.\\

\noindent{\bf Applications.} The CP decomposition has already been established as an advanced tool for blind
signal separation in vastly diverse branches of signal processing
and machine learning \cite{Acar2008,Kolda08,Morup11,Anandkumar2014tensor,SmolaNIPS15,Tresp2015learning,Sidiro-Lath2016}.
It is also routinely used in exploratory data analysis, where the rank-1 terms capture
essential properties of dynamically complex datasets, 
while in wireless communication systems, signals transmitted by different users correspond to
rank-1 terms in the case of line-of-sight propagation  and therefore admit  analysis in the CP format.
Another potential application is in harmonic retrieval and direction of arrival problems, where real or complex exponentials have rank-1 structures,
for which the use of CP decomposition is quite  natural \cite{Sidiropoulos00Bro,sid_carat,sorensenVDM}.

\section{The Tucker Tensor Format}
\label{sect:Tucker}

\begin{figure}
(a)  {\small{Standard block diagrams of Tucker (top) and Tucker-CP (bottom) decompositions for a 3rd-order tensor}}\\
\vspace{-0.1cm}
\begin{center}
\includegraphics[width=11.305cm]{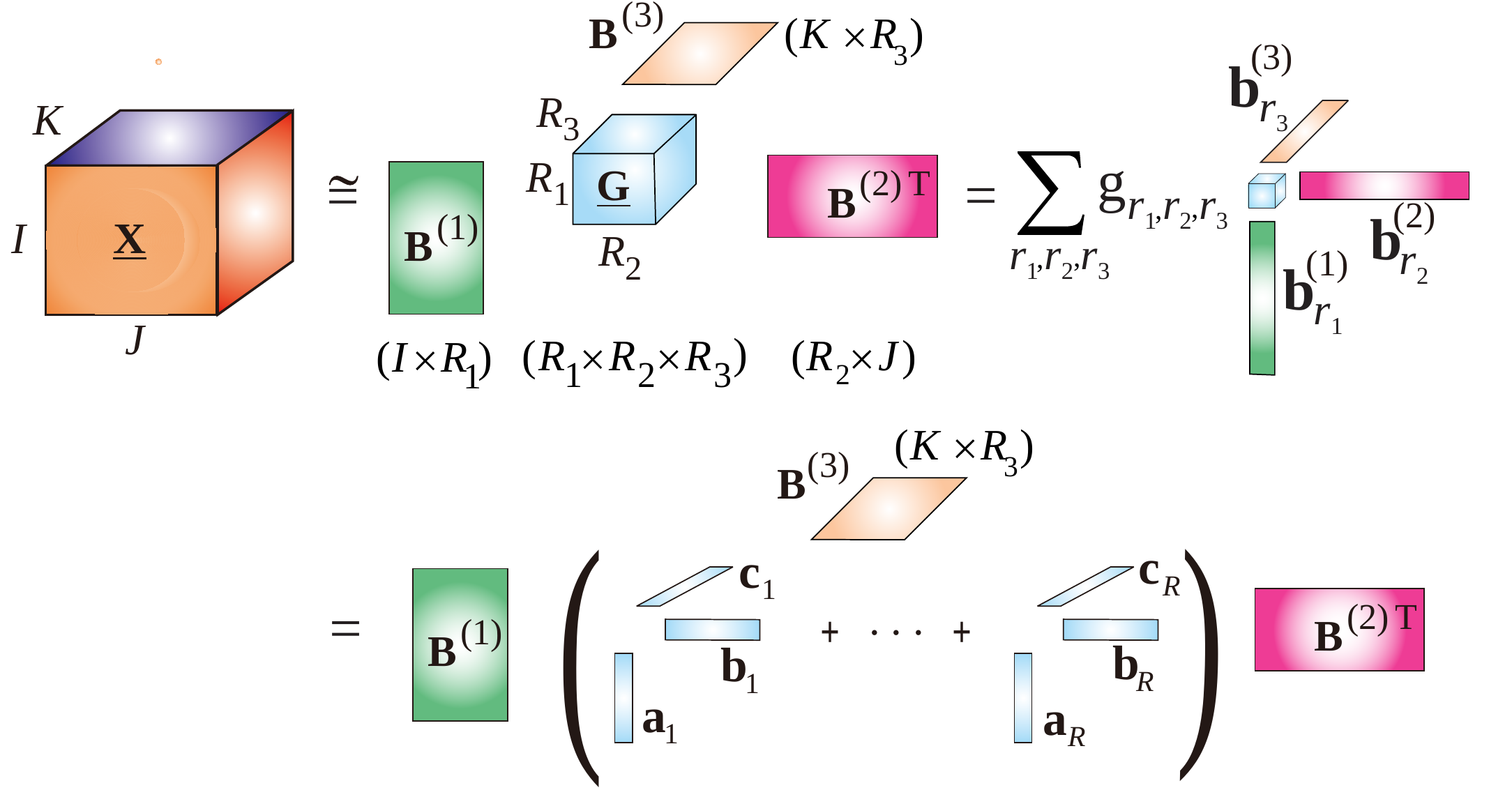}
\end{center}
(b) {\small{The TN diagram for the  Tucker and Tucker/CP decompositions of a 4th-order tensor}}  \\
\vspace{-0.1cm}
\begin{center}
  \includegraphics[width=10.83cm]{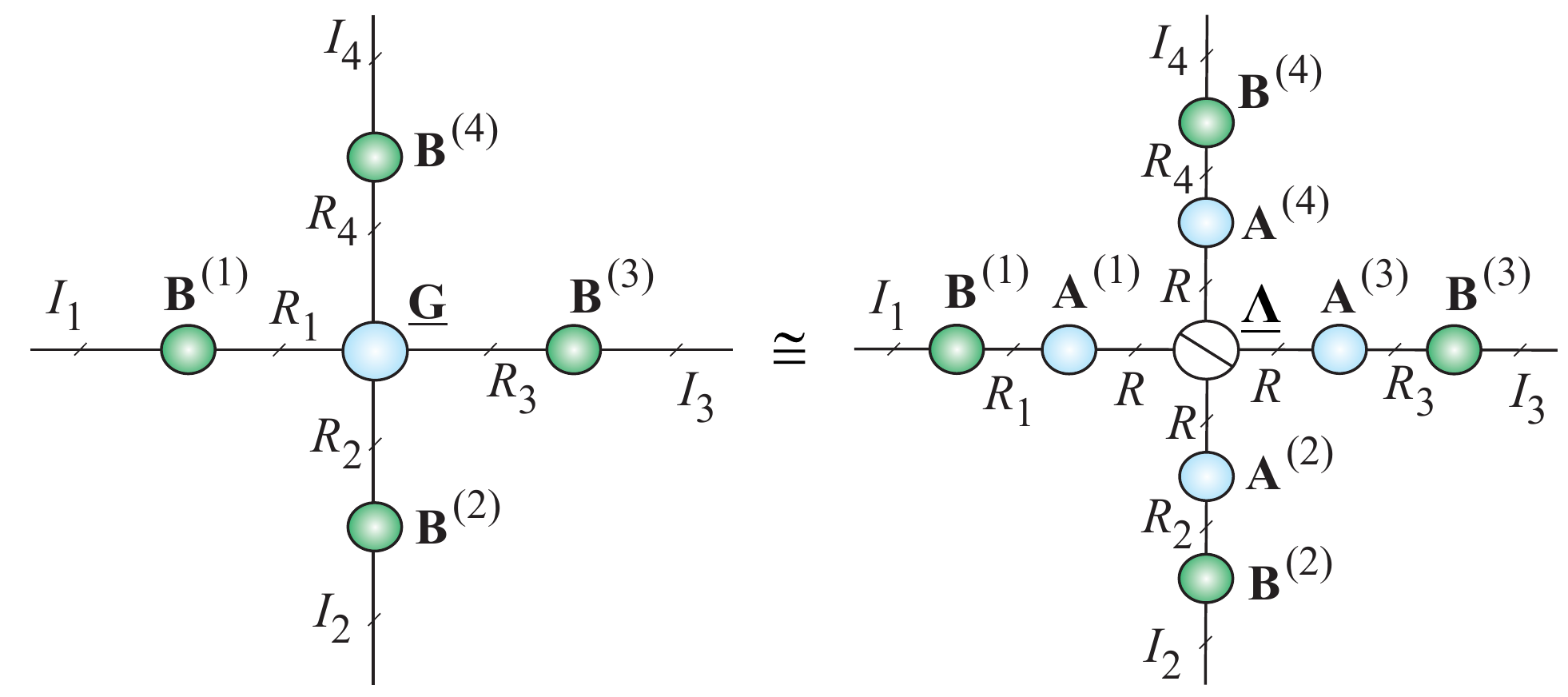}
 \end{center}
\caption{{\small Illustration of the Tucker and Tucker-CP decompositions, where the objective
is to compute the factor matrices, $\bB^{(n)}$, and the core tensor, $\underline \bG$. (a) Tucker decomposition of a 3rd-order tensor, $\underline \bX \cong \underline \bG \times_1 \bB^{(1)} \times_2 \bB^{(2)} \times_3 \bB^{(3)}$.  In some applications,  the core tensor can be further approximately factorized using the CP decomposition as $\underline \bG \cong \sum_{r=1}^R \ba_r \circ \bb_r \circ \bc_r$ (bottom diagram), or alternatively using TT/HT decompositions.  (b) Graphical representation of the Tucker-CP decomposition
for a 4th-order tensor,  $\underline \bX \cong \underline \bG \times_1 \bB^{(1)} \times_2 \bB^{(2)} \times_3 \bB^{(3)} \times_4 \bB^{(4)} =\llbracket \underline \bG; \bB^{(1)},  \bB^{(2)}, \bB^{(3)}, \bB^{(4)} \rrbracket \cong (\underline {\mbi \Lambda} \times_1 \bA^{(1)} \times_2 \bA^{(2)} \times_3 \bA^{(3)} \times_4 \bA^{(4)})  \times_1 \bB^{(1)} \times_2 \bB^{(2)}  \times_3 \bB^{(3)}  \times_4 \bB^{(4)} = \llbracket \underline {\mbi \Lambda}; \; \bB^{(1)} \bA^{(1)}, \; \bB^{(2)} \bA^{(2)}, \; \bB^{(3)} \bA^{(3)}, \; \bB^{(4)} \bA^{(4)} \rrbracket$.}}
\label{Fig:Tucker}
\end{figure}

Compared to the CP decomposition, the Tucker decomposition provides a more general factorization of an $N$th-order tensor into a relatively small size core tensor and factor matrices, and can be expressed  as follows:
\be
  \underline \bX
 & \cong & \sum\limits_{r_1 = 1}^{R_1} {
 {\cdots
 \sum\limits_{r_N = 1}^{R_N}{g_{r_1 r_2 \cdots r_N} \, \left(\; \bb^{(1)}_{r_1} \circ \bb^{(2)}_{r_2} \circ \cdots \circ \bb^{(N)}_{r_N} \right)}}} \notag \\
& = & \underline \bG \times_1 \bB^{(1)} \times_2 \bB^{(2)} \cdots \times_N \bB^{(N)} \notag \\
& = & \llbracket\underline \bG; \bB^{(1)},\bB^{(2)},\ldots,\bB^{(N)}\rrbracket,
 \label{GeneralTDModel}
 \ee
 where $\underline \bX \in \Real^{I_1\times I_2 \times \cdots \times I_N}$ is the given data tensor,  $\underline \bG  \in \Real^{R_1 \times R_2 \times \cdots \times R_N}$  is the core tensor, and $\bB^{(n)} =[\bb_1^{(n)},\bb_2^{(n)},\ldots,\bb_{R_n}^{(n)}] \in \Real^{I_n\times R_n}$ are the mode-$n$ factor (component)  matrices, $n=1,2,\ldots,N$ (see Figure \ref{Fig:Tucker}). The core tensor (typically,  $R_n << I_n$) models a potentially complex pattern of mutual interaction between the vectors  in  different modes. The model in (\ref{GeneralTDModel}) is  often referred to as the Tucker-$N$ model.

The CP and Tucker decompositions have  long history. For  recent  surveys   and  more detailed
 information we refer to  \cite{Kolda08,Grasedyck-rev,comon2014tensors,Cich-Lath,Sidiro-Lath2016}.

 Using the properties of the Kronecker tensor product, the Tucker-$N$ decomposition in (\ref{GeneralTDModel}) can be expressed in an equivalent  matrix and vector form as
\begin{align}
\bX_{(n)} &\cong \bB^{(n)} \bG_{(n)} (\bB^{(1)} \otimes_L \cdots \otimes_L \bB^{(n-1)} \otimes_L \bB^{(n+1)} \otimes_L \cdots \otimes_L \bB^{(N)})^{\text{T}} \notag \\
&= \bB^{(n)} \bG_{(n)} (\bB^{(N)} \otimes \cdots \otimes \bB^{(n+1)} \otimes \bB^{(n-1)} \otimes \cdots \otimes \bB^{(1)})^{\text{T}},\\
\bX_{<n>} &\cong  (\bB^{(1)} \otimes_L \cdots \otimes_L \bB^{(n)}) \; \bG_{<n>} (\bB^{(n+1)} \otimes_L  \cdots \otimes_L \bB^{(N)})^{\text{T}} \notag \\
&= (\bB^{(n)} \otimes \cdots \otimes \bB^{(1)}) \; \bG_{<n>} (\bB^{(N)} \otimes \bB^{(N-1)} \otimes \cdots \otimes \bB^{(n+1)})^{\text{T}}, \notag\\ \\
\text{vec}(\underline \bX) &\cong [\bB^{(1)} \otimes_L \bB^{(2)} \otimes_L\cdots \otimes_L \bB^{(N)}] \; \text{vec}(\underline \bG)\notag \\
&= [\bB^{(N)} \otimes \bB^{(N-1)} \otimes \cdots \otimes \bB^{(1)}] \; \text{vec}(\underline \bG),
\label{Tucker-kron2}
\end{align}
 where the multi-indices are ordered in a reverse lexicographic order (little-endian).

Table \ref{tab:CP_Tucker3}  and Table \ref{tab:CP_TuckerN} summarize fundamental
 mathematical representations of CP and Tucker decompositions for 3rd-order and $N$th-order tensors.

The Tucker decomposition is said to be in an \emph{independent Tucker format} if all the factor matrices,  $\bB^{(n)}$, are full column rank, while a Tucker format is termed an \emph{orthonormal format},  if in addition,  all the factor matrices, $\bB^{(n)}=\bU^{(n)}$, are  orthogonal. The standard Tucker model often has orthogonal factor matrices.\\

\noindent{\bf Multilinear rank.} The multilinear rank of an $N$th-order tensor $\underline \bX \in \Real^{I_1 \times I_2 \times \cdots \times I_N}$ corresponds to the $N$-tuple
 $(R_1,R_2,\ldots,R_N)$  consisting of the dimensions of the different subspaces.
If the Tucker decomposition (\ref{GeneralTDModel}) holds exactly it is mathematically defined as
\be
\mbox{rank}_{ML} (\underline \bX)= \{\mbox{rank} (\bX_{(1)}), \mbox{rank} (\bX_{(2)}), \ldots, \mbox{rank} (\bX_{(N)})\},
\ee
with $\bX_{(n)} \in \Real^{I_n \times I_1 \cdots I_{n-1} I_{n+1} \cdots I_N}$ for $n=1,2,\ldots,N$. 
Rank of the Tucker decompositions can be determined using information criteria \cite{Yokota_MDL}, or through the number of dominant eigenvalues when an approximation accuracy of the decomposition or a noise level is given (see Algorithm \ref{alg:Tucker2}).

 The independent Tucker format has the following important properties if the equality in (\ref{GeneralTDModel}) holds exactly
  (see, e.g., \cite{Jiang-Tucker} and references therein):

 \begin{enumerate}

 \item   The tensor (CP) rank  of any tensor, $\underline \bX = \llbracket\underline \bG; \bB^{(1)},\bB^{(2)},\ldots,\bB^{(N)}\rrbracket \in \Real^{I_1\times I_2 \times \cdots \times I_N}$, and the rank of its core tensor,  $\underline \bG  \in \Real^{R_1 \times R_2 \times \cdots \times R_N}$,  are exactly the same, i.e.,
     \be
     \mbox{rank}_{CP} (\underline{\bX}) = \mbox{rank}_{CP} (\underline{\bG}).
     \ee

     \item If a tensor, $\underline \bX \in \Real^{I_1 \times I_2 \times \cdots \times I_N}
     =\llbracket\underline \bG; \bB^{(1)},\bB^{(2)},\ldots,\bB^{(N)}\rrbracket $, admits an independent Tucker format with
     multilinear rank $\{R_1,R_2, \ldots, R_N\}$, then
     \be
     R_n \leq \prod_{p\neq n}^N R_p \quad \forall n.
     \ee
     Moreover, without loss of generality, under the assumption $R_1 \leq R_2 \leq \cdots \leq R_N$, we have
     \be
      R_1 \leq \mbox{rank}_{CP}(\underline \bX) \leq R_2 R_3 \cdots R_N.
      \ee

    \item If a data  tensor is symmetric and admits an independent Tucker format, $\underline \bX = \llbracket \underline \bG; \bB,\bB,\ldots,\bB \rrbracket \in \Real^{I\times I \times \cdots \times I} $, then its core tensor,
     $\underline \bG \in \Real^{R \times R \times  \cdots \times R}$, is also symmetric, with $\mbox{rank}_{CP} (\underline \bX)=\mbox{rank}_{CP} (\underline \bG$).

  \item  For the orthonormal Tucker format, that is, $\underline \bX
     =\llbracket\underline \bG; \bU^{(1)},\bU^{(2)},\ldots,\bU^{(N)}\rrbracket \in \Real^{I_1\times I_2 \times \cdots \times I_N} $,  with $\bU^{(n) \text{T}} \,\bU^{(n)} = \bI, \, \forall n$, the  Frobenius norms and the Schatten $p$-norms{\footnote{The Schatten $p$-norm of an $N$th-order tensor $\underline \bX$ is defined as the average of the Schatten norms of mode-$n$ unfoldings, i.e.,
     $\|\underline \bX\|_{{\cal{S}}_p}= (1/N) \sum_{n=1}^N \|\bX_{(n)}\|_{{\cal{S}}_p}$  and $\|\bX\|_{{\cal{S}}_p} = (\sum_r \sigma_r^p)^{1/p}$, where $\sigma_r$ is the $r$th singular value of the matrix $\bX$. For $p=1$, the Schatten norm of a matrix $\bX$ is called the nuclear norm or the trace norm,  while for  $p=0$ the Schatten norm is the rank of $\bX$, which can be replaced by the surrogate function $\log \det(\bX \bX^{\text{T}} +\varepsilon \bI),\;\;\varepsilon>0$.}} of the data tensor, $\underline \bX$ ,  and its core tensor, $\underline \bG$, are  equal, i.e.,
     \be
     \|\underline \bX\|_F &=& \|\underline \bG\|_F, \notag \\
     \|\underline \bX\|_{{\cal {S}}_p} &=& \|\underline \bG\|_{{\cal {S}}_p}, \quad  1 \leq p < \infty. \notag
     \ee
Thus, the computation of the Frobenius
norms can be performed with an ${\cal{O}} (R^N)$   complexity $(R=\max\{R_1,\ldots,R_N\})$, instead of the usual order ${\cal{O}} (I^N)$ complexity (typically $R \ll I$).
 \end{enumerate}

\begin{table}
\vspace{-2.5cm}
\caption{Different forms of CP and Tucker representations of a 3rd-order tensor $\underline \bX \in \Real^{I \times J \times K}$, where $\;\mbi \lambda=[\lambda_1,\lambda_2, \ldots, \lambda_R]^{\text{T}}$, and $\;\mbi \Lambda= \mbox{diag}\{\lambda_1,\lambda_2, \ldots, \lambda_R\}$.}
 \centering
 {
 \shadingbox{
   \begin{tabular*}{1.00\textwidth}[t]{@{\extracolsep{\fill}}@{\hspace{2ex}}l@{\hspace{-4em}}l} \hline  &  \\
{ \raisebox{0mm}[0mm][4mm]{CP Decomposition}}& \hspace{-3.2cm} Tucker Decomposition \\ \hline
 &   \\
%
  \multicolumn{2}{c}{Scalar representation}\\ & \\
 \begin{tabular*}{.5\linewidth}[t]{@{\extracolsep{\fill}}cc@{}}
     $ x_{ijk} = \sum\limits_{r=1}^{R} {\lambda_r \,a_{i\,r} \, b_{j\,r} \, c_{k\,r}} $  & \hspace{0.4cm} 
$ \quad x_{ijk} =\sum\limits_{r_1=1}^{R_1} \sum\limits_{r_2=1}^{R_2} \sum\limits_{r_3=1}^{R_3} {g_{r_1 \, r_2 \, r_3} \; a_{i\,r_1} \, b_{j\,r_2} \, c_{k\,r_3}}  $
\end{tabular*}
                   \\  &   \\ \hline \\
 \multicolumn{2}{c}{Tensor representation, outer products}\\
 & \\
\begin{tabular*}{.9\linewidth}[t]{@{\extracolsep{\fill}}cc@{}}
 $\;\;\underline \bX  = \sum\limits_{r=1}^{R} {\lambda_r \; \ba_{r} \circ \bb_{r} \circ \bc_{r}} $  & \hspace{0.4cm}
$\;\; \quad \underline \bX  = \sum\limits_{r_1=1}^{R_1} \sum\limits_{r_2=1}^{R_2} \sum\limits_{r_3=1}^{R_3} {g_{r_1 \,r_2 \, r_3} \; \ba_{r_1} \circ \bb_{r_2} \circ \bc_{r_3}} $
\end{tabular*}
\\
 &  \\ \hline
  & \\
 \multicolumn{2}{c}{Tensor representation, multilinear products}\\
  & \\
   $ \;\; \underline \bX = \underline {\mbi{\Lambda}} \times_1 \bA \times_2 \bB \times_3 \bC $  &  \hspace{-3.2cm}
 $ \underline \bX = \underline \bG \times_1 \bA \times_2 \bB \times_3 \bC $ \\
 $ \;\; \quad = \llbracket \underline {\mbi{\Lambda}}; \; \bA, \; \bB, \; \bC \rrbracket$  &  \hspace{-3.2cm}
 $ \quad  =  \llbracket \underline \bG; \;  \bA, \; \bB, \; \bC \rrbracket $
 \\  &   \\ \hline
     & \\
 \multicolumn{2}{c}{ Matrix representations}\\
    &  \\
       $\;\;\bX_{(1)} = \bA \; \mbi\Lambda \; (\bB \odot_L \bC)^{\text{T}} $  &  \hspace{-3.2cm}
       $\bX_{(1)} = \bA \;\bG_{(1)} \; (\bB \otimes_L \bC)^{\text{T}} $ \\
      {$\;\;\bX_{(2)} = \bB \;  \mbi\Lambda \; (\bA \odot_L \bC)^{\text{T}} $} &  \hspace{-3.2cm}
      {$\bX_{(2)} = \bB \;\bG_{(2)} \; (\bA \otimes_L  \bC)^{\text{T}} $} \\
       {$\;\;\bX_{(3)} = \bC \;  \mbi\Lambda \; (\bA \odot_L \bB)^{\text{T}} $}
        &  \hspace{-3.2cm} {$\bX_{(3)} =\bC \;\bG_{(3)} \; (\bA \otimes_L  \bB)^{\text{T}} $ } \\
       &  \\ \hline
    & \\
    \multicolumn{2}{c}{Vector representation} \\
                    & \\
                       $\;\;\text{vec}(\underline \bX ) =  (\bA \odot_L \bB \odot_L \bC) \mbi \lambda$
                      &  \hspace{-3.4cm}
                       {$\;\; \text{vec}(\underline \bX ) = (\bA \otimes_L  \bB \otimes_L  \bC) \;
                        \text{vec}(\underline \bG ) $} \\  &  \hspace{0.1cm}   \\ \hline
                      & \\
       & \\
  \multicolumn{2}{c}{Matrix slices $\bX_k = \bX(:,:,k)$} \\
    &  \\
    $\bX_k =  \bA \, {\rm diag}(\lambda_1 \, c_{k,1},\ldots,\lambda_R \, c_{k,R}) \, \bB^{\text{T}}$  &
    \hspace{-3.4cm}
      $\;\; \bX_k= \bA\, \left(\sum\limits_{r_3=1}^{R_3} c_{k r_3} \bG(:,:,r_3)\right) \, \bB^{\text{T}} $ \\
      &  \\[-6pt] \hline
    \end{tabular*}
    }}
\label{tab:CP_Tucker3}
\end{table}


\begin{table} 
\caption{Different forms of CP and Tucker representations
of an $N$th-order tensor $\underline \bX \in \Real^{I_1\times I_2\times\cdots\times I_N}$.}
 \centering
 {\small\tabcolsep 3pt
 \shadingbox{
\begin{tabular}{ll}
\hline\\[-8pt]
\multicolumn{1}{c}{CP}   &    \multicolumn{1}{c}{Tucker}\\
\hline
\multicolumn{2}{c}{Scalar product} \\ [1.2pc]
	$\displaystyle
	x_{i_1,\ldots,i_N} = \sum_{r=1}^R \lambda_r \; b^{(1)}_{i_1,r} \cdots b^{(N)}_{i_N,r}
	$
	&
	\hspace{-0.5cm}$\displaystyle
	x_{i_1,\ldots,i_N} = \sum_{r_1=1}^{R_1} \cdots\sum_{r_N=1}^{R_N} g_{r_1,\ldots,r_N} \; b^{(1)}_{i_1,r_1} \cdots b^{(N)}_{i_N,r_N}
	$
\\ [1.2pc]
\hline
\multicolumn{2}{c}{Outer product} \\ [1.2pc]
	$\displaystyle
	\underline \bX = \sum_{r=1}^R \lambda_r \; \bb^{(1)}_r \circ \cdots \circ \bb^{(N)}_r
	$
	&\qquad
	\hspace{-0.4cm} $\displaystyle
	\underline \bX = \sum_{r_1=1}^{R_1} \cdots \sum_{r_N=1}^{R_N} g_{r_1,\ldots,r_N} \; \bb^{(1)}_{r_1} \circ\cdots \circ \bb^{(N)}_{r_N}
	$
\\ [1.2pc]
\hline
\multicolumn{2}{c}{Multilinear product} \\ [1.2pc]
	$\displaystyle
	\underline \bX = \underline{\mbi{\Lambda}} \times_1 \bB^{(1)} \times_2 \bB^{(2)} \cdots \times_N \bB^{(N)}
	$
	&\qquad
	$\displaystyle
	\underline \bX = \underline \bG \times_1 \bB^{(1)} \times_2 \bB^{(2)} \cdots \times_N \bB^{(N)}
	$
\\[1.2pc]
	$\displaystyle
	\underline \bX = \left\llbracket \underline{\mbi{\Lambda}}; \bB^{(1)}, \bB^{(2)}, \ldots, \bB^{(N)} \right\rrbracket
	$
	&\qquad
	$\displaystyle
	\underline \bX = \left\llbracket \underline \bG; \bB^{(1)},\bB^{(2)}, \ldots,\bB^{(N)} \right\rrbracket
	$
\\[1.2pc]
\hline
\multicolumn{2}{c}{Vectorization} \\ [1.2pc]
	$\displaystyle
	\text{vec}(\underline \bX) = \left( \bigodot_{n=N}^1 \bB^{(n)} \right) \mbi\lambda
	$
	&\qquad
	$\displaystyle
	\text{vec}(\underline \bX) = \left( \bigotimes_{n=N}^1 \bB^{(n)} \right) \text{vec}(\underline \bG)
	$
\\ [1.2pc]
\hline
\multicolumn{2}{c}{Matricization} \\ [1.2pc]
	$\displaystyle
	\bX_{(n)} = \bB^{(n)} \mbi\Lambda
\left( \bigodot_{m=N,\; m\neq n}^1\bB^{(m)} \right)^\text{T}
	$
	&\qquad
	$\displaystyle
	\bX_{(n)} = \bB^{(n)} \bG_{(n)}
\left( \bigotimes_{{m=N,\;m\neq n}}^1 \bB^{(m)} \right)^\text{T}
	$
\\ [1.2pc]
	$ \displaystyle
	\bX_{<n>} = ( \bigodot_{m=n}^1 \bB^{(m)} ) \mbi\Lambda
	( \bigodot_{m=N}^{n+1} \bB^{(m)} )^\text{T},
	$
	&
	\hspace{-0.2cm}$ \qquad \displaystyle
	\bX_{<n>} = ( \bigotimes_{m=n}^1 \bB^{(m)} ) \bG_{<n>}
	( \bigotimes_{m=N}^{n+1} \bB^{(m)} )^\text{T}
	$
\\ [1.2pc]
\hline
\multicolumn{2}{c}{Slice representation} \\ [1.2pc]
	$\displaystyle
	\underline \bX (:,:,k_3) = \bB^{(1)} \widetilde{\bD}_{k_3} \bB^{(2)\;\text{T}}
	$
	&
	$\displaystyle
	\underline \bX (:,:,k_3) = \bB^{(1)} \widetilde{\bG}_{k_3} \bB^{(2) \;\text{T}}, \;k_3=\overline{i_3 i_4 \cdots i_N}
	$
\\ [1.2pc]
 $\widetilde{\bD}_{k_3} = \text{diag}(\tilde{d}_{11},\ldots,\tilde{d}_{RR}) \in \Real^{R\times R}$
& \hspace{-0.5cm} with entries
	$\;\tilde{d}_{rr}=\lambda_r b^{(3)}_{i_3,r}\cdots b^{(N)}_{i_N,r}$
\\ [1.2pc]
\multicolumn{2}{l}{
	 $\displaystyle \widetilde{\bG}_{k_3} =
	\sum_{r_3}\cdots\sum_{r_N} b^{(3)}_{i_3,r_3}\cdots b^{(N)}_{i_N,r_N} \bG_{:,:,r_3,\ldots,r_N}$
	is the sum of frontal slices.}
\\ [0.5pc]
\hline
\end{tabular}
}}
\label{tab:CP_TuckerN}
\end{table}

Note that the CP decomposition can be considered as a special case of the Tucker decomposition, whereby the cube core tensor has nonzero elements only on the main diagonal (see Figure \ref{Fig:CPD1}). In contrast to the CP decomposition, the unconstrained Tucker decomposition is not  unique. However, constraints imposed on all factor matrices and/or core tensor can reduce the indeterminacies inherent in CA to only column-wise permutation and scaling, thus yielding a unique core tensor
and factor matrices  \cite{Zhou-Cichocki-MBSS}.

The Tucker-$N$ model, in which $(N-K)$ factor matrices are identity matrices is called the Tucker-$(K,N)$ model. In the simplest scenario, for a 3rd-order tensor $\underline \bX \in \Real^{I \times J\times K} $, the Tucker-(2,3) or simply Tucker-2 model, can be described as{\footnote{For a 3rd-order tensor, the Tucker-2 model is equivalent to the TT model. The case where the factor matrices and the core tensor are non-negative is referred to as the  NTD-2 (Nonnegative Tucker-2 decomposition).}}
\be
\underline \bX \cong \underline \bG \times_1 \bA \times_2 \bB \times_3 \bI =\underline \bG \times_1 \bA \times_2 \bB,
\ee
or in an equivalent  matrix form
\be
\bX_k =\bA \bG_k \bB^{\text{T}}, \qquad (k=1,2,\ldots,K),
\ee
where    $\bX_k =\underline \bX(:,:,k) \in \Real^{I \times J}$ and $\bG_k =\underline \bG(:,:,k) \in \Real^{R_1 \times R_2}$ are  respectively  the frontal slices of the data tensor $\underline \bX$ and the core tensor $\underline \bG \in \Real^{R_1 \times R_2 \times R_3}$, and $\bA \in \Real^{I \times R_1}$,  $\bB \in \Real^{J \times R_2}$.

\noindent{\bf Generalized Tucker format and its links to TTNS model.} For high-order tensors, $\underline \bX \in \Real^{I_{1,1} \times \cdots \times I_{1,K_1} \times I_{2,1} \times \cdots \times I_{N, K_N}}$, the Tucker-$N$ format can be naturally  generalized  by replacing the
factor matrices, $\bB^{(n)} \in \Real^{I_n \times R_n}$, by higher-order  tensors
$\underline \bB^{(n)} \in \Real^{I_{n,1} \times I_{n,2}\times \cdots  \times I_{n, K_n} \times R_n}$, to give
\be
\underline \bX
\cong \llbracket\underline \bG; \underline \bB^{(1)}, \underline \bB^{(2)}, \ldots,  \underline\bB^{(N)}\rrbracket,
 \label{GenTF2}
 \ee
where the entries of the data tensor are computed as
\be \underline \bX(\bi_1, \ldots,\bi_N) =\sum_{r_1=1}^{R_1}  \cdots \sum_{r_N=1}^{R_N} \underline \bG(r_1, \ldots, r_N) \notag
 \underline \bB^{(1)}(\bi_1,r_1)  \cdots \bB^{(N)}(\bi_N,r_N),
 \ee
and $\bi_n = (i_{n,1} i_{n,2} \ldots i_{n,K_n})$  \cite{Lee-TTfund1}.

Furthermore, the nested (hierarchical) form of such a generalized Tucker decomposition leads to the Tree Tensor  Networks State (TTNS) model \cite{DMRG2013}
(see Figure  \ref{Fig:HT8} and Figure  \ref{Fig:TTNS}), with possibly a varying order of cores,  which can be formulated as 
\be
\underline \bX & = &  \llbracket\underline \bG_1; \; \bB^{(1)},  \bB^{(2)}, \ldots,\bB^{(N_1)}\rrbracket  \notag \\
\underline \bG_1 & = &  \llbracket\underline \bG_2; \; \underline \bA^{(1,2)}, \underline \bA^{(2,2)}, \ldots, \underline\bA^{(N_2,2)}\rrbracket. \notag\\
\cdots \notag \\
\underline \bG_P & = & \llbracket\underline \bG_{P+1}; \; \underline \bA^{(1,P+1)}, \underline \bA^{(2,P+1)}, \ldots, \underline \bA^{(N_{P+1},P+1)}\rrbracket,
 \label{GenTFN3}
 \ee
where $\underline \bG_{p} \in \Real^{R_1^{(p)}\times R_2^{(p)} \times \cdots \times R_{N_p}^{(p)}}$ and
$ \underline \bA^{(n_p,p)} \in \Real^{R_{l_{n_p}}^{(p-1)} \times \cdots \times R_{m_{n_p}}^{(p-1)} \times R_{n_p}^{(p)}}$,
with $p=2,\ldots,P+1$.

 Note that some factor tensors, $\underline \bA^{(n,1)}$  and/or  $\underline \bA^{(n_p,p)}$, can be identity tensors which yield an irregular structure, possibly with a varying order of tensors.
This follows from the simple observation that a mode-$n$ product may have, e.g., the following form
 \be
 \underline \bX \times_n  \bB^{(n)} = \llbracket \underline \bX; \bI_1, \ldots,  \bI_{I_{n-1}},
  \bB^{(n)}, \bI_{I_{n+1}}, \ldots,  \bI_{I_N} \rrbracket. \notag
 \ee
 The efficiency of this representation strongly relies on an appropriate choice of the tree structure.
It is usually assumed that the tree structure of TTNS is given or assumed {\it a priori},
and  recent efforts  aim to find an optimal
tree structure from a subset of tensor entries and without any {\it a priori} knowledge of the
tree structure. This is achieved using so-called
 rank-adaptive cross-approximation techniques which
approximate a tensor by hierarchical tensor formats \cite{Ballani2014adaptiveHT,Ballani2014review}.\\


\noindent {\bf Operations in the Tucker format.} If large-scale data tensors admit an  exact or approximate
representation in their Tucker formats, then most mathematical  operations can be performed more efficiently using the so obtained much smaller core tensors and factor matrices. Consider the $N$th-order tensors $ \underline \bX$ and $ \underline \bY$ in the Tucker format, given by
 \begin{equation}
 \underline \bX = \llbracket \underline \bG_X; \bX^{(1)}, \ldots, \bX^{(N)}\rrbracket \quad \mbox{and} \quad
 \underline \bY = \llbracket \underline \bG_Y; \bY^{(1)}, \ldots, \bY^{(N)} \rrbracket,
 \end{equation}
for which the respective multilinear ranks are $\{R_1,R_2,\ldots,R_N\}$ and $\{Q_1,Q_2,\ldots,Q_N\}$,
 then the following mathematical operations can be performed directly in the Tucker
format{\footnote{Similar operations can be performed in the CP format, assuming that the core tensors are diagonal.}},
which admits a significant reduction in computational costs  \cite{Phan_BTDLxR,Phan_ICASSP15,Lee-TTfund1}:
\begin{itemize}
\item {\bf The addition} of two Tucker tensors of the same order and sizes  
\be
\underline \bX +  \underline \bY  = \llbracket \underline \bG_X \oplus \underline \bG_Y; \; [\bX^{(1)}, \bY^{(1)}], \ldots,  [\bX^{(N)},  \bY^{(N)}] \rrbracket,\phantom{XXXX}
\ee
where $\oplus$ denotes a direct sum of two tensors, and $[\bX^{(n)},\bY^{(n)}] \in \Real^{I_n \times (R_n + Q_n)}$, $\;\;\bX^{(n)} \in \Real^{I_n \times R_n}$ and $\bY^{(n)} \in \Real^{I_n \times Q_n}, \; \forall n$.

\item {\bf The Kronecker product} of two Tucker tensors of arbitrary orders and sizes 
\be
\underline \bX \otimes  \underline \bY  = \llbracket \underline \bG_X \otimes \underline \bG_Y; \; \bX^{(1)} \otimes \bY^{(1)}, \ldots,
\bX^{(N)} \otimes \bY^{(N)}\rrbracket. \quad
\ee

\item {\bf The Hadamard} or element-wise product of two Tucker tensors of  the same order and the same  sizes
\be
\underline \bX \*  \underline \bY  = \llbracket \underline \bG_X \otimes \underline \bG_Y; \; \bX^{(1)} \odot_1 \bY^{(1)}, \ldots, \bX^{(N)} \odot_1 \bY^{(N)}\rrbracket,\qquad
\ee
where $\odot_1$ denotes the mode-1 Khatri--Rao product, also called  the transposed Khatri--Rao product or row-wise Kronecker product.

\item {\bf The inner product} of two Tucker tensors of the same order  and sizes can be reduced to the inner product of two smaller tensors by exploiting the Kronecker product structure in the vectorized form, as follows
\begin{align}
\langle \underline \bX,  \underline \bY \rangle &= \mbox{vec}(\underline \bX)^{\text{T}}  \; \mbox{vec}(\underline \bY)   \\
&= \mbox{vec}(\underline \bG_X)^{\text{T}} \left(\bigotimes_{n=1}^N \bX^{(n)\;\text{T}} \right) \left(\bigotimes_{n=1}^N \bY^{(n)} \right) \mbox{vec}(\underline \bG_Y)  \notag \\
&= \mbox{vec}(\underline \bG_X)^{\text{T}} \left(\bigotimes_{n=1}^N \bX^{(n)\,\text{T}} \; \bY^{(n)} \right) \; \mbox{vec}(\underline \bG_Y)  \notag \\
&=  \langle \llbracket \underline \bG_X; (\bX^{(1)\,\text{T}} \; \bY^{(1)}), \ldots, (\bX^{(N)\,\text{T}} \; \bY^{(N)}) \rrbracket , \underline \bG_Y \rangle. \notag
\end{align}
\item {\bf The Frobenius norm} can be computed in a particularly simple way if the factor matrices are orthogonal, since then all products $\bX^{(n)\,\text{T}} \, \bX^{(n)}, \; \forall n$, become the identity matrices, so that
\begin{align}
\|\bX\|_F &= \langle \underline \bX,  \underline \bX \rangle \notag \\
&= \mbox{vec}\left( \llbracket \underline \bG_X; (\bX^{(1)\,\text{T}} \; \bX^{(1)}), \ldots, (\bX^{(N)\,\text{T}} \; \bX^{(N)}) \rrbracket  \right)^{\text{T}}
\mbox{vec}(\underline \bG_X)  \notag \\
&= \mbox{vec}(\underline \bG_X)^{\text{T}} \; \mbox{vec}(\underline \bG_X) = \|\underline \bG_X\|_F.
\end{align}

\item {\bf The $N$-D discrete convolution} of  tensors $\underline \bX \in \Real^{I_1 \times \cdots \times I_N}$
and $\underline \bY \in \Real^{J_1 \times \cdots \times J_N}$
in their Tucker formats
can be expressed as
\be
\underline \bZ &=& \underline \bX \ast \underline \bY = \llbracket \underline \bG_Z; \bZ^{(1)}, \ldots, \bZ^{(N)} \rrbracket  \\
&\in& \Real^{(I_1+J_1-1) \times \cdots \times (I_N+J_N-1)}. \notag
 \ee
If  $\{R_1,R_2,\ldots,R_N\}$  is the multilinear rank of $\underline \bX$ and $\{Q_1,Q_2,\ldots,Q_N\}$  the multilinear rank $\underline \bY$, then the core tensor $\underline \bG_Z = \underline \bG_X  \otimes \underline \bG_Y \in \Real^{R_1 Q_1 \times \cdots \times R_N Q_N} $  and the factor matrices
 \be
  \bZ^{(n)} = \bX^{(n)} \boxdot_1  \bY^{(n)} \in \Real^{(I_n+J_n-1) \times R_n Q_n},
 \ee
 where $\bZ^{(n)}(:,s_n) = \bX^{(n)}(:,r_n) \ast  \bY^{(n)}(:,q_n) \in \Real^{(I_n+J_n-1)}$ for $s_n= \overline{r_n q_n} = 1,2, \ldots, R_n Q_n$.\\

 \item {\bf Super Fast discrete Fourier transform} (MATLAB functions
 fftn$(\underline \bX)$ and fft$(\bX^{(n)}, [], 1)$) of a tensor in the Tucker format
\be
 {\cal {F}} (\underline \bX) = \llbracket \underline \bG_X;  {\cal {F}} (\bX^{(1)}), \ldots, {\cal {F}} (\bX^{(N)})\rrbracket.
 \ee
Note that if the data tensor admits low multilinear rank approximation, then performing the FFT on factor matrices of relatively small size  $\bX^{(n)} \in \Real^{I_n \times R_n}$, instead of
a large-scale data tensor, decreases
considerably computational complexity.
This approach is referred to as the super fast Fourier transform in Tucker format.

\end{itemize}

\section{Higher Order SVD (HOSVD) for Large-Scale Problems}
\label{sect:HOSVD}

\begin{figure}
(a)
\vspace{-0.1cm}
\begin{center}
\includegraphics[width=11.3905cm]{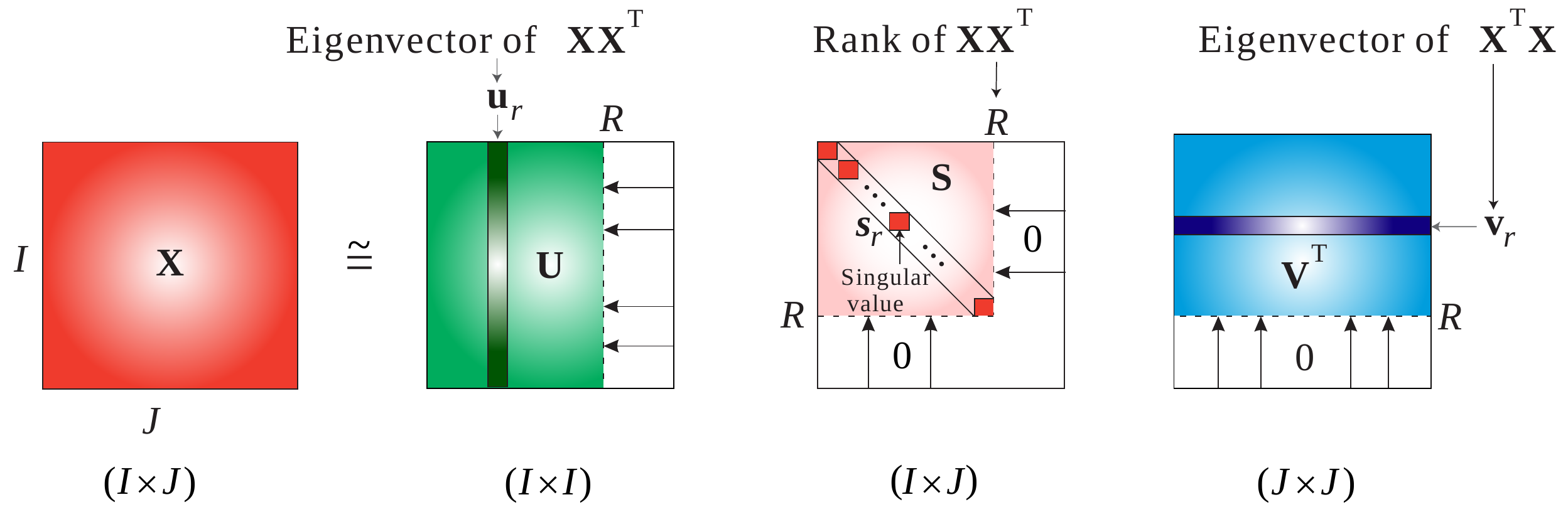}
\end{center}
(b)
\vspace{-0.1cm}
\begin{center}
\includegraphics[width=10.545cm]{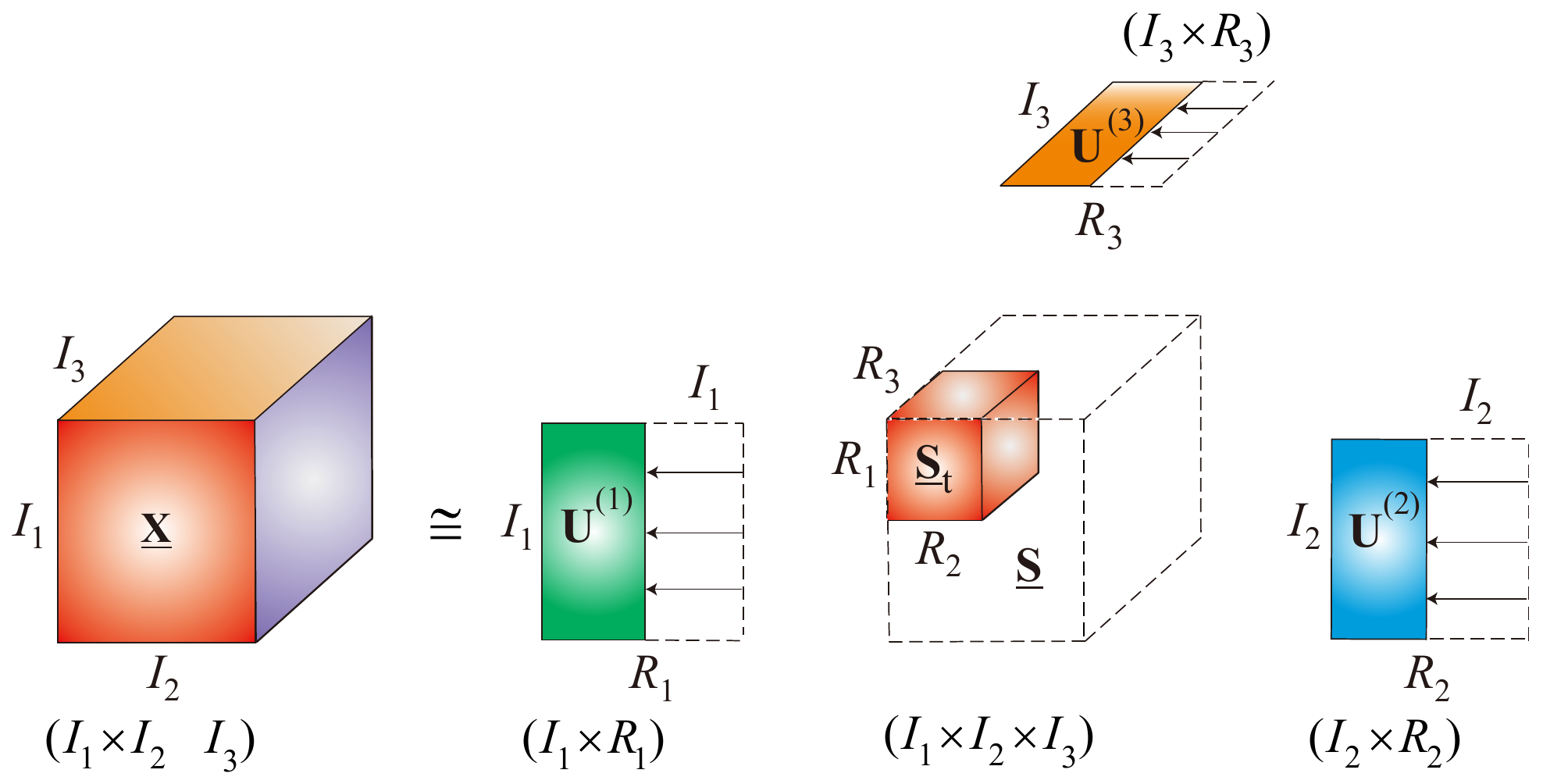}
\end{center}
(c)
\vspace{-0.1cm}
\begin{center}
\includegraphics[width=6.46cm]{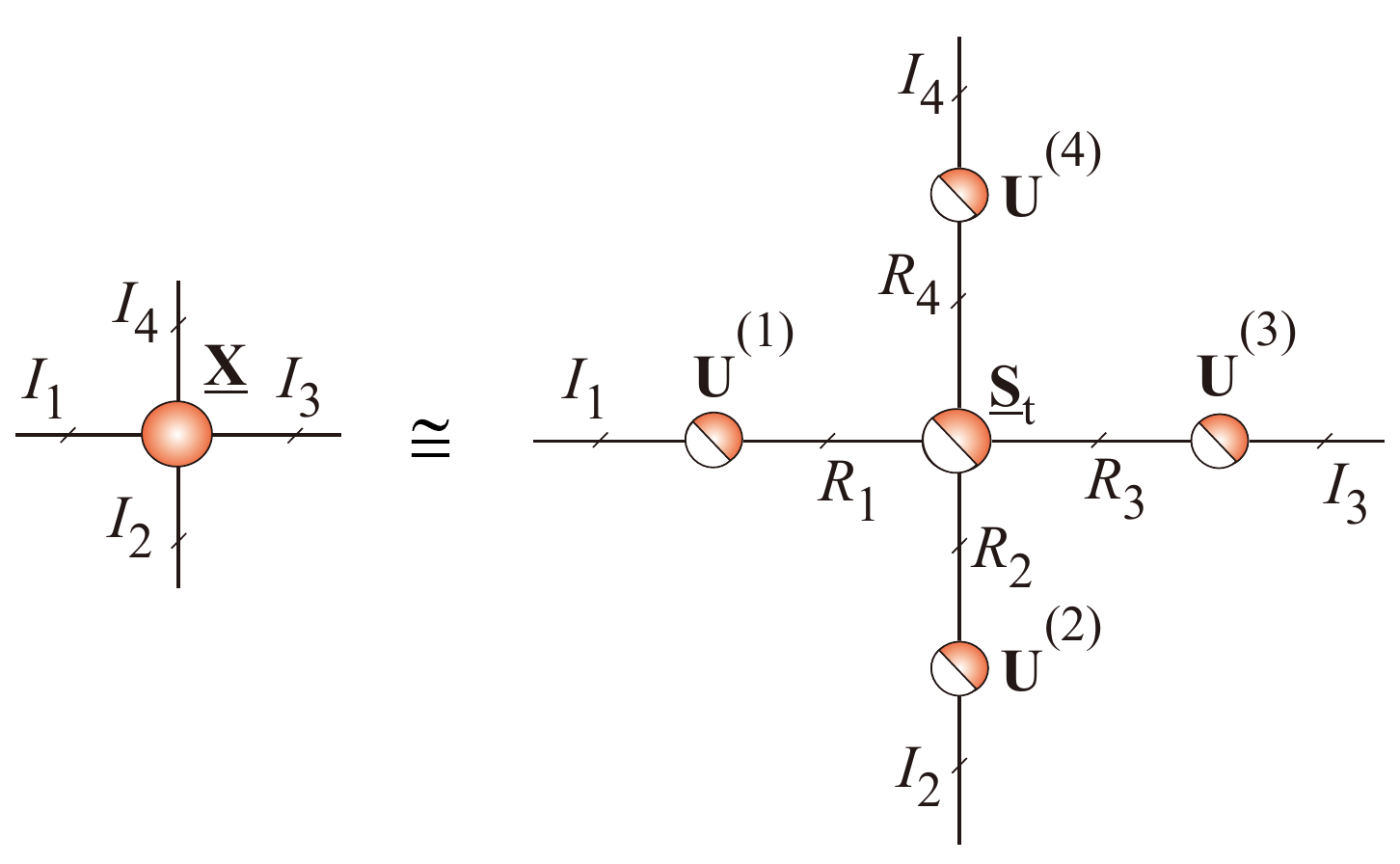}
\end{center}
\caption{{\small Graphical illustration of the truncated SVD and HOSVD. (a) The exact and truncated standard matrix SVD, $\bX \cong \bU \bS \bV^{\text{T}}$. (b) The  truncated (approximative) HOSVD  for a 3rd-order tensor calculated  as
  $\underline \bX \cong \underline \bS_t \times_1 \bU^{(1)} \times_2 \bU^{(2)}  \times_3 \bU^{(3)}$.  (c) Tensor network notation for the HOSVD of a 4th-order tensor $\underline \bX \cong \underline \bS_t \times_1 \bU^{(1)} \times_2 \bU^{(2)}  \times_3 \bU^{(3)} \times_4 \bU^{(4)}$. All the factor matrices,
  $\bU^{(n)} \in \Real^{I_n \times R_n}$, and the core tensor, $\underline \bS_t = \underline \bG \in \Real^{R_1 \times \cdots \times R_N}$, are orthogonal.}}
\label{Fig:HOSVD}
\end{figure}

  The  MultiLinear Singular Value Decomposition (MLSVD),
   also called  the higher-order SVD (HOSVD), can be considered as a
   special form of the constrained Tucker decomposition \cite{HOSVD2000,Lathauwer_HOOI},
in which all factor matrices, $\bB^{(n)}=\bU^{(n)} \in \Real^{I_n \times I_n}$, are orthogonal and the
 core tensor,
$\underline \bG = \underline \bS \in \Real^{I_1 \times I_2 \times \cdots \times I_N}$, is all-orthogonal (see Figure \ref{Fig:HOSVD}).

The orthogonality properties of the core tensor are defined through the following conditions:

\begin{enumerate}

\item {\it All orthogonality.} The slices in each mode are mutually orthogonal, e.g., for a 3rd-order tensor and its lateral slices
\be
\langle\bS_{:,k,:} \bS_{:,l,:}\rangle=0, \quad \mbox{for} \quad k \neq l,
\ee

\item {\it Pseudo-diagonality.} The Frobenius norms of slices in each mode are decreasing with the increase  in
the running index, e.g., for a 3rd-order tensor and its lateral slices
\be
\|\bS_{:,k,:}\|_F \geq  \|\bS_{:,l,:}\|_F, \quad k \geq l.
\ee
These norms play a  role similar  to singular values in  standard matrix SVD.

\end{enumerate}

In practice, the orthogonal matrices $\bU^{(n)} \in \Real^{I_n \times R_n}$, with $R_n \leq I_n$, can be  computed by applying both the randomized  and standard truncated SVD to the unfolded mode-$n$ matrices, $\bX_{(n)} \cong \bU^{(n)} \bS_n \bV^{(n) \text{T}}\in \Real^{I_n \times I_1 \cdots I_{n-1} I_{n+1} \cdots I_N}$.
After obtaining the orthogonal matrices $\bU^{(n)}$  of left singular vectors of $\bX_{(n)}$, for each $n$, the core tensor $\underline \bG =\underline \bS$ can be computed as
\be
\underline \bS = \underline \bX  \times_1 \bU^{(1)\;\text{T}} \times_2 \bU^{(2)\;\text{T}} \cdots  \times_N \bU^{(N)\;\text{T}},
\label{HOSVD:core}
\ee
so that
\be
\underline \bX = \underline \bS  \times_1 \bU^{(1)}  \times_2 \bU^{(2)}   \cdots  \times_N \bU^{(N)}.
\label{HOSVD}
\ee
Due to the orthogonality of the core tensor $\underline \bS$, its slices are also mutually orthogonal.

 Analogous to the standard truncated SVD, a large-scale data tensor, $\underline \bX$, can be approximated by discarding the multilinear
 singular vectors and slices of the core tensor corresponding to small multilinear singular values.  Figure \ref{Fig:HOSVD} and  Algorithm \ref{alg:STHOSVD} outline the truncated HOSVD,
 for which any optimized  matrix SVD procedure can be applied.

 For large-scale tensors, the unfolding matrices, $\bX_{(n)} \in \Real^{I_n \times I_{\bar n}}$ ($I_{\bar n}=I_1 \cdots I_n I_{n+1} \cdots I_N$) may become  prohibitively large (with $I_{\bar n} \gg I_n$), easily exceeding the memory of  standard computers.
  Using a direct and simple divide-and-conquer approach,  the truncated SVD of an  unfolding matrix, $\mathbf{X}_{(n)}=\bU^{(n)} \bS_n
  \bV^{(n) \text{T}}$,  can be  partitioned  into $Q$
slices, as $\mathbf{X}_{(n)} = [
  \mathbf{X}_{1,n}, \mathbf{X}_{2,n}, \ldots , \mathbf{X}_{Q,n}]=\bU^{(n)} \bS_n[
  \mathbf{V}_{1,n}^{\text{T}}, \mathbf{V}_{2,n}^{\text{T}}, \ldots, \mathbf{V}_{Q,n}^{\text{T}}
]$. Next, the orthogonal matrices $\bU^{(n)}$ and the diagonal matrices $\bS_n$ can be obtained from the eigenvalue decompositions $\mathbf{X}_{(n)}\mathbf{X}_{(n)}^{\text{T}}=
\bU^{(n)} \bS_{n}^2 \bU^{(n) \text{T}} =\sum_q\mathbf{X}_{q,n}\mathbf{X}_{q,n}^{\text{T}} \in \Real^{I_n\times I_n}$,
allowing for the terms $\mathbf{V}_{q,n}=\mathbf{X}_{q,n}^{\text{T}} \bU^{(n)} \bS_n^{-1}$ to be computed separately.
This enables us to optimize the  size of the $q$th slice $\bX_{q,n} \in \Real^{I_n \times (I_{\bar n}/Q)}$ so as  to match the available computer memory.                            Such a simple approach to compute matrices $\bU^{(n)}$ and/or $\bV^{(n)}$ does not require loading the entire unfolding matrices
at once into  computer memory; instead the  access to the datasets is sequential.
For current standard sizes of computer memory,  the dimension $I_n$ is typically less than 10,000, while  there is no limit on the dimension $I_{\bar n} = \prod_{k\neq n} I_k$. 

 For very large-scale and low-rank matrices,  instead of the standard truncated SVD approach, we can alternatively apply
the randomized SVD algorithm,  which reduces the  original data matrix $\bX$  to a  relatively small matrix by random sketching, i.e. through multiplication
 with a random  sampling matrix $\mbi \Omega$ (see Algorithm \ref{alg:rSVD}). 
Note that we explicitly allow the rank of the data matrix $\bX$ to be overestimated (that is, $\tilde R= R+P$, where $R$ is a true but unknown rank and $P$ is the over-sampling parameter)  because it is easier to obtain more accurate approximation of this form. Performance of randomized SVD can be further improved by integrating multiple random sketches, that is, by multiplying a data matrix $\bX$ by a set of random matrices $\mbi \Omega_p$  for $p=1,2,\ldots,P$ and  integrating leading low-dimensional subspaces  by applying a Monte Carlo integration method \cite{mrSVD16}.

 Using special  random sampling matrices, for instance, a sub-sampled random Fourier transform,
 substantial gain in the execution time can be achieved, together with the asymptotic complexity of ${\cal {O}}(IJ \log(R))$.
 Unfortunately, this approach  is  not accurate enough for matrices  for which the singular values decay slowly \cite{siam_probLowRank}.

 \begin{algorithm}[t!]
{\small
\caption{\textbf{Sequentially Truncated HOSVD \cite{vannieuwenhoven2012}}}
\label{alg:STHOSVD}
\SetKwFunction{svd}{svd}
\SetKwFunction{tsvd}{truncated\_svd}
\SetKwFunction{reshape}{reshape}
\begin{algorithmic}[1]
\REQUIRE $N$th-order tensor $\underline \bX \in \Real^{I_1 \times I_2 \times \cdots \times I_N}$ and
 approximation \\ accuracy $\varepsilon$
\ENSURE HOSVD in the Tucker format $\underline {\hat\bX} = \llbracket \underline \bS; \bU^{(1)},\ldots, \underline \bU^{(N)}\rrbracket $,  \\ such that  $\|\underline \bX - \underline {\hat\bX}\|_F \leq \varepsilon$
\STATE $\underline \bS \leftarrow \underline \bX$
\FOR{$n=1$ to $N$}
    \STATE $[\bU^{(n)}, \bS,\bV] =  \tsvd(\bS_{(n)},\frac{\varepsilon}{\sqrt {N}})$
    \STATE  $\underline \bS \leftarrow \bV \bS$
\ENDFOR
\STATE $\underline \bS \leftarrow \reshape(\underline \bS,[R_1, \ldots, R_N])$
\RETURN  Core tensor $\underline \bS$ and orthogonal factor matrices\\ $\bU^{(n)} \in \Real^{I_n \times R_n} $.
\end{algorithmic}
}
\end{algorithm}

\begin{algorithm}[t!]
\caption{\textbf{Randomized SVD (rSVD) for large-scale and low-rank matrices
with single sketch \cite{siam_probLowRank}}}
\label{alg:rSVD}
{\small
\begin{algorithmic}[1]
\REQUIRE A matrix $\bX  \in \Real^{I \times J}$, desired or estimated rank $R$, and \\ oversampling parameter $P$ or overestimated rank $\widetilde R =R+P$, \\ exponent of the power method $q$ ($q=0$ or $q=1$)
\ENSURE An approximate  rank-$\widetilde R$ SVD, $\bX \cong \bU \bS \bV^{\text{T}}$, \\
i.e., orthogonal matrices $\bU \in \Real^{I \times\widetilde R}$, $\bV \in \Real^{J \times \widetilde R}$\\
and diagonal matrix of singular values $\bS \in \Real^{\widetilde R \times \widetilde R}$
\STATE Draw a random Gaussian  matrix $\mbi \Omega \in \Real^{J \times \widetilde R}$,
\STATE Form the sample matrix $\bY = (\bX \bX^{\text{T}})^q \;  \bX \mbi \Omega \in \Real^{I \times \widetilde R}$
\STATE Compute a QR decomposition $\bY = \bQ \bR $
\STATE Form the  matrix $\bA = \bQ^{\text{T}} \bX \in \Real^{\widetilde R \times J}$
\STATE Compute the SVD of the small matrix $\bA$ as $ \bA= \widehat \bU \bS \bV^{\text{T}}$
\STATE Form the matrix  $\bU = \bQ \widehat \bU $.
\end{algorithmic}
}
\end{algorithm}

 The truncated HOSVD can be optimized  and implemented in several alternative
 ways. For example,  if $R_n \ll I_n$, the truncated tensor $\underline \bZ \leftarrow \underline \bX \times_1 \bU^{(1)\text{T}}$ yields
 a smaller unfolding matrix $\bZ_{(2)} \in \Real^{I_2 \times R_1I_3\cdots I_N}$, so that the multiplication $\bZ_{(2)} \bZ_{(2)}^{\text{T}}$ can be faster in the next iterations \cite{vannieuwenhoven2012,Ballard-Tucker}.

 Furthermore, since the unfolding matrices  $\bY_{(n)}^{\text{T}}$ are  typically very ``tall and skinny'',  a huge-scale  truncated  SVD and other
 constrained low-rank matrix factorizations can be computed efficiently  based on the Hadoop / MapReduce paradigm \cite{TallSVD2014,TallQR2011,TallNMF2014}.

\emph{Low multilinear rank approximation is always well-posed}, however,
 in contrast to the standard truncated SVD for matrices, \emph{the truncated HOSVD does not yield the best multilinear rank approximation}, but  satisfies the quasi-best approximation property \cite{HOSVD2000} 
 \be
 \|\underline \bX - \llbracket \underline \bS; \bU^{(1)}, \ldots, \bU^{(N)}\rrbracket\| \leq \sqrt{N} \| \underline \bX - \underline \bX_{\mbox{Best}} \|,
 \ee
 where $\underline \bX_{\mbox{Best}}$ is the best multilinear rank approximation of $\underline \bX$, for a specific tensor norm $\|\cdot\|$.

  When it comes to the problem of finding the best approximation, the ALS type algorithm called the Higher Order  Orthogonal Iteration
  (HOOI) exhibits both the advantages and drawbacks of ALS algorithms for CP decomposition. 
  For the HOOI algorithms,
  see Algorithm \ref{alg:HOOI} and Algorithm  \ref{alg:RandHOOI}.
%
For more  sophisticated  algorithms for Tucker decompositions  with orthogonality and  nonnegativity constraints, suitable for large-scale data tensors,
see \cite{PhanHALS2011,Zhou-Cichocki-SP,TallSVD2014,UKang2016}.

\begin{table} 
\vspace{-2.0cm}
 \setlength{\tabcolsep}{2pt}
\renewcommand{\arraystretch}{2}
\centering
\caption{Basic multiway component analysis (MWCA)/Low-Rank Tensor Approximations (LRTA) and related  multiway dimensionality reduction models.
The symbol $\underline \bX \in \Real^{I_1 \times I_2 \times \cdots \times I_N}$ denotes a noisy data tensor, while
$\bY= \underline \bG \times_1 \bB^{(1)} \times_2 \bB^{(2)}  \cdots  \times_N \bB^{(N)}$ is the general constrained Tucker model with the latent factor  matrices $\bB^{(n)} \in \Real^{I_n \times R_n}$ and the core tensor $\underline \bG \in \Real^{R_1 \times R_2 \times \cdots \times R_N}$. In the special case of a CP decomposition,
the core tensor is  diagonal, $\underline \bG = \underline {\mbi \Lambda} \in \Real^{R \times \cdots \times R}$, so that $\underline \bY = \sum_{r=1}^R \lambda_r (\bb_r^{(1)} \circ \bb_r^{(2)} \circ \cdots \circ \bb_r^{(N)})$. }
\label{tab:MWCA}
\centerline{
 {\tabsize \small
 \shadingbox{
{\begin{tabular*}{1.00\textwidth}[t]{@{\extracolsep{\fill}}@{\hspace{1ex}}l|l}
\hline 
 Cost Function  & Constraints \\ \hline
\minitab[l]{Multilinear (sparse) PCA (MPCA) \\
$\max_{\bu_r^{(n)}} \; \underline \bX  \bar \times_1 \bu_r^{(1)} \bar \times_2 \bu_r^{(2)} \cdots  \bar \times_N \bu_r^{(N)}+\gamma \sum_{n=1}^N \|\bu_r^{(n)}\|_1$}  &
\minitab[l]{$\bu_r^{(n) \, \text{T}}\, \bu_r^{(n)} =1, \; \forall (n,r)$ \\
$\bu_r^{(n) \,\text{T}} \, \bu_q^{(n)} =0 \;\; \mbox{for} \;\; r\neq q$}\\  \hline
\minitab[l]{HOSVD/HOOI \\
$\min_{\bU^{(n)}} \; \|\underline \bX - \underline \bG \times_1 \bU^{(1)} \times_2 \bU^{(2)} \cdots  \times_N \bU^{(N)}\|_F^2$}  &
$\bU^{(n) \,\text{T}} \, \bU^{(n)} =\bI_{R_n}, \;
\forall n$ \\  \hline
%
%
%
\minitab[l]{ Multilinear ICA  \\ $\min_{\bB^{(n)}} \;\| \underline \bX - \underline \bG \times_1 \bB^{(1)} \times_2 \bB^{(2)} \cdots  \times_N \bB^{(N)}\|_F^2$}  &\minitab[l] {Vectors of $\bB^{(n)}$ statistically \\as independent as possible} \\  \hline
\minitab[l]{Nonnegative CP/Tucker decomposition \\ (NTF/NTD) \cite{NMF-book} \\
$\min_{\bB^{(n)}} \;\| \underline \bX - \underline \bG \times_1 \bB^{(1)}  \cdots  \times_N \bB^{(N)}\|_F^2$ \\
$+ \gamma \sum_{n=1}^N  \sum_{r_n=1}^{R_n} \|\bb^{(n)}_{r_n}\|_{1}$}
& \minitab[l]{Entries of $\underline \bG$ and $\bB^{(n)}, \; \forall n$ \\
are nonnegative} \\  \hline
%
%
\minitab[l]{Sparse CP/Tucker decomposition \\ $\min_{\bB^{(n)}} \;\| \underline \bX - \underline \bG \times_1 \bB^{(1)}  \cdots  \times_N \bB^{(N)}\|_F^2$ \\
 $+\gamma \sum_{n=1}^N \sum_{r_n=1}^{R_n} \|\bb^{(n)}_{r_n}\|_{1} $ } &\minitab[l]{Sparsity constraints \\imposed on $\bB^{(n)}$} \\  \hline 
\minitab[l]{ Smooth CP/Tucker decomposition \\(SmCP/SmTD) \cite{Yokota-SmCA} \\ $\min_{\bB^{(n)}} \;\| \underline \bX - \underline {\mbi \Lambda} \times_1 \bB^{(1)}  \cdots  \times_N \bB^{(N)}\|_F^2$ \\
$+ \gamma \sum_{n=1}^N \sum_{r=1}^R \|\bL \bb^{(n)}_r\|_{2}$ } & \minitab[l] {Smoothness imposed\\ on vectors $\bb_r^{(n)}$ \\of $\bB^{(n)} \in \Real^{I_n \times R}$, \; $\forall n$\\
via a difference operator $\bL$} \\
\hline 
\end{tabular*}}
}
}}
\end{table}

\begin{algorithm}[t]
{
\caption{\textbf{Higher Order Orthogonal Iteration (HOOI) \cite{Lathauwer_HOOI,Ballard-Tucker}}}
\label{alg:HOOI}
\begin{algorithmic}[1]
\REQUIRE $N$th-order tensor $\underline \bX \in \Real^{I_1 \times I_2 \times \cdots \times I_N}$
(usually  in Tucker/HOSVD format)
\ENSURE Improved  Tucker approximation   using ALS approach, with orthogonal factor matrices $\bU^{(n)}$
\STATE  Initialization  via the standard HOSVD (see Algorithm \ref{alg:STHOSVD})
\REPEAT
\FOR{$n=1$ to $N$}
     \STATE  $\underline \bZ \leftarrow \underline \bX \times_{p \neq n} \{\bU^{(p)\,\text{T}} \}$
    \STATE $\bC \leftarrow \bZ_{(n)} \bZ_{(n)}^{\text{T}}\in \Real^{R \times R}$
    \STATE  $\bU^{(n)} \leftarrow$ leading $R_n$ eigenvectors of $\bC$
\ENDFOR
\STATE $\underline \bG \leftarrow \underline \bZ \times_N \bU^{(N)\,\text{T}}$
\UNTIL the cost function $(\|\underline \bX\|^2_F - \|\underline \bG\|^2_F)$ ceases to decrease
\RETURN  $\llbracket \underline \bG; \bU^{(1)}, \bU^{(2)}, \ldots ,  \bU^{(N)}\rrbracket$
\end{algorithmic}
}
\end{algorithm}

\begin{algorithm}[h!]
{
\caption{\textbf{HOOI using randomization for large-scale data \cite{ZhouIP15}}}
\label{alg:RandHOOI}
\begin{algorithmic}[1]
\REQUIRE $N$th-order tensor $\underline \bX \in \Real^{I_1 \times I_2 \times \cdots \times I_N}$
 and multilinear  rank $\{R_1,R_2, \ldots, R_N\}$
\ENSURE Approximative representation of a tensor in Tucker  format, \\with orthogonal factor matrices $\bU^{(n)} \in \Real^{I_n \times R_n}$
\STATE Initialize factor matrices $\bU^{(n)}$ as random Gaussian matrices \\
Repeat steps (2)-(6) only two times:
\FOR{$n=1$ to $N$}
 \STATE  $\underline \bZ = \underline \bX \times_{p \neq n} \{ \bU^{(p)\, \text{T}}\}$
 \STATE Compute $\tilde \bZ^{(n)} =\bZ_{(n)} \mbi \Omega^{(n)} \in \Real^{I_n \times R_n}$,  where $\mbi \Omega^{(n)} \in \Real^{\prod_{p \neq n} R_p \times R_n}$ \\
 is a random matrix  drawn from Gaussian distribution
    \STATE Compute $\bU^{(n)}$ as an orthonormal basis of  $\tilde \bZ^{(n)}$, e.g., by using QR decomposition
\ENDFOR
\STATE Construct the core tensor as\\
$\underline \bG = \underline \bX \times_1 \bU^{(1)\;\text{T}} \times_2 \bU^{(2)\;\text{T}} \cdots  \times_N  \bU^{(N)\;\text{T}}$
\RETURN    $\underline \bX \cong \llbracket \underline \bG; \bU^{(1)}, \bU^{(2)}, \ldots , \underline \bU^{(N)}\rrbracket$
\end{algorithmic}
}
\end{algorithm}

\begin{algorithm}[h!]
{
\caption{\textbf{Tucker decomposition with constrained factor matrices via 2-way CA /LRMF }}
\label{alg:Tucker-LRMF}
\begin{algorithmic}[1]
\REQUIRE $N$th-order tensor $\underline \bX \in \Real^{I_1 \times I_2 \times \cdots \times I_N}$,
  multilinear rank $\{R_1,\ldots, R_N\}$~and  desired constraints imposed
 on factor matrices $\bB^{(n)} \in \Real^{I_n \times R_n}$
\ENSURE Tucker decomposition with constrained factor matrices $\bB^{(n)}$ \\
using LRMF and a simple unfolding   approach
\STATE  Initialize randomly or via standard HOSVD (see Algorithm \ref{alg:STHOSVD})
\FOR{$n=1$ to $N$}
     \STATE  Compute specific LRMF or 2-way  CA (e.g., RPCA, ICA, NMF) of unfolding
     $\bX^{\text{T}}_{(n)} \cong \bA^{(n)} \bB^{(n)\;\text{T}} $ or  $\bX_{(n)} \cong \bB^{(n)} \bA^{(n)\;\text{T}} $
\ENDFOR
\STATE  Compute core tensor
$\underline \bG =\underline \bX \times_1 [\bB^{(1)}]^{\dagger} \times_2 [\bB^{(2)}]^{\dagger} \cdots \times_N [\bB^{(N)]^{\dagger}}$
\RETURN  Constrained Tucker decomposition
$\underline \bX \cong \llbracket \underline \bG,\bB^{(1)},\ldots , \underline \bB^{(N)}\rrbracket$
\end{algorithmic}
}
\end{algorithm}

\begin{figure}
(a) {\small{Sequential computation}}
\begin{center}
\includegraphics[width=9.5cm]{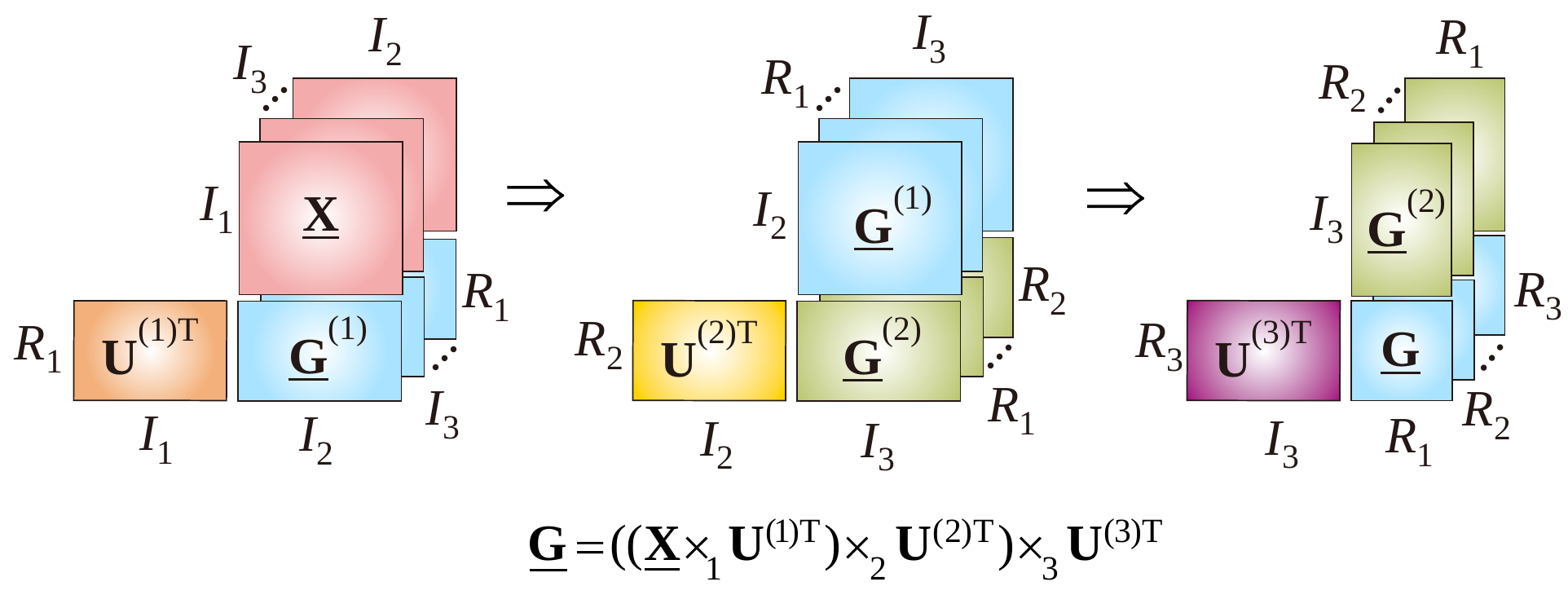}
\end{center}
(b) {\small{Fast matrix-by-matrix approach}}
\begin{center}
\includegraphics[width=9.5cm]{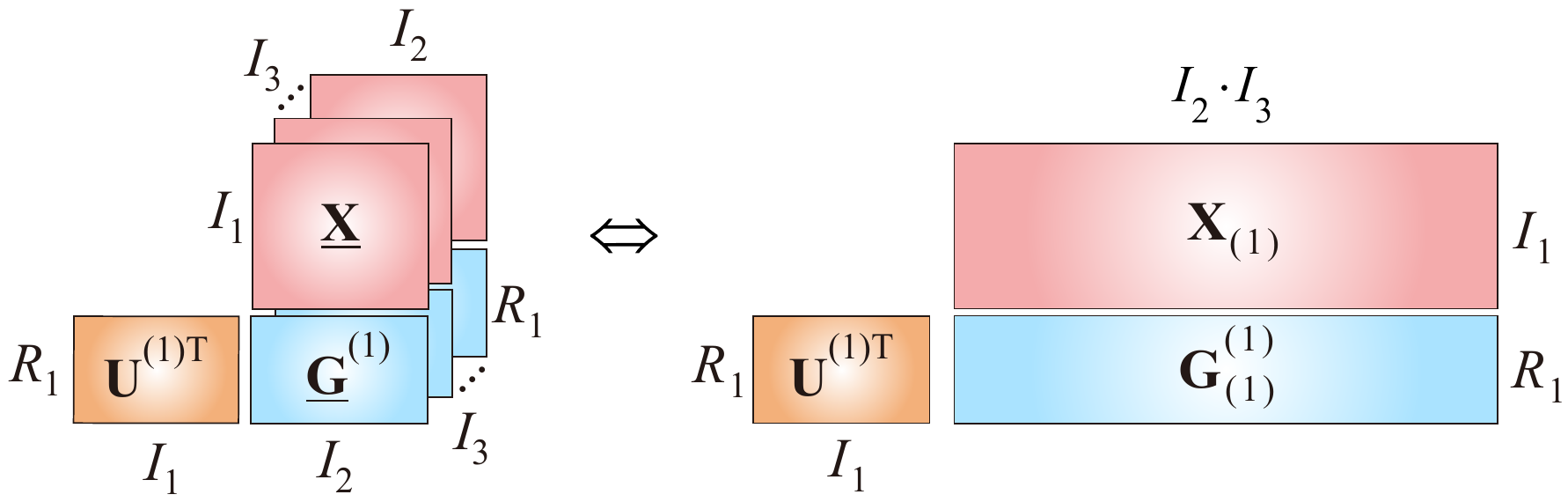}\\
\end{center}
(c) {\small{Divide-and-conquer approach}}
\begin{center}
\includegraphics[width=9.7cm]{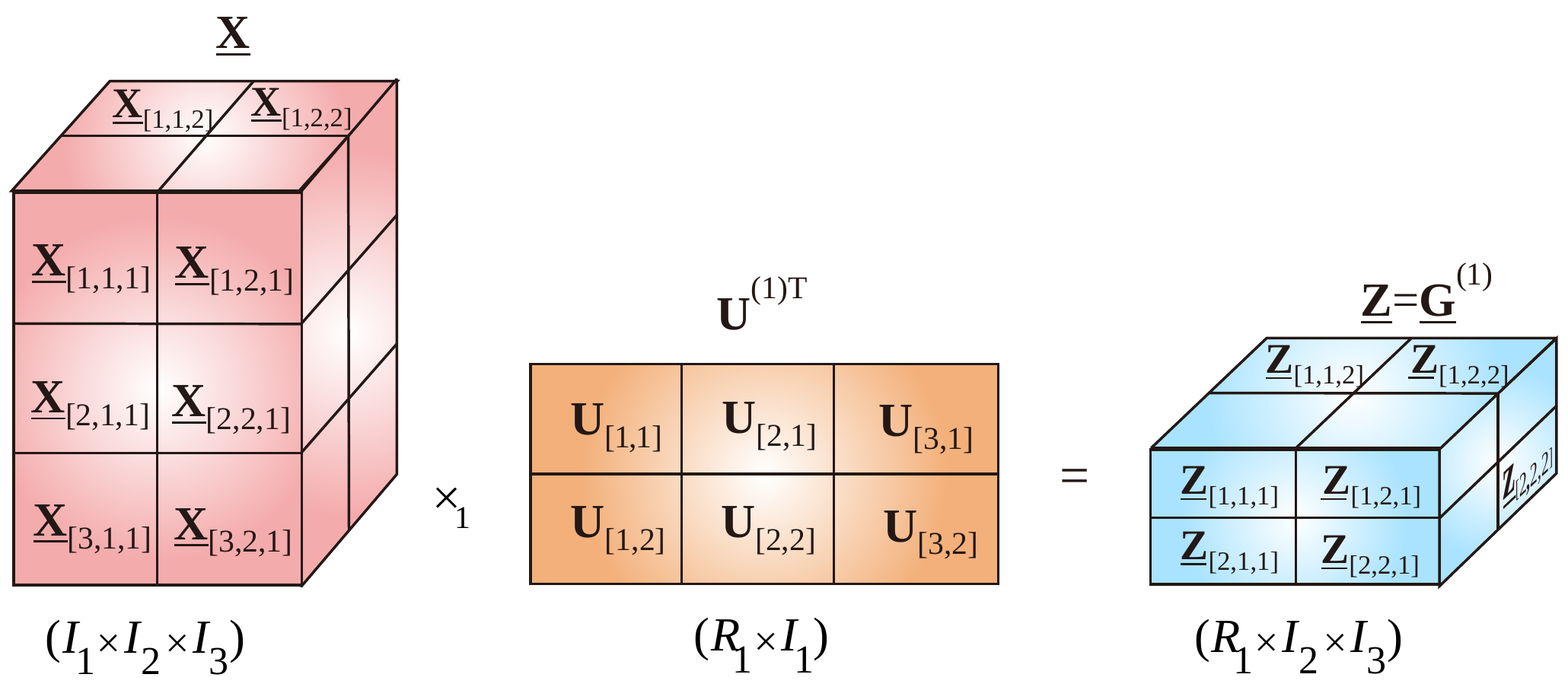}
\end{center}
\caption{{\small Computation of a multilinear (Tucker) product  for  large-scale  HOSVD.  (a) Standard sequential computing of multilinear products (TTM) $\underline \bG = \underline \bS= (((\underline \bX \times_1 \bU^{(1)\text{T}}) \times_2  \bU^{(2)\text{T}}) \times_3 \bU^{(3)\text{T}})$. (b) Distributed implementation through fast matrix-by-matrix multiplications. (c) An alternative method for  large-scale problems using the ``divide and conquer'' approach, whereby  a data tensor,  $\underline \bX$,  and  factor matrices, $\bU^{(n)\text{T}}$,  are partitioned into suitable small blocks: Subtensors $\underline \bX_{[k_1,k_2,k_3]} $ and block  matrices $\bU^{(1)\text{T}}_{[k_1,p_1]}$. The blocks of a tensor,  $\underline \bZ=\underline \bG^{(1)}=\underline \bX \times_1 \bU^{(1)\text{T}}$,  are computed as $\underline \bZ_{[q_1,k_2,k_3]}= \sum_{k_1=1}^{K_1} \bX_{[k_1,k_2,k_3]} \times_1 \bU^{(1)\text{T}}_{[k_1,q_1]}$ (see Eq.~(\ref{outcore-prod}) for a general case).}}
\label{Fig:outcore}
\end{figure}

When a data tensor $\underline \bX$ is very large and cannot be stored in  computer memory,
another challenge  is to compute a core tensor $\underline \bG =\underline \bS$  directly, using the formula (\ref{HOSVD:core}).
Such computation  is performed sequentially by fast matrix-by-matrix
multiplications{\footnote{Efficient and parallel (state of the art) algorithms  for multiplications  of such very large-scale matrices are proposed in \cite{TTM2015,Ballard15MbyMimproved}.}}, as illustrated in Figure \ref{Fig:outcore}(a)  and~(b). 

We have shown that for very large-scale problems, it is useful to divide a data tensor $\underline \bX$ into small blocks $\underline \bX_{[k_1,k_2,\ldots, k_N]} $.
   In a similar way, we can partition  the orthogonal factor matrices $\bU^{(n)\text{T}}$ into the corresponding blocks of matrices $\bU^{(n)\text{T}}_{[k_n,p_n]}$, as illustrated in Figure \ref{Fig:outcore}(c) for  3rd-order tensors \cite{Wang-out-core05,Suter13}. For example, the blocks within the resulting tensor $\underline \bG^{(n)}$  can be computed
   sequentially or in parallel, as follows:
\be
\underline \bG^{(n)}_{[k_1,k_2,\ldots,q_n,\ldots,k_N]}= \sum_{k_n=1}^{K_n} \bX_{[k_1,k_2,\ldots,k_n,\ldots,k_N]} \times_n \bU^{(n) \;\text{T}}_{[k_n,q_n]}.
\label{outcore-prod}
\ee

\noindent{\bf Applications.}
We have shown that the Tucker/HOSVD decomposition may be considered as a multilinear extension of PCA \cite{Kroonenberg}; it therefore generalizes signal subspace techniques and finds application in areas including multilinear blind source separation, classification, feature extraction, and subspace-based harmonic retrieval \cite{Vasilescu,Haardt08,Phan2010TF,MSLSurvey2011}.
In this way, a low multilinear rank approximation achieved through
Tucker decomposition may yield higher Signal-to-Noise Ratio (SNR) than the SNR for the original raw data tensor, which also makes Tucker decomposition a  natural tool for signal compression and  enhancement. 

It was recently shown that  HOSVD can also perform simultaneous subspace selection (data compression) and K-means clustering,
both unsupervised learning tasks \cite{Huang2008,Papalexakis2013}. This is important, as  a combination of these methods can both identify and classify   ``relevant'' data, and in this way not only reveal desired information  but also simplify feature  extraction.\\

\begin{figure}[t]
\centering
\includegraphics[width=10.6cm]{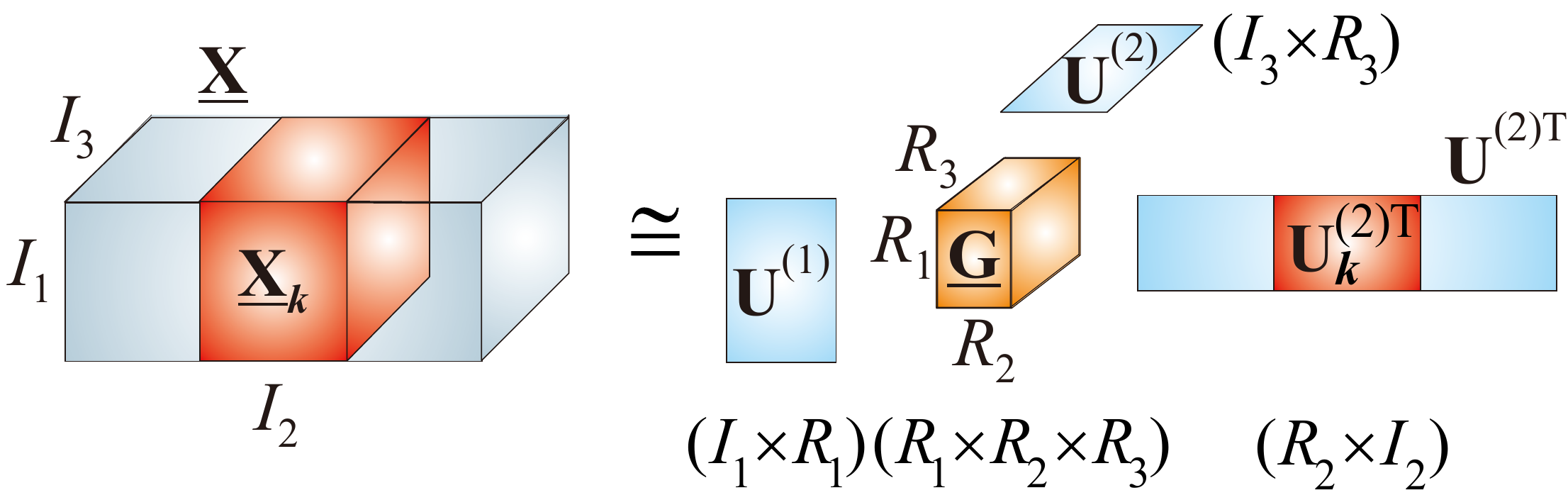}
\caption{Conceptual model for performing  the  HOSVD for a very  large-scale 3rd-order data tensor. This is achieved by dividing the  tensor into blocks
  $\underline \bX_k \cong \underline \bG \times_1 \bU^{(1)} \times_2 \bU_k^{(2)} \times_3 \bU^{(3)}$, $(k=1,2\ldots,K)$. It assumed that the data tensor $\bX \in \Real^{I_1 \times I_2 \times I_3}$ is sampled by  sliding  the block $\underline \bX_k$  from  left to  right (with an overlapping sliding window).
The model can be used for  anomaly detection by fixing the core tensor and some factor matrices while monitoring the
changes along one or more specific modes (in our case mode two). Tensor decomposition is then first performed for a sampled
(pre-selected) small block, followed by the analysis of changes in specific smaller--dimensional factor matrices  $\bU^{(n)}$.}
\label{Fig:blocks}
\end{figure}

\noindent {\bf Anomaly detection using HOSVD.} Anomaly detection refers to the discrimination of some specific patterns, signals, outliers or features   that do not conform to certain expected behaviors, trends  or properties \cite{anomaly2009,anomaly-tensor}. While such analysis can be performed in different domains, it is most frequently based on spectral  methods such as  PCA, whereby  high dimensional data are projected onto a lower-dimensional subspace in which the anomalies may be identified more easier.
The main assumption within
such approaches is that  the normal and abnormal patterns, which may be difficult to distinguish in the original space,  appear significantly different in the projected subspace.
When considering very large datasets, since the basic Tucker decomposition model
 generalizes PCA and SVD, it offers a natural framework
 for anomaly detection via HOSVD,  as illustrated in Figure  \ref{Fig:blocks}. To handle the exceedingly large dimensionality,
we may first  compute tensor decompositions for sampled
(pre-selected) small blocks of the original large-scale 3rd-order tensor, followed by the analysis of changes in specific factor matrices  $\bU^{(n)}$.  A simpler form is straightforwardly obtained by fixing the core tensor and some factor matrices while monitoring the changes along one or more specific modes, as the block tensor moves from left to right as shown in Figure \ref{Fig:blocks}.

\section{Tensor Sketching Using Tucker Model}
\label{sect:Sketching}

The notion of sketches refers to replacing the original huge matrix or tensor by a new matrix or tensor of a significantly smaller size  or compactness,
 but which  approximates well the original matrix/tensor.
Finding such sketches in  an efficient way is important for the analysis of big data,
as a computer processor (and memory) is often incapable of handling  the whole data-set
 in a feasible amount of time. 
 For these reasons,  the computation is often spread among a set of processors
which for standard ``all-in-one'' SVD algorithms, are unfeasible.

Given a very large-scale tensor $\underline \bX$ , a useful approach
is to compute  a sketch tensor $\underline \bZ$
or set of sketch tensors $\underline \bZ_n$ that are of significantly smaller sizes than the original one.

There exist  several  matrix and tensor  sketching approaches:
sparsification, random projections, fiber subset selections, iterative sketching techniques and distributed  sketching.
 We review the main sketching approaches which are promising for tensors.

\smallbreak
\noindent {\bf 1. Sparsification}  generates a sparser version of the tensor which, in general, can be stored more efficiently and admit faster multiplications by factor matrices. This is achieved by decreasing the number on non-zero entries and quantizing or rounding up entries.
A simple technique is element-wise sparsification which zeroes out all sufficiently small elements (below some threshold) of a  data tensor, keeps all sufficiently large elements, and randomly samples the remaining elements  of the tensor with sample probabilities proportional to the square of their magnitudes \cite{Drineas15}.

\smallbreak
\noindent {\bf 2. Random Projection} based sketching randomly combines fibers  of a data tensor
in all or selected modes, and is related to the concept of
a randomized subspace embedding, which is used to solve a variety
of numerical linear algebra problems (see  \cite{Tropp16} and references therein).

%
\smallbreak
\noindent {\bf 3. Fiber subset selection},  also called tensor cross approximation (TCA),
 finds a small subset of fibers
which approximates the entire data tensor. For the matrix case, this problem
is known as the Column/Row Subset Selection  or CUR Problem which has been thoroughly
investigated  and  for which there exist several algorithms with almost matching lower bounds
\cite{RandomizedLowRank_Michael_2011,desai2016improved,ghashami2016efficient}.

\section{Tensor Sketching via Multiple Random Projections}

The random projection framework  has been developed for computing structured low-rank
approximations of a data tensor from  (random) linear projections of much lower dimensions
than the data tensor itself \cite{caiafa2015stable,Tropp16}.
Such techniques have many potential applications in large-scale numerical
multilinear algebra  and optimization problems.

\begin{figure}
(a)
\vspace{-0.1cm}
\begin{center}
\includegraphics[width=9.86 cm]{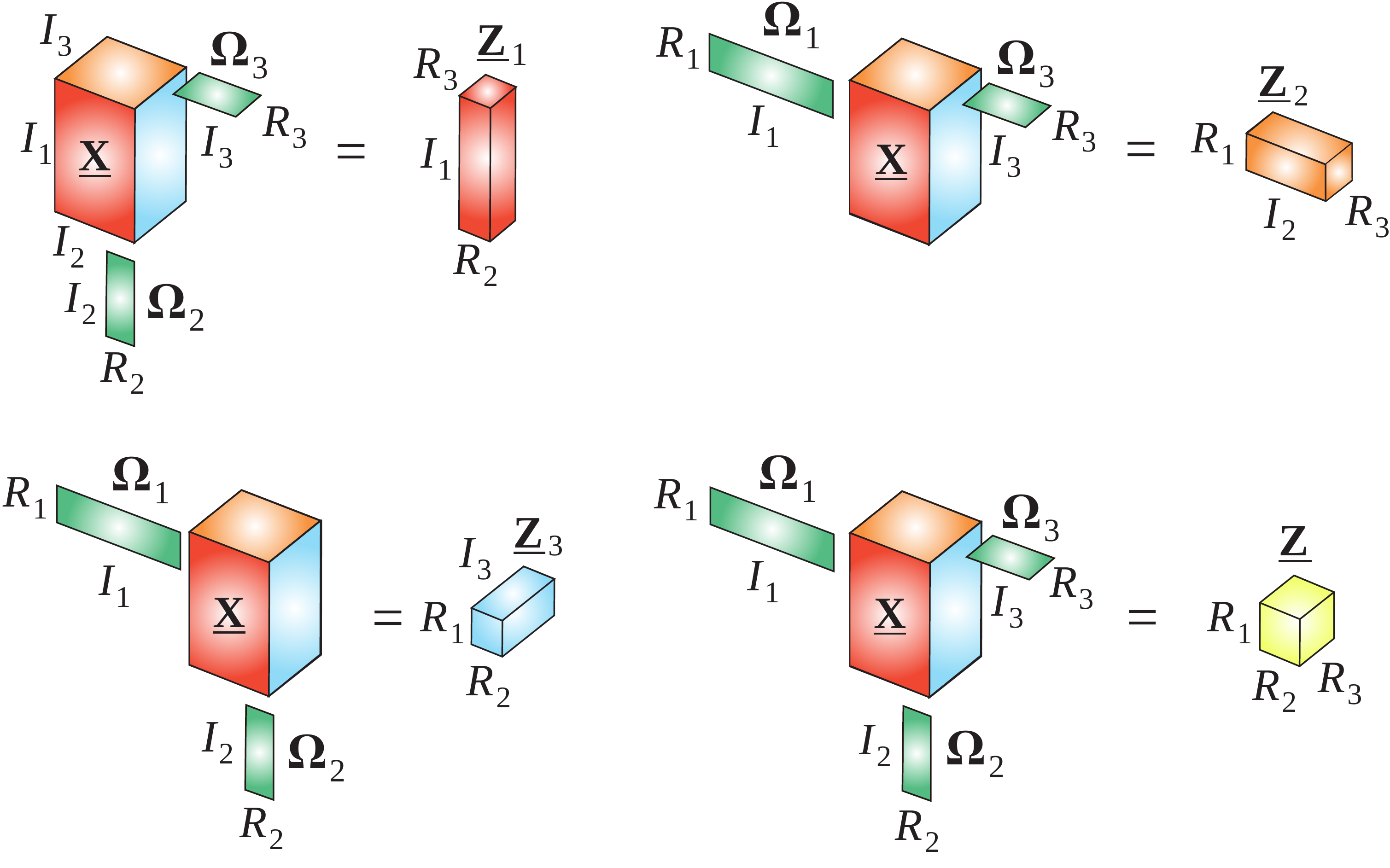}
\end{center}
(b)
\vspace{-0.1cm}
\begin{center}
\includegraphics[width=9.86 cm]{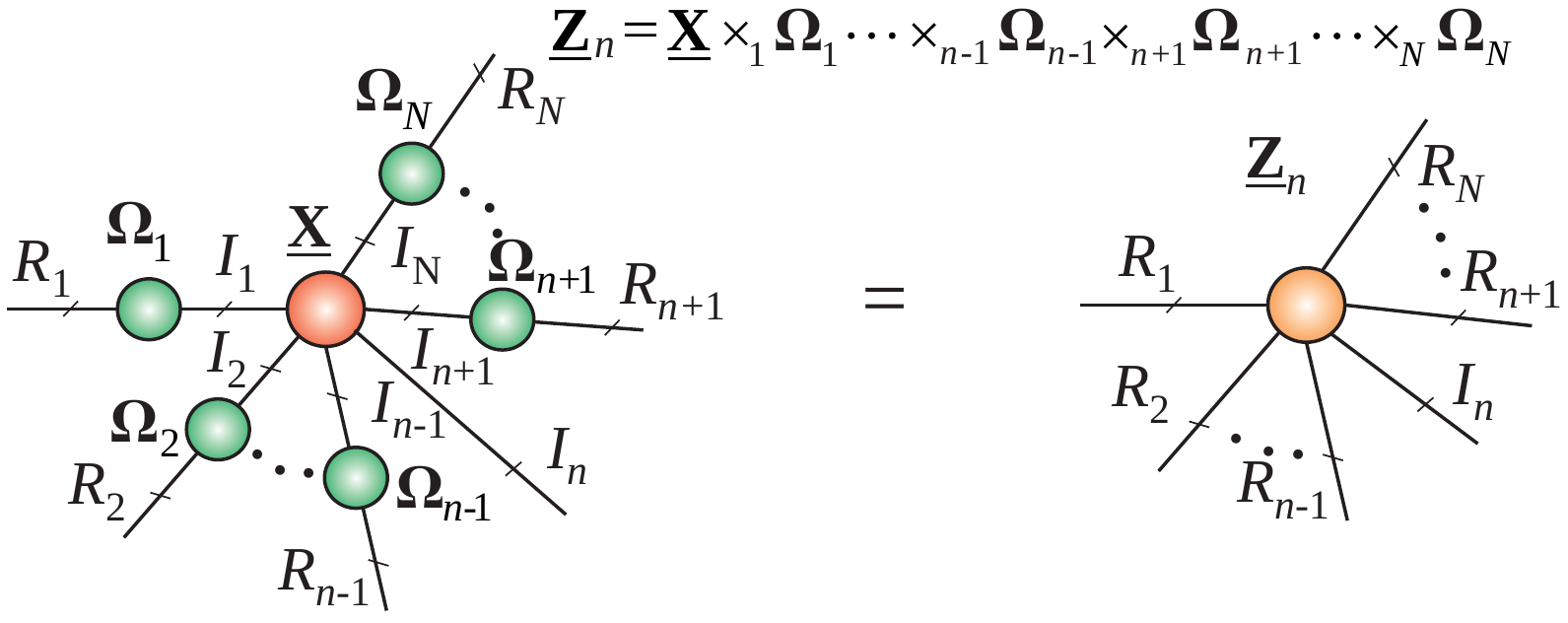} \\
\includegraphics[width=8.16 cm]{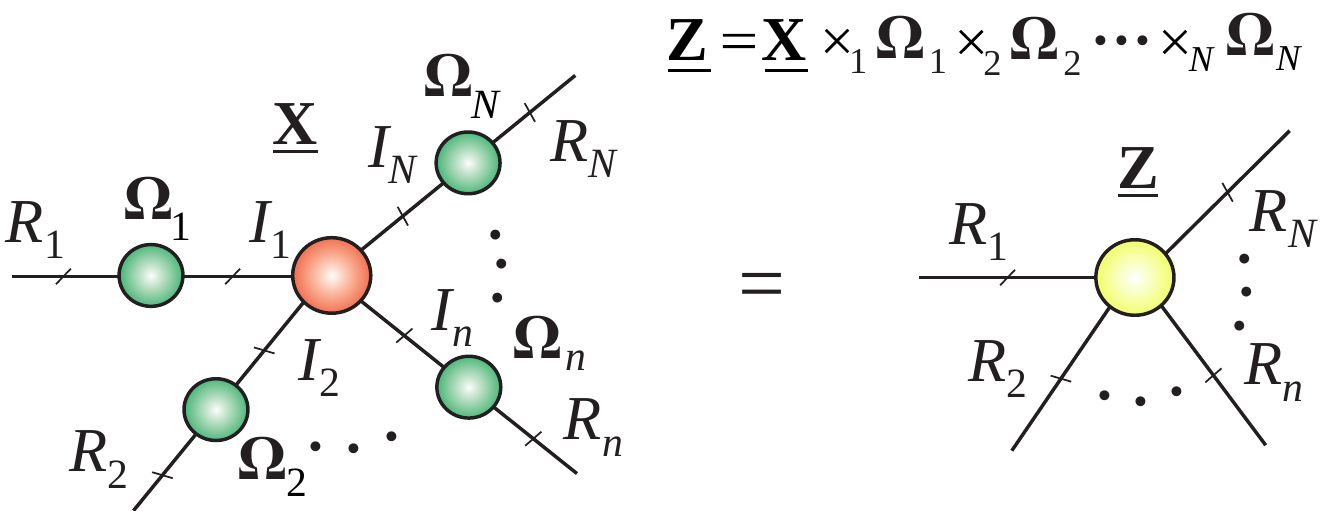}
\end{center}
\caption{Illustration of  tensor sketching using  random projections of a data tensor. (a) Sketches of a 3rd-order tensor   $\underline \bX \in \Real^{I_1 \times I_2 \times I_3}$ given by
$\underline \bZ_1 = \underline \bX \times_2 \mbi \Omega_2 \times_3 \mbi \Omega_3 \in \Real^{I_1 \times R_2 \times R_3}$, $\underline \bZ_2 = \underline \bX \times_1 \mbi \Omega_1 \times_3 \mbi \Omega_3 \in \Real^{R_1 \times I_2 \times R_3}$, $\underline \bZ_3 = \underline \bX \times_1 \mbi \Omega_1 \times_2 \mbi \Omega_2 \in \Real^{R_1 \times R_2 \times I_3}$, and $\underline \bZ = \underline \bX \times_1 \mbi \Omega_1 \times_2 \mbi \Omega_2 \times_3 \mbi \Omega_3 \in \Real^{R_1 \times R_2 \times R_3}$. (b) Sketches for an  $N$th-order tensor  $\underline \bX \in \Real^{I_1 \times \cdots \times I_N}$.}
\label{Fig:Sketching}
\end{figure}

Notice that for an $N$th-order tensor $\underline \bX \in \Real^{I_1 \times I_2 \times \cdots \times I_N}$, we can compute the following sketches
\be
\underline \bZ  = \underline \bX \times_1 \mbi \Omega_1 \times_{2} \mbi \Omega_{2}  \cdots \times _N \mbi \Omega_N
\ee
and
\begin{equation}
\underline \bZ_{\,n}  = \underline \bX \times_1 \mbi \Omega_1   \cdots \times_{n-1} \mbi \Omega_{n-1} \times_{n+1} \mbi \Omega_{n+1}   \cdots \times _N \mbi \Omega_N,
\end{equation}
for $n=,1,2,\ldots,N$, where $\mbi \Omega_n \in \Real^{R_n \times I_n}$ are statistically independent random matrices with $R_n \ll I_n$, usually   called test (or sensing) matrices.

A sketch can be implemented using test matrices drawn from various  distributions. The
choice of a distribution leads to some tradeoffs \cite{Tropp16}, especially regarding
 (i) the costs of randomization, computation, and communication to generate the test matrices;
(ii) the storage costs for the test matrices and the sketch; (iii) the arithmetic costs for
sketching and updates; (iv) the numerical stability of reconstruction algorithms; and (v) the
quality of a priori error bounds. The most important distributions of random test matrices include: 
\begin{itemize}
\item  {\bf Gaussian random projections} which  generate random matrices with standard normal distribution. Such matrices  usually   provide excellent performance in practical scenarios  and accurate a priori error bounds.

\item {\bf Random matrices with  orthonormal columns}  that span uniformly distributed random subspaces of dimensions $R_n$.
Such  matrices behave similar to Gaussian case, but  usually
exhibit even better numerical stability, especially when $R_n$ are large.

\item {\bf Rademacher  and super-sparse Rademacher random projections} that have independent Rademacher entries which take the values $\pm 1$ with equal probability.
Their properties are similar to  standard normal test matrices, but
exhibit some improvements in the cost of storage  and computational complexity. In a special case,  we may use ultra sparse Rademacher test matrices, whereby in  each  column of a test matrix independent Rademacher random variables are placed only in very few  uniformly random locations determined by a
sampling parameter $s$; the remaining entries are set to zero.  In an extreme case of
maximum sparsity, $s=1$, and each column of a test matrix  has exactly only one nonzero entry.

\item {\bf Subsampled randomized Fourier transforms} based on test matrices take
the following form
\be
\mbi \Omega_n = \bP_n \bF_n \bD_n,
\ee
where $\bD_n$ are diagonal square matrices with independent
Rademacher entries,  $\bF_n$ are discrete cosine transform (DCT)
 or discrete Fourier transform (DFT) matrices, and entries of the matrix $\bP_n$ are
 drawn   at random from a uniform distribution.
\end{itemize}
 {\bf Example.} The concept of  tensor sketching via random projections is illustrated in Figure
\ref{Fig:Sketching} for a 3rd-order tensor and for a general case of $N$th-order tensors.
  For a 3rd-order tensor with volume (number of entries) $I_1 I_2 I_3$ we have  four possible sketches which are  subtensors of much
   smaller sizes, e.g., $I_1  R_2  R_3$, with $R_n \ll I_n$, if the sketching is performed along mode-2 and mode-3, or $R_1  R_2  R_3$, if the sketching is performed along all three modes (Figure \ref{Fig:Sketching}(a) bottom right). From these subtensors we can reconstruct any huge  tensor if it has low a multilinear rank
(lower than $\{R_1,R_2,\ldots,R_n\}$).

In more general scenario, it can be  shown \cite{caiafa2015stable} that the $N$th order tensor data tensor $\underline \bX$  with sufficiently low-multilinear rank
can be reconstructed perfectly from the sketch tensors
$\underline \bZ_n $, for $n=1,2,\ldots,N$, as follows
\be
\underline {\hat{\bX}}  = \underline \bZ \times_1 \bB^{(1)}  \times_2 \bB^{(2)} \cdots \times _N \bB^{(N)},
\ee
where $\bB^{(n)} = [\underline \bZ_n]_{(n)} \bZ_{(n)}^{\dag}$ for $n=1,2, \ldots,N$
(for more detail see the next section).

\section{Matrix/Tensor Cross-Approximation (MCA/TCA)}

Huge-scale matrices can be factorized using the Matrix Cross-Approximation (MCA) method, which is also known under the names of  Pseudo-Skeleton or  CUR matrix decompositions \cite{Goreinov:1997,Goreinov:1997b,mahoney2008tensor,Mahoney:2009,oseledets2010tt,bebendorf2011adaptive,Bebendorf2015,Khoromskij16efficient}.
The main idea behind the MCA  is  to provide reduced dimensionality of data through a linear combination
of only a few ``meaningful'' components, which are exact replicas of columns and rows of the original data matrix. Such an approach is based on the fundamental assumption that  large datasets are highly redundant and can therefore  be approximated by low-rank matrices, which significantly  reduces  computational complexity at the cost of a marginal loss of information.
%

\begin{figure}[t]
\centering
\includegraphics[width=10.6 cm]{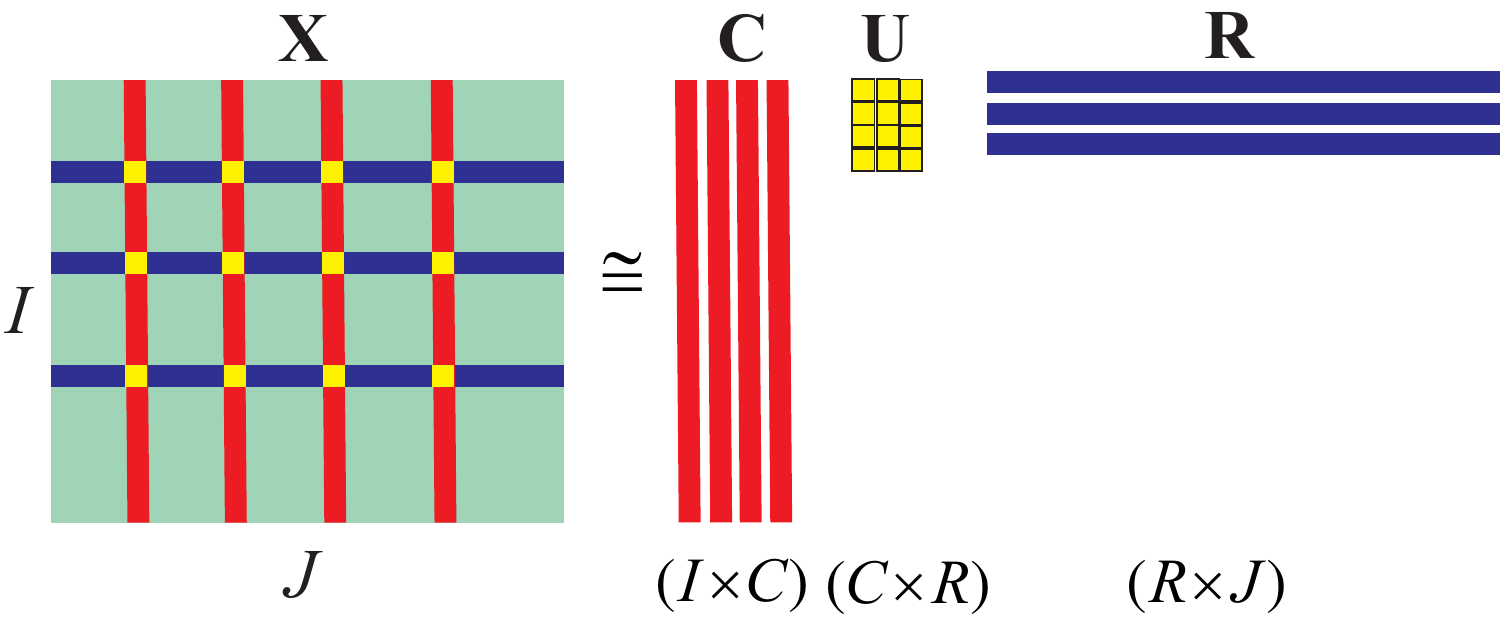}
\caption{Principle of the matrix cross-approximation which decomposes  a huge matrix $\bX$ into a product of three matrices, whereby only a small-size core matrix $\bU$ needs to be computed.}
\label{Fig:CUR}
\end{figure}

The MCA method factorizes a data matrix $\bX \in \Real^{ I \times J}$ as \cite{Goreinov:1997b,Goreinov:1997}
(see Figure \ref{Fig:CUR})
\be
\bX =  \bC \bU \bR + \bE,
\label{CUR1}
\ee
where $\bC \in \Real^{I \times C}$ is a matrix
constructed from $C$ suitably selected columns of the data matrix $\bX$, matrix $\bR \in
\Real^{R \times J}$ consists of $R$ appropriately selected rows of $\bX$, and matrix $\bU \in
\Real^{C \times R}$ is calculated so as to minimize the norm of the error $\bE \in \Real^{I
\times J}$.

\begin{figure}[t]
\centering
\includegraphics[width=11.75 cm]{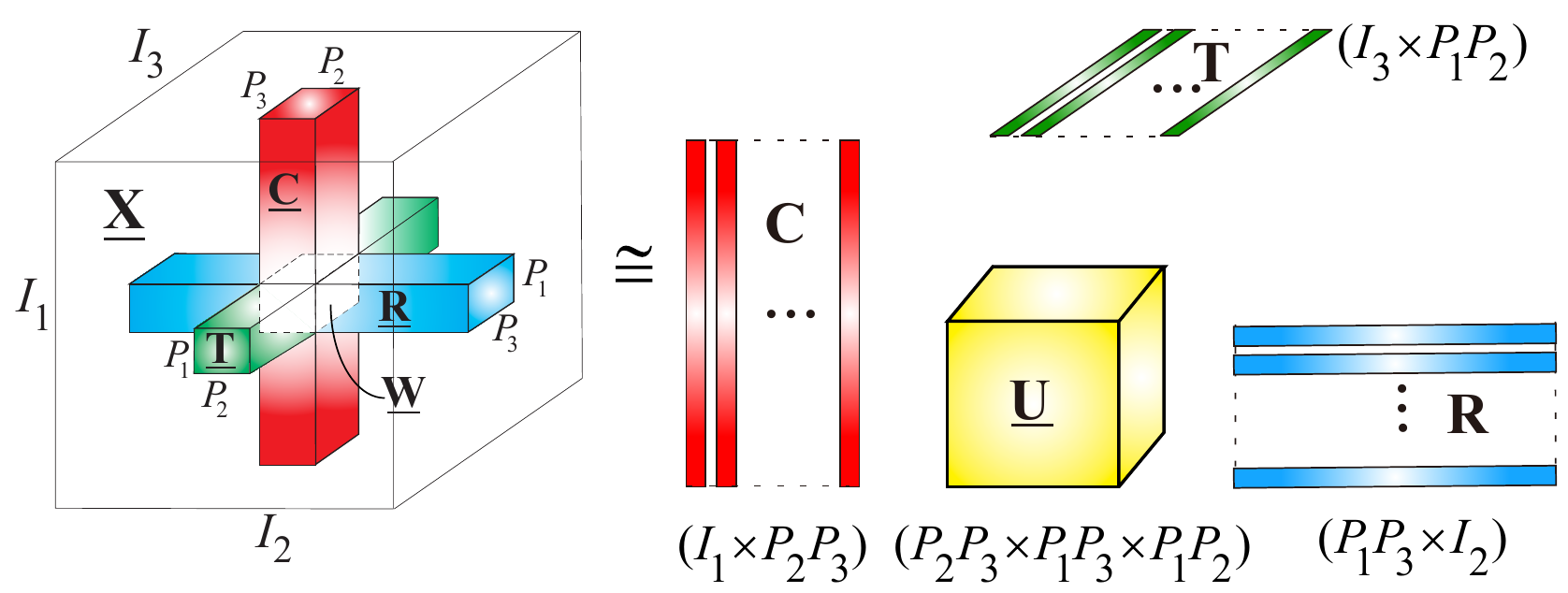}
\caption{The principle of the tensor cross-approximation (TCA) algorithm, illustrated for a large-scale  3rd-order tensor  $\underline \bX \cong \underline \bU \times_1 \bC \times_2 \bR \times_3 \bT
=\llbracket \underline \bU; \bC, \bR, \bT \rrbracket$, where $\underline \bU = \underline \bW \times_1\mathbf{W}_{(1)}^{\dagger} \times_2 \mathbf{W}_{(2)}^{\dagger} \times_3 \mathbf{W}_{(3)}^{\dagger} =\llbracket \underline \bW;\mathbf{W}_{(1)}^{\dagger}, \mathbf{W}_{(2)}^{\dagger},  \mathbf{W}_{(3)}^{\dagger}\rrbracket \in \Real^{P_2 P_3 \times P_1 P_3 \times P_1 P_2}$ and $\underline \bW \in \Real^{P_1 \times P_2 \times P_3}$. For simplicity of illustration, we assume that the selected fibers are permuted, so as to become  clustered  as subtensors, $\underline \bC \in \Real^{I_1 \times P_2 \times P_3}$, $\underline \bR \in \Real^{P_1 \times I_2 \times P_3}$ and $\underline \bT \in \Real^{P_1 \times P_2 \times I_3}$.}
\label{Fig:CUR-tensor}
\end{figure}

\begin{figure}
(a)
\begin{center}
\includegraphics[width=7.74 cm]{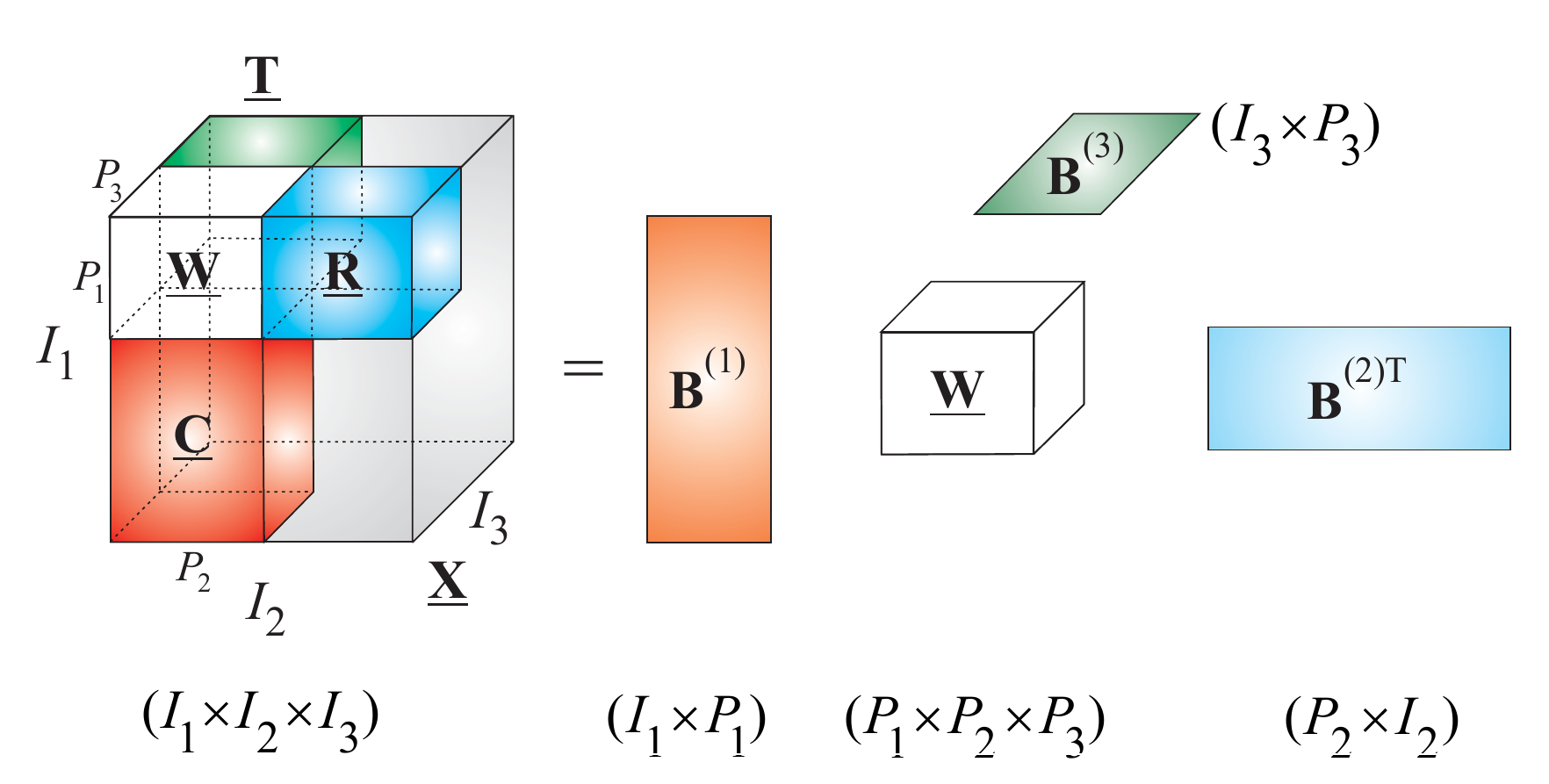}\\
\end{center}
(b)
\begin{center}
\includegraphics[width=7.29cm]{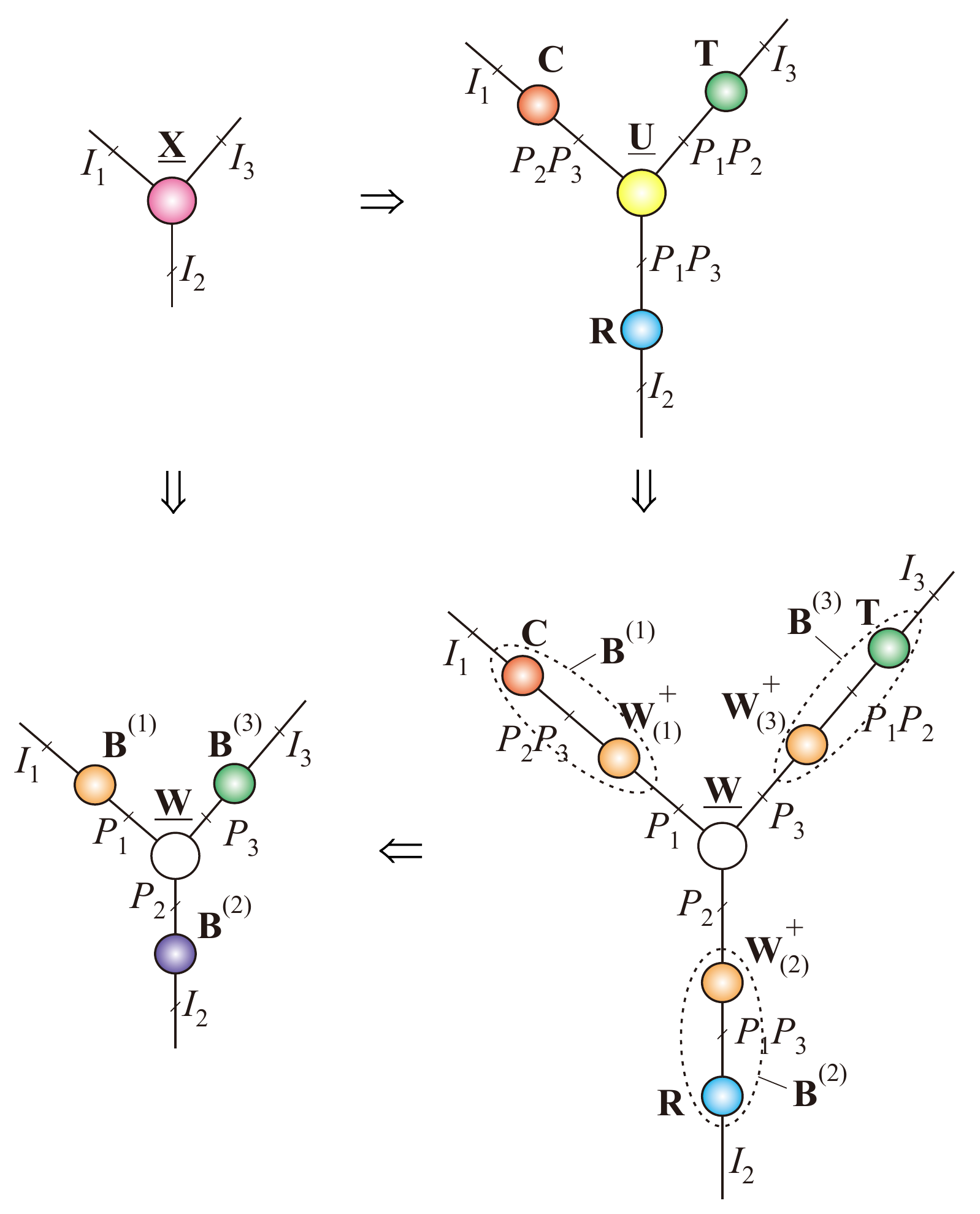}
\end{center}
\caption{The Tucker decomposition of a low multilinear rank 3rd-order tensor using the
cross-approximation approach.
(a) Standard block diagram.
 (b) Transformation from  the
 TCA  in the Tucker format, $\bX \cong \underline \bU \times_1 \bC \times_2 \bR \times_3 \bT$, into a standard Tucker representation, $\underline \bX \cong \underline \bW \times_1 \bB^{(1)} \times_2 \bB^{(2)} \times_3 \bB^{(3)}
=\llbracket \underline \bW; \bC \bW_{(1)}^{\dagger}, \bR \bW_{(2)}^{\dagger}, \bT  \bW_{(3)}^{\dagger}\rrbracket$, with a prescribed core tensor $\underline \bW$.}
\label{Fig:CUR-tensor2}
\end{figure}

A  simple modification of this formula, whereby the matrix $\bU$ is absorbed into either $\bC$ or $\bR$, yields the so-called CR matrix factorization or Column/Row Subset selection:
 \be
\bX \cong \bC  \tilde \bR = \tilde \bC \bR
\label{CUR2}
\ee
for which the  bases can be either the columns, $\bC$,  or rows, $\bR$, while $\tilde \bR = \bU \bR$ and $\tilde \bC= \bC \bU$.

For dimensionality reduction, $C \ll J$ and $R \ll I$,
and the  columns and rows of $\bX$ should be chosen optimally, in the sense of providing a high ``statistical leverage'' and  the best low-rank fit to the data matrix, while at the same time minimizing the cost function $\|\bE\|_F^2$.
For a given set of columns, $\bC$, and rows, $\bR$, the optimal choice for the core matrix is $\bU= \bC^{\dagger} \bX (\bR^{\dagger})^{\text{T}}$. This requires access to all the entries of  $\bX$ and is not practical or feasible for large-scale data. In such cases, a pragmatic choice for the core matrix would be $\bU =\bW^{\dagger}$, where the matrix  $\bW \in \Real^{R \times C}$  is composed from the intersections of the selected rows and columns. It should be noted that for rank$(\bX)\le  \min\{C,R\}$ the cross-approximation is exact.  For the general case, it has been proven that when the intersection submatrix $\bW$ is of maximum volume\footnote{The volume of a square submatrix $\bW$ is defined as $|\det(\bW)|$.}, the matrix cross-approximation is close to the optimal SVD solution. 
 The problem of finding  a submatrix with maximum volume  has exponential complexity, however, suboptimal matrices can be found using fast greedy algorithms
\cite{WangCUR13,OseledetsCA15,Oseledets-RakhubaCA15,MahoneyCUR15}.

The concept of MCA  can be generalized to tensor cross-approximation (TCA) (see Figure \ref{Fig:CUR-tensor}) through several approaches, including:
\begin{itemize}

\item Applying the MCA decomposition to a matricized version of the tensor data \cite{mahoney2008tensor};

\item Operating directly on fibers of  a data tensor which admits a low-rank Tucker approximation, an approach termed the Fiber Sampling Tucker Decomposition (FSTD) \cite{Caiafa-Cichocki-CUR,Caiafa2012-NC,caiafa2015stable}.
\end{itemize}

Real-life structured data often admit good low-multilinear rank approximations, and the FSTD provides such a low-rank Tucker decomposition which is practical as it is directly expressed in terms of a relatively small number of fibers  of the data tensor.


For example, for a  3rd-order tensor, $\underline \bX \in \Real^{I_1\times I_2 \times I_3}$, for which an exact rank-$(R_1,R_2,R_3)$ Tucker representation exists, the FSTD selects $P_n \geq R_n$, $n=1,2,3$, indices in each mode; this determines an intersection subtensor, $\underline \bW \in \Real ^{P_1\times P_2 \times P_3}$, so that the following exact Tucker representation can be obtained (see Figure \ref{Fig:CUR-tensor2})
\begin{equation}\label{NTensorCUR}
    \underline \bX = \llbracket \underline \bU; \bC, \bR, \bT \rrbracket,
 \end{equation}
where the core tensor is computed as $\underline \bU=\underline \bG=\llbracket \underline \bW;\mathbf{W}_{(1)}^{\dagger}, \mathbf{W}_{(2)}^{\dagger},  \mathbf{W}_{(3)}^{\dagger}\rrbracket$, while the factor matrices, $\bC \in \Real^{I_1 \times P_2 P_3}, \bR \in \Real^{I_2 \times P_1 P_3}, \bT \in \Real^{I_3 \times P_1 P_2}$, contain the  fibers which are the respective subsets of the columns $\underline \bC$, rows  $\underline \bR$ and tubes  $\underline \bT$. An equivalent Tucker representation is then given by
\begin{equation}\label{remark1}
 \underline \bX=\llbracket \underline \bW; \mathbf{C}\mathbf{W}_{(1)}^{\dagger},\mathbf{R}
 \mathbf{W}_{(2)}^{\dagger},\mathbf{T}\mathbf{W}_{(3)}^{\dagger} \rrbracket.
\end{equation}
Observe that for $N=2$,  the TCA  model simplifies into the MCA  for a matrix case, $\bX=\bC \bU \bR$,
for which the core matrix is $\mathbf{U}=\llbracket \mathbf{W}; \mathbf{W}_{(1)}^{\dagger}, \mathbf{W}_{(2)}^{\dagger}\rrbracket=\mathbf{W}^{\dagger}\mathbf{W}\mathbf{W}^{\dagger}=\mathbf{W}^{\dagger}$.


For a general case of an $N$th-order tensor, we can show  \cite{Caiafa-Cichocki-CUR} that
 a tensor, $\underline \bX \in \Real^{I_1 \times I_2 \times \cdots \times I_N}$, with a low multilinear rank $\{R_1,R_2,\ldots,R_N\}$, where $R_n \leq I_n, \, \forall n$,  can be fully reconstructed via the TCA FSTD, $\underline \bX = \llbracket \underline \bU; \bC^{(1)},\bC^{(2)},\ldots,\bC^{(N)}\rrbracket$, using only $N$  factor  matrices $\bC^{(n)} \in \Real^{I_n \times P_n}$ $(n=1,2,\ldots, N)$, built up from the fibers of the data and core tensors, $\underline \bU=\underline \bG=\llbracket \underline \bW; \mathbf{W}_{(1)}^{\dagger}, \mathbf{W}_{(2)}^{\dagger}, \ldots, \mathbf{W}_{(N)}^{\dagger}\rrbracket
$, under the condition that the subtensor $\underline \bW \in \Real^{P_1 \times P_2  \times \cdots \times P_N}$ with $P_n \geq R_n, \; \forall n$,  has the multilinear rank $\{R_1,R_2,\ldots,R_N\}$.

The selection of a minimum number of suitable fibers depends upon a chosen optimization criterion. A strategy which requires access to only a small subset of entries of a data tensor, achieved by selecting the entries with maximum modulus within each single fiber,  is given in \cite{Caiafa-Cichocki-CUR}. These entries are selected sequentially using a deflation approach, thus  making the tensor cross-approximation FSTD algorithm  suitable for the approximation of very large-scale but relatively low-order tensors (including tensors with missing fibers or entries).

It should be noted that an alternative efficient way to estimate subtensors $\underline \bW, \underline \bC, \underline \bR$ and $\underline \bT$  is to apply random projections as follows
\be
\underline \bW &=& \underline \bZ = \underline \bX \times_1 \mbi \Omega_1 \times_2 \mbi \Omega_2 \times_3 \mbi \Omega_3 \in \Real^{P_1 \times P_2 \times P_3}, \notag \\
\underline \bC &=& \underline \bZ_1 =\underline \bX  \times_2 \mbi \Omega_2 \times_3 \mbi \Omega_3
\in \Real^{I_1 \times P_2 \times P_3}, \notag \\
\underline \bR &=& \underline \bZ_2 = \underline \bX \times_1 \mbi \Omega_1 \times_3 \mbi \Omega_3
 \in \Real^{P_1 \times I_2 \times P_3}, \notag \\
\underline \bT &=& \underline \bZ_3 =\underline \bX \times_1 \mbi \Omega_1 \times_2 \mbi \Omega_2
 \in \Real^{P_1 \times P_2 \times I_3},
\ee
where $\mbi \Omega_n \in \Real^{P_n \times I_n}$ with $P_n \ge R_n$ for  $n=1,2,3$  are independent random matrices.
We explicitly assume that the multilinear rank $\{P_1,P_2, \ldots, P_N\}$ of  approximated tensor to be somewhat larger than a true multilinear rank $\{R_1,R_2, \ldots, R_N\}$
of target tensor, because it is easier to obtain an accurate approximation in this form.

\section{Multiway Component Analysis (MWCA)}
\label{sect:MBSS}

\subsection{Multilinear Component Analysis Using Constrained Tucker Decomposition}

The great success of 2-way component analyses (PCA, ICA, NMF, SCA) is
largely due to  the existence of very efficient algorithms for their computation and  the possibility to extract  components with a desired  physical meaning, provided by the various
flexible constraints exploited in these methods.
Without these constraints, matrix factorizations would be less useful in practice, as the components would have only mathematical
but not physical meaning.

\begin{figure}[t]
 \centering
 \includegraphics[width=11.95cm]{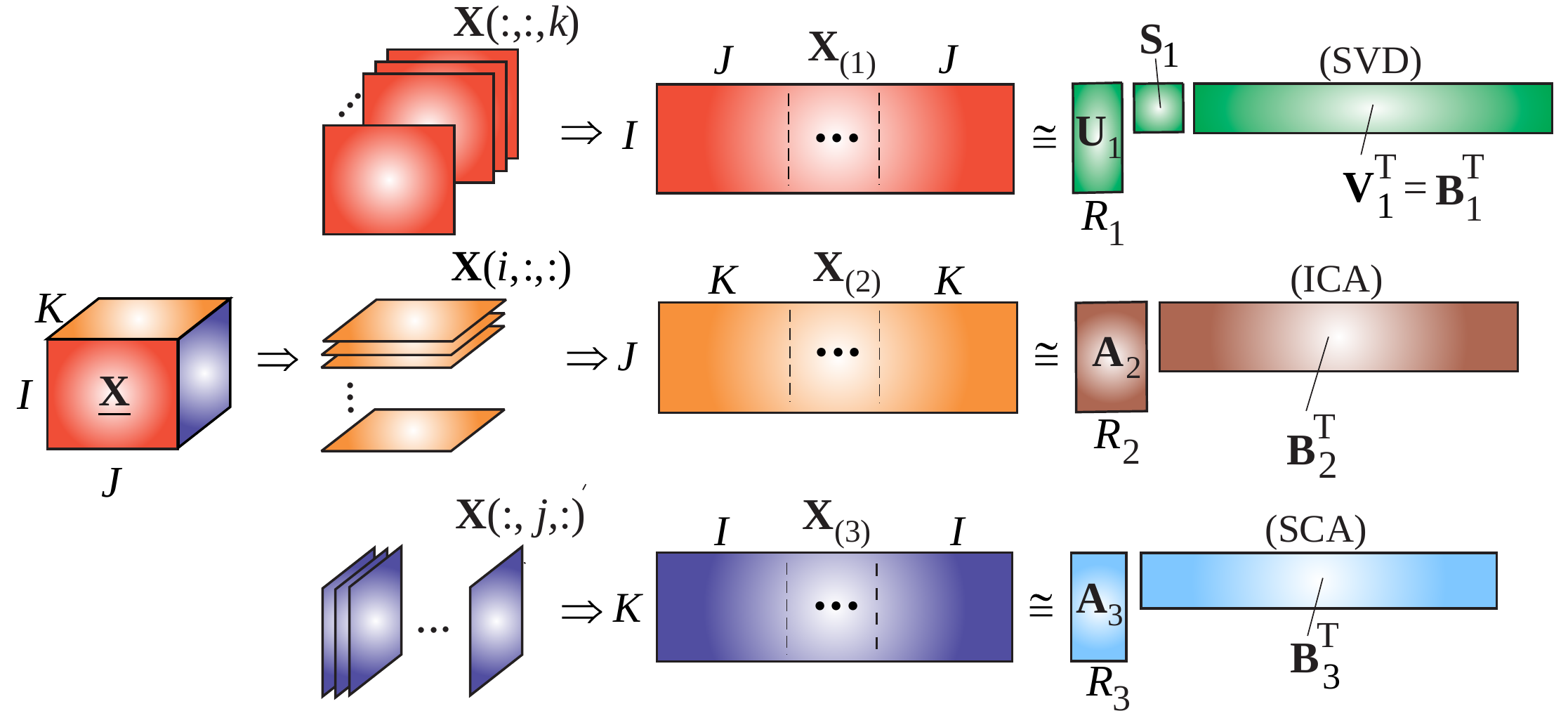}
 \caption{Multiway Component Analysis (MWCA) for a third-order tensor via constrained matrix factorizations, assuming that the components are:
 orthogonal in the first mode, statistically independent in the second mode and sparse
in the third mode.}
 \label{Fig:MWCA3}
\end{figure}

Similarly, to exploit the full potential of  tensor factorization/decompositions, it is a prerequisite  to impose suitable constraints on  the desired components.
In fact, there is much more flexibility for tensors, since different constraints can be imposed on the matrix factorizations in every mode $n$ a  matricized tensor $\bX_{(n)}$ (see Algorithm \ref{alg:Tucker-LRMF} and Figure \ref{Fig:MWCA3}).

Such physically meaningful representation through flexible mode-wise  constraints underpins  the concept of multiway component analysis (MWCA).
The  Tucker representation of MWCA naturally accommodates such diversities in different modes. 
Besides the orthogonality, alternative constraints in the Tucker format include statistical independence, sparsity, smoothness and nonnegativity \cite{Vasilescu,NMF-book,Zhou-Cichocki-MBSS,Cich-Lath} (see Table~\ref{tab:MWCA}).

The multiway component analysis (MWCA) based on the Tucker-$N$ model can be computed directly in two or three steps:
\begin{enumerate}
\item For each mode $n$ $(n=1,2, \ldots, N)$ perform model reduction and matricization of data tensors sequentially, then apply a  suitable set of  2-way CA/BSS algorithms to the so reduced unfolding matrices, $\tilde \bX_{(n)}$. In each mode, we can  apply different  constraints and a different 2-way CA algorithms.

\item  Compute the core tensor using, e.g., the inversion formula,
 $ \hat{\underline \bG} = \underline \bX  \times_{1}\matn[1]{B}{}^\dagger
 \times_{2}\matn[2]{B}{}^\dagger
  \cdots \times_{N}\matn[N]{B}{}^\dagger$. This step is quite important because core tensors often model the  complex links among the multiple components in different modes.

\item Optionally,  perform fine tuning of factor matrices and the core tensor by the ALS minimization of a suitable cost function, e.g., $\|\underline \bX - \llbracket \underline \bG;\bB^{(1)}, \ldots, \bB^{(N)}\rrbracket\|^2_F$, subject to specific  imposed constraints.
\end{enumerate}

\section{Analysis of Coupled Multi-block Matrix/Tensors --  Linked Multiway Component Analysis (LMWCA)}
\label{sect:LMCA}

We have shown that TDs provide natural extensions of  blind source separation (BSS)
and 2-way (matrix) Component Analysis  to multi-way component analysis (MWCA) methods.

\begin{figure}[t!]
(a)
\begin{center}
 \includegraphics[width=5.9cm]{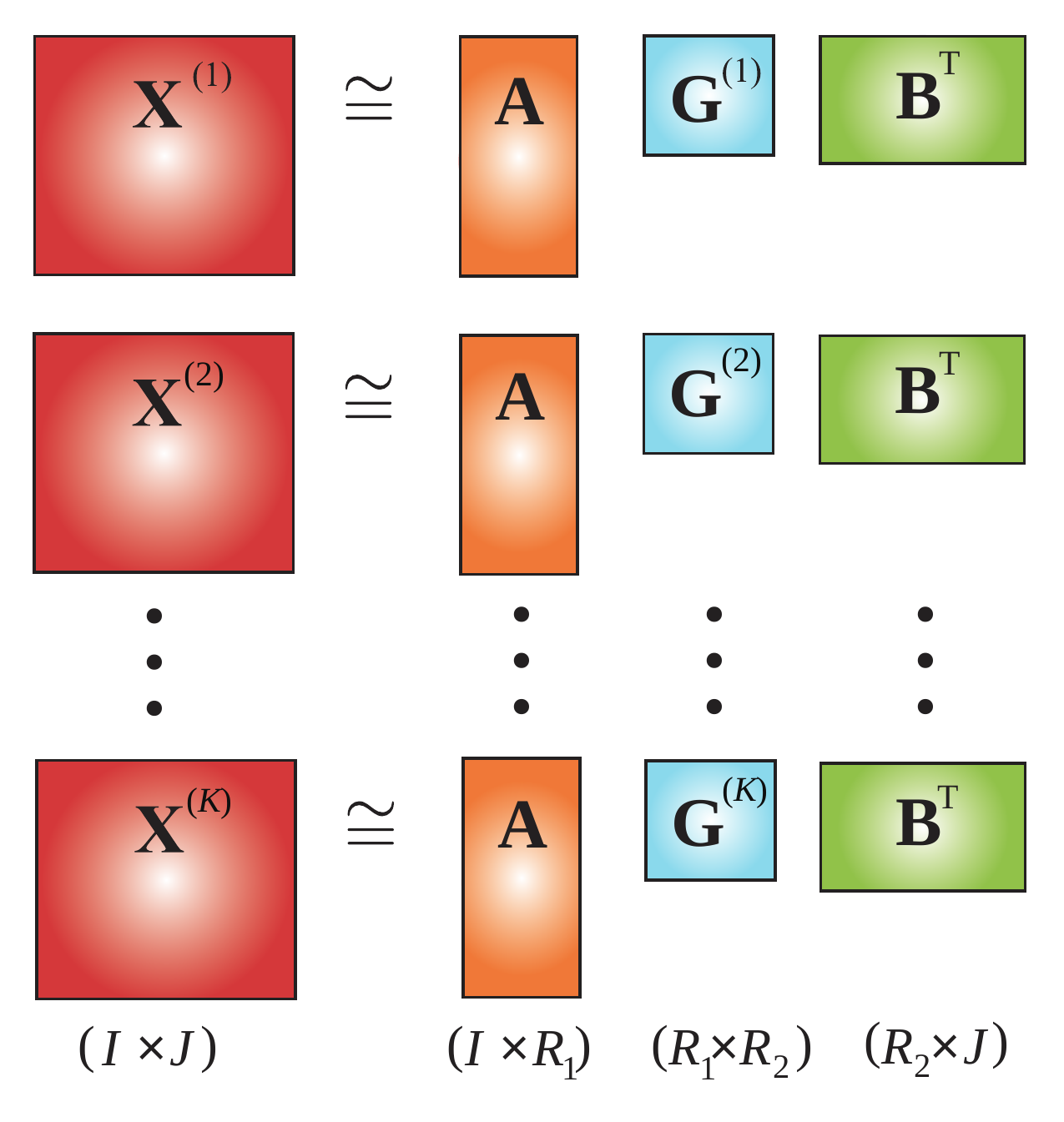}
 \end{center}
 (b)
 \begin{center}
 \includegraphics[width=7.5cm]{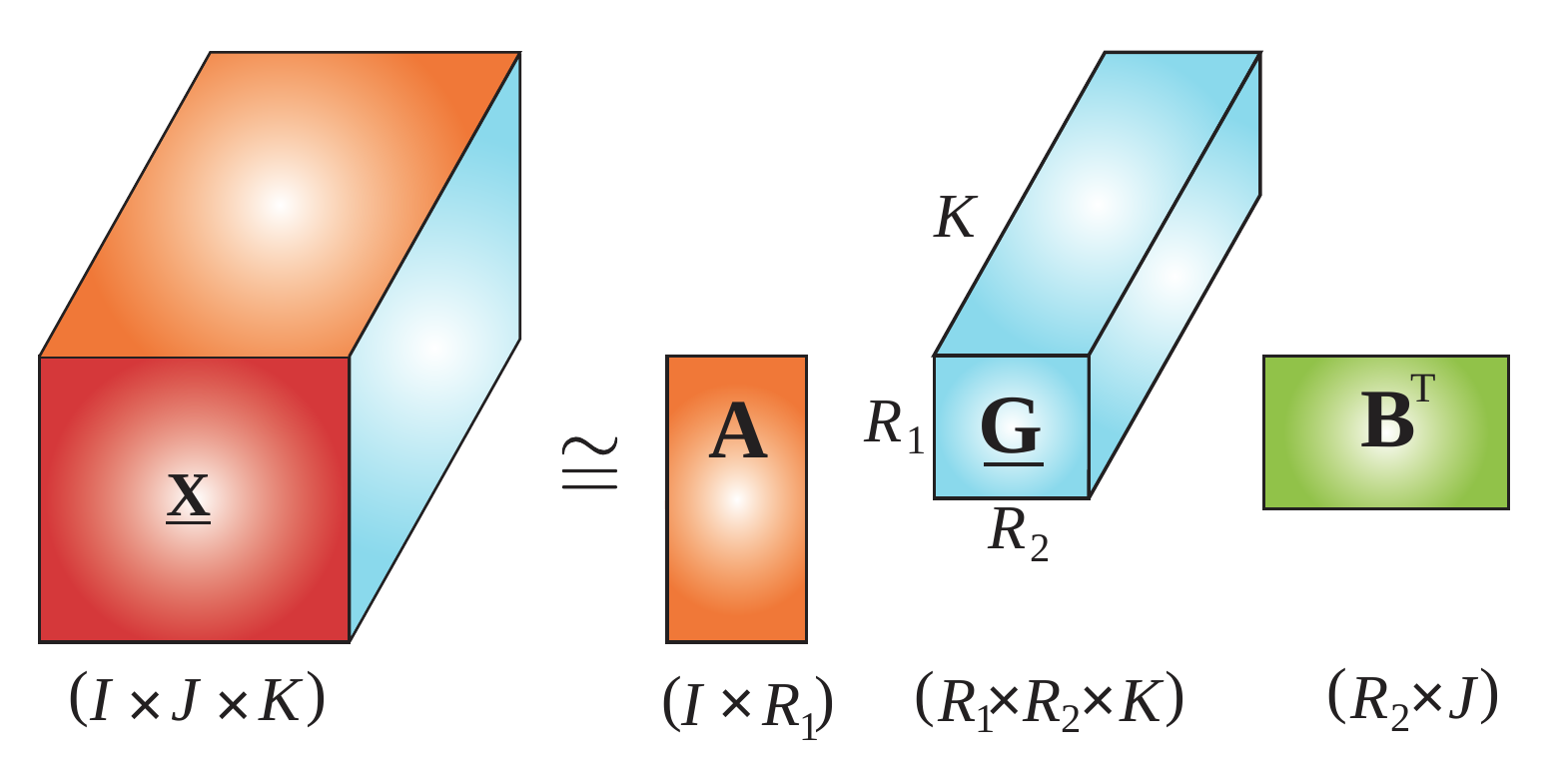}
 \end{center}
 \caption{Concept of the Population Value Decomposition (PVD). (a) Principle of simultaneous multi-block matrix factorizations. (b)  Equivalent representation of the PVD as the  constrained or unconstrained Tucker-2 decomposition, $\underline \bX \cong \underline \bG \times_1 \bA \times_2 \bB$.
  The objective is to find the common factor matrices,  $\bA$, $\bB$ and the core tensor, $\underline \bG \in \Real^{R_1 \times R_2 \times K}$.}
 \label{Fig:PVD}
\end{figure}

 In addition, TDs are suitable  for the coupled multiway analysis of multi-block datasets,
 possibly with missing values and corrupted by noise.
 To illustrate the simplest scenario for multi-block analysis, consider the block matrices, $\bX^{(k)} \in \Real^{I \times J}$,
 which  need to be approximately jointly factorized as
 \be
 \bX^{(k)} \cong \bA \bG^{(k)} \bB^{\text{T}}, \quad (k=1,2, \ldots,K),
 \ee
 where   $\bA \in \Real^{I \times R_1}$ and $\bB \in \Real^{J \times R_2}$ are common factor matrices and $\bG^{(k)} \in \Real^{R_1 \times R_2}$
 are reduced-size matrices, while the number of data matrices $K$ can be huge (hundreds of millions or more  matrices).
 Such a simple model is referred to as the Population Value Decomposition (PVD) \cite{PVD2011}. Note that the PVD is equivalent to the
  unconstrained or constrained Tucker-2 model, as illustrated in Figure  \ref{Fig:PVD}. In a special case  with  square diagonal matrices, $\bG^{(k)}$, the model is equivalent to the CP decomposition and is related to joint matrix diagonalization \cite{deLathauwer-JAD,WASOBI,chabriel2014JAD}.
  Furthermore, if $\bA = \bB$ then the PVD model is equivalent to  the
  RESCAL model \cite{nickel2016review}.

  Observe that the PVD/Tucker-2 model is quite general and flexible, since any  high-order  tensor,
  $\underline \bX \in \Real^{I_1 \times I_2 \times \cdots \times I_N}$ (with $N>3$), can be  reshaped and optionally permuted into a  ``skinny and tall''
  3rd-order tensor, $\underline {\widetilde \bX} \in \Real^{J \; \times \; J\; \times \; K}$, with e.g., $I=I_1, \; J=I_2$ and $K=  I_3 I_4 \cdots I_N$, for which PVD/Tucker-2 Algorithm \ref{alg:Tucker2} can be applied.

   As previously mentioned, various constraints, including  sparsity, nonnegativity or smoothness
   can be imposed on the factor matrices,  $\bA$ and $\bB$, to obtain physically meaningful and unique components.

\begin{algorithm}[t!]
\caption{\textbf{Population Value Decomposition (PVD) with orthogonality constraints}}
\label{alg:PVD}
{\small
\begin{algorithmic}[1]
\REQUIRE A set  of matrices $\bX_k \in \Real^{I \times J}$, for $k=1, \ldots,K$ (typically,  $\;K \gg \max\{I,J\}$)
\ENSURE Factor matrices $\bA \in \Real^{I \times R_1}$, $\bB \in \Real^{J \times R_2}$  and $\bG_k \in \Real^{R_1 \times R_2}$,
  with~orthogonality constraints $\bA^{\text{T}} \bA =\bI_{R_1}$ and
$\bB^{\text{T}} \bB =\bI_{R_2}$
\FOR {$k=1$ to $K$}
\STATE Perform truncated SVD, $\bX_k = \bU_k \bS_k \bV_k^{\text{T}}$, using $R$ largest singular \\values
\ENDFOR
\STATE Construct short and wide matrices:\\ $\bU =[\bU_1 \bS_1, \ldots, \bU_K \bS_K] \in \Real^{I \times K R}$ and
$\bV =[\bV_1 \bS_1, \ldots, \bV_K \bS_K] \in \Real^{J \times K R}$
\STATE Perform SVD  (or QR) for the  matrices $\bU$ and $\bV$ \\
Obtain common  orthogonal matrices  $\bA$ and $\bB$ as left-singular \\
matrices of $\bU$ and $\bV$, respectively
\FOR {$k=1$ to $K$}
\STATE Compute $\bG_k = \bA^{\text{T}} \bX_k \bB$
\ENDFOR
\end{algorithmic}
}
\end{algorithm}
\begin{algorithm}[t!]
\caption{\textbf{Orthogonal Tucker-2 decomposition with a prescribed approximation accuracy \cite{Phan_TT_part1}}}
\label{alg:Tucker2}
{\small
\begin{algorithmic}[1]
\REQUIRE A 3rd-order tensor  $\underline \bX \in \Real^{I \times J \times K}$ (typically, $K \gg \max \{I,J\}$) \\and  estimation accuracy $\varepsilon$
\ENSURE A set of orthogonal matrices $\bA \in \Real^{I \times R_1}$, $\bB \in \Real^{J \times R_2}$  and core tensor
$\underline \bG \in \Real^{R_1 \times R_2 \times K}$,
  which satisfies the constraint $\|\underline \bX - \underline \bG \times_1 \bA \times \bB\|_F^2 \leq \varepsilon^2$ , s.t, $\bA^{\text{T}} \bA =\bI_{R_1}$ and
$\bB^{\text{T}} \bB =\bI_{R_2}$.
\STATE Initialize $\bA=\bI_{I} \in \Real^{I \times I}$, $R_1=I$
\WHILE {not converged or iteration limit is not reached}
    \STATE Compute the tensor $\underline \bZ^{(1)} = \underline \bX \times_1 \bA^{\text{T}}
    \in \Real^{R_1 \times J \times K}$
      \STATE Compute EVD of a small matrix $\bQ_1 =\bZ_{(2)}^{(1)} \bZ_{(2)}^{(1)\;\text{T}}
      \in \Real^{J \times J}$ as
       $\bQ_1 =\bB \; \diag \left(\lambda_1,\cdots, \lambda_{R_2}\right) \; \bB^{\text{T}} $, such that $\sum_{r_2=1}^{R_2} \lambda_{r_2}
    \geq \|\underline \bX\|_F^2 -\varepsilon^2 \geq \sum_{r_2=1}^{R_2-1} \lambda_{r_2}$
     \STATE Compute tensor $\underline \bZ^{(2)} = \underline \bX \times_2 \bB^{\text{T}}
     \in \Real^{I \times R_2 \times K}$
      \STATE Compute EVD of a small matrix $\bQ_2 =\bZ_{(1)}^{(2)} \bZ_{(1)}^{(2)\;\text{T}}
      \in \Real^{I \times I}$ as
    $\bQ_2 =\bA \; \diag\left(\lambda_1,\ldots, \lambda_{R_1}\right) \; \bA^{\text{T}} $, such that $\sum_{r_1=1}^{R_1} \lambda_{r_1}
    \geq \|\underline \bX\|_F^2 -\varepsilon^2 \geq \sum_{r_1=1}^{R_1-1} \lambda_{r_1}$
\ENDWHILE
\STATE Compute the core tensor $\underline \bG = \underline \bX \times_1 \bA^{\text{T}} \times_2 \bB^{\text{T}}$
\RETURN  $\bA, \bB$ and $ \underline \bG$.
\end{algorithmic}
}
\end{algorithm}

 A simple SVD/QR based algorithm for the PVD with orthogonality constraints  is presented in Algorithm \ref{alg:PVD} \cite{PVD2011,TallSVD2014,PVD2016}.
 However, it should be noted that this algorithm does not provide an optimal solution
 in the sense of the absolute minimum of the cost function, $\sum_{k=1}^K \|\bX_{k} -\bA \bG_k \bB^{\text{T}}\|_F^2$, and for data corrupted by Gaussian noise, better performance can be achieved using the HOOI-2 given in  Algorithm \ref{alg:HOOI},  for $N=3$. An improved PVD algorithm referred to as Tucker-2 algorithm is given in Algorithm \ref{alg:Tucker2} \cite{Phan_TT_part1}.

\begin{figure}[t!]
 \centering
 \includegraphics[width=9.9cm]{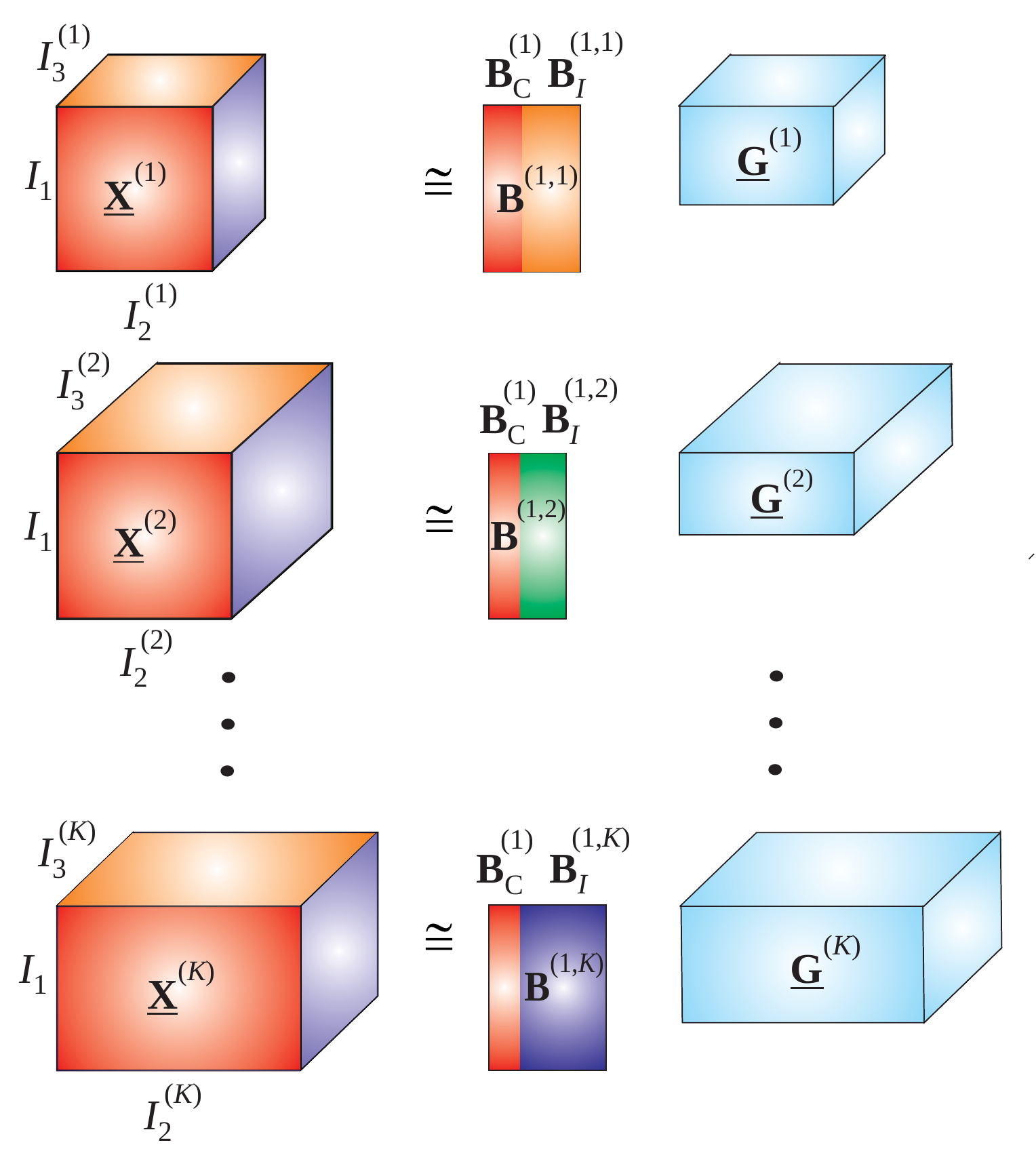}
 \caption{Linked Multiway Component Analysis (LMWCA) for
  coupled  3rd-order data tensors $\underline \bX^{(1)},\ldots, \underline \bX^{(K)}$; these can have different dimensions in every mode, except for the mode-1 for which the size is $I_1$ for all $\underline \bX^{(k)}$.
  Linked Tucker-1 decompositions are then performed in the form $\underline \bX^{(k)} \cong \underline \bG^{(k)} \times_1 \bB^{(1,k)}$, where partially correlated  factor matrices are $\bB^{(1,k)} = [\bB^{(1)}_C, \bB^{(1,k)}_I] \in \Real^{I_1 \times R_k}$, $(k=1,2,\ldots,K)$. The objective is to find the common components, $\bB^{(1)}_C \in \Real^{I_1 \times C}$, and individual components, $\bB^{(1,k)}_I \in \Real^{I_1 \times (R_k-C)}$, where $C \leq \min\{R_1,\ldots,R_K\}$ is the number of  common components in mode-1.}
 \label{Fig:LMCA}
\end{figure}

\noindent \textbf{Linked MWCA.} Consider the analysis of multi-modal high-dimensional data collected under the same or very similar conditions, for example, a set of EEG and MEG or EEG and fMRI signals recorded for different subjects over many trials and under the same experimental configurations and mental tasks. Such data share
some common latent (hidden) components but can also have their own independent features. As a result, it is advantageous and natural to
analyze such data in a linked way instead of treating them independently.
%
%
In such a scenario, the PVD  model can be  generalized to multi-block matrix/tensor datasets \cite{Cichocki-SICE,Zhou-PAMI,Zhou-PIEEE}.

\begin{figure}
(a)
\begin{center}
  \includegraphics[width=11.2cm]{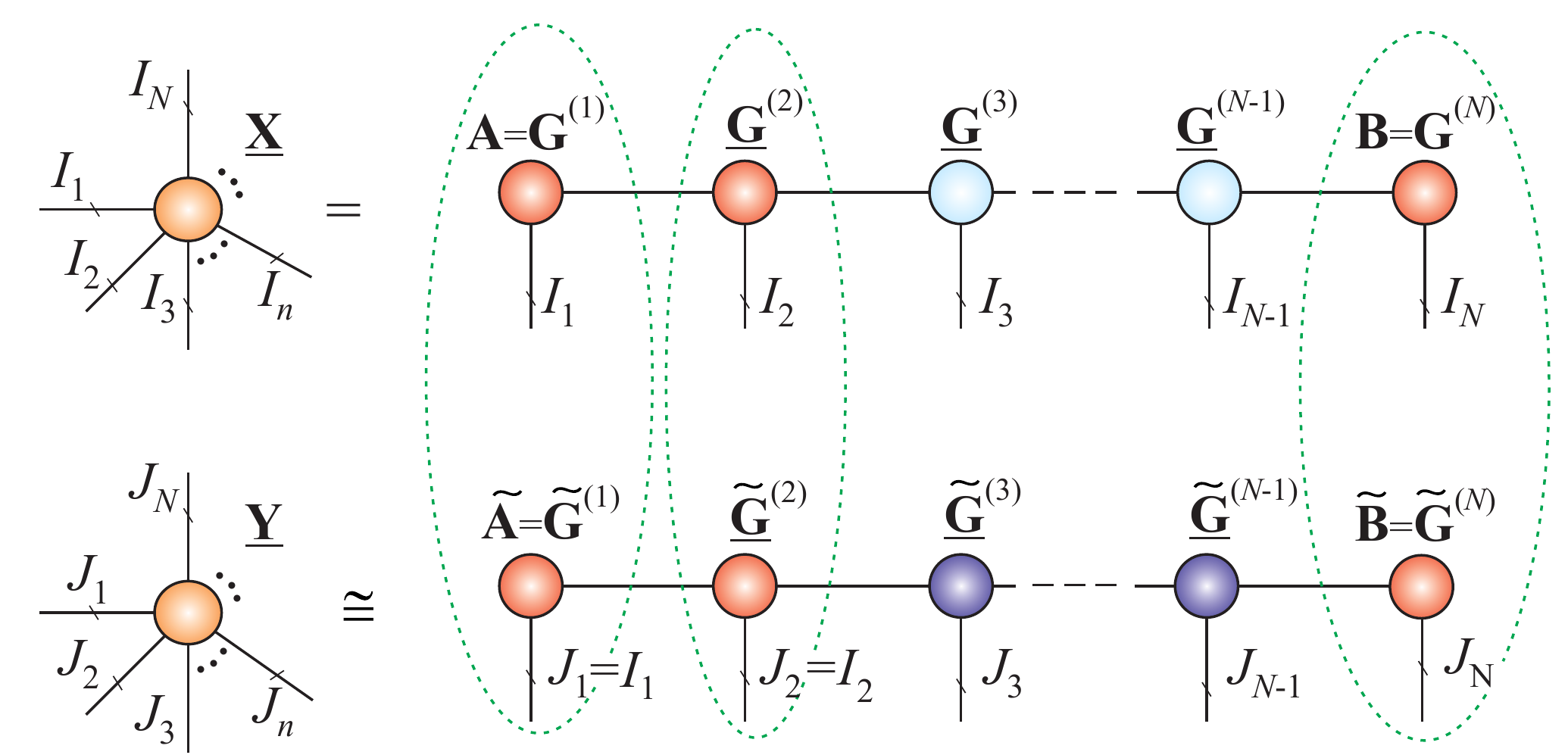}
\end{center}
\vspace{0.8cm}
(b)
\begin{center}
 \includegraphics[width=9.0cm]{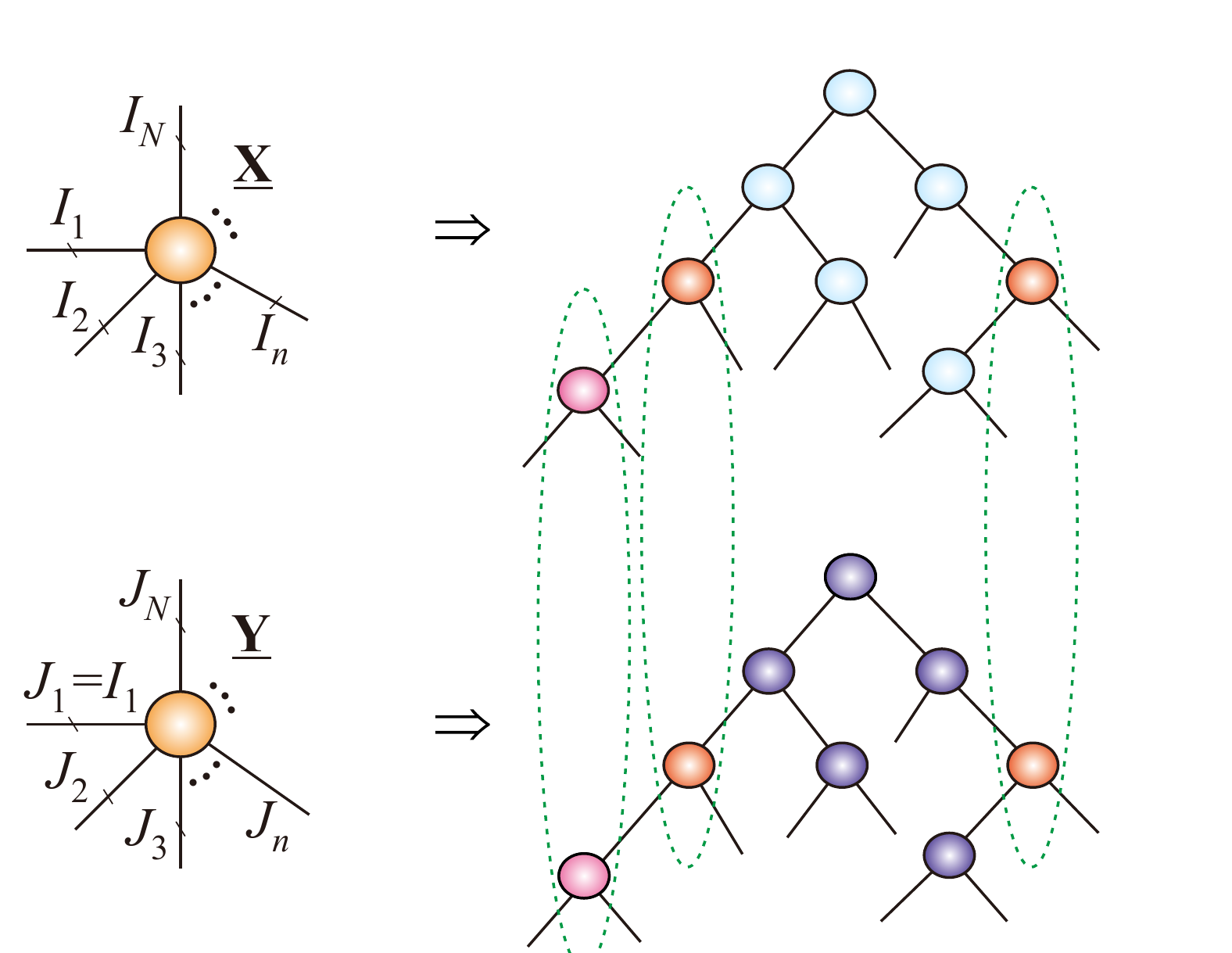}
 \end{center}
 \caption{Conceptual models of generalized Linked  Multiway Component Analysis (LMWCA) applied to  the cores of high-order TNs. The objective
 is to find  a  suitable tensor decomposition which yields the maximum number of cores that are as much  correlated as possible. (a) Linked Tensor Train (TT) networks.
(b) Linked Hierarchical Tucker (HT) networks with the correlated cores indicated by ellipses  in broken lines.}
 \label{Fig:LMCA2}
\end{figure}

The linked multiway component analysis (LMWCA) for multi-block tensor data can therefore be
formulated as a set of approximate simultaneous (joint) Tucker-$(1,N)$ decompositions of a set of data tensors,
$ \underline \bX^{(k)} \in  \Real^{I_1^{(k)} \times I_2^{(k)} \times \cdots \times I_N^{(k)}}$, with $I_1^{(k)}=I_1$ for $k=1,2,\ldots,K$, in the form (see Figure   \ref{Fig:LMCA})
\be
   \underline \bX^{(k)} = \underline \bG^{(k)}  \times_{1} \bB^{(1,k)}, \quad (k=1,2,\ldots K)
 \label{LMBSS1}
\ee
where each factor (component) matrix, $\bB^{(1,k)}=[\bB^{(1)}_C, \; \bB^{(1,k)}_I]
\in \Real^{I_1 \times R_k}$, comprises  two sets of components: (1) Components $\bB^{(1)}_C \in \Real^{I_1 \times C} $ (with $0 \leq C \leq R_k$), $\;\forall k$,
  which are common for  all the available blocks and correspond to  identical or maximally correlated  components, and (2) components
  $\bB^{(1,k)}_I \in \Real^{I_1 \times (R_k-C)}$, which are  different independent processes for each block, $k$, these can be, for example, latent variables
   independent of excitations or stimuli/tasks. The objective is therefore to estimate the common (strongly correlated) components, $\bB_C^{(1)}$,
   and statistically independent (individual) components, $\bB^{(1,k)}_I $ \cite{Cichocki-SICE}.

%
 If $\bB^{(n,k)}= \bB^{(n)}_C \in \Real^{I_n \times R_n}$ for a specific mode $n$ (in our case $n=1$), and under the additional assumption that the block tensors are of the same order and size,  the problem simplifies into generalized Common Component Analysis or tensor Population Value Decomposition (PVD)  and can be solved by concatenating all data tensors along one mode, followed by constrained Tucker
 or CP decompositions \cite{Phan2010TF}.

In a more general scenario, when $C_n<R_n$,  we can  unfold  each data tensor $\underline \bX^{(k)}$ in the common mode, and perform a set of simultaneous
 matrix factorizations, e.g., $\bX_{(1)}^{(k)}\cong \bB^{(1)}_C\matn[1,k]{A}_C + \matn[1,k]{B}_I\matn[1,k]{A}_I$, through solving the constrained optimization problems
 \begin{equation}
   \label{LWCA}
   \begin{split}
   \min \; & \sum_{k=1}^{K}\|\bX_{(1)}^{(k)}-\bB^{(1)}_C\matn[1,k]{A}_C-\matn[1,k]{B}_I\matn[1,k]{A}_I\|_F \\
    & + P(\bB^{(1)}_C), \;\;
    s.t. \;\; \bB_C^{(1)\;\text{T}}  \matn[1,k]{B}_I=\mathbf{0} \; \forall k,
   \end{split}
 \end{equation}
 where the symbol $P$ denotes the penalty terms which impose additional constraints on the common components, $\bB_C^{(1)}$, in order to extract as many common components as possible. In the special case of orthogonality constraints, the problem can be transformed into a generalized eigenvalue problem.
  The key point is to assume that common factor submatrices, $\bB_C^{(1)}$, are present in all data blocks and hence reflect structurally complex  latent (hidden) and intrinsic links between the data blocks. \emph{In practice, the number of common components,  $C$,  is unknown and should be estimated} \cite{Zhou-PAMI}.

The linked multiway component analysis (LMWCA) model complements currently available  techniques for group
 component analysis and  feature extraction from multi-block datasets,  and is a natural extension of  group ICA, PVD, and  CCA/PLS methods (see \cite{Cichocki-SICE,Qibin-HOPLS,Zhou-PAMI,Zhou-PIEEE} and references therein).
%
Moreover, the concept of LMWCA can be generalized to tensor networks, as illustrated in Figure  \ref{Fig:LMCA2}.

\section{Nonlinear Tensor Decompositions -- Infinite Tucker}

The Infinite Tucker model and  its modification, the
Distributed Infinite Tucker (DinTucker), generalize the standard Tucker decomposition to infinitely dimensional feature spaces using  kernel and Bayesian approaches \cite{Ta-2013,InfTucker,DINTucker}.

 Consider the classic Tucker-$N$ model of an $N$th-order tensor $\underline \bX \in \Real^{I_1 \times \cdots \times I_N}$, given by
\be
\underline \bX &=& \underline \bG \times_1 \bB^{(1)} \times_2 \bB^{(2)} \cdots \times_N \bB^{(N)} \notag \\
&=& \llbracket  \underline \bG; \bB^{(1)}, \bB^{(2)}, \ldots, \bB^{(N)} \rrbracket
\ee
in its vectorized version
\be
\mbox{vec}(\underline \bX) =  (\bB^{(1)} \otimes_L \cdots \otimes_L \bB^{(N)}) \; \mbox{vec}(\underline \bG).\notag
\ee
Furthermore, assume that the noisy data tensor is modeled as
\be \underline \bY = \underline \bX  + \underline \bE,
\ee
 where $\underline \bE$ represents the tensor of additive Gaussian noise.
Using the Bayesian framework and tensor-variate Gaussian processes (TGP)  for
Tucker decomposition, a standard  normal prior can be assigned over each entry,  $g_{r_1,r_2, \ldots , r_N}$, of an $N$th-order core tensor,
$\underline \bG \in \Real^{R_1 \times \cdots \times R_N}$, in order to marginalize out $\underline \bG$ and express the probability density function of tensor $\underline \bX$ \cite{chu2009probabilistic,InfTucker,DINTucker} in the form
\be
&& p \left(\underline \bX \,| \, \bB^{(1)}, \ldots, \bB^{(N)}\right) = {\cal{N}} \left(\mbox{vec}(\underline \bX); \0,
\bC^{(1)} \otimes_L \cdots \otimes_L \bC^{(N)}\right) \notag \\
&& = \frac{\exp\left( -\frac{1}{2} \| \llbracket  \underline \bX; (\bC^{(1)})^{-1/2}, \ldots, (\bC^{(N)})^{-1/2} \rrbracket \|_F^2 \right)}
{(2 \pi)^{I/2} \prod_{n=1}^N |\bC^{(n)}|^{-I/(2I_n)}}
\label{eq:InfTucker}
\ee
where $I= \prod_n I_n$ and $\bC^{(n)} = \bB^{(n)} \; \bB^{(n)\;\text{T}} \in \Real^{I_n \times I_n} $ for $n=1,2 ,\ldots, N$.

In order to model unknown, complex, and potentially nonlinear interactions between the latent factors, each row,
 $\bar \bb_{i_n}^{(n)} \in \Real^{1 \times R_n}$, within $\bB^{(n)}$, is replaced by a nonlinear feature transformation $\Phi ( \bar \bb_{i_n}^{(n)})$
 using the kernel trick \cite{zhao2013kernelization}, whereby the nonlinear covariance matrix
 $\bC^{(n)} = k( \bB^{(n)}, \; \bB^{(n)}) $ replaces the standard covariance matrix,
  $\bB^{(n)}  \bB^{(n)\;\text{T}}$.
Using such a nonlinear feature mapping, the original Tucker factorization is performed in an infinite feature space, while Eq.~(\ref{eq:InfTucker})
defines a Gaussian process (GP) on a tensor, called the Tensor-variate GP (TGP), where the inputs come from a set of factor matrices
$\{\bB^{(1)}, \ldots, \bB^{(N)}\}=\{\bB^{(n)}\}$.

For a noisy data tensor $\underline \bY$,  the joint probability density function is given by
\be
 p ( \underline \bY, \underline \bX, \{ \bB^{(n)} \}) =  p(\{\bB^{(n)}\}) \; p(\underline \bX \, | \, \{\bB^{(n)}\}) \; p(\underline \bY | \underline \bX ).
 \ee
To improve scalability,  the observed noisy tensor $\underline \bY$  can be split into   $K$  subtensors
$\{\underline \bY_1, \ldots, \underline \bY_K\}$, whereby each subtensor $\underline \bY_k$  is sampled from its own GP based  model with factor  matrices,
$\{\tilde \bB^{(n)}_k\} =\{ \tilde \bB^{(1)}_k, \ldots, \tilde \bB^{(N)}_k\}$. The factor matrices can then be merged via a prior distribution
\be
 p ( \{ \tilde \bB_k^{(n)} \} | \{ \bB^{(n)} \}) &=& \prod_{n=1}^N  p(\tilde \bB_k^{(n)} | \bB^{(n)}) \notag \\
&=& \prod_{n=1}^N {\cal{N}} (\mbox{vec} (\tilde \bB_k^{(n)}) | \mbox{vec}(\bB^{(n)})), \lambda \bI),
 \ee
 where $\lambda >0 $ is a variance parameter which controls the similarity between the corresponding factor matrices.
 The above model is referred to as DinTucker  \cite{DINTucker}.

The full covariance matrix, $\bC^{(1)} \otimes \cdots \otimes \bC^{(N)} \in \Real^{\prod_n I_n \times \prod_n I_n}$,
may have a prohibitively large size  and can be extremely sparse. For such cases, an
alternative nonlinear tensor decomposition model has been recently developed, which does not, either explicitly or implicitly, exploit the Kronecker structure of covariance matrices \cite{cichocki2015log}.
Within this model, for each tensor entry, $x_{i_1,\ldots,i_N}=x_{\bi}$,
with $\bi= (i_1,i_2,\ldots,i_N)$,  an input vector $\bb_{\bi}$  is constructed by concatenating the corresponding row vectors
of factor (latent) matrices, $\bB^{(n)}$, for all $N$ modes, as
\be
\bb_{\bi} =[\bar \bb^{(1)}_{i_1}, \ldots, \bar \bb^{(N)}_{i_N}] \in \Real^{1 \times \sum_{n=1}^N R_n}.
\ee
We can formalize an (unknown) nonlinear transformation as
\be
x_{\bi}= f(\bb_{\bi}) =f([\bar \bb^{(1)}_{i_1}, \ldots, \bar \bb^{(N)}_{i_N}])
\ee
for which a zero-mean multivariate Gaussian distribution is
determined by
$\bB_{\cal {S}} = \{\bb_{{\bi}_1}, \ldots, \bb_{{\bi}_M}\}$ and  $\boldf_{\cal {S}} =
\{f(\bb_{{\bi}_1}), \ldots, f(\bb_{{\bi}_M})\}$. This allows us to construct the following
probability function
\be
p \left(\boldf_{\cal {S}} | \{\bB^{(n)}\} \right) = {\cal {N}} \left( \boldf_{\cal {S}} | \0, \; k (\bB_{\cal {S}},\; \bB_{\cal {S}}) \right),
\label{eq:InfTucker4}
\ee
where $k(\cdot,\cdot)$ is a nonlinear covariance function  which can be expressed as $k(\bb_{\bi},\bb_{\bj}) = k(([\bar \bb^{(1)}_{i_1}, \ldots, \bar \bb^{(N)}_{i_N}]),([\bar \bb^{(1)}_{j_1}, \ldots, \bar \bb^{(N)}_{j_N}]))$ and ${\cal{S}} = [\bi_1,\ldots,\bi_M]$.

In order to assign a standard normal prior over the factor matrices,  $\{\bB^{(n)}\}$,
we assume that for selected entries, $\bx = [x_{\bi_1},\ldots, x_{\bi_M}]$, of a tensor $\underline \bX$,  the  noisy entries, $\by=[y_{\bi_1},\ldots, y_{\bi_M}]$, of the observed tensor  $\underline \bY$, are sampled from the following joint probability model
\be
&& p (\by, \bx, \{\bB^{(n)}\})  \\
 &&\quad= \prod_{n=1}^N
  {\cal{N}} (\mbox{vec} (\bB^{(n)}) | \0, \bI ) \; {\cal{N}} (\bx | \0, k(\bB_{\cal {S}},\; \bB_{\cal {S}})) \; {\cal{N}}(\by | \bx, \beta^{-1} \bI), \notag
 \ee
where $\beta$ represents noise variance.

These nonlinear and  probabilistic models can be potentially applied for  data tensors or function-related tensors comprising large number of entries,
typically with millions of non-zero entries
and billions of zero entries. Even if  only nonzero entries are used, exact inference of the
 above  nonlinear tensor decomposition models may still be    intractable.
To alleviate this problem, a  distributed variational  inference algorithm  has been developed,
 which  is based on sparse GP, together with an efficient
MapReduce framework which uses a small set of inducing
points to break up the dependencies between random function
values \cite{DINTucker,Titsias2009}.

\chapter{Tensor Train Decompositions: Graphical Interpretations and Algorithms}
\chaptermark{Tensor Train Decompositions}
\label{chap:TT}

\vspace{0cm}

Efficient implementation of the various operations in tensor train
(TT) formats requires compact
and easy-to-understand  mathematical and graphical representations
\cite{CichockiSISA,Cichocki2014optim}.
To this end, we next present mathematical formulations of the TT decompositions and demonstrate  their advantages in both theoretical and practical scenarios.


\section{Tensor Train Decomposition -- Matrix Product State}
\label{sect:TT}

\begin{figure} 
(a)\\
\vspace{-0.1cm}
\begin{center}
\includegraphics[width=11.99cm]{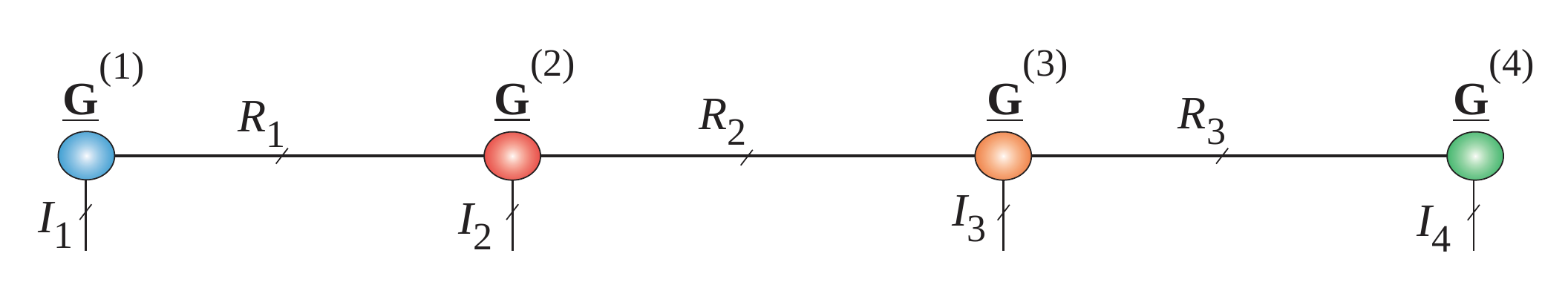}\\
\end{center}
\vspace{-0.1cm}
\begin{center}
\includegraphics[width=11.8cm]{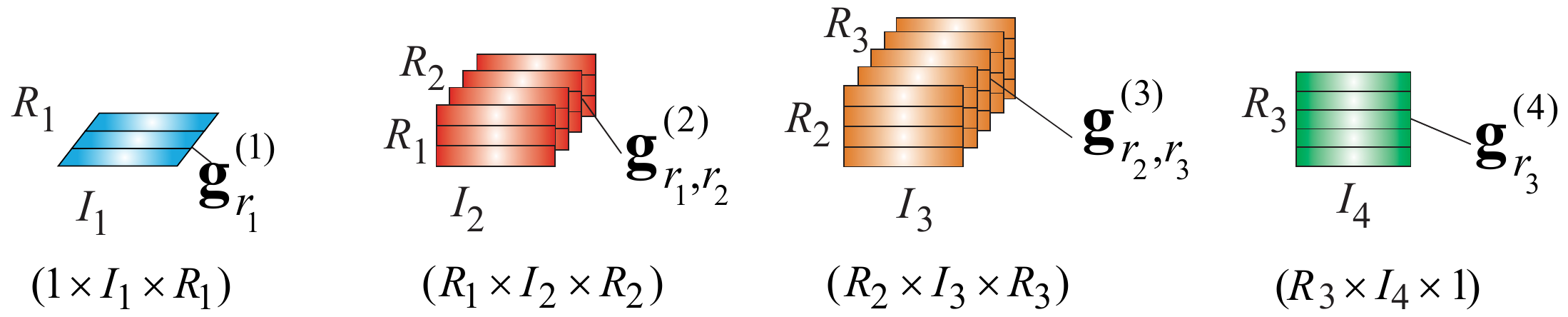}\\
\end{center}
\vspace{0.5cm}
(b)\\
\begin{center}
\includegraphics[width=11.5cm]{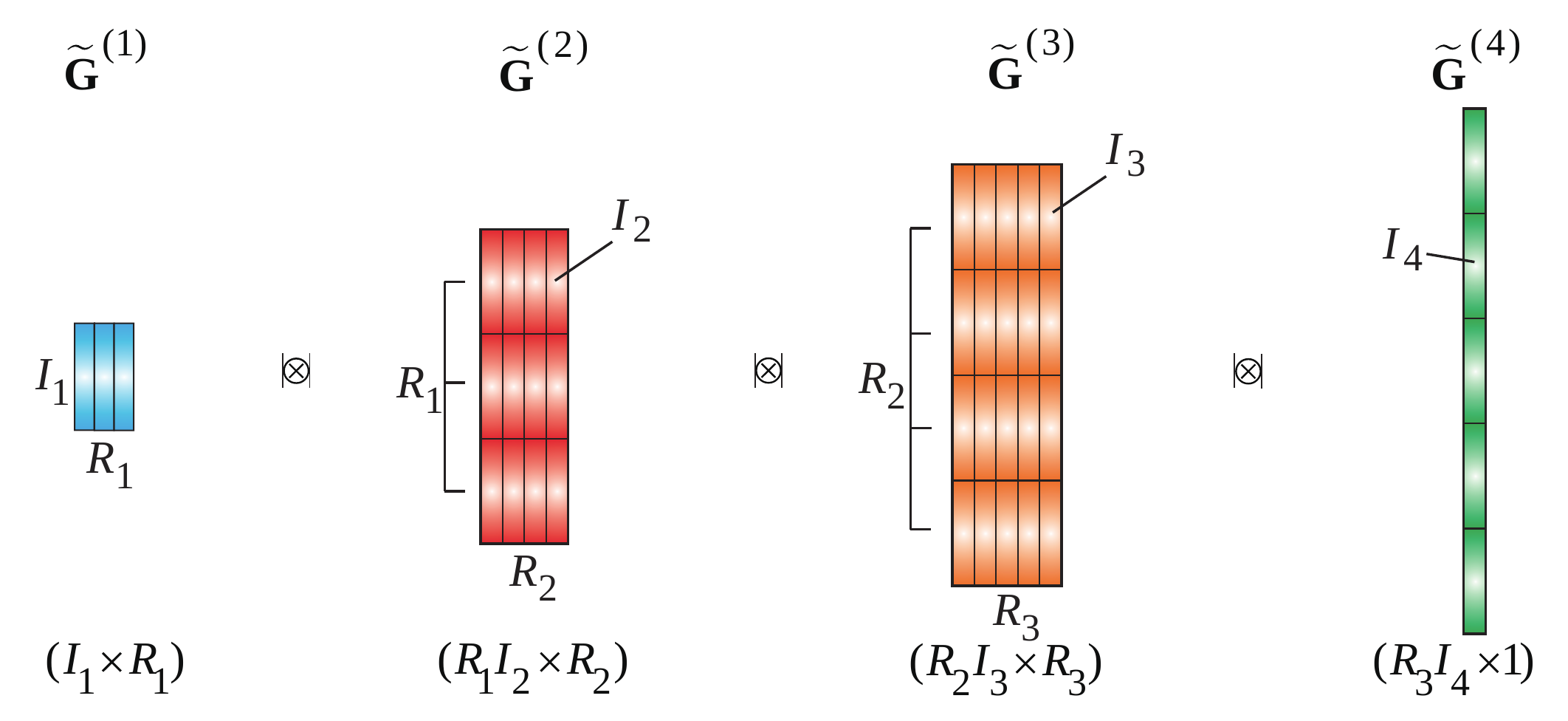}
\end{center}
\caption{{\small TT decomposition of a 4th-order  tensor, $\underline \bX$,  for which the  TT rank is $R_1=3, \; R_2=4, \;R_3=5$. (a) (Upper panel) Representation of the TT via a multilinear product of the cores,  $\underline \bX \cong \underline  \bG^{(1)} \, \times^1 \, \underline  \bG^{(2)} \, \times^1 \, \underline  \bG^{(3)} \times^1 \, \underline \bG^{(4)}  = \llangle \underline  \bG^{(1)}, \underline \bG^{(2)},\underline \bG^{(3)}, \underline \bG^{(4)} \rrangle$, and (lower panel) an equivalent representation via the outer product of
mode-2 fibers  (sum of rank-1 tensors) in the form,  $\underline \bX \cong \sum_{r_1=1}^{R_1} \sum_{r_2=1}^{R_2}  \sum_{r_{3}=1}^{R_{3}}  \sum_{r_{4}=1}^{R_{4}} (\bg^{(1)}_{\,r_1} \; \circ \; \bg^{(2)}_{\,r_1, \,r_2} \; \circ  \; \bg^{(3)}_{\,r_{2},\,r_{3}} \; \circ \; \bg^{(4)}_{\,r_{3}})$.   (b) TT decomposition  in a vectorized form  represented via strong Kronecker products of block matrices, $\bx \cong \widetilde \bG^{(1)} \skron \widetilde \bG^{(2)} \skron \widetilde \bG^{(3)} \skron  \widetilde \bG^{(4)}  \in \Real^{I_1 I_2 I_3 I_4}$, where the block matrices are defined as $\widetilde \bG^{(n)} \in \Real^{R_{n-1} I_{n} \times R_{n}}$, with block vectors $\bg^{(n)}_{r_{n-1},\,r_n} \in \Real^{I_n \times 1}$, $n=1,\ldots,4$ and $R_0=R_4=1$.}}
\label{Fig:TT=MPS4}
\end{figure}
\begin{table}
\vspace{-1.5cm}
 \centering
\caption{Equivalent representations of the  Tensor Train decomposition (MPS  with open boundary  conditions)  approximating an $N$th-order tensor $\underline \bX \in \Real^{I_1 \times I_2 \times \cdots \times I_N}$. It is assumed that the TT rank is $\brr_{TT}= \{R_1,R_2,\ldots,R_{N-1}\}$, with $R_0=R_N=1$.}
{
\shadingbox{
    \begin{tabular*}{1.0\textwidth}[t]{@{\extracolsep{\fill}}@{\hspace{2ex}}l@{\hspace{1em}}l}
     \hline  &  \\[-6pt]
 \multicolumn{1}{c}{Tensor representation:  Multilinear products of TT-cores}\\
  & \\
   $ \underline \bX = \underline \bG^{(1)} \times^1 \; \underline \bG^{(2)} \times^1 \; \cdots \; \times^1 \;\underline \bG^{(N)} \in \Real^{I_1 \times I_2 \times \cdots \times I_N}$  
 \\ & \\
 with the 3rd-order cores $\underline \bG^{(n)} \in \Real^{ R_{n-1} \times I_n \times R_n}$, $\; (n=1,2,\ldots,N)$ 
 \\  &   \\[-6pt] \hline & \\[-6pt]

  \multicolumn{1}{c}{Tensor representation: Outer products}\\
  & \\
   $ \underline \bX  = \displaystyle{\sum_{r_1,\,r_2,\ldots,r_{N-1}=1}^{R_1,R_2, \ldots,
 R_{N-1}}} \;\; \bg^{(1)}_{\;1,r_1} \; \circ \; \bg^{(2)}_{\;r_1,\, r_2}  \circ \cdots
 \circ \; \bg^{(N-1)}_{\;r_{N-2}, \, r_{N-1}} \; \circ \; \bg^{(N)}_{\;r_{N-1},\,1} $ 
   \\ & \\
 where $\bg^{(n)}_{\;r_{n-1}, \,r_n}=\underline \bG^{(n)}(r_{n-1},:,r_n)  \in \Real^{I_n}$ are fiber vectors.
 \\  &   \\[-6pt] \hline & \\[-6pt]
  \multicolumn{1}{c}{Vector representation: Strong Kronecker products}
   \\ & \\
   $ \bx  =  \widetilde \bG^{(1)} \; \skron \;  \widetilde \bG^{(2)} \; \skron \cdots \skron\;  \widetilde \bG^{(N)} \in \Real^{I_1 I_2 \cdots I_N}$, $\quad$ where
    \\ & \\
 $\widetilde\bG^{(n)} \in \Real^{ R_{n-1} I_n \times  R_n}$ are block matrices with blocks $\bg^{(n)}_{r_{n-1}, r_n} \in \Real^{I_n}$
  \\   &  \\[-6pt] \hline & \\[-6pt]
    \multicolumn{1}{c}{Scalar representation}\\
  & \\
   $ x_{\;i_1,i_2,\ldots,i_N}  = \displaystyle{\sum_{r_1,r_2,\ldots,r_{N-1}=1}^{R_1,R_2, \ldots,
 R_{N-1}}} \!\!  g^{(1)}_{\;1,\, i_1, \,r_1} \; g^{(2)}_{\;r_1, \,i_2, \,r_2} \;
   \cdots \; g^{(N-1)}_{\;r_{N-2},\, i_{N-1}, \,r_{N-1}}
 g^{(N)}_{\;r_{N-1}, \,i_N,1}$
  \\ & \\
where $g^{(n)}_{\;r_{n-1}, \,i_n, \,r_n}$ are entries of a 3rd-order core $\underline \bG^{(n)} \in \Real^{ R_{n-1} \times I_n \times R_n}$ 
  \\  &   \\[-6pt] \hline & \\[-6pt]
  \multicolumn{1}{c}{Slice (MPS) representation}\\
  & \\
   $ x_{\;i_1,\,i_2,\ldots,i_N} =  \bG^{(1)}_{i_1} \; \bG^{(2)}_{i_2} \; \cdots  \; \bG^{(N)}_{i_N}$, $\quad$ where
  \\ & \\
  $\bG^{(n)}_{i_n} = \underline \bG^{(n)}(:,i_n,:) \in \Real^{R_{n-1} \times R_n}$ are lateral slices of $\underline \bG^{(n)} \in \Real^{ R_{n-1} \times I_n \times R_n}$
\\  &  \\[-6pt] \hline  
    \end{tabular*}
    }}
\label{table:MPS}
\end{table}

\begin{table}
\vspace{-2.5cm}
 \centering
\caption{Equivalent representations of the  Tensor Chain (TC) decomposition (MPS  with periodic boundary  conditions)  approximating an $N$th-order tensor $\underline \bX \in \Real^{I_1 \times I_2 \times \cdots \times I_N}$. It is assumed that the TC rank is $\brr_{TC}= \{R_1,R_2,\ldots,R_{N-1},R_N\}$.}
{
\shadingbox{
    \begin{tabular*}{1.00\textwidth}[t]{@{\extracolsep{\fill}}@{\hspace{2ex}}l@{\hspace{1em}}l}
    \hline  &  \\[-6pt]
 \multicolumn{1}{c}{Tensor representation: Trace of multilinear products of cores}\\
  & \\
   $ \underline \bX = \mbox{Tr} \; (\underline \bG^{(1)} \times^1 \; \underline \bG^{(2)} \times^1 \; \cdots \; \times^1 \;\underline \bG^{(N)} )
   \in \Real^{I_1 \times I_2 \times \cdots \times I_N}$  
 \\ & \\
 with the 3rd-order cores $\underline \bG^{(n)} \in \Real^{ R_{n-1} \times I_n \times R_n}$,\; $R_0=R_N$,\;
 $n=1,2,\ldots,N$ 
 \\  &   \\[-6pt] \hline & \\[-6pt]
  \multicolumn{1}{c}{Tensor/Vector representation: Outer/Kronecker products}\\
  & \\
   $ \underline \bX  = \displaystyle{\sum_{r_1,\,r_2,\ldots,r_{N}=1}^{R_1,R_2, \ldots,
 R_{N}}} \;\; \bg^{(1)}_{\;r_N,\,r_1} \; \circ \; \bg^{(2)}_{\;r_1,\, r_2} \; \circ \cdots \circ \; \bg^{(N)}_{\;r_{N-1},\,r_N} \in \Real^{I_1 \times I_2 \times \cdots \times I_N}$ 
   \\ & \\
    $  \bx  = \displaystyle{\sum_{r_1,\,r_2,\ldots,r_{N}=1}^{R_1,R_2, \ldots,
 R_{N}}} \;\; \bg^{(1)}_{\;r_N,\,r_1} \; \otimes_L \; \bg^{(2)}_{\;r_1,\, r_2} \; \otimes_L \cdots \otimes_L \; \bg^{(N)}_{\;r_{N-1},\,r_N} \in \Real^{I_1 I_2 \cdots I_N}$ 
   \\ & \\
 where $\bg^{(n)}_{\;r_{n-1}, \; r_n} \in \Real^{I_n}$ are fiber vectors  within $\underline \bG^{(n)}(r_{n-1},\, :, \, r_n) \in \Real^{I_n}$
 \\  &   \\[-6pt] \hline & \\[-6pt]
  \multicolumn{1}{c}{Vector representation: Strong Kronecker products}
   \\ & \\
   $ \bx  = \displaystyle\sum_{r_N=1}^{R_N}(\widetilde \bG^{(1)}_{\,r_N} \; \skron \;  \widetilde \bG^{(2)} \; \skron \cdots \skron\;  \widetilde \bG^{(N-1)} \skron\;  \widetilde \bG^{(N)}_{\,r_N}) \in \Real^{I_1 \, I_2 \cdots I_N}$  $\quad$ where
    \\ & \\
 $\widetilde\bG^{(n)} \in \Real^{ R_{n-1} \, I_n \times  R_n}$ are block matrices with blocks $\bg^{(n)}_{\;r_{n-1}, \,r_n} \in \Real^{I_n}$,
\\ & \\
 $\widetilde\bG^{(1)}_{\,r_N} \in \Real^{ I_1 \times  R_1}$ is a matrix with blocks (columns) $\bg^{(1)}_{\;r_{N},\, r_1} \in \Real^{I_1}$,
\\ & \\
$\widetilde\bG^{(N)}_{\,r_N} \in \Real^{ R_{N-1} \, I_N \times  1}$ is a block vector with blocks $\bg^{(N)}_{\; r_{N-1},\, r_N} \in \Real^{I_N}$
  \\   &  \\[-6pt] \hline & \\[-6pt]
    \multicolumn{1}{c}{Scalar representations}
 \\ & \\
   $ x_{\;i_1,\,i_2,\ldots,i_N} =  \tr(\bG^{(1)}_{i_1} \, \bG^{(2)}_{i_2}  \cdots   \bG^{(N)}_{i_N})
   =\!\!\displaystyle\sum_{r_N=1}^{R_N}\!\! (\bg^{(1)\,\text{T}}_{r_N,\, i_1, \, :} \, \bG^{(2)}_{i_2}  \cdots   \bG^{(N-1)}_{i_{N-1}} \bg^{(N)}_{:, \, i_N, \, r_N})$
  \\ & \\
 where
 $\bg^{(1)}_{r_N, \, i_1, \,:} = \underline \bG^{(1)}(r_N,\,i_1,\,:) \in \Real^{R_1}$,
 $ \bg^{(N)}_{:, \,i_N,\,r_N}\! = \underline \bG^{(N)}(:, \, i_N, \,r_N) \in \Real^{R_{N-1}}$
\\  &  \\[-6pt] \hline  
    \end{tabular*}
    }}
\label{table:MPS_PBC}
\end{table}
\begin{figure} 
(a)
\begin{center}
\includegraphics[width=4.8cm]{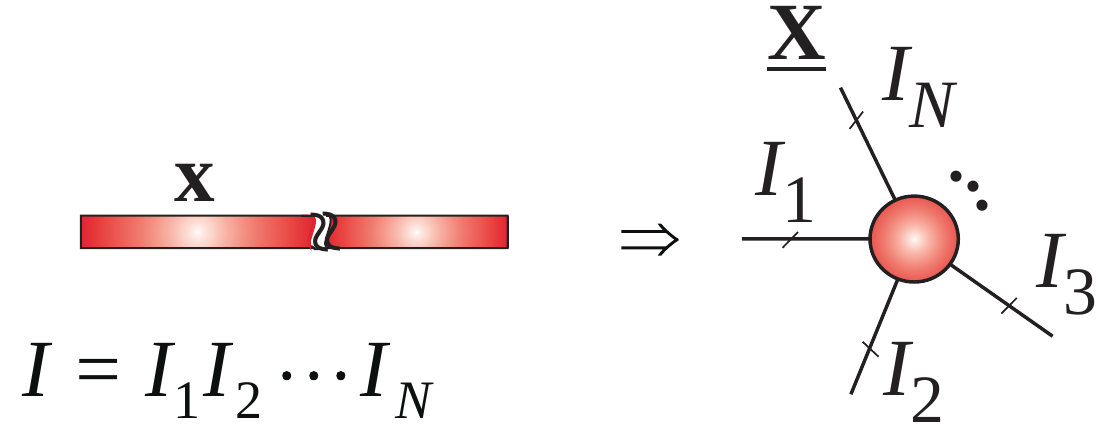}
\end{center}
(b)
\begin{center}
\includegraphics[width=11.8cm]{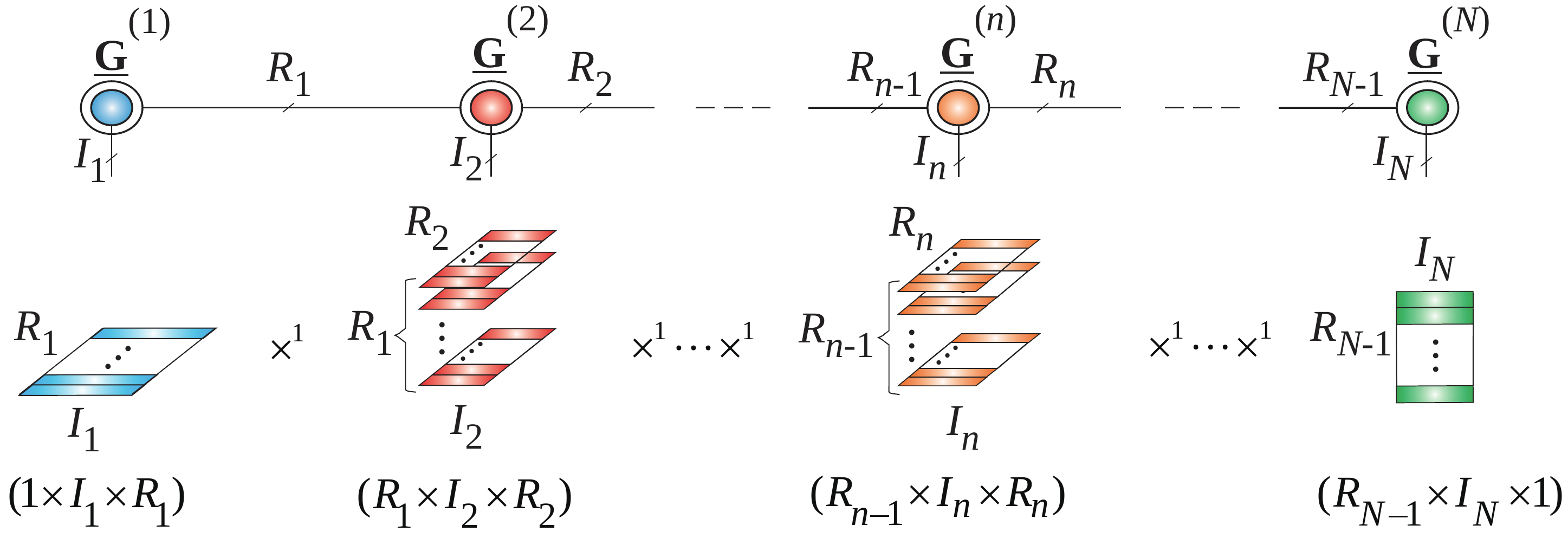}\\
\end{center}
(c)
\begin{center}
\includegraphics[width=11.8cm]{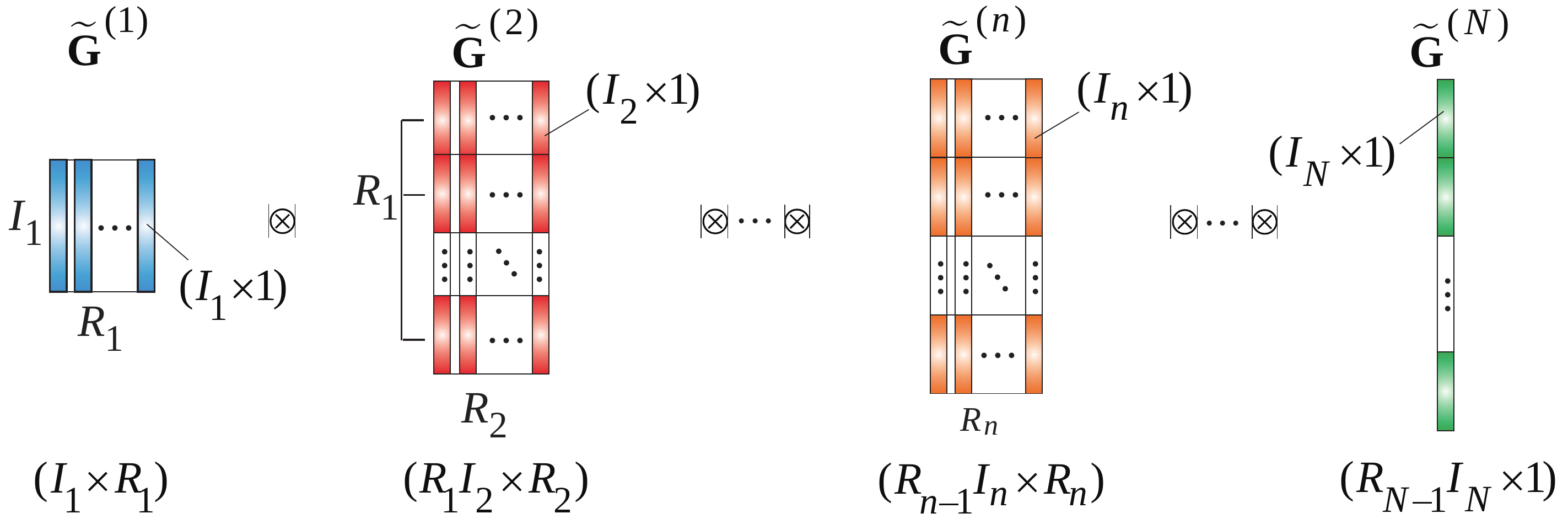}
\end{center}
\caption{TT/MPS decomposition of an $N$th-order data tensor, $\underline \bX$,  for which the  TT rank is $\{R_1,R_2,\ldots,R_{N-1}\}$. (a) Tensorization of a huge-scale vector, $\bx \in \Real^I$, into an $N$th-order tensor, $\underline \bX \in \Real^{I_1 \times I_2 \times \cdots \times I_N}$.
 (b) The data tensor can be represented exactly or approximately via  a tensor train (TT/MPS), consisting of 3rd-order cores in the form $\underline \bX \cong \underline  \bG^{(1)} \, \times^1 \, \underline  \bG^{(2)} \, \times^1 \, \cdots \times^1 \, \underline \bG^{(N)}  = \llangle \underline  \bG^{(1)}, \underline \bG^{(2)}, \ldots, \underline \bG^{(N)} \rrangle$, where $\underline \bG^{(n)} \in \Real^{R_{n-1} \times I_n \times R_n}$ for $n=1,2,\ldots,N$ with $R_0=R_N=1$. (c) Equivalently,  using the strong Kronecker products, the TT tensor can be expressed   in a vectorized form, $\bx \cong \widetilde \bG^{(1)} \skron \widetilde \bG^{(2)} \skron \cdots \skron \widetilde \bG^{(N)}  \in \Real^{I_1 I_2 \cdots I_N}$, where the block  matrices are defined as $\widetilde \bG^{(n)}  \in \Real^{R_{n-1} I_{n} \times R_{n}}$, with blocks $\bg^{(n)}_{r_{n-1},\,r_n} \in \Real^{I_n \times 1}$.}
\label{Fig:TT=MPS=MPO}
\end{figure}

The   tensor train (TT/MPS) representation of an $N$th-order data tensor, $\underline \bX \in \Real^{I_1 \times I_2  \times \cdots \times I_N}$, can be described  in
several equivalent forms (see Figures \ref{Fig:TT=MPS4}, \ref{Fig:TT=MPS=MPO} and  Table \ref{table:MPS}) listed below:

\begin{enumerate}

\item The entry-wise scalar form, given by 
\be
x_{i_1,i_2,\ldots,i_N} &\cong& \sum_{r_1,\,r_2,\ldots,r_{N-1}=1}^{R_1,R_2, \ldots,
 R_{N-1}} g^{(1)}_{1,\,i_1,\,r_1} \; g^{(2)}_{r_1,\,i_2,\,r_2}  \cdots g^{(N)}_{r_{N-1},\,i_N,1}.
\notag \\
\ee
\item The slice representation (see Figure \ref{Fig:TTouter}) in the form
\be
x_{i_1,i_2,\ldots,i_N} \cong \bG^{(1)}_{i_1} \; \bG^{(2)}_{i_2}
 \cdots  \bG^{(N)}_{i_N},
\ee
where the slice matrices  are defined as
\begin{equation}
\bG^{(n)}_{i_n}= \underline \bG^{(n)}(:, i_n, :) \in \Real^{R_{n-1} \times R_n}, \quad i_n=1,2,\ldots,I_n \notag
\end{equation}
 with $\bG^{(n)}_{i_n}$ being the $i_n$th  lateral slice of the  core
 $\underline \bG^{(n)} \in \Real^{R_{n-1} \times I_n \times R_n}$, $n=1,2,\ldots,N$ and $R_0=R_N=1$.

\item  The (global) tensor form, based on multilinear products (contraction) of cores (see Figure  \ref{Fig:TT=MPS4}(a)) given by
\be
\underline \bX &\cong& \underline \bG^{(1)}\times^1 \underline \bG^{(2)} \times^1
\cdots \times^1 \underline \bG^{(N-1)} \times^1 \underline \bG^{(N)} \nonumber \\
  &=& \llangle \underline \bG^{(1)}, \underline \bG^{(2)},  \ldots, \underline \bG^{(N-1)}, \underline \bG^{(N)} \rrangle,
\ee
where the 3rd-order cores{\footnote{Note that the cores  $\underline \bG^{(1)}$ and $\underline \bG^{(N)}$ are now two-dimensional arrays (matrices), but for a uniform representation, we assume that these matrices are treated as 3rd-order cores of sizes $1 \times I_1 \times R_1$ and $R_{N-1} \times I_N \times 1$, respectively.}} $\underline \bG^{(n)} \in \Real^{R_{n-1} \times I_n \times R_n}$, $n=1,2, \ldots, N$ and $R_0=R_N=1$ (see also Figure \ref{Fig:TT=MPS=MPO}(b)).

\item  The tensor form, expressed as a sum of rank-1 tensors  (see Figure \ref{Fig:TT=MPS4}(a))
\begin{equation}
\underline \bX \cong \sum_{r_1,\,r_2,\ldots,r_{N-1}=1}^{R_1,R_2, \ldots,
 R_{N-1}} \bg^{(1)}_{\,1,\,r_1} \; \circ \; \bg^{(2)}_{\,r_1,\, r_2} \; \circ \; \cdots
\circ \; \bg^{(N-1)}_{\,r_{N-2},\,r_{N-1}} \circ \; \bg^{(N)}_{\,r_{N-1},\,1},
\label{TT-outerprod}
\end{equation}
where
$\bg^{(n)}_{r_{n-1},r_{n}} = \underline \bG^{(n)}(r_{n-1}, \, :, \, r_n) \in \Real^{I_n}$ are mode-2 fibers,
$n=1,2,\ldots, N$ and $R_0=R_N=1$.

\item  A  vector form, expressed by Kronecker products of the fibers
\begin{align}
\bx
&\cong \sum_{r_1,r_2,\ldots,r_{N-1}=1}^{R_1,R_2, \ldots,
 R_{N-1}} \bg^{(1)}_{\,1,\,r_1} \; \otimes_L \; \bg^{(2)}_{\,r_1, \,r_2} \; \otimes_L \;\notag\\
&\quad  \cdots \otimes_L \; \bg^{(N-1)}_{\,r_{N-2},\,r_{N-1}}
\otimes_L \; \bg^{(N)}_{\,r_{N-1},\,1},\label{TT-Kron1}
\end{align}
where $\bx= \mbox{vec} (\bX) \in \Real^{I_1 I_2 \cdots I_N}$. 

\item  An alternative vector form, produced by strong Kronecker products of block matrices 
(see Figure \ref{Fig:TT=MPS4}(b)) and Figure \ref{Fig:TT=MPS=MPO}(c)), given by
\begin{equation}
\bx \cong  \widetilde \bG^{(1)} \; \skron \;  \widetilde \bG^{(2)} \; \skron \cdots \skron \widetilde \bG^{(N)}, \\
\end{equation}
where
the block matrices $\widetilde \bG^{(n)}  \in \Real^{R_{n-1} I_n \times R_n}$, for $n=1,2,\ldots,N$, consist of blocks $\bg^{(n)}_{r_{n-1},r_n} \in \Real^{I_n \times 1}$, $n=1,2, \ldots, N$, with $R_0=R_N=1$, and the symbol  $\skron$ denotes the strong Kronecker product.

\end{enumerate}

Analogous relationships can be established for Tensor Chain (i.e., MPS with PBC (see Figure \ref{Fig:TTouter}(b))
and summarized  in Table \ref{table:MPS_PBC}.

\section{Matrix  TT Decomposition -- Matrix Product Operator}

\begin{figure}
(a)\\
\vspace{-0.1cm}
\begin{center}
\includegraphics[width=11.99cm]{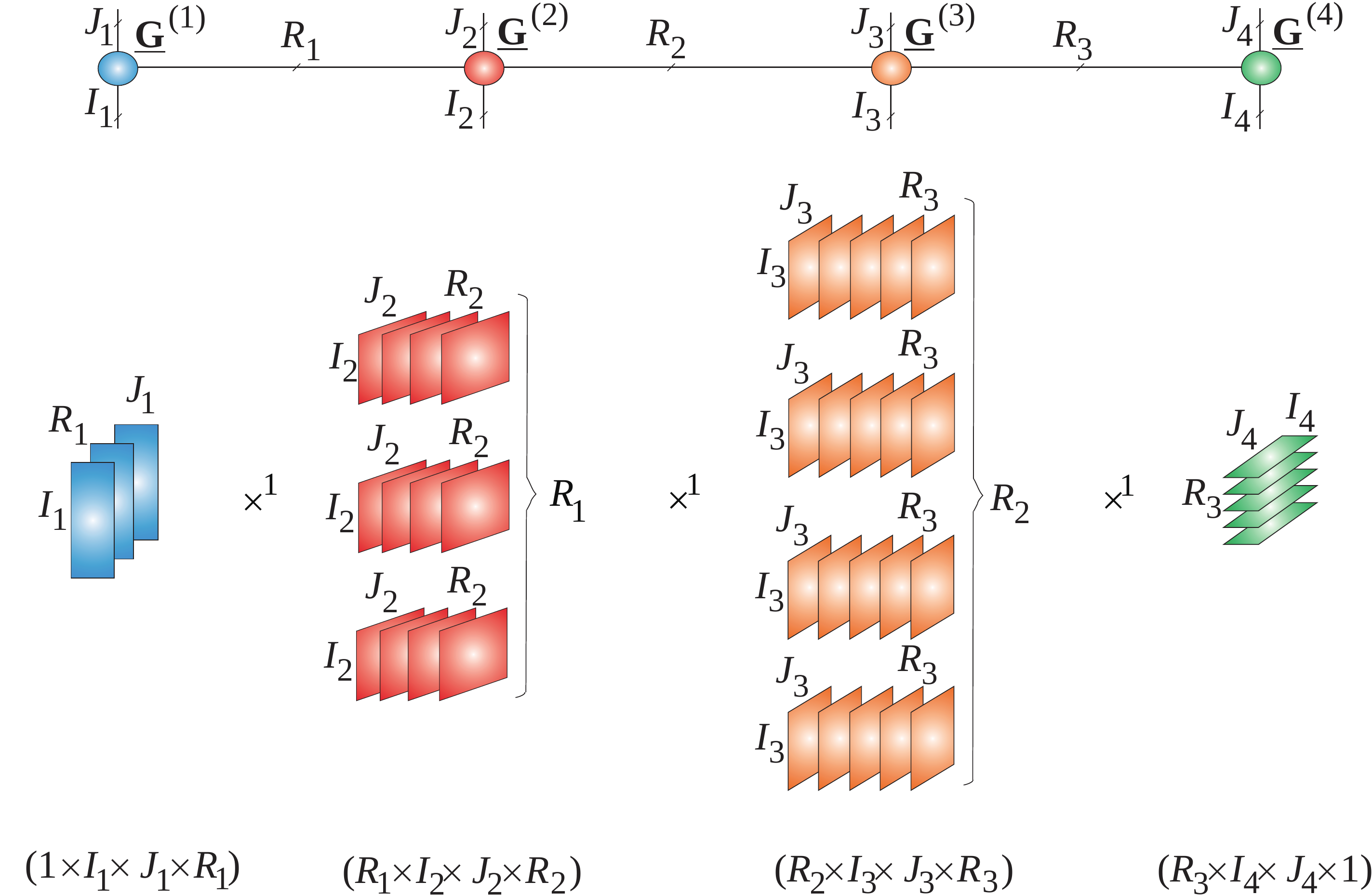}\\
\end{center}
(b)\\
\begin{center}
\includegraphics[width=11.6cm]{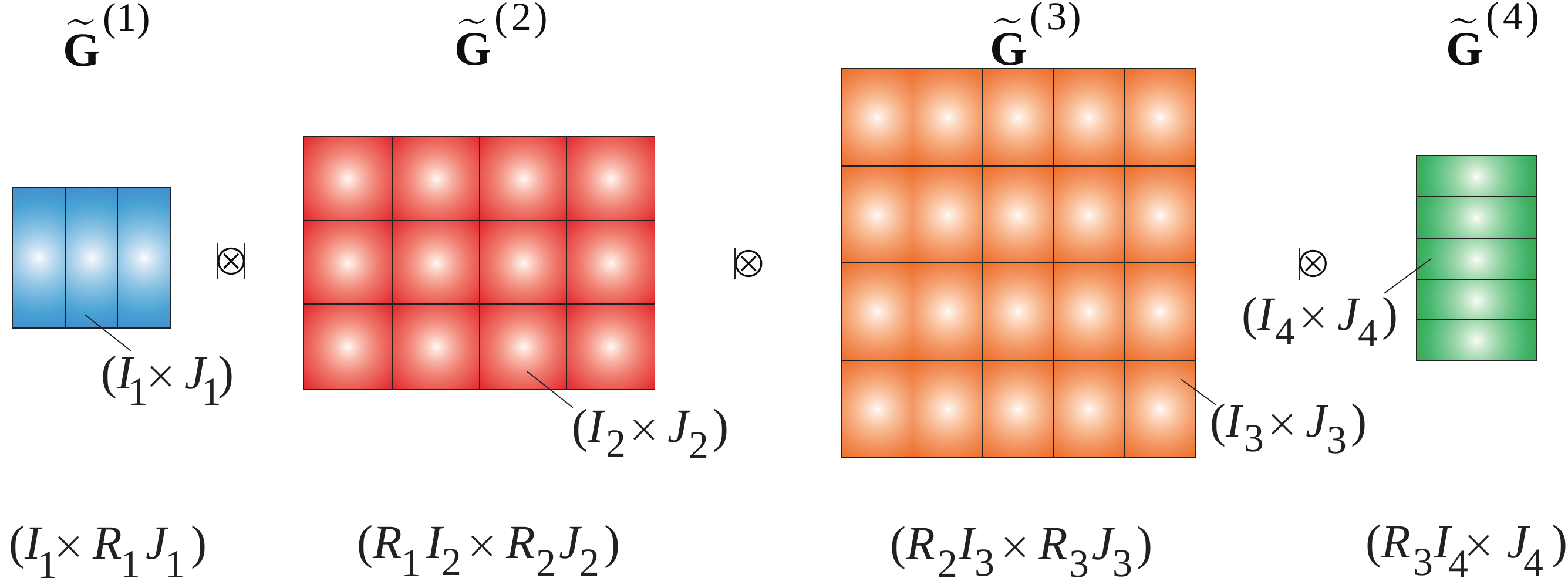}\\
\end{center}
\caption{TT/MPO decomposition of a matrix, $\bX \in \Real^{I \times J}$, reshaped as an  8th-order tensor,  $\underline \bX \in \Real^{I_1 \times J_1 \times \cdots \times I_4 \times J_4} $, where $I=I_1 I_2 I_3 I_4$ and $J=J_1 J_2 J_3 J_4$. (a) Basic TT representation via multilinear products (tensor contractions) of cores
$\underline \bX = \underline \bG^{(1)} \times^1 \underline \bG^{(2)} \times^1 \underline \bG^{(3)} \times^1  \underline \bG^{(4)}$, with $\underline \bG^{(n)} \in \Real^{R_{n-1} \times I_n \times R_n}$ for $R_1=3, R_2=4, R_3=5, R_0=R_4=1$.
(b) Representation of a matrix or a matricized tensor via  strong Kronecker products of block matrices, in the form
$\bX = \widetilde \bG^{(1)} \skron \widetilde \bG^{(2)} \skron \widetilde \bG^{(3)} \skron \widetilde \bG^{(4)} \in \Real^{I_1 I_2 I_3 I_4 \; \times  \; J_1 J_2 J_3 J_4} $.}
\label{Fig:TT=MPO4}
\end{figure}

\begin{figure}
(a)
\begin{center}
\includegraphics[width=5.9cm]{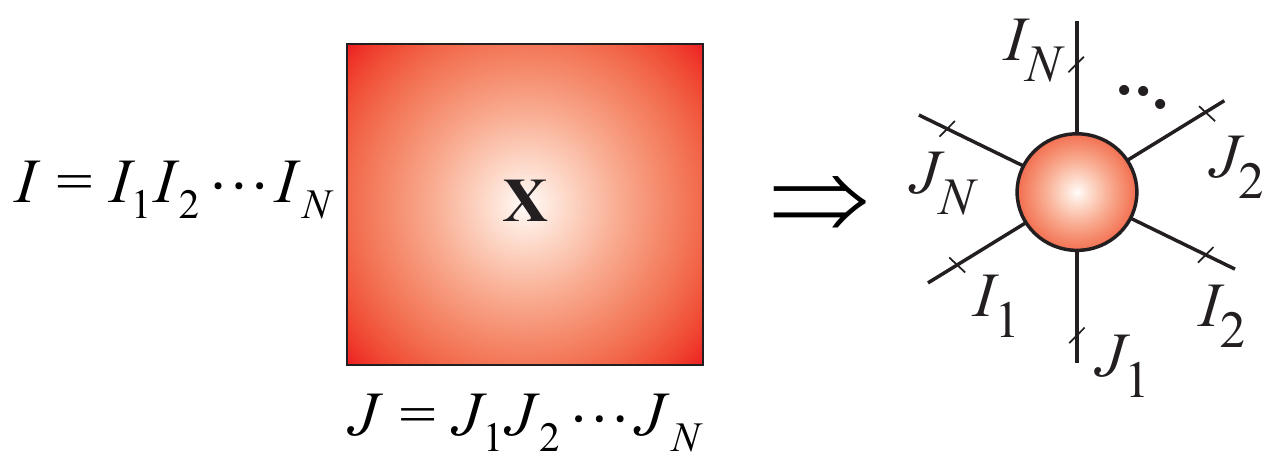}
\end{center}
(b)
\begin{center}
\includegraphics[width=11.8cm]{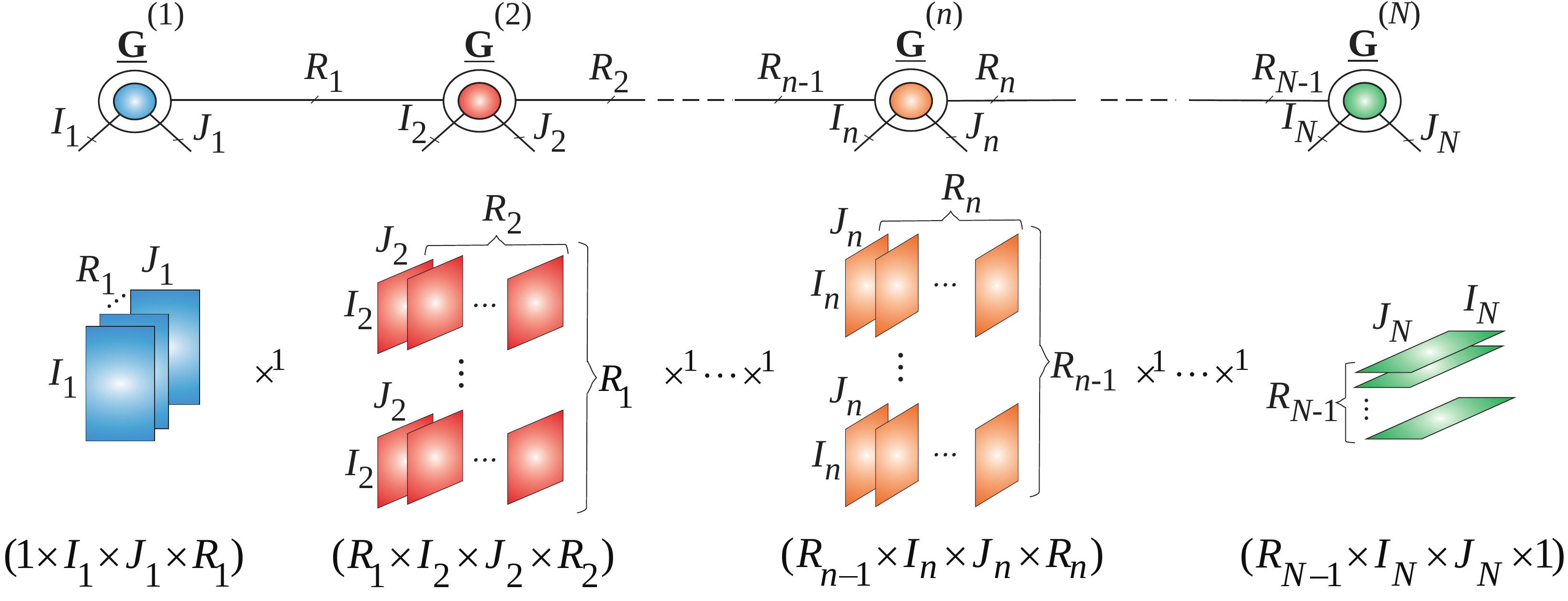}\\
\end{center}
(c)
\begin{center}
\includegraphics[width=11.8cm]{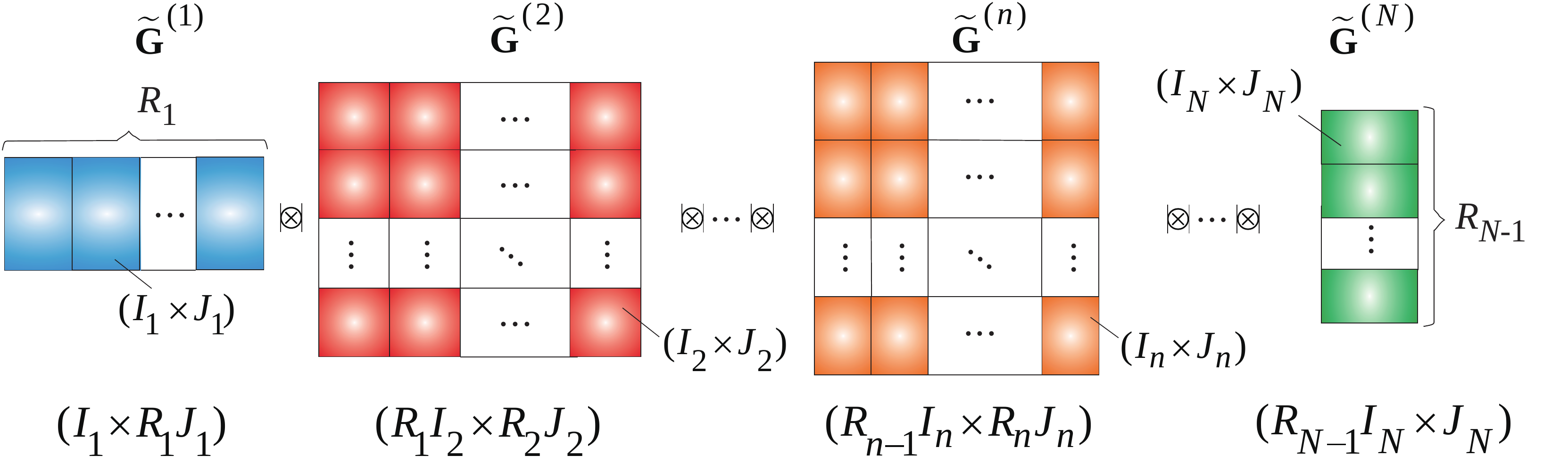}
\end{center}
\caption{Representations of huge matrices by ``linked'' block matrices. (a) Tensorization of a huge-scale  matrix,  $\bX \in \Real^{I \times J}$, into a $2N$th-order tensor $\underline \bX \in \Real^{I_1 \times J_2 \times \cdots \times I_N \times J_N}$. (b) The  TT/MPO decomposition of a huge matrix, $\bX$,  expressed by 4th-order cores, $\underline \bG^{(n)} \in \Real^{R_{n-1} \times I_n \times J_n \times R_n}$.
(c) Alternative graphical representation of a matrix, $\bX \in \Real^{I_1 I_2 \cdots I_N \; \times \; J_1 J_2 \cdots J_N}$,  via strong Kronecker products of block matrices $\widetilde \bG^{(n)} \in \Real^{R_{n-1} \, I_n \; \times \; R_n \, J_n}$ for $n=1,2,\ldots,N$ with $R_0=R_N=1$.}
\label{Fig:TT=MPS=MPO2}
\end{figure}

The matrix tensor train, also called the Matrix Product Operator (MPO) with  open boundary conditions (TT/MPO), is an important TN model which
 first represents huge-scale structured matrices, $\bX \in \Real^{I \times J}$, as  $2N$th-order tensors,  $\underline \bX \in \Real^{I_1 \times J_1 \times I_2 \times J_2 \times \cdots I_N \times J_N}$, where $I=I_1 I_2 \cdots I_N$ and $J=J_1 J_2 \cdots J_N$
(see Figures \ref{Fig:TT=MPO4}, \ref{Fig:TT=MPS=MPO2} and Table \ref{table:MPO}).
Then, the matrix TT/MPO  converts such a $2N$th-order tensor into a chain (train) of 4th-order cores{\footnote{The cores  $\underline \bG^{(1)}$ and $\underline \bG^{(N)}$
 are in fact three-dimensional arrays, however for uniform representation, we treat them as   4th-order cores of sizes $1 \times I_1 \times J_1 \times R_1$ and $R_{N-1} \times I_N \times J_N \times 1$.}}.
 It should be noted that the matrix TT decomposition is equivalent to the vector TT, created by merging all index pairs $(i_n,j_n)$ into a single index ranging  from 1 to  $I_n J_n$, in a reverse lexicographic order.

Similarly to the vector TT decomposition, a large scale $2N$th-order tensor, $\underline \bX \in \Real^{I_1 \times J_1 \times I_2 \times J_2 \times \cdots \times I_N \times J_N}$,
can be represented in a TT/MPO format via the following mathematical representations:
\begin{enumerate}

\item The scalar (entry-wise) form
 \be
x_{i_1,j_1, \ldots,i_{N},j_N} &\cong& \sum_{r_1=1}^{R_1} \sum_{r_2=1}^{R_2}\cdots \sum_{r_{N-1}=1}^{R_{N-1}} g^{(1)}_{\,1,\,i_1,j_1,r_1} \; g^{(2)}_{\,r_1,\,i_2,\,j_2,\,r_2}  \notag \\
&&\cdots g^{(N-1)}_{\,r_{N-2},\,i_{N-1},\,j_{N-1},\,r_{N-1}} \; g^{(N)}_{\,r_{N-1},\,i_{N},\,j_{N},\,1}.
\ee

 \item The slice representation
\begin{equation}
x_{i_1,j_1, \ldots,i_{N},j_N} \cong  \bG^{(1)}_{i_1, j_1} \; \bG^{(2)}_{i_2, j_2} \cdots \bG^{(N)}_{i_N, j_N}, \\
\end{equation}
where $\bG^{(n)}_{i_n, j_n} = \underline \bG^{(n)}(:, \,i_n,\,j_n,\,:) \in \Real^{R_{n-1} \times R_n}$ are slices of the cores
$\underline \bG^{(n)} \in \Real^{R_{n-1} \times I_n \times J_n \times R_n}$, $\;\; n=1,2,\ldots,N$ and $R_0 =R_N=1$.

\item  The  compact tensor form based on multilinear products (Figure \ref{Fig:TT=MPS=MPO2}(b))
\be
\underline \bX & \cong & \underline \bG^{(1)} \; \times^1 \; \underline \bG^{(2)} \; \times^1 \; \cdots \times^1 \; \underline \bG^{(N)} \nonumber \\
  &=& \llangle \underline \bG^{(1)}, \underline \bG^{(2)}, \dots, \underline \bG^{(N)} \rrangle,
\ee
where the TT-cores are defined as $\underline \bG^{(n)} \in \Real^{R_{n-1} \times I_{n} \times J_{n} \times R_{n}}$, $\;n=1,2,\ldots,N$ and $R_0=R_N=1$.

\item  A  matrix form, based on strong Kronecker products of block matrices (Figures  \ref{Fig:TT=MPO4}(b) and \ref{Fig:TT=MPS=MPO2}(c))
\begin{equation}
\bX \cong  \widetilde \bG^{(1)} \; \skron \; \widetilde \bG^{(2)} \; \skron \cdots \skron\; \widetilde\bG^{(N)} \in \Real^{I_1 \cdots I_N  \; \times
\; J_1 \cdots J_N},
\end{equation}
where
$\widetilde \bG^{(n)} \in \Real^{R_{n-1} I_{n} \times R_{n} J_{n}}$ are block  matrices with blocks  $\bG^{(n)}_{r_{n-1},r_n} \in \Real^{I_{n} \times J_{n}}$ and the number of blocks  is $R_{n-1}\times R_{n}$.
 In a special case, when the TT ranks $R_n=1, \; \forall n$, the strong Kronecker products simplify into standard (left) Kronecker products.

\end{enumerate}

\begin{table}
\vspace{-2.5cm}
\centering
\caption{Equivalent forms of the matrix Tensor Train  decomposition (MPO with open boundary  conditions) for a $2N$th-order tensor $\underline \bX \in \Real^{I_1 \times J_1 \times I_2 \times J_2 \times \cdots \times  I_N \times J_N }$. It is assumed that the TT rank is $\{R_1,R_2,\ldots,R_{N-1}\}$, with $R_0=R_N=1$.}
 {
\shadingbox{
    \begin{tabular*}{1.00\textwidth}[t]{@{\extracolsep{\fill}}@{\hspace{2ex}}l@{\hspace{1em}}l}
    \hline  &  \\[-6pt]
 \multicolumn{1}{c}{Tensor representation: Multilinear products (tensor contractions)}\\
  & \\
 $ \underline \bX = \underline \bG^{(1)} \times^1 \; \underline \bG^{(2)} \times^1  \cdots \times^1 \; \underline \bG^{(N-1)} \times^1 \; \underline \bG^{(N)}$
 \\ & \\
  \hspace{2em}
 with 4th-order cores $\underline \bG^{(n)} \in \Real^{ R_{n-1} \times I_n \times J_n \times R_n}$, $\quad(n=1,2,\ldots,N)$
 \\  &   \\[-6pt] \hline & \\[-6pt]
  \multicolumn{1}{c}{Tensor representation: Outer products}\\
  & \\
 $  \underline \bX  = \displaystyle{\sum_{r_1,r_2,\ldots,r_{N-1}=1}^{R_1,R_2, \ldots,
 R_{N-1}}} \;\; \bG^{(1)}_{\;1,\,r_1} \; \circ \; \bG^{(2)}_{\;r_1,\, r_2}  \circ \cdots
 \circ \; \bG^{(N-1)}_{\;r_{N-2}, \,r_{N-1}} \; \circ \; \bG^{(N)}_{\;r_{N-1}, \,1}$
   \\ & \\
 \hspace{1em}
 where $\bG^{(n)}_{\;r_{n-1}, \,r_n} \in \Real^{I_n \times J_n} $ are blocks of  $\widetilde \bG^{(n)}  \in \Real^{ R_{n-1} I_n  \times R_n J_n}$
 \\  &   \\[-6pt] \hline & \\[-6pt]
  \multicolumn{1}{c}{Matrix representation: Strong Kronecker products}
   \\ & \\
 $ \bX =   \widetilde \bG^{(1)} \; \skron \; \widetilde \bG^{(2)} \; \skron \cdots \;\skron\; \widetilde\bG^{(N)} \in \Real^{I_1  \cdots I_N \; \times \; J_1  \cdots J_N}$
  \\ & \\
where $\widetilde \bG^{(n)} \in \Real^{ R_{n-1} I_n  \times R_n J_n}$
 are block matrices with blocks\\[4pt]
 \quad\quad $\underline \bG^{(n)}(r_{n-1},:,:,r_n)$
 %
  \\   &  \\[-6pt] \hline & \\[-6pt]
    \multicolumn{1}{c}{Scalar  representation}\\
  & \\
 $x_{i_1,j_1, i_2, j_2, \ldots, i_N, j_N}  = \displaystyle{\sum_{r_1,r_2,\ldots,r_{N-1}=1}^{R_1,R_2, \ldots,
 R_{N-1}}} \;\; g^{(1)}_{\;1,i_1, \,j_1, \,r_1} \; g^{(2)}_{\;r_1,\, i_2, \,j_2,\, r_2} \cdots \;
  g^{(N)}_{\;r_{N-1}, \, i_{N}, \,j_{N}, 1}$
  \\ & \\
 where $g^{(n)}_{\;r_{n-1}, \,i_n,\, j_n,\, r_n}$ are entries of a 4th-order core\\[4pt]
  \quad\quad $\underline \bG^{(n)} \in \Real^{ R_{n-1} \times I_n \times J_n \times R_n}$
  \\  &   \\[-6pt] \hline & \\[-6pt]
  \multicolumn{1}{c}{Slice (MPS) representation}\\
  & \\
  $x_{i_1, j_1, i_2, j_2,\ldots, i_N, j_N}  =  \bG^{(1)}_{i_1, j_1} \; \bG^{(2)}_{i_2, j_2} \cdots \bG^{(N)}_{i_N, j_N}$  $\quad$ where
  \\ & \\
 $\bG^{(n)}_{i_n, \,j_n} =\underline \bG^{(n)}(:,i_n,\,j_n,:) \in \Real^{R_{n-1} \times R_n}$ are slices of  $\underline \bG^{(n)} \in \Real^{ R_{n-1} \times I_n \times J_n \times R_n}$
\\  &  \\ \hline  
    \end{tabular*}\label{table:MPO}
    }}
 \vspace{12pt}
\end{table}

The strong Kronecker product representation of a TT is probably the most comprehensive  and useful form for displaying
tensor trains  in  their vector/matrix form, since it allows us to perform many  operations  using relatively small block matrices.

\noindent \textbf{Example.} For two matrices (in the TT format) expressed via the strong Kronecker products, $\bA = \tilde \bA^{(1)} \skron \tilde\bA^{(2)} \skron \cdots \skron \tilde\bA^{(N)} $ and $\;\;\bB = \tilde\bB^{(1)} \skron \tilde\bB^{(2)} \skron \cdots \skron \tilde\bB^{(N)} $, their Kronecker product can be efficiently computed  as
$\bA \otimes_L \bB = \tilde \bA^{(1)} \skron \cdots \skron \tilde \bA^{(N)} \skron \tilde \bB^{(1)} \skron \cdots \skron \tilde \bB^{(N)}$.
Furthermore, if the matrices $\bA$ and $\bB$ have the same mode sizes{\footnote{Note that,  wile original matrices $\bA \in \Real^{I_1 \cdots I_N \times J_1 \cdots J_N}$ and $\bB \in \Real^{I_1 \cdots I_N \times J_1 \cdots J_N}$ must have the same mode sizes, the corresponding core tenors, $\underline \bA^{(n)} =\in \Real^{R^A_{n-1} \times I_n \times J_n \times R^A_n}$ and $\underline \bB^{(n)}  =\in \Real^{R^B_{n-1} \times I_n \times J_n \times R^B_n}$, may have arbitrary mode sizes.}}, then their linear combination, $\bC = \alpha \bA + \beta \bB$ can be  compactly expressed as \cite{OseledetsTT11,Kazeev_Toeplitz13,Kazeev2013LRT}
\begin{equation}
\bC=[\tilde\bA^{(1)} \; \tilde\bB^{(1)}] \skron \left[\begin{array}{@{}c@{\hspace{1ex}}c@{}}\tilde\bA^{(2)}& \0 \\ \0 &\tilde\bB^{(2)}\end{array}\right] \skron
 \cdots \skron
\left[\begin{array}{@{}c@{\hspace{1ex}}c@{}}\tilde\bA^{(N-1)}&0\\0&\tilde\bB^{(N-1)}\end{array}\right] \skron \left[\begin{matrix} \alpha \tilde\bA^{(N)} \\ \beta \tilde\bB^{(N)} \end{matrix} \right]. \notag
\end{equation}

\noindent Consider its reshaped   tensor $\underline \bC = \llangle \underline \bC^{(1)},  \underline \bC^{(2)}, \ldots, \underline \bC^{(N)}\rrangle$ in the TT format; then its  cores $\underline \bC^{(n)} \in \Real^{R_{n-1} \times I_n \times J_n \times R_n}$, $n=1,2, \ldots,N$ can be expressed through their unfolding matrices,
$\bC^{(n)}_{<n>} \in \Real^{R_{n-1} I_n \times R_n J_n}$,  or equivalently by the lateral slices, $\bC^{(n)}_{i_n,j_n} \in \Real^{R_{n-1} \times R_n}$,  as follows
\be
\label{eq:TTadd2}
\bC_{i_n,j_n}^{(n)} &=& \left[\begin{matrix}\bA_{i_n,j_n}^{(n)}&\0 \\ \0 &\bB_{i_n,j_n}^{(n)}\end{matrix}\right], \;\; n=2,3,\ldots, N-1,
\ee
while  for the  border cores
\be
\bC_{i_1,j_1}^{(1)} =  \left[\bA_{i_1,j_1}^{(1)} \; \bB_{i_1,j_1}^{(1)} \right],   \qquad
\bC_{i_N,j_N}^{(N)} &=& \left[\begin{matrix} \alpha \; \bA_{i_N,j_N}^{(N)} \; \\ \; \beta \; \bB_{i_N,j_N}^{(N)} \end{matrix} \right]
\label{eq:TTadd22}
\ee
for $i_n=1,2,\ldots, I_n, \; j_n=1,2,\ldots,J_N, \; n=1,2, \ldots, N$.

Note that the various mathematical and graphical representations of TT/MPS and TT/MPO  can be used interchangeably for different purposes or applications.
With these representations, all basic mathematical
operations in TT format   can be performed on the constitutive  block matrices, even
without the need to explicitly construct core tensors \cite{OseledetsTT11,Dolgovth}.\\

\noindent{\bf Remark.} In the  TT/MPO paradigm, compression of large
matrices is not performed  by global  (standard) low-rank matrix approximations, but by low-rank approximations of
block-matrices (submatrices)  arranged in a hierarchical (linked) fashion.
However, to achieve a low-rank  TT and consequently a good compression ratio, ranks of  all  the corresponding unfolding matrices of  a specific structured data tensor
must be low, i.e., their singular values
must rapidly decrease to zero. While this is true for many structured matrices, unfortunately in general, this assumption does not hold. 

\section{Links Between CP, BTD Formats and TT/TC Formats}

It is important to note that any specific TN format can be converted into the TT format. This very useful  property is next illustrated for two
simple but important cases which establish links between the CP and TT and the BTD and TT formats.

 \begin{figure}
(a)
\begin{center}
\includegraphics[width=9.69cm]{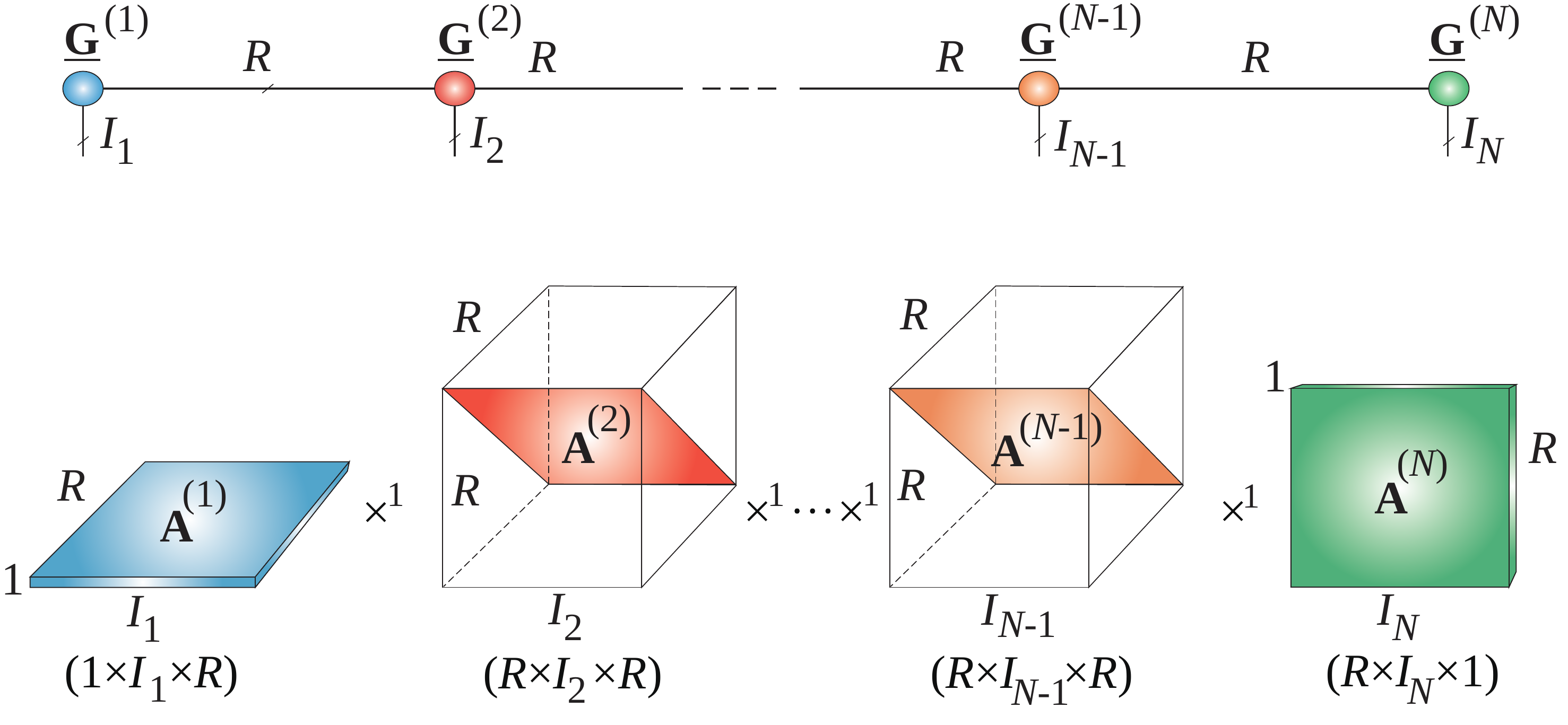}
\end{center}
(b)
\begin{center}
\includegraphics[width=11.8cm]{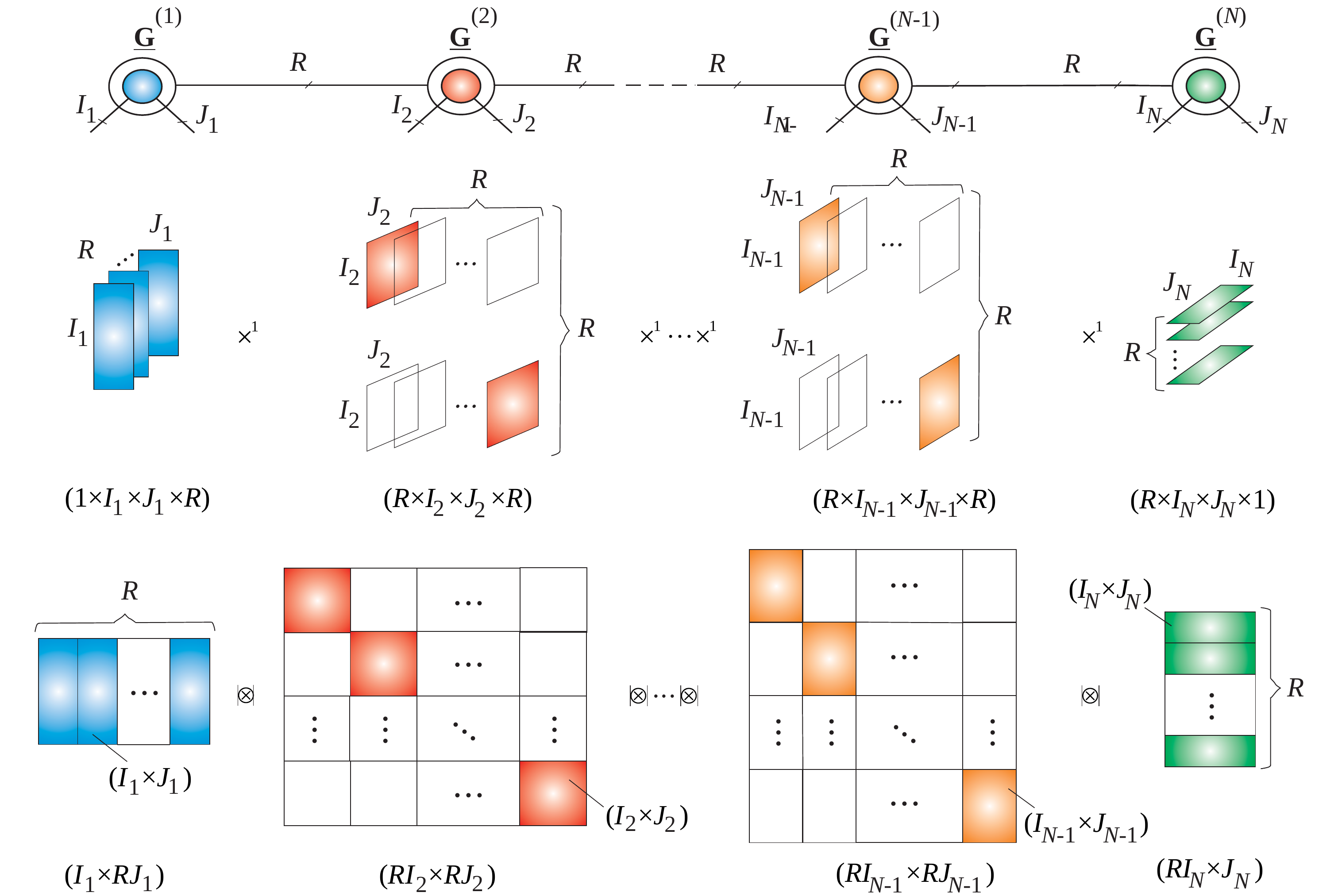}
\end{center}
\caption{Links between the  TT format and other tensor network formats. (a) Representation of the CP decomposition for an $N$th-order tensor, $\underline \bX = \underline \bI \times_1 \bA^{(1)} \times_2 \bA^{(2)}  \cdots \times_N \bA^{(N)}$, in the TT format.
 (b) Representation of the BTD model given by Eqs.~(\ref{HOPTA_AN}) and~(\ref{HOPTA_SK}) in the TT/MPO format. Observe that the TT-cores  are very sparse and the TT ranks are $\{R,R, \ldots, R\}$. Similar relationships can be established straightforwardly for the TC format.}
\label{Fig:CPTT}
\end{figure}

\begin{enumerate}
\item A tensor in  the CP format, given by
\be
\underline \bX  = \sum_{r=1}^R \ba_r^{(1)} \circ \ba_r^{(2)} \circ \cdots  \circ \ba_r^{(N)}, 
\ee
 can be  straightforwardly converted into the TT/MPS format as follows.
 Since each of the $R$ rank-1 tensors can be represented in the TT format of TT rank $(1,1,\ldots,1)$,
 using formulas (\ref{eq:TTadd2}) and (\ref{eq:TTadd22}), we have
 \be
\underline \bX  &=&  
\sum_{r=1}^R \llangle  \ba^{(1)\text{T}}_r, \ba^{(2)\text{T}}_r, \ldots, \ba^{(N)\text{T}}_r \rrangle \\ \notag
&=& \llangle  \underline \bG^{(1)}, \underline \bG^{(2)}, \ldots,  \underline \bG^{(N-1)}, \underline \bG^{(N)}\rrangle,
\ee
where the  TT-cores  $\underline  \bG^{(n)} \in \Real^{R \times I_n \times R}$ have
 diagonal lateral   slices $\bG^{(n)}(:, i_n, :) = \bG^{(n)}_{i_n} =\diag(a_{i_n,1},a_{i_n,2}, \ldots, a_{i_n,R}) \in \Real^{R \times R}$ for
 $n=2,3,\ldots,N-1$ and $\bG^{(1)}= \bA^{(1)} \in \Real^{I_1 \times R}$ and $\bG^{(N)}= \bA^{(N)\;\text{T}} \in \Real^{R \times I_N}$
 (see Figure  \ref{Fig:CPTT}(a)).

 \item A more general Block Term Decomposition (BTD)  for a $2N$th-order data tensor
\be
\underline \bX  = \sum_{r=1}^R (\bA_r^{(1)} \circ \bA_r^{(2)} \circ \cdots  \circ \bA_r^{(N)}) \in \Real^{I_1 \times J_1 \times \cdots \times I_N \times J_N}
\label{HOPTA_AN}
\ee
with full rank matrices, $\bA_r^{(n)} \in \Real^{I_n \times J_n},\; \forall r$,
  can be converted into a matrix  TT/MPO format, as illustrated in Figure  \ref{Fig:CPTT}(b).

  Note that (\ref{HOPTA_AN})  can be
  expressed in a matricized (unfolding) form  via strong Kronecker products of  block diagonal matrices (see  formulas (4.11)), 
 given by
 \be \label{HOPTA_SK}
 \bX &=& \sum_{r=1}^R (\bA_r^{(1)} \otimes_L \bA_r^{(2)} \otimes_L \cdots  \otimes_L \bA_r^{(N)}) \\
 &=&  \widetilde \bG^{(1)} \skron \widetilde \bG^{(2)} \skron \cdots \skron  \widetilde \bG^{(N)}
 \in \Real^{I_1 \cdots I_N  \; \times  \; J_1 \cdots \times J_N }, \notag
\ee
 with the TT rank, $R_n=R$ for $ n=1,2,\ldots N-1$, and the  block diagonal matrices,
 $ \widetilde \bG^{(n)} =\diag(\bA_1^{(n)},\bA_2^{(n)},\ldots,  \bA_R^{(n)}) \in \Real^{R I_{n} \times R J_n}$,  for $n=2,3,\ldots,N-1$, while
 $\widetilde \bG^{(1)} =[\bA_1^{(1)},\bA_2^{(1)},\ldots,  \bA_R^{(1)}] \in \Real^{I_1 \times RJ_1}$ is a row block matrix,  and
$ \widetilde \bG^{(N)} = \left[\begin{matrix} \bA_1^{(N)} \\ \vdots \\  \bA_R^{(N)} \end{matrix} \right] \in \Real^{R I_N  \times J_N}$  a column block matrix
(see Figure \ref{Fig:CPTT}(b)).

\end{enumerate}
Several algorithms exist for decompositions in the form  (\ref{HOPTA_AN}) and  (\ref{HOPTA_SK})
\cite{VisaSP-09,Batselier2014,Batselier2016}. In this way, TT/MPO decompositions for  huge-scale structured matrices can be  constructed indirectly.

 \section{Quantized Tensor Train (QTT) -- Blessing of Dimensionality}
\label{sect:QTT}

The procedure of creating a higher-order tensor from lower-order original data is referred to as tensorization,
 while in a special case where each mode has a very small size 2, 3 or 4, it is referred to as quantization.  In addition to vectors and matrices,  lower-order tensors can also be  reshaped  into higher-order tensors. By virtue of quantization, low-rank TN approximations with high compression ratios
  can be obtained, which is not possible to achieve with original raw data formats.  \cite{Oseledets10,Khoromskij-SC}.
%

\begin{figure}[t!]
\begin{center}
\includegraphics[width=11.2cm]{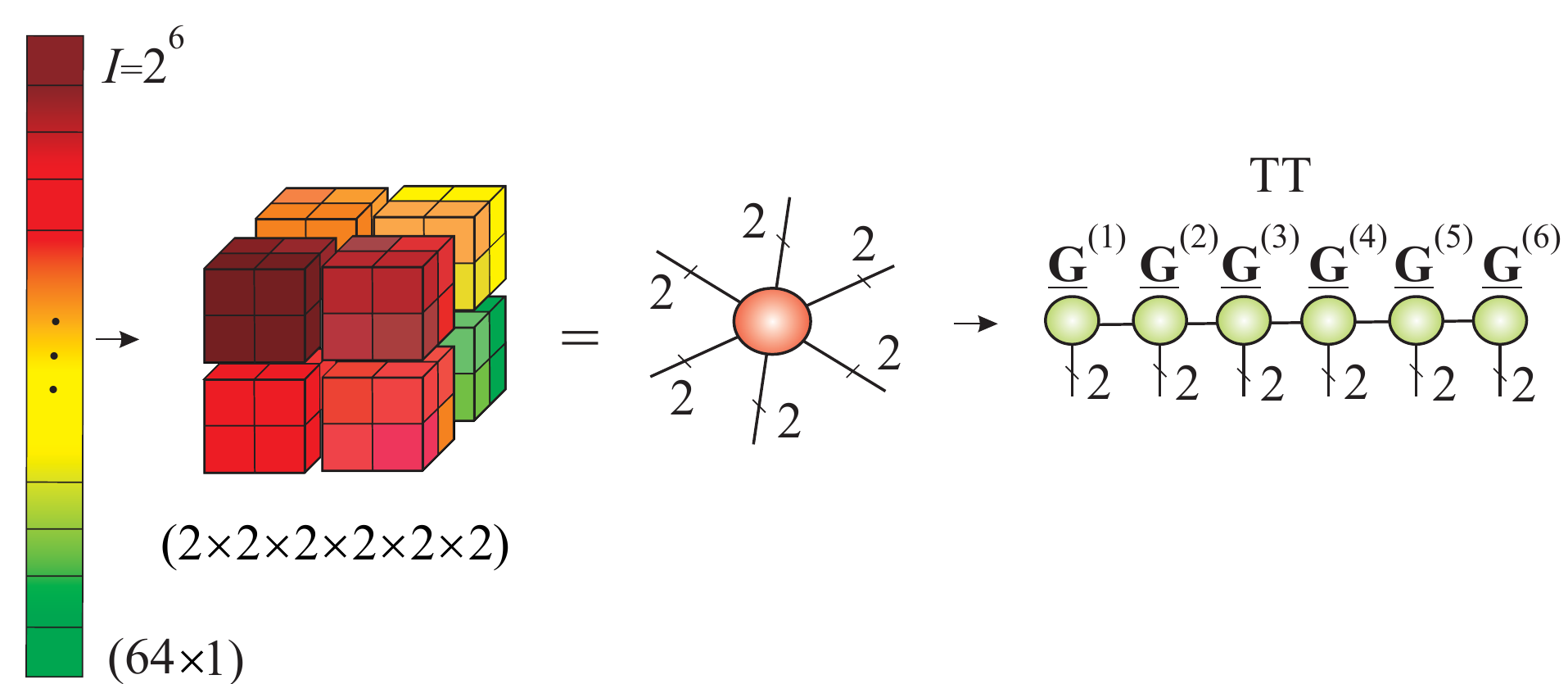}
\end{center}
\caption{Concept of tensorization/quantization of a large-scale vector
into a higher-order quantized tensor.
In order to achieve a good compression ratio, we need to apply a suitable  tensor decomposition such as the  quantized TT (QTT) using 3rd-order cores, $\underline \bX = \underline \bG^{(1)} \times^1 \underline \bG^{(2)} \times^1 \cdots \times^1  \underline \bG^{(6)}$.}
\label{Fig:Tensorization}
\end{figure}

Therefore, \emph{the quantization can be considered as a special form of tensorization where size of each mode is very small, typically 2  or 3}.
  The concept of quantized tensor networks (QTN) was first proposed in  \cite{Oseledets10} and  \cite{Khoromskij-SC},
whereby  low-size 3rd-order cores are sparsely interconnected via tensor contractions.
The so obtained model often provides an efficient,  highly compressed, and low-rank representation of a data tensor and helps to mitigate the curse of dimensionality, as illustrated below.\\

\noindent \textbf{Example.} The  quantization   of a huge vector, $\bx \in \Real^I$, $I = 2^K$, can be achieved through reshaping to give a
$(2 \times 2  \times \cdots \times 2)$  tensor $\underline \bX$ of order $K$,   as illustrated in Figure \ref{Fig:Tensorization}.
For structured data such a quantized tensor, $\underline \bX$, often admits low-rank TN approximation, so that a good compression of a huge vector $\bx$ can be achieved by enforcing the maximum possible low-rank structure on the tensor network. Even more generally,  an $N$th-order tensor, $\underline \bX \in \Real^{I_1 \times \cdots \times I_N}$, with  $I_n=q^{K_n}$,  can be quantized  in all modes simultaneously  to yield a
$(q \times q \times \cdots \times q)$  quantized tensor  of  higher-order and with small value of $q$. \\

\noindent \textbf{Example.} Since large-scale tensors (even of low-order) cannot be loaded directly into the computer memory, our approach to this problem is to represent the huge-scale data by tensor networks in a distributed and compressed TT format, so as  to avoid the explicit requirement for unfeasible large computer memory.

\begin{figure}
(a)\\
\begin{center}
\includegraphics[width=7.3cm]{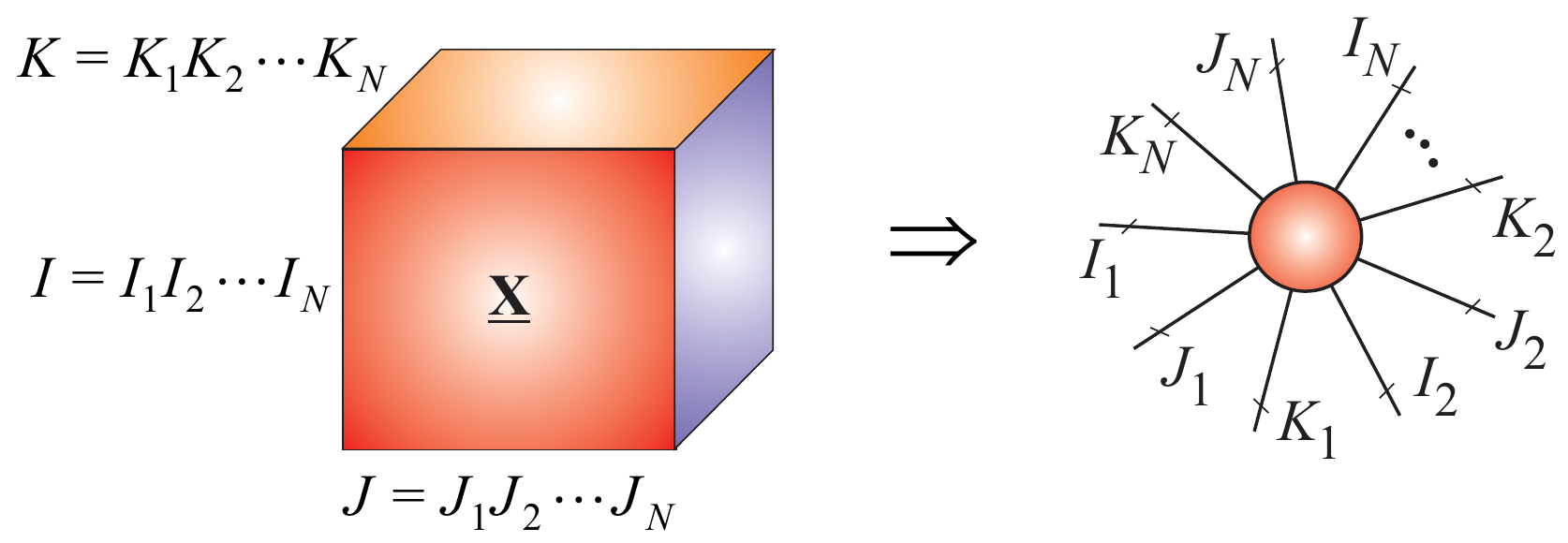}
\end{center}
(b)
\begin{center}
\includegraphics[width=11.8cm]{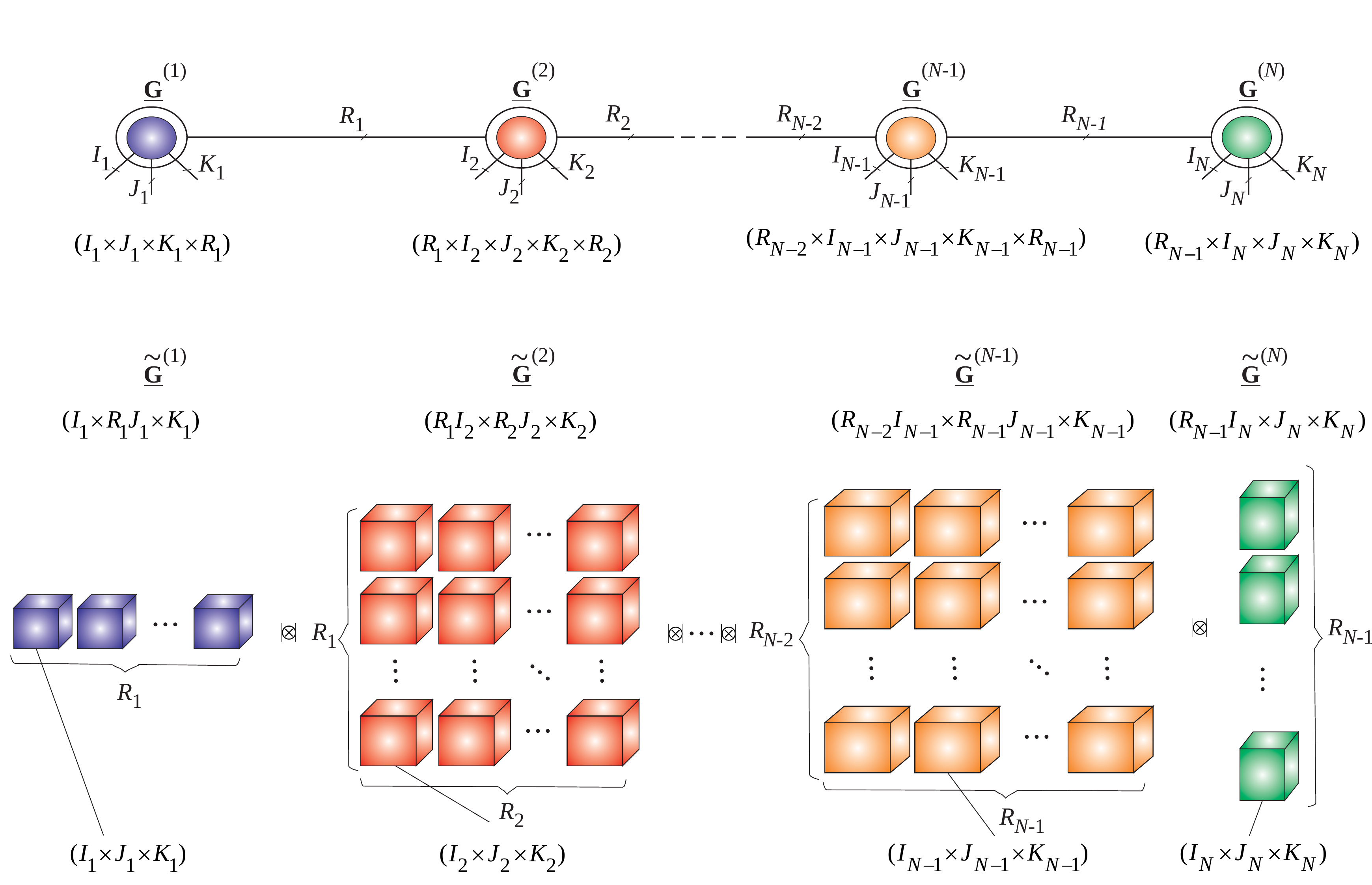}
\end{center}
\caption{Tensorization/quantization  of a huge-scale 3rd-order tensor into a higher order tensor and its TT representation. (a) Example  of tensorization/quantization of a 3rd-order tensor, $\underline \bX \in \Real^{ I \times J \times K}$, into a $3N$th-order tensor, assuming that the mode sizes can be factorized as, $I =I_1 I_2 \cdots I_N$, $J =J_1 J_2 \cdots J_N $ and $K =K_1 K_2 \cdots K_N$. (b) Decomposition of the high-order tensor via a generalized Tensor Train and its representation by the strong Kronecker product of block tensors as $\underline \bX \cong \underline {\widetilde \bG}^{(1)} \; \skron \; \underline {\widetilde \bG}^{(2)} \;  \skron  \cdots \skron \; \underline {\widetilde \bG}^{(N)}  \in \Real^{I_1 \cdots I_N \times  J_1 \cdots J_N  \times  K_1 \cdots K_N}$,
where each block  $\underline {\widetilde \bG}^{(n)} \in \Real^{R_{n-1} I_{n} \times R_n J_{n} \times K_{n}}$ is  also a 3rd-order tensor of size $(I_{n} \times J_{n} \times K_{n})$, for $n=1,2,\ldots,N$ with $R_0=R_N=1$. In the special case when $J=K=1$, the model simplifies into the standard TT/MPS model.} 
\label{Fig:TPS}
\end{figure}

 In the example shown in Figure \ref{Fig:TPS}, the tensor train of a huge 3rd-order tensor is expressed by the strong Kronecker products of block tensors with relatively small 3rd-order tensor blocks. The QTT
 is  mathematically represented in a distributed form via strong Kronecker products of  block 5th-order tensors.
Recall that  the strong Kronecker product of two block core tensors,  $\underline {\widetilde\bG}^{(n)} \in \Real^{R_{n-1} I_{n} \times R_n J_{n} \times K_{n}}$ and
 $\underline {\widetilde\bG}^{(n+1)} \in \Real^{R_n I_{n+1} \times R_{n+1} J_{n+1} \times K_{n+1}}$, is defined as the block tensor,
$\underline \bC = \underline {\widetilde\bG}^{(n)} \skron  \underline {\widetilde\bG}^{(n+1)} \in \Real^{R_{n-1}I_n I_{n+1} \times R_{n+1} J_n J_{n+1} \times K_n K_{n+1}}$, with 3rd-order tensor blocks, $\underline \bC_{r_{n-1},r_{n+1}}=\sum_{r_n=1}^{R_n} \underline \bG^{(n)}_{r_{n-1},r_n} \otimes_L  \underline \bG^{(n+1)}_{r_n,r_{n+1}} \in$ \\ $ \Real^{I_n I_{n+1} \times J_n J_{n+1} \times K_n K_{n+1}}$, where $\underline \bG^{(n)}_{r_{n-1},r_n} \in \Real^{I_{n} \times J_{n} \times K_{n}}$ and  $\underline\bG^{(n+1)}_{r_{n},r_{n+1}} \in \Real^{I_{n+1} \times J_{n+1} \times K_{n+1}}$ are the block tensors of $\underline {\widetilde\bG}^{(n)}$ and $\underline {\widetilde\bG}^{(n+1)}$, respectively.\\

 In practice, a fine ($q=2,3,4$ ) quantization is desirable to create as many virtual (additional) modes as possible, thus allowing for the  implementation of efficient low-rank tensor approximations.
For example, the binary encoding ($q=2$) reshapes an $N$th-order tensor with $(2^{K_1} \times 2^{K_2} \times \cdots \times 2^{K_N})$ elements into a tensor of order $(K_1+K_2 +\cdots + K_N)$, with the same number of elements.
In other words, the idea is to quantize each of the $n$ ``physical'' modes
(dimensions) by replacing them with $K_n$ ``virtual'' modes, provided that the corresponding mode sizes, $I_n$, are factorized as $I_n= I_{n,1} I_{n,2} \cdots I_{n, K_n}$. This, in turn, corresponds to reshaping the $n$th mode of size $I_n$ into $K_n$ modes of sizes $I_{n,1}, I_{n,2}, \ldots, I_{n,K_n}$.

The TT decomposition applied to quantized tensors is
referred to as the QTT, Quantics-TT or Quantized-TT, and was first introduced as a compression scheme for large-scale matrices  \cite{Oseledets10}, and also independently  for  more general settings.

\minrowclearance 2ex
\begin{table}[t!]
\caption{Storage complexities of  tensor decomposition models for an $N$th-order  tensor, $\underline \bX \in \Real^{I_1 \times I_2 \times \cdots \times I_N}$,
for which the original storage complexity is ${\cal{O}}(I^N)$, where $I=\max\{I_1,I_2,\ldots,I_N\}$, while $R$ is the upper bound on the ranks of tensor decompositions considered, that is,  $R=\max\{R_1,R_2,\ldots, R_{N-1} \}$ or $R=\max\{R_1,R_2,\ldots, R_N \}$.}
 \centering
  {\shadingbox{
\begin{tabular*}{0.99\linewidth}[t]{@{\extracolsep{\fill}}@{\hspace{2ex}}ll} \hline
1.   Full (raw) tensor format  &  ${\cal{O}}(I^N)$ \\
2.  CP   &  ${\cal{O}}(NIR)$ \\
3. Tucker   & ${\cal{O}}(NIR +R^N)$ \\
4. TT/MPS   & ${\cal{O}}(NIR^2)$ \\
5. TT/MPO   & ${\cal{O}}(NI^2R^2)$ \\
6.  Quantized TT/MPS (QTT)  &${\cal{O}}(N R^2 \log_q(I))$ \\
7.  QTT+Tucker  & ${\cal{O}}(N R^2 \log_q(I)+N R^3)$ \\
8. Hierarchical Tucker (HT) & ${\cal{O}}(N I R + N R^3)$  \\
\hline
    \end{tabular*}
   }}
\label{table_complexity}
\end{table}
\minrowclearance 0ex

\begin{figure}
(a)
\vspace{-0.1cm}
\begin{center}
\includegraphics[width=7.31cm]{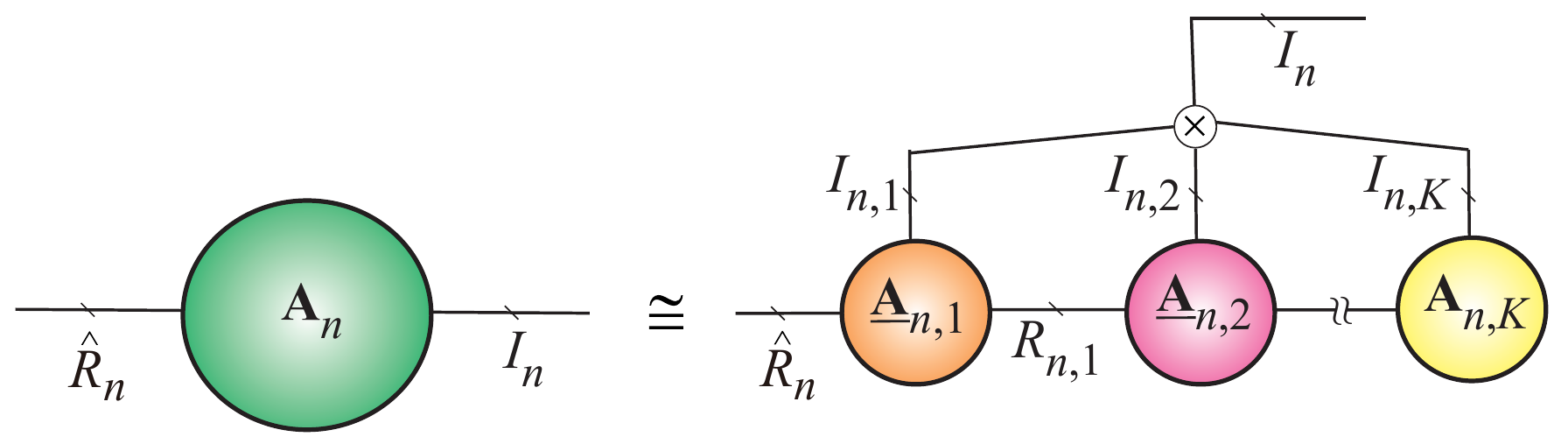}
\end{center}
(b)
\vspace{-0.1cm}
\begin{center}
\includegraphics[width=8.4915cm]{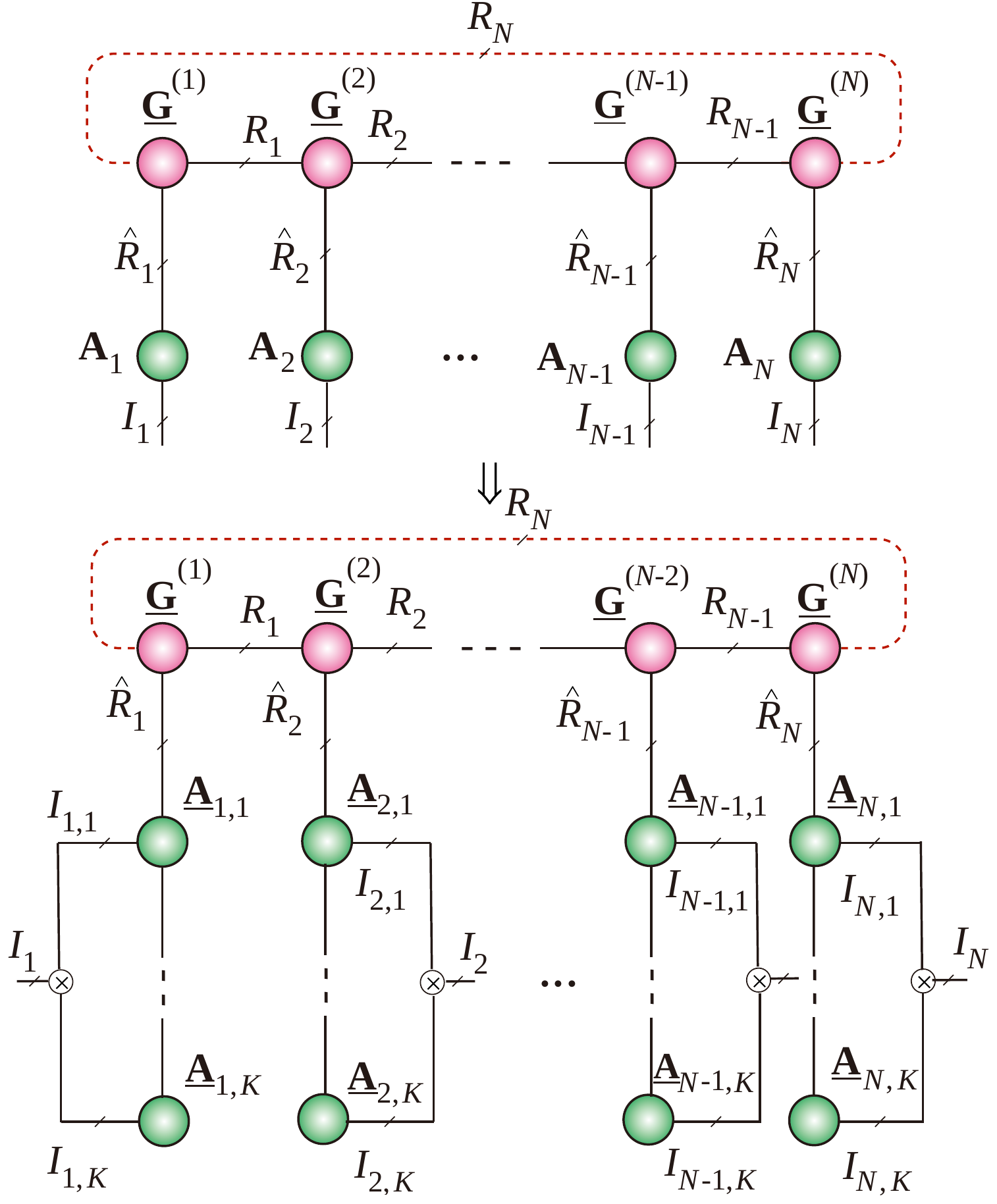}
\end{center}
\caption{The QTT-Tucker or alternatively QTC-Tucker (Quantized Tensor-Chain-Tucker) format. (a) Distributed representation of a  matrix $\bA_n \in \Real^{I_n \times \hat R_n}$ with a very large value of $I_n$ via QTT, by tensorization to a high-order quantized tensor, followed by QTT decomposition. (b) Distributed representation of a large-scale   Tucker-$N$  model,
$\underline \bX \cong \underline \bG \times_1 \bA_1 \times \bA_2  \cdots \times_N \bA_N$, via  a quantized TC model in
which the core tensor $\underline \bG \in \Real^{\hat R_1 \times \hat R_2 \times \cdots \times \hat R_N}$ and optionally all large-scale factor matrices $\bA_n$ ($n=1,2,\ldots,N$) are represented by MPS models (for more detail see  \cite{QTT_Tucker}).}
\label{Fig:QTT-Tucker}
\end{figure}

The attractive properties of QTT are:
\begin{enumerate}

\item Not only QTT ranks are typically small (usually, below 20) but they are also almost  independent\footnote{At least uniformly bounded.} of the data  size (even for  $I=2^{50}$), thus providing a logarithmic (sub-linear) reduction of storage requirements from  ${\cal{O}}(I^N)$ to  ${\cal{O}}(N R^2 \log_q(I))$ which is referred to as super-compression
\cite{Khoromskij-SC,QTT_Laplace,Kazeev_Toeplitz13,QTT_Tucker,dolgovEIG2013}.
Comparisons of the storage complexity of various tensor formats are given in Table \ref{table_complexity}.

\item Compared to the TT decomposition (without quantization), \emph{the QTT format often represents deep structures in the data by introducing
 ``virtual'' dimensions or modes}.  For data which exhibit high degrees of structure, the high compressibility of the QTT approximation is a consequence of the better separability properties of the quantized tensor.

\item The fact  that the QTT ranks are often moderate or even low\footnote{The TT/QTT ranks are constant or
growing linearly with respect to the tensor order $N$ and are constant or growing logarithmically with respect to the dimension of tensor modes $I$.}
offers unique advantages in the context of big data analytics (see \cite{Khoromskij-SC,Khoromskij-TT,Kazeev_Toeplitz13} and references therein),
together with  high efficiency
of multilinear algebra within the TT/QTT algorithms which rests upon the well-posedness of
the  low-rank TT  approximations.

\end{enumerate}
The ranks of the  QTT format often grow dramatically with data size,
but with a linear increase in the approximation accuracy. To overcome this
problem, Dolgov and Khoromskij proposed the QTT-Tucker format
\cite{QTT_Tucker} (see Figure~\ref{Fig:QTT-Tucker}), which
 exploits the TT approximation not  only for the Tucker core tensor, but also   for the  factor matrices.
This model naturally admits distributed computation, and often  yields  bounded ranks, thus avoiding the curse of dimensionality.

The TT/QTT tensor networks  have already found  application in very large-scale problems in scientific computing, such  as in eigenanalysis, super-fast Fourier transforms, and in solving huge systems of large linear equations  (see \cite{QTT_Tucker,Huckle2013,dolgovEIG2013,wahls2014learning,KressnerEIG2014,KressnerEVD2016} and references therein).

\section{Basic Operations in TT Formats}
\label{sect:TT-oper}

For big tensors in their TT formats, basic mathematical operations, such as the addition, inner product, computation of tensor norms, Hadamard and Kronecker product, and matrix-by-vector and matrix-by-matrix multiplications can be very efficiently performed  using  block (slice)  matrices of individual
(relatively small size) core tensors.\\

\noindent Consider two $N$th-order tensors in the TT format
\be
\underline \bX &=& \llangle \underline \bX^{(1)}, \underline \bX^{(2)}, \ldots , \underline \bX^{(N)} \rrangle \in \Real^{I_1 \times I_2 \times \cdots \times I_N} \nonumber \\
\underline \bY &=& \llangle \underline \bY^{(1)}, \underline \bY^{(2)}, \ldots , \underline \bY^{(N)}\rrangle \in \Real^{I_1 \times I_2 \times \cdots \times I_N}, \nonumber
\ee
for which the TT ranks are $\brr_X =\{R_1,R_2, \ldots, R_{N-1}\}$ and $\brr_Y =\{\tilde R_1, \tilde R_2, \ldots, \tilde R_{N-1}\}$.  The following operations can  then be  performed  directly in the TT formats.        \\

\noindent{\bf Tensor addition.} The sum of two tensors
\be \underline \bZ =\underline \bX + \underline \bY = \llangle \underline \bZ^{(1)}, \underline \bZ^{(2)}, \ldots , \underline \bZ^{(N)}\rrangle \in \Real^{I_1 \times I_2 \times \cdots \times I_N}
\ee
 has the TT rank $\brr_Z = \brr_X+\brr_Y$ and can be  expressed via lateral slices of the cores $\underline \bZ \in \Real^{R_{n-1} \times I_n \times R_n}$ as
\be
\bZ_{i_n}^{(n)} &=& \left[\begin{matrix}\bX_{i_n}^{(n)}&\0 \\ \0 &\bY_{i_n}^{(n)}\end{matrix}\right], \;\; n=2,3,\ldots, N-1.
\label{eq:TTadd2ad}
\ee
For the border cores, we have
\be
\bZ_{i_1}^{(1)} =  \left[\bX_{i_1}^{(1)} \; \bY_{i_1}^{(1)} \right],   \qquad
\bZ_{i_N}^{(N)} &=& \left[\begin{matrix} \bX_{i_N}^{(N)} \; \\ \;  \bY_{i_N}^{(N)} \end{matrix} \right]
\ee
for $i_n=1,2,\ldots, I_n, \;\ n=1,2, \ldots, N$.\\

\noindent{\bf Hadamard product.} The computation of the Hadamard  (element-wise) product, $\underline \bZ = \underline \bX \* \underline \bY$, of two tensors, $\bX$ and $\bY$, of the same order and the same size
  can be  performed very efficiently in the TT format by expressing  the slices
of the cores, $\underline \bZ \in \Real^{R_{n-1} \times I_n \times R_n}$, as
\be
\bZ_{i_n}^{(n)} = \bX_{i_n}^{(n)} \otimes \bY_{i_n}^{(n)}, \quad n=1, \ldots, N, \;\; i_n=1,\ldots,I_n.
\ee
This increases the TT ranks for the tensor $\underline \bZ$ to at most $R_n \tilde R_n$, $n=1,2, \ldots, N$,
but the associated computational complexity can be  reduced from being exponential in $N$, ${\cal {O}}(I^N)$,
to being linear in both $I$ and $N$, ${\cal {O}}(I N (R \tilde R)^2))$.\\

 \noindent{\bf Super fast Fourier transform} of a tensor in the TT format (MATLAB functions:  fftn$(\underline \bX)$ and  fft$(\underline \bX^{(n)}, [], 2)$) can be computed as
\be
 {\cal {F}} (\underline \bX) &=& \llangle   {\cal {F}} (\underline \bX^{(1)}), {\cal {F}} (\underline \bX^{(2)}), \ldots, {\cal {F}} (\underline \bX^{(N)})\rrangle \notag \\
 &=& {\cal {F}} (\underline \bX^{(1)}) \times^1 {\cal {F}} (\underline \bX^{(2)})  \times^1  \cdots  \times^1  {\cal {F}} (\underline \bX^{(N)}).
 \label{eq:DFT_TT}
 \ee
 It should be emphasized that performing computation of the FFT on relatively  small core tensors $\underline \bX^{(n)} \in \Real^{R_{n-1} \times I_n \times R_n}$ \emph{reduces dramatically computational complexity under condition that a data tensor admits low-rank TT approximation}.  This approach is referred to as the
 super fast Fourier transform (SFFT) in TT format. Wavelets, DCT, and other linear integral transformations admit a similar form to the SFFT in (\ref{eq:DFT_TT}),
 for example, for the wavelet transform in the TT format, we have
 \be
 {\cal {W}} (\underline \bX) &=& \llangle   {\cal {W}} (\underline \bX^{(1)}), {\cal {W}} (\underline \bX^{(2)}), \ldots, {\cal {W}} (\underline \bX^{(N)})\rrangle \notag \\
 &=& {\cal {W}} (\underline \bX^{(1)}) \times^1 {\cal {W}} (\underline \bX^{(2)})  \times^1  \cdots  \times^1  {\cal {W}} (\underline \bX^{(N)}).
 \label{eq:WT_TT}
 \ee

\noindent{\bf The N-D discrete convolution in a TT format} of  tensors  $\underline \bX \in \Real^{I_1 \times \cdots \times I_N}$
 with TT rank $\{R_1, R_2,\ldots, R_{N-1}\}$ and $\underline \bY \in \Real^{J_1 \times \cdots \times J_N}$
  with TT rank $\{Q_1,Q_2, \ldots, Q_{N-1}\}$
can be computed as
\be
\underline \bZ &=& \underline \bX \ast \underline \bY  \\
&=& \llangle \underline \bZ^{(1)}, \underline \bZ^{(2)}, \ldots, \underline \bZ^{(N)} \rrangle
\in \Real^{(I_1+J_1-1) \times (I_2+J_2-1) \times \cdots \times (I_N+J_N-1)}, \notag
 \ee
  with the  TT-cores given by
 \be
 \underline \bZ^{(n)} =\underline \bX^{(n)} \boxdot_2  \underline \bY^{(n)}
 \in \Real^{(R_{n-1} Q_{n-1}) \times (I_n+J_n-1) \times (R_n Q_n)},
 \ee
or, equivalently, using the standard convolution  $\underline \bZ^{(n)}(s_{n-1},:,s_n) = \underline \bX^{(n)}(r_{n-1}, :,r_n) \ast
\underline \bY^{(n)}(q_{n-1}, :, q_n) \in \Real^{(I_n+J_n-1)}$ for $s_n=1,2, \ldots, R_n Q_n$  and  $n=1,2, \ldots, N$, $\;R_0=R_N=1$.\\

\noindent{\bf Inner product.}  The  computation of the inner (scalar, dot)  product of two $N$th-order tensors,
$\underline \bX = \llangle \underline \bX^{(1)}, \underline \bX^{(2)}, \ldots , \underline \bX^{(N)} \rrangle \in \Real^{I_1 \times I_2 \times \cdots \times I_N}$  and
 $\bY = \llangle \underline \bY^{(1)}, \underline \bY^{(2)}, \ldots , \underline \bY^{(N)}\rrangle \in \Real^{I_1 \times I_2 \times \cdots \times I_N}$, is given by
\be
\langle \underline \bX, \underline \bY \rangle &=&  \langle \mbox{vec}(\underline \bX), \mbox{vec}(\underline \bY) \rangle \\
&=& \sum_{i_1=1}^{I_1} \cdots \sum_{i_N=1}^{I_N} x_{i_1 \dots i_n}  \; y_{i_1 \cdots i_N} \notag
\ee
 and has the complexity of ${\cal{O}}(I^N)$ in the raw tensor format.
 In TT formats, the inner product
 can be computed with the reduced  complexity of only
 ${\cal{O}}(NI(R^2 \tilde R + R \tilde R^2))$ when the
 inner product is  calculated by moving TT-cores from left to right and performing calculations on relatively small matrices,
 $\bS_n = \underline \bX^{(n)} \times^{1,2}_{1,2} (\underline \bY^{(n)} \times_1 \bS_{n-1}) \in \Real^{R_n \times \widetilde R_n}$ for $n=1,2, \ldots,N$. The results are then sequentially multiplied by the next core $\underline \bY^{(n+1)}$ (see  Algorithm \ref{alg:TT-inner-prod}).\\

\noindent{\bf Computation of the Frobenius norm.} In a similar way, we can efficiently compute the Frobenius  norm of a tensor, $\|\underline \bX\|_F =\sqrt{\langle \underline \bX, \underline \bX \rangle}$, in the TT format. For the so-called
  $n$-orthogonal{\footnote{An $N$th-order tensor $\underline \bX =\llangle \underline \bX^{(1)},\underline \bX^{(2)} \ldots, \underline \bX^{(N)} \rrangle$
   in the TT format is  called $n$-orthogonal if all the cores to the left of the core $\underline \bX^{(n)}$  are left-orthogonalized  and all the cores to the right of the
    core $\bX^{(n)}$ are  right-orthogonalized (see Part 2 for more detail).}}
      TT format, it is easy to show that
      \be
      \|\underline \bX\|_F =  \|\underline \bX^{(n)}\|_F.
      \ee

  \begin{algorithm}[t]
{\small
\caption{\textbf{Inner product of two large-scale tensors in the TT Format \cite{OseledetsTT11,Dolgovth}}}
\label{alg:TT-inner-prod}
 \begin{algorithmic}[1] 
\REQUIRE $N$th-order tensors, $\underline \bX = \llangle \underline \bX^{(1)}, \underline \bX^{(2)}, \ldots , \underline \bX^{(N)} \rrangle \in \Real^{I_1 \times I_2 \times \cdots \times I_N}$  \\ and
 $\bY = \llangle \underline \bY^{(1)}, \underline \bY^{(2)}, \ldots , \underline \bY^{(N)}\rrangle \in \Real^{I_1 \times I_2 \times \cdots \times I_N}$ in TT formats, with \\ TT-cores  $\underline \bX \in \Real^{R_{n-1} \times I_n \times R_{n}}$ and $\underline \bY \in \Real^{\widetilde R_{n-1} \times I_n \times \widetilde R_{n}}$ \\
 and $R_0=\widetilde R_0 = R_N =\widetilde R_N =1$
\ENSURE Inner product $\langle\underline \bX, \underline \bY\rangle = \mbox{vec}(\underline \bX)^{\text{T}} \mbox{vec}(\underline \bY)$
\STATE Initialization $\bS_0=1$
\FOR{$n=1$ to $N$}
    \STATE $\bZ^{(n)}_{(1)} =\bS_{n-1} \bY^{(n)}_{(1)} \in \Real^{R_{n-1} \times I_n \widetilde R_{n}}$
    \STATE $\bS_{n} =\bX^{(n)\;\text{T}}_{<2>} \bZ^{(n)}_{<2>} \in \Real^{R_{n} \times \widetilde R_n}$
\ENDFOR
\RETURN Scalar $\langle \underline \bX, \underline \bY \rangle = \bS_{N}  \in \Real^{R_{N} \times \widetilde R_{N}}=\Real$, with $R_N=\widetilde R_N=1$
\end{algorithmic}
}
\end{algorithm}

\begin{figure}
(a)
\begin{center}
\includegraphics[width=4.5cm,height=6cm]{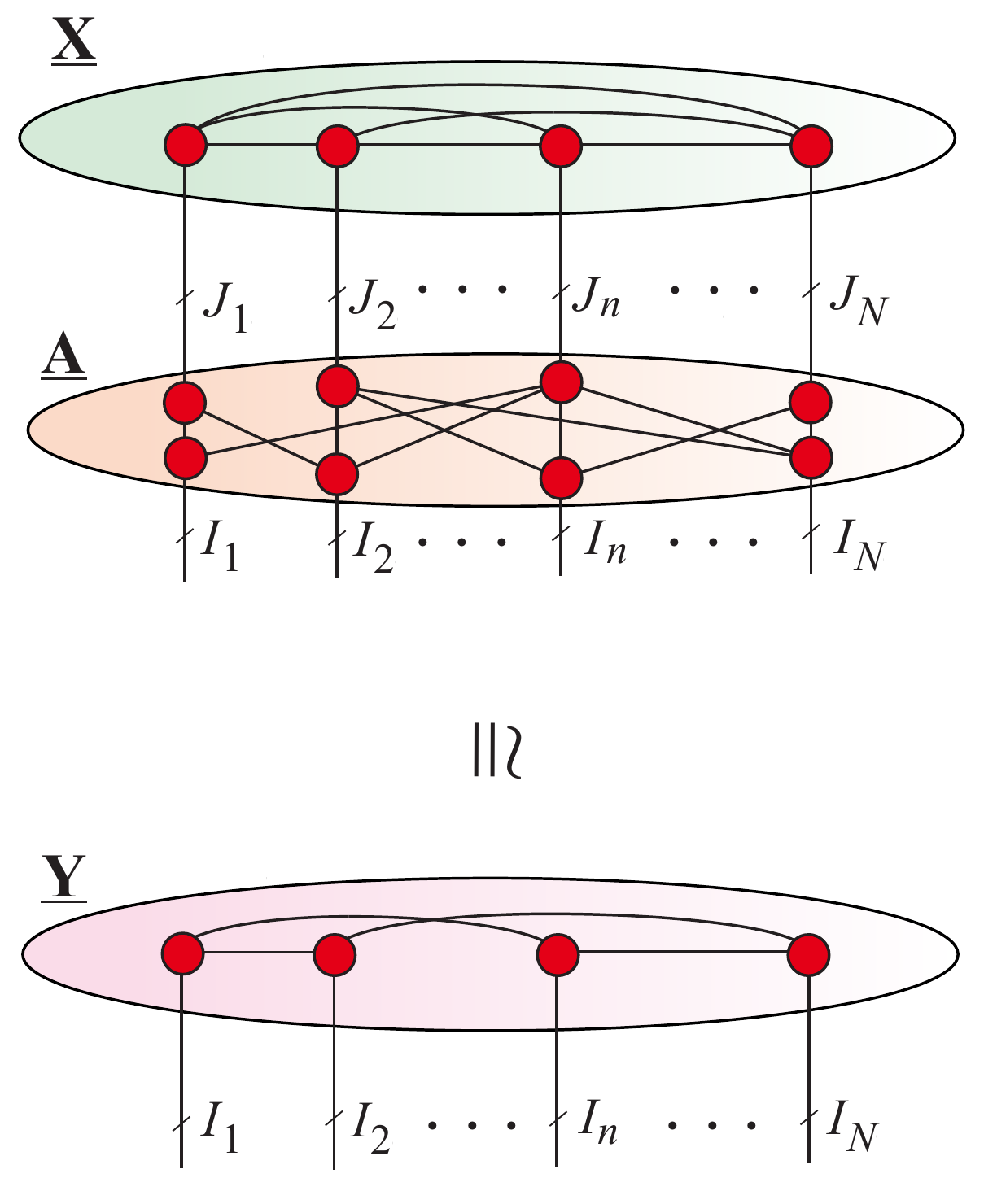} \hspace{0.2cm} \includegraphics[width=6.8cm,height=6cm]{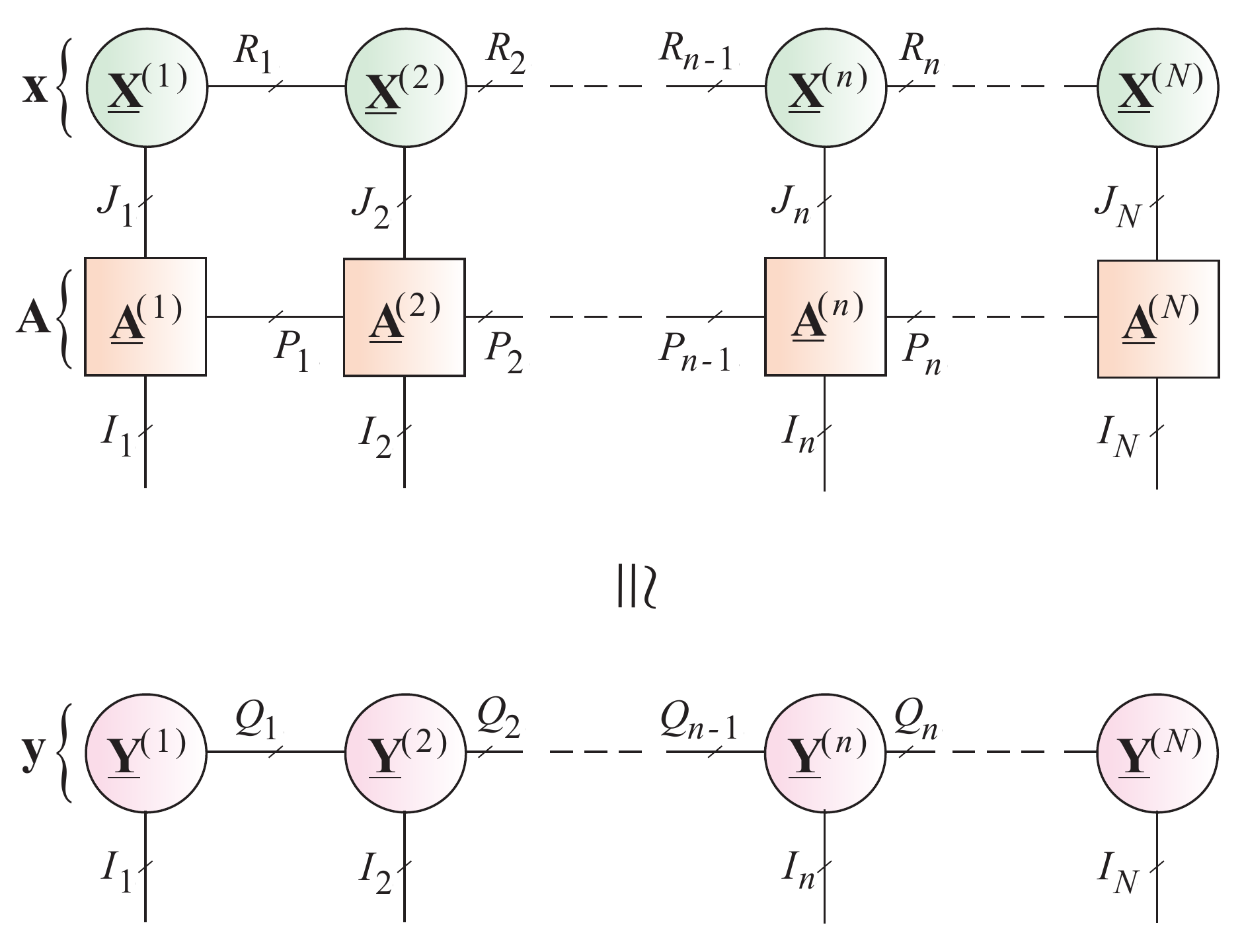}
\end{center}
\vspace{0.5cm}
(b)
\begin{center}
\includegraphics[width=4.7cm,height=6cm]{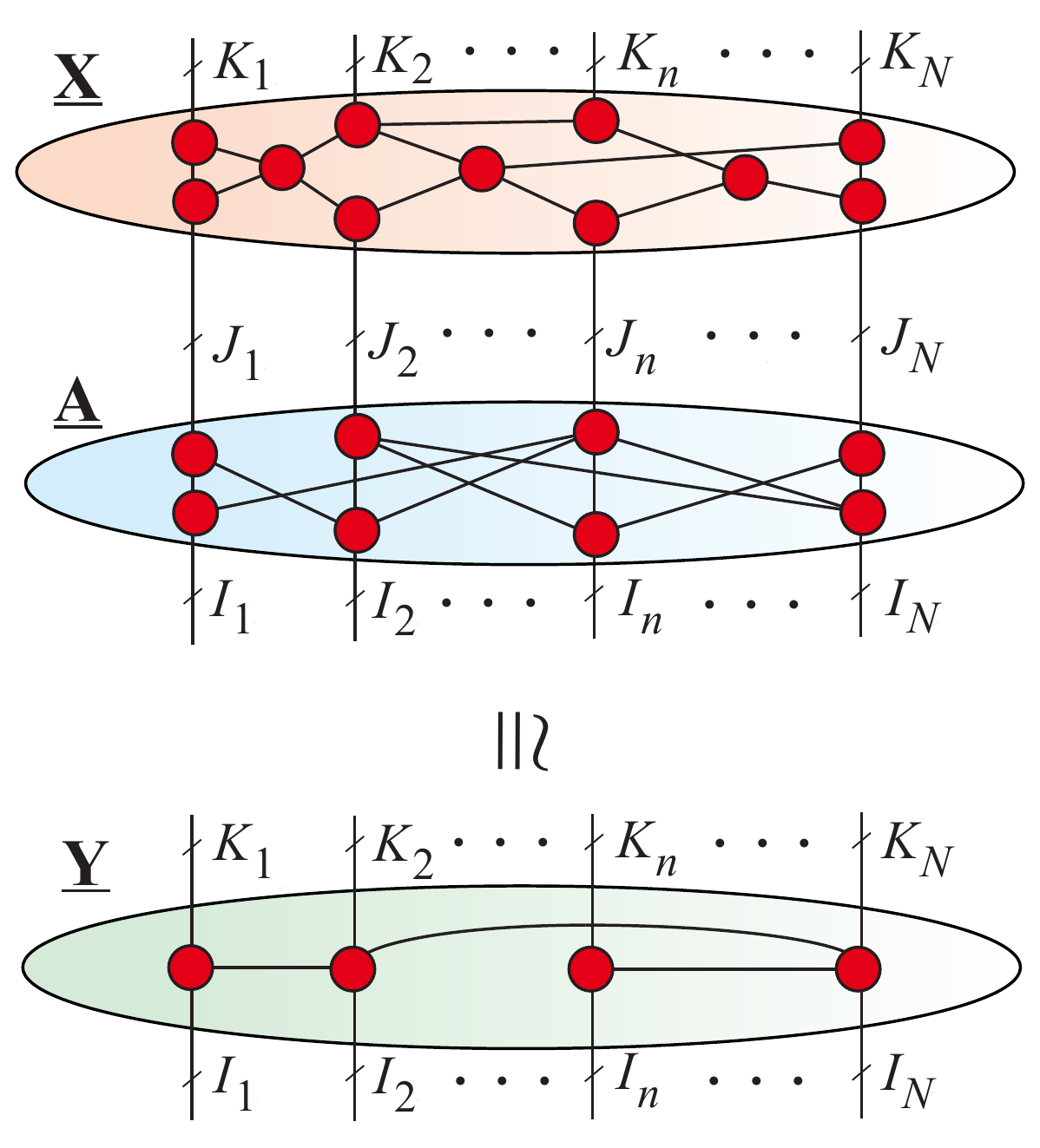} \hspace{0.4cm} \includegraphics[width=6.5cm,height=6cm]{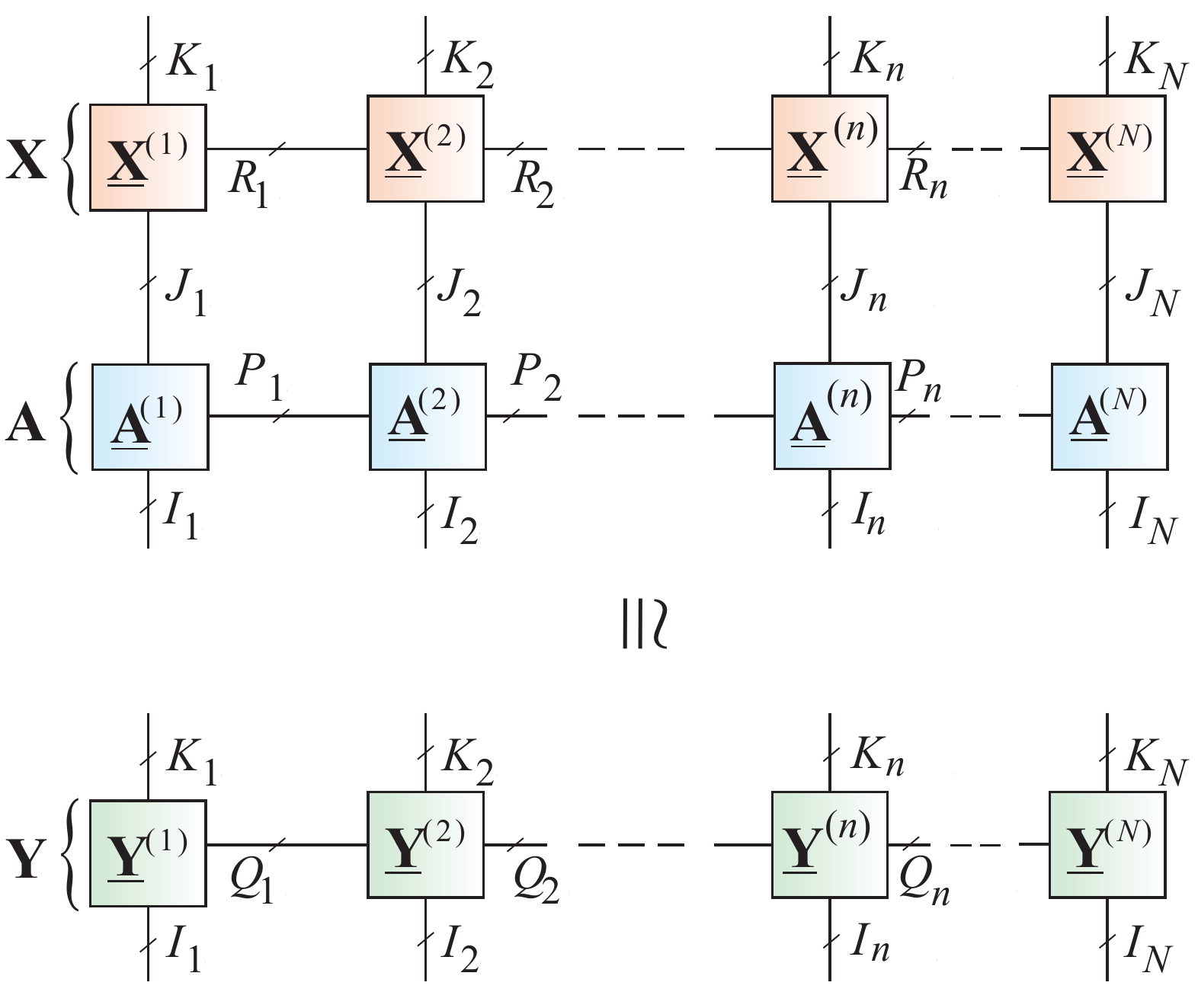}
\end{center}
\caption{Linear systems represented by arbitrary  tensor networks ({\it{left}}) and  TT networks ({\it{right}}) for  (a) $\bA \bx \cong \by$ and (b) $\bA \bX \cong \bY$.}
\label{Fig:YAX}
\vspace{18pt}
\end{figure}

\begin{figure}[t]
(a)\vspace{-0.1cm}
\begin{center}
\includegraphics[width=4.5cm]{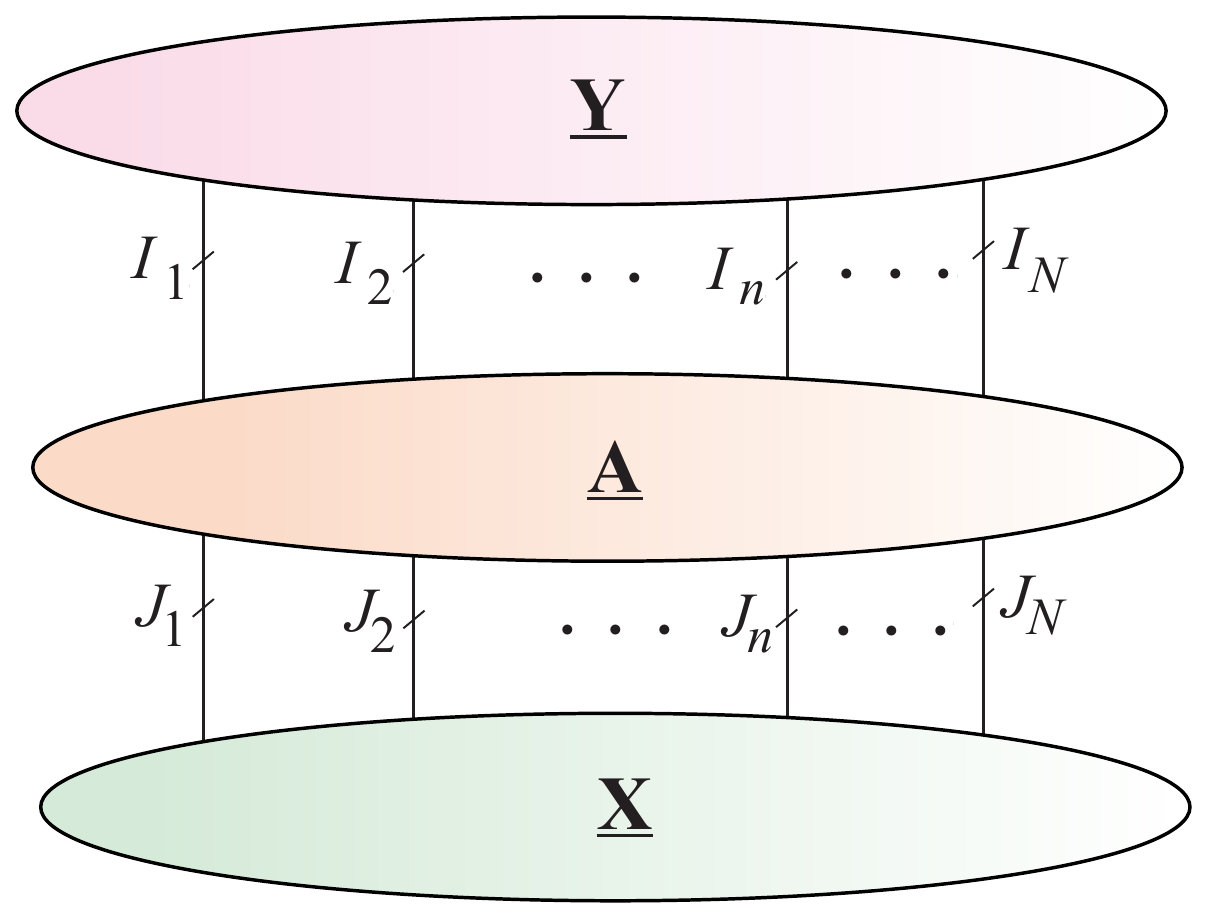} \hspace{0.2cm}
\includegraphics[width=6.7cm]{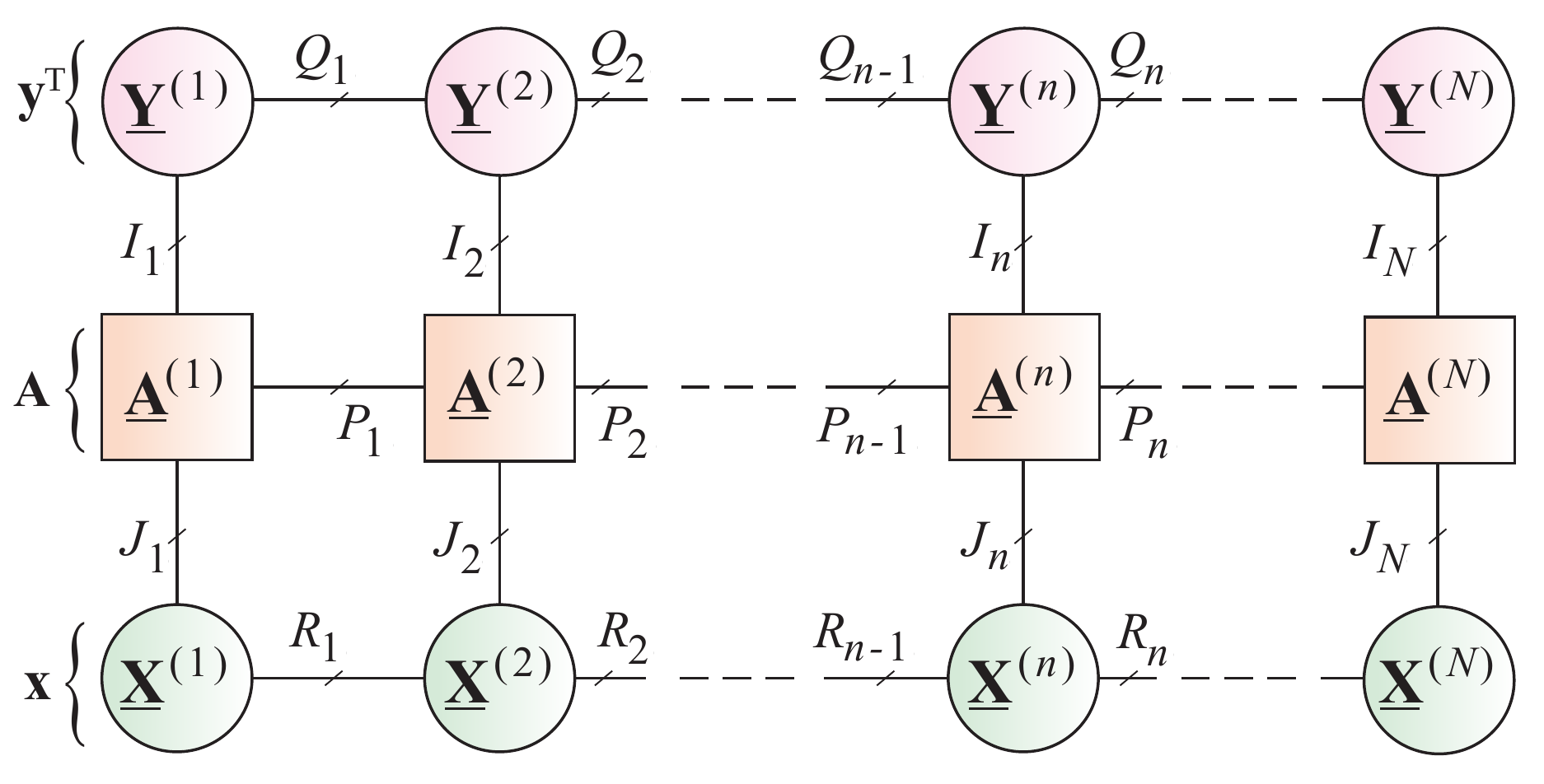}
\end{center}
(b)\vspace{-0.1cm}
\begin{center}
\includegraphics[width=4.4cm]{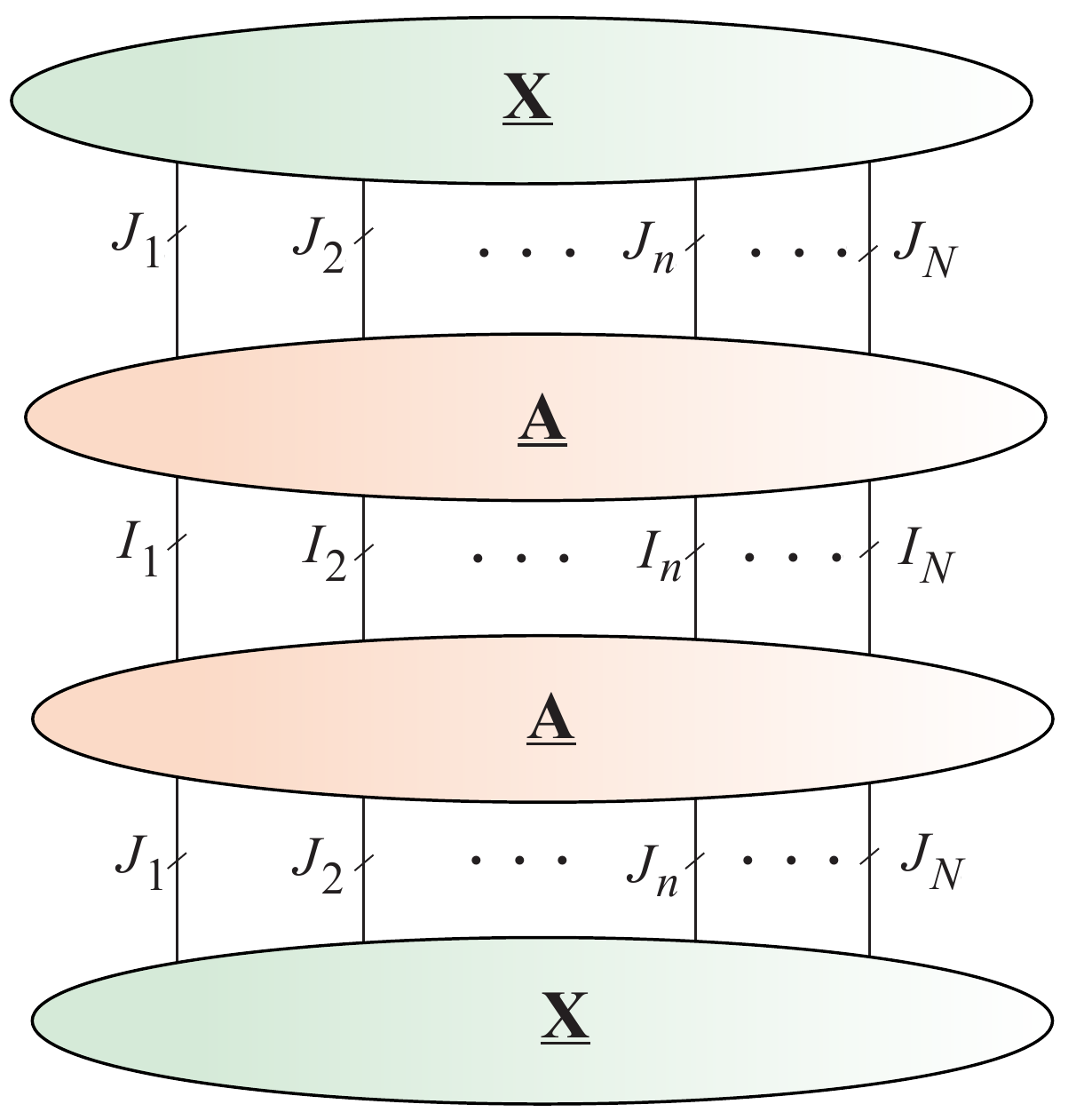} \hspace{0.2cm}
\includegraphics[width=6.6cm]{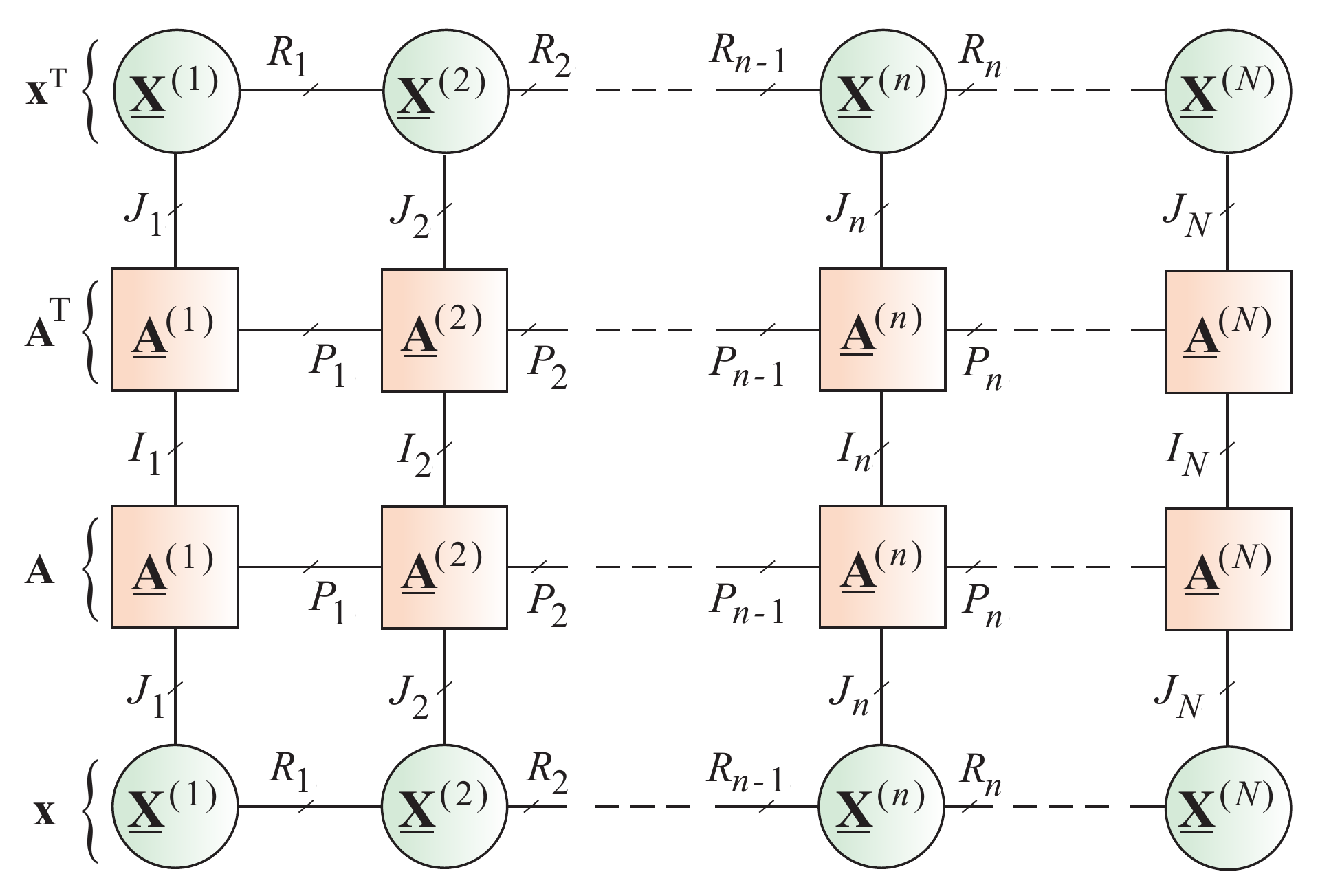}
\end{center}
\caption{Representation of typical cost functions by arbitrary TNs and  by TT networks: (a) $J_1(\bx)=\by^{\text{T}} \bA \bx$ and (b) $J_2(\bx)= \bx^{\text{T}} \bA^{\text{T}} \bA \bx$.
Note that tensors
$\underline \bA$, $\underline \bX$ and $\underline \bY$ can be, in general,
 approximated by any TNs that  provide good low-rank representations.}
\label{Fig:ATAX=Y}
\end{figure}

\noindent{\bf Matrix-by-vector multiplication.} Consider a huge-scale  matrix equation
(see Figure \ref{Fig:YAX} and Figure \ref{Fig:ATAX=Y})
\be
\bA \bx = \by,
\label{Axy}
\ee
where $\bA \in \Real^{I \times J}$, $\;\;\bx \in \Real^J$ and $\by \in \Real^I$ are represented approximately in the TT format, with $I=I_1 I_2 \cdots I_N$ and $J=J_1 J_2 \cdots J_N$. As shown in Figure \ref{Fig:YAX}(a), the cores are defined as
$\underline \bA^{(n)} \in \Real^{P_{n-1} \times I_n \times J_n \times P_n}$,
$\; \underline \bX^{(n)} \in \Real^{R_{n-1} \times J_n \times R_n}$ and
$\;\underline \bY^{(n)} \in \Real^{Q_{n-1} \times I_n  \times Q_n}$.

Upon representing the entries of the matrix $\bA$ and vectors $\bx$ and $\by$ in their tensorized forms, given by
\be
 \underline \bA  &=& \displaystyle{\sum_{p_1,p_2,\ldots,p_{N-1}=1}^{P_1,P_2, \ldots,
 P_{N-1}}} \;\; \bA^{(1)}_{\;1,p_1} \; \circ \; \bA^{(2)}_{\;p_1, p_2}  \circ \cdots
 \circ \; \bA^{(N)}_{\;p_{N-1},1} \notag\\
\underline \bX  &=& \displaystyle{\sum_{r_1,r_2,\ldots,r_{N-1}=1}^{R_1,R_2, \ldots,
 R_{N-1}}} \;\; \bx^{(1)}_{\;r_1} \; \circ \; \bx^{(2)}_{\;r_1, r_2}  \circ \cdots
 \circ \; \bx^{(N)}_{\;r_{N-1}}\\
 \underline \bY & =& \displaystyle{\sum_{q_1,q_2,\ldots,q_{N-1}=1}^{Q_1,Q_2, \ldots,
 Q_{N-1}}} \;\; \by^{(1)}_{\;q_1} \; \circ \; \by^{(2)}_{\;q_1, q_2}  \circ \cdots
 \circ \; \by^{(N)}_{\;q_{N-1}}, \notag
 \label{Ax=y-outer}
\ee
we arrive at a simple formula for the tubes of the tensor $\underline \bY$, in the form
\be
 \by^{(n)}_{q_{n-1},q_n} = \by^{(n)}_{\overline{r_{n-1} \, p_{n-1}},\;\overline{r_n \, p_n}}
= \bA^{(n)}_{p_{n-1},\,p_n}  \; \bx^{(n)}_{r_{n-1},\,r_n}
 \in \Real^{I_n},  \nonumber
\ee
with $Q_n=P_n \,R_n$ for $n=1,2,\ldots,N$.

Furthermore, by representing the  matrix $\bA$ and vectors $\bx$, $\by$ via the strong Kronecker products
\be
\bA &=& \tilde \bA^{(1)} \skron \tilde \bA^{(2)} \skron \cdots \skron \tilde \bA^{(N)} \notag\\
\bx &=& \tilde \bX^{(1)} \skron \tilde \bX^{(2)} \skron \cdots \skron \tilde \bX^{(N)}\\
\by &=& \tilde \bY^{(1)} \skron \tilde \bY^{(2)} \skron \cdots \skron \tilde \bY^{(N)}, \notag
\ee
with $\tilde \bA^{(n)} \in \Real^{P_{n-1}  I_n \times J_n P_n}, \;\;
\tilde \bX^{(n)} \in \Real^{R_{n-1}  J_n \times R_n}$ and $
\tilde \bY^{(n)} \in \Real^{Q_{n-1} I_n  \times Q_n}$,
we can establish a simple  relationship
\begin{equation}
\tilde \bY^{(n)} =\tilde \bA^{(n)} \; \Cprod \; \tilde \bX^{(n)} \in
\Real^{R_{n-1} \,P_{n-1} \,I_n \times R_{n} \,P_{n}}, \quad  n=1,\ldots,N,
\label{ACproduct}
\end{equation}
where the operator $\Cprod$ represents the C (Core) product of two block matrices.

\begin{figure}[t]
\centering
\includegraphics[width=10.5cm]{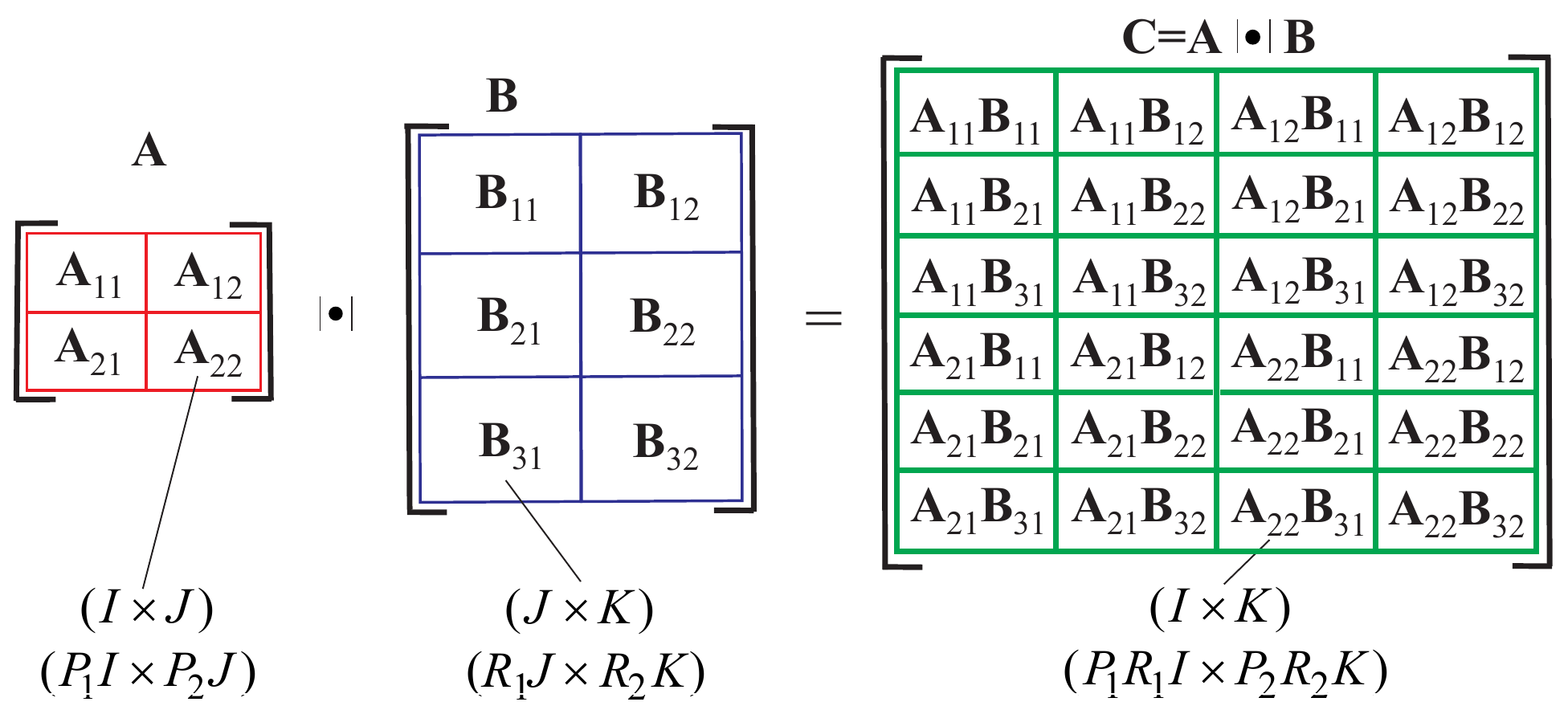}
\caption{Graphical illustration of the C product of two block matrices.}
\label{Fig:SKPAC-product}
\end{figure}

The  C product of a block  matrix $\bA^{(n)} \in \Real^{P_{n-1} I_n \times P_n J_n}$
 with blocks $\bA^{(n)}_{p_{n-1},p_n} \in \Real^{I_n \times J_n}$, and a block matrix $\bB^{(n)} \in \Real^{R_{n-1} J_n \times R_n K_n}$, with  blocks $\bB^{(n)}_{r_{n-1},r_n} \in \Real^{J_n \times K_n}$,  is defined as $\bC^{(n)}= \bA^{(n)} \; \Cprod \; \bB^{(n)} \in \Real^{Q_{n-1} I_n \times Q_n K_n}$, the  blocks of which are given by $\bC^{(n)}_{q_{n-1},q_n} = \bA^{(n)}_{p_{n-1},p_n} \bB^{(n)}_{r_{n-1},r_n} \in \Real^{I_n \times K_n}$, where $\;q_n =\overline{p_n r_n}$, as illustrated in Figure \ref{Fig:SKPAC-product}. 

\begin{table}
\vspace{-1.5cm}
\centering
\caption{Basic operations on tensors in TT formats, where
$\underline \bX= \underline \bX^{(1)} \times^1 \underline \bX^{(2)} \times^1 \cdots \times^1 \underline \bX^{(N)}
\in \Real^{I_1 \times I_2 \times \cdots \times I_N}$,
$\underline \bY= \underline \bY^{(1)} \times^1 \underline \bY^{(2)} \times^1 \cdots \times^1 \underline \bY^{(N)} \in \Real^{J_1 \times J_2 \times \cdots \times J_N}$,
and $\underline \bZ= \underline \bZ^{(1)} \times^1 \underline \bZ^{(2)} \times^1 \cdots \times^1 \underline \bZ^{(N)} \in \Real^{K_1 \times K_2 \times \cdots \times K_N}$.}
\vspace{.5pc}
\centering
 {\tabsize
 \normalsize \shadingbox{
\begin{tabular}{l}
\hline \\ 
Operation \hspace{2.5cm} TT-cores   \\[0.6pc]
\hline \\ 
{
$\underline \bZ=\underline \bX+\underline \bY
		 = \left(\underline \bX^{(1)} \oplus_2 \underline \bY^{(1)}\right) \times^1
		\left(\underline \bX^{(2)} \oplus_2 \underline \bY^{(2)}\right) \times^1\cdots\times^1
		\left(\underline \bX^{(N)} \oplus_2 \underline\bY^{(N)}\right) $} \\[1.2pc]
 $\underline \bZ^{(n)} = \underline \bX^{(n)} \oplus_2 \underline \bY^{(n)}$,
with TT core slices $\bZ^{(n)}_{i_n} = \bX^{(n)}_{i_n} \oplus \bY^{(n)}_{i_n}, \;(I_n=J_n=K_n,\; \forall n)$
\\[1.0pc]
\hline
\\[0.2pc]
{
$\underline \bZ=\underline \bX \oplus \underline \bY
		 = \left(\underline \bX^{(1)} \oplus \underline \bY^{(1)}\right) \times^1
		\left(\underline \bX^{(2)} \oplus \underline \bY^{(2)}\right) \times^1\cdots\times^1
		\left(\underline \bX^{(N)} \oplus \underline\bY^{(N)}\right) $}
\\[1.0pc]
\hline
\\[0.2pc]
{
$\underline \bZ=\underline \bX\circledast\underline \bY
		 = \left(\underline \bX^{(1)} \odot_2 \underline \bY^{(1)}\right) \times^1
		\left(\underline \bX^{(2)} \odot_2 \underline \bY^{(2)}\right) \times^1\cdots\times^1
		\left(\underline \bX^{(N)} \odot_2 \underline \bY^{(N)}\right) $}
\\[1pc]
 $\underline \bZ^{(n)}=\underline \bX^{(n)} \odot_2 \underline \bY^{(n)}$,
 with TT core slices $\bZ^{(n)}_{i_n} = \bX^{(n)}_{i_n} \otimes \bY^{(n)}_{i_n} , \;(I_n=J_n=K_n,\; \forall n)$
 \\[1.0pc]
 \hline
 \\[0.0pc]
{
$\underline \bZ=\underline \bX \otimes \underline \bY
		 = \left(\underline \bX^{(1)} \otimes \underline \bY^{(1)}\right) \times^1
		\left(\underline \bX^{(2)} \otimes \underline \bY^{(2)}\right) \times^1\cdots\times^1
		\left(\underline \bX^{(N)} \otimes \underline \bY^{(N)}\right) $}
\\[1pc]
 $\underline \bZ^{(n)}=\underline \bX^{(n)} \otimes \underline \bY^{(n)}$,
 with TT core slices
  $\bZ^{(n)}_{k_n} = \bX^{(n)}_{i_n} \otimes \bY^{(n)}_{j_n}$\quad ($k_n=\overline{i_nj_n}$)
\\[1.0pc]
\hline
\\[0.0pc]
{
$\underline \bZ=\underline \bX \ast \underline \bY
		 = (\underline \bX^{(1)} \boxdot_2 \underline \bY^{(1)}) \times^1 \cdots \times^1 (\underline \bX^{(N)} \boxdot_2 \underline \bY^{(N)})$}
\\[1pc]
$\underline \bZ^{(n)} =\underline \bX^{(n)} \boxdot_2  \underline \bY^{(n)}
 \in \Real^{(R_{n-1} Q_{n-1}) \times (I_n+J_n-1) \times (R_n Q_n)}$,
with vectors \\[0.5pc]
 $\underline \bZ^{(n)}(s_{n-1},:,s_n) = \underline \bX^{(n)}(r_{n-1}, :,r_n) \ast
\underline \bY^{(n)}(q_{n-1}, :, q_n) \in \Real^{(I_n+J_n-1)}$ \\[0.5pc]
for $s_n=1,2, \ldots, R_n Q_n$
 and  $n=1,2, \ldots,N$, $\;R_0=R_N=1$.
\\[1.0pc]
\hline
\\[0.0pc]
{
$\underline \bZ = \underline \bX \times_n \bA=\underline \bX^{(1)} \times^1 \cdots \times^1 \underline \bX^{(n-1)}
		\times^1 \left(\underline \bX^{(n)}\times_2\bA \right) \times^1
		\underline \bX^{(n+1)}  \times^1 \cdots \times^1 \underline \bX^{(N)}$}
\\[1.0pc]
\hline
\\[0.0pc]

$z=\left\langle \underline \bX, \underline \bY \right\rangle =
 \bZ^{(1)} \times^1  \bZ^{(2)} \times^1 \cdots \times^1  \bZ^{(N)}=  \bZ^{(1)}   \bZ^{(2)}  \cdots   \bZ^{(N)} $\\[1pc]
$ \bZ^{(n)} = \left( \underline \bX^{(n)} \odot_2 \underline \bY^{(n)} \right)
\overline{\times}_2 {\bf{1}}_{I_n}= \sum_{i_n}\bX^{(n)}_{i_n}\otimes \bY^{(n)}_{i_n}  \;\;(I_n=J_n,\; \forall n)$ \\[1pc]
\hline
\end{tabular}\label{tab:basic_operations_TT1}
}}
\vspace{18pt}
\end{table}

\begin{table} 
\vspace{-1.5cm}
\centering
\caption{Basic operations in the TT format expressed via the strong Kronecker and C products of block matrices, where
$\bA = \widetilde{\bA}^{(1)} \skron \widetilde{\bA}^{(2)} \skron \cdots \skron \widetilde{\bA}^{(N)}$,
$\bB= \widetilde{\bB}^{(1)} \skron \widetilde{\bB}^{(2)} \skron \cdots \skron \widetilde{\bB}^{(N)}$,
$\bx = \widetilde{\bX}^{(1)} \skron \widetilde{\bX}^{(2)} \skron \cdots \skron \widetilde{\bX}^{(N)}$,
$\by = \widetilde{\bY}^{(1)} \skron \widetilde{\bY}^{(2)} \skron \cdots \skron \widetilde{\bY}^{(N)}$ and the  block matrices
$\widetilde \bA^{(n)} \in \Real^{R^A_{n-1} I_{n} \times J_n R^A_n}$, $\widetilde \bB^{(n)} \in \Real^{R^B_{n-1} J_{n} \times K_n R^B_n}$,
$\widetilde \bX^{(n)} \in \Real^{R^x_{n-1} I_{n} \times  R^x_n}$, $\widetilde \bY^{(n)} \in \Real^{R^y_{n-1} I_{n} \times  R^y_n}$.}
\vspace{.5pc}
\centering
 {\small  \shadingbox{
\begin{tabular*}{1.00\linewidth}{lll}
\hline\\[-0.8pc]
Operation \hspace{2.5cm} Block matrices of TT-cores & & \\
\hline\\[-0.8pc]
%
\multicolumn{3}{l}{
$\bZ = \bA + \bB $} \\
	$\quad = \left[\begin{array}{@{}c@{\hspace{1ex}}c@{}} \widetilde{\bA}^{(1)}&\widetilde{\bB}^{(1)}\end{array}\right]
	\skron \left[\begin{array}{@{}c@{\hspace{1ex}}c@{}} \widetilde{\bA}^{(2)}&\bf{0}\\ \bf{0}&\widetilde{\bB}^{(2)}\end{array}\right]
	\skron\cdots \skron \left[\begin{array}{@{}c@{\hspace{1ex}}c@{}}\widetilde{\bA}^{(N-1)}&\bf{0}\\ \bf{0}&\widetilde{\bB}^{(N-1)}\end{array}\right]
	\skron \begin{bmatrix} \widetilde{\bA}^{(N)}\\ \widetilde{\bB}^{(N)}\end{bmatrix}$
\\[1.5pc]
\hline
\\[0pc]
\multicolumn{3}{l}{
$\bZ = \bA \otimes \bB
	= \widetilde{\bA}^{(1)} \skron \cdots \skron \widetilde{\bA}^{(N)}
		\skron \widetilde{\bB}^{(1)}  \skron \cdots \skron \widetilde{\bB}^{(N)}$}
\\[1pc]
\hline
\\[0pc]									
\multicolumn{3}{l}{	
$z=\bx^{\text{T}}\by=\left\langle \bx, \by\right\rangle
	= \left( \widetilde{\bX}^{(1)} \Cprod \widetilde{\bY}^{(1)} \right) \skron\cdots\skron
	\left( \widetilde{\bX}^{(N)} \Cprod \widetilde{\bY}^{(N)} \right)	
	$}
\\[1pc]
	$\widetilde{\bZ}^{(n)}=\widetilde{\bX}^{(n)} \Cprod
		\widetilde{\bY}^{(n)} \in \Real^{R_{n-1}^x R^y_{n-1} \times R_{n}^x R^y_{n}}$, with core slices $\bZ^{(n)}= \sum_{i_n}\bX^{(n)}_{i_n}\otimes \bY^{(n)}_{i_n}$
\\[1pc]
\hline
\\[0pc]
\multicolumn{3}{l}{			
$\bz= \bA \bx = \left( \widetilde{\bA}^{(1)} \Cprod	\widetilde{\bX}^{(1)} \right) \skron\cdots\skron
	\left( \widetilde{\bA}^{(N)} \Cprod \widetilde{\bX}^{(N)} \right) $}
\\[1pc]
	 $\widetilde{\bZ}^{(n)}=\widetilde{\bA}^{(n)} \times^1 \widetilde{\bX}^{(n)}$,
with blocks (vectors) \\[0.5pc]
$\bz^{(n)}_{s_{n-1}, s_n}= \bA^{(n)}_{r^A_{n-1}, r^A_n} \bx^{(n)}_{r^x_{n-1},r^x_n}$ \quad ($s_n = \overline{r^A_n \, r^x_n})$
\\[1pc]
\hline
\\[0pc]
\multicolumn{3}{l}{			
$\bZ= \bA \bB
	= \left( \widetilde{\bA}^{(1)} \Cprod
		\widetilde{\bB}^{(1)} \right) \skron\cdots\skron
	\left( \widetilde{\bA}^{(N)} \Cprod \widetilde{\bB}^{(N)} \right)	
	$}\\[1pc]
$\widetilde{\bZ}^{(n)}=\widetilde{\bA}^{(n)} \Cprod \widetilde{\bB}^{(n)}$,
with blocks \\[0.5pc]
$\bZ^{(n)}_{s_{n-1},\,s_n}= \bA^{(n)}_{r^A_{n-1}, r^A_n} \bB^{(n)}_{r^B_{n-1}, r^B_n}$ \quad ($s_n = \overline{r^A_n \, r^B_n}$)
\\[1pc]
\hline
\\[0.1pc]
\multicolumn{3}{l}{			
$z = \bx^{\text{T}} \bA \bx = \langle \bx,\bA \bx\rangle$}\\
	$\quad = \left( \widetilde{\bX}^{(1)} \Cprod \widetilde{\bA}^{(1)} \Cprod \widetilde{\bX}^{(1)} \right) \skron\cdots\skron
	\left( \widetilde{\bX}^{(N)} \Cprod \widetilde{\bA}^{(N)} \Cprod	\widetilde{\bX}^{(N)} \right)	
	$ \\[1pc]
	$\widetilde{\bZ}^{(n)}=\widetilde{\bX}^{(n)} \Cprod \widetilde{\bA}^{(n)} \Cprod \widetilde{\bX}^{(n)} \in \Real^{R_{n-1}^x R^A_{n-1} R^x_{n-1} \times R_{n}^x R_n^A  R^x_{n}}$, with blocks (entries) \\[0.4pc]
$z^{(n)}_{s_{n-1},s_n}=
		\left\langle \bx^{(n)}_{r^x_{n-1},r^x_n},
		\bA^{(n)}_{r^A_{n-1},r^A_n}\bx^{(n)}_{r^y_{n-1},r^y_n}
		\right\rangle $ \quad ($s_n = \overline{r^x_n \,r^A_n \, r^y_n})$
\\[0.2pc]
\hline
\end{tabular*}
}}
\label{tab:TT_Strong_Kron}
\end{table}
\begin{algorithm}[t]
\caption{\textbf{Computation of a Matrix-by-Vector {Product}~in the TT Format}}
\label{alg:Ax}
\begin{algorithmic}[1]
\REQUIRE Matrix $\bA \in \Real^{I \times J}$ and vector $\bx \in \Real^{J}$  in  their respective TT format\\
 $\bA = \llangle \underline \bA^{(1)}, \underline \bA^{(2)}, \ldots , \underline \bA^{(N)} \rrangle \in \Real^{I_1 \times J_1 \times I_2 \times J_2 \times \cdots \times I_N \times J_N}$,\\ and $\underline \bX = \llangle \underline \bX^{(1)}, \underline \bX^{(2)}, \ldots , \underline \bX^{(N)}\rrangle \in \Real^{J_1 \times J_2 \times \cdots \times J_N}$,\\
 with TT-cores  $\underline \bX^{(n)} \in \Real^{R_{n-1} \times J_n \times R_{n}}$ and $\underline \bA^{(n)} \in \Real^{R^A_{n-1} \times I_n \times  I_n \times R^A_{n}}$
\ENSURE Matrix by vector product  $\by =\bA \bx$  in the TT format $\bY = \llangle \underline \bY^{(1)}, \underline \bY^{(2)}, \ldots , \underline \bY^{(N)}\rrangle \in \Real^{I_1 \times I_2 \times \cdots \times I_N}$, with cores  $\underline \bY^{(n)} \in \Real^{R^Y_{n-1} \times J_n \times R^Y_{n}}$
\FOR{$n=1$ to $N$}
\FOR {$i_n=1$ to $I_n$}
    \STATE $\bY^{(n)}_{i_n} =\sum_{j_n=1}^{J_n} \left( \bA^{(n)}_{i_n,j_n}
    \otimes_L \bX^{(n)}_{j_n}  \right)$
\ENDFOR
\ENDFOR
\RETURN  $\by \in \Real^{I_1 I_2 \cdots I_N}$ in the TT format $\underline \bY = \llangle \underline \bY^{(1)}, \underline \bY^{(2)}, \ldots , \underline \bY^{(N)}\rrangle$
\end{algorithmic}
\end{algorithm}

Note that, equivalently to Eq.~(\ref{ACproduct}), for $\bA \bx =\by$,  we can use a slice representation, given by
\be
\bY^{(n)}_{i_n} = \sum_{j_n=1}^{J_n} (\bA^{(n)}_{i_n,j_n} \otimes_L \bX^{(n)}_{j_n}),
\ee
which can be implemented by fast matrix-by matrix multiplication algorithms (see  Algorithm \ref{alg:Ax}). In practice, for very large scale data, we usually perform TT core contractions (MPO-MPS product) approximately, with reduced TT ranks, e.g.,  via the ``zip-up'' method proposed by \cite{Stoudenmire2010}.

In a similar way, the matrix equation
 \be
\bY \cong \bA \bX, 
\label{AXY2}
\ee
where $\bA \in \Real^{I \times J}$, $\;\;\bX  \in \Real^{J \times K}$,  $\bY \in \Real^{I \times K}$, with $I=I_1 I_2 \cdots I_N$, $J=J_1 J_2 \cdots J_N$ and
$K=K_1 K_2 \cdots K_N$, can be represented  in TT formats.  This is illustrated in  Figure \ref{Fig:YAX}(b) for the corresponding TT-cores
  defined as
\be
&&\underline \bA^{(n)} \in \Real^{P_{n-1} \times I_n \times J_n \times P_n} \notag \\
&&\underline \bX^{(n)} \in \Real^{R_{n-1} \times J_n \times K_n \times R_n} \notag \\
&&\underline \bY^{(n)} \in \Real^{Q_{n-1} \times I_n  \times K_n \times Q_n}. \notag
\ee
It is straightforward to show that when the matrices, $\bA \in \Real^{I \times J}$ and $\bX \in \Real^{J \times K}$, are represented in  their TT formats, they can be
 expressed via a strong Kronecker product of block matrices as
$\bA =  \tilde \bA^{(1)} \skron \tilde \bA^{(2)} \skron \cdots \skron \tilde \bA^{(N)}$ and
$\bX = \tilde \bX^{(1)} \skron  \tilde \bX^{(2)} \skron \cdots \skron \tilde \bX^{(N)}$, where the factor matrices are $\tilde \bA^{(n)} \in \Real^{P_{n-1} \,I_n \times J_n \,P_n}$ and $\tilde \bX^{(n)} \in \Real^{R_{n-1} \, J_n \times K_n \, R_n}$. Then, the matrix $\bY = \bA \bX$  can also be expressed via the strong Kronecker products,
  $\bY =  \tilde \bY^{(1)} \skron  \cdots \skron \tilde \bY^{(N)}$,
  where
  $\tilde \bY^{(n)} = \tilde \bA^{(n)} \; \Cprod  \; \tilde \bX^{(n)} \in \Real^{Q_{n-1} \,I_n \times K_n \,Q_n}$, $(n=1,2, \ldots, N)$, with blocks
  $\tilde \bY_{q_{n-1}, \,q_n}^{(n)} = \tilde \bA^{(n)}_{p_{n-1},\,p_n}   \tilde \bX^{(n)}_{r_{n-1},\,r_n} $, where $Q_n= R_n \, P_n,\;  q_n = \overline{p_n r_n}, \; \forall n$.

 Similarly, a quadratic form, $z=\bx^{\text{T}} \bA \bx$, for a huge symmetric
 matrix $\bA$, can be computed by first computing (in TT formats), a vector $\by = \bA \bx$, followed by the inner product $\bx^{\text{T}} \by$.

 Basic operations in the TT format are summarized in Table~\ref{tab:basic_operations_TT1}, while  Table \ref{tab:TT_Strong_Kron} presents
 these operations expressed via strong Kronecker and C products of block matrices of
 TT-cores.
 For more advanced and sophisticated operations in TT/QTT formats, see \cite{Kazeev_Toeplitz13,Kazeev2013LRT,Lee-TTfund1}.

\section{Algorithms for TT Decompositions}

 We have shown that a major advantage of the TT decomposition is  the existence of efficient algorithms for  an exact representation of higher-order tensors and/or their low-rank approximate representations with a prescribed accuracy. Similarly to the quasi-best approximation  property of the HOSVD,  the TT approximation
 $\widehat {\underline \bX} =\llangle \underline {\widehat \bX}^{(1)}, \underline {\widehat \bX}^{(2)}, \ldots, \underline {\widehat \bX}^{(N)} \rrangle
 \in \Real^{I_1 \times I_2 \times \cdots \times I_N}$ (with core tensors denoted by $\bX^{(n)}=\underline \bG^{(n)}$), obtained by the TT-SVD algorithm,  satisfies the following inequality
 \be
 \|\underline \bX - \widehat {\underline \bX}\|^2_2 \leq \sum_{n=1}^{N-1} \sum_{j=R_{n+1}}^{I_n} \sigma^2_j(\bX_{<n>}),
 \ee
 where the $\ell_2$-norm of a tensor is defined via its vectorization and $\sigma_j(\bX_{<n>})$ denotes the $j$th
 largest singular value of the unfolding matrix $\bX_{<n>}$ \cite{OseledetsTT11}.

\begin{figure}[t!]
\begin{center}
\includegraphics[width=10.5cm]{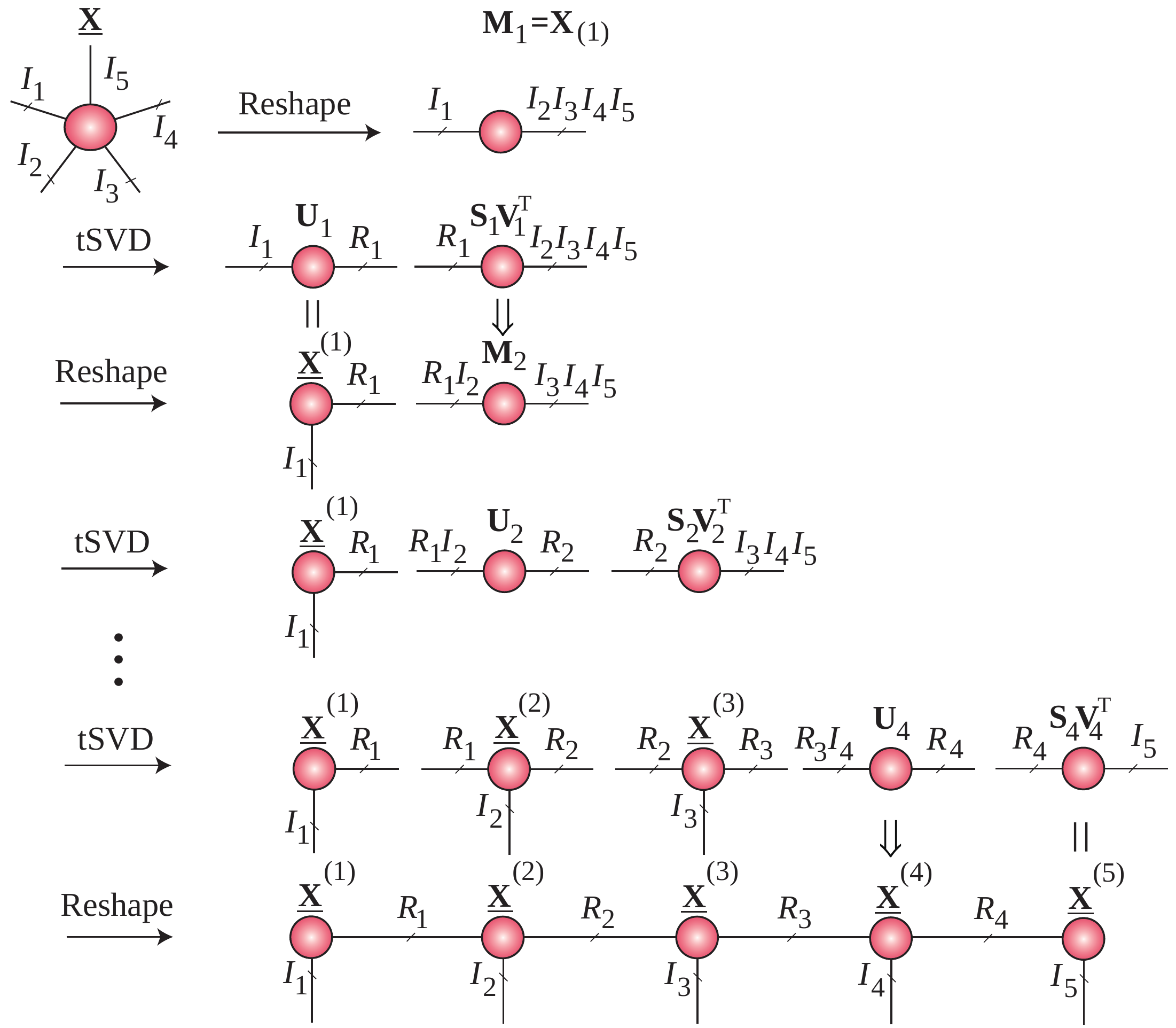}
\end{center}
\caption{The  TT-SVD algorithm  for a 5th-order  data
tensor using truncated SVD. Instead of
the  SVD, any alternative LRMF algorithm can be employed, such as randomized SVD, RPCA, CUR/CA, NMF, SCA, ICA. Top panel: A 6th-order tensor $\underline \bX$ of size $I_1 \times I_2 \times \cdots \times I_5$ is first reshaped into a long matrix $\bM_1$ of size $I_1 \times I_2 \cdots I_5$. Second panel: The tSVD is performed to produce low-rank matrix factorization, with $I_1 \times R_1$ factor matrix $\bU_1$ and the $R_1 \times I_2 \cdots I_5$ matrix $\bS_1 \bV_1^{\text{T}}$, so that $\bM_1 \cong \bU_1 \bS_1 \bV_1^{\text{T}}$. Third panel: the matrix $\bU_1$ becomes the first core core $\underline \bX^{(1)} \in \Real^{1 \times I_1 \times R_1}$, while the matrix $\bS_1 \bV_1^{\text{T}}$ is reshaped into the $R_1 I_2 \times I_3 I_4 I_5$ matrix $\bM_2$. Remaining panels: Perform tSVD to yield $\bM_2 \cong \bU_2 \bS_2 \bV_2^{\text{T}}$, reshape $\bU_2$ into an $R_1 \times I_2 \times R_2$ core $\underline \bX^{(2)}$ and repeat the procedure until all the five cores are extracted (bottom panel). The same procedure applies to higher order tensors of any order.}
\label{Fig:AlgSVDMPS}
\end{figure}

 The two basic approaches to perform efficiently TT decompositions are based on: (1) low-rank matrix factorizations (LRMF), and (2) constrained Tucker-2 decompositions.

\section{Sequential SVD/LRMF Algorithms}

The most important  algorithm for the TT decomposition is the TT-SVD algorithm
(see Algorithm \ref{alg:TT-SVD}) \cite{Vidal03,OseledetsTT09}, which  applies the truncated SVD sequentially to the unfolding matrices, as illustrated  in Figure  \ref{Fig:AlgSVDMPS}.
Instead of SVD,  alternative and efficient LRMF algorithms can be used \cite{Zorin-TT}, see also Algorithm \ref{alg:TT-LRMA}). For example, in \cite{oseledets2010tt}  a new approximate  formula for TT decomposition is proposed, where  an $N$th-order data tensor $\underline \bX$ is interpolated using a special form of cross-approximation.
In fact, the TT-Cross-Approximation is analogous to the TT-SVD algorithm, but uses adaptive
cross-approximation instead  of the computationally more expensive SVD.
 The complexity of the cross-approximation algorithms scales linearly with the order
  $N$  of a data tensor.

\begin{algorithm}[t!]
{\small
\caption{\textbf{TT-SVD Decomposition using truncated SVD (tSVD)
or  randomized SVD (rSVD) \cite{Vidal03,OseledetsTT11}}}
\label{alg:TT-SVD}
\begin{algorithmic}[1]
\REQUIRE $N$th-order tensor $\underline \bX \in \Real^{I_1 \times I_2 \times \cdots \times I_N}$ and
 approximation accuracy $\varepsilon$
\ENSURE Approximative representation of a tensor in the TT format $\underline {\widehat \bX} = \llangle \underline {\widehat \bX}^{(1)},
\underline {\widehat \bX}^{(2)}, \ldots , \underline {\widehat \bX}^{(N)} \rrangle $, such that
$\|\underline \bX - \underline {\widehat \bX}\|_F \leq \varepsilon$
\STATE Unfolding of tensor $\underline \bX$ in mode-1 $\bM_1 =\bX_{(1)}$
\STATE Initialization $R_0=1$
\FOR{$n=1$ to $N-1$}
    \STATE Perform tSVD $[\bU_n, \bS_n, \bV_n] = \mbox{tSVD} (\bM_n,\varepsilon/\sqrt{N-1})$
    \STATE Estimate $n$th TT rank $R_n= \mbox{size}(\bU_n,2)$
    \STATE Reshape orthogonal matrix $\bU_n$ into a 3rd-order core
    $\underline {\widehat\bX}^{(n)} =\mbox{reshape}(\bU_n,[R_{n-1}, I_n, R_n])$
    \STATE Reshape the matrix $\bV_n$ into a  matrix $\bM_{n+1} =\mbox{reshape}\left(\bS_n \bV_n^{\text{T}}, [R_{n}I_{n+1}, \prod_{p=n+2}^N I_p]\right)$
\ENDFOR
\STATE Construct the last core as $\underline {\widehat\bX}^{(N)} =\mbox{reshape}(\bM_N,[R_{N-1},I_N,1])$
\RETURN  $\llangle \underline {\widehat \bX}^{(1)}, \underline {\widehat \bX}^{(2)}, \ldots , \underline {\widehat \bX}^{(N)} \rrangle $.
\end{algorithmic}
}
\end{algorithm}

\begin{algorithm}[h!]
{
\caption{\textbf{TT Decomposition using any efficient LRMF}}
\label{alg:TT-LRMA}
\begin{algorithmic}[1]
\REQUIRE Tensor $\underline \bX \in \Real^{I_1 \times I_2 \times \cdots \times I_N}$ and the approximation accuracy $\varepsilon$
\ENSURE Approximate tensor representation in the TT format
 $\underline {\widehat \bX} \cong \llangle \underline {\widehat \bX}^{(1)}, \underline  {\widehat \bX}^{(2)}, \ldots , \underline  {\widehat \bX}^{(N)}\rrangle $
\STATE Initialization $R_0=1$
\STATE Unfolding of tensor $\underline \bX$ in mode-1 as $\bM_1 =\bX_{(1)}$
\FOR{$n=1$ to $N-1$}
    \STATE Perform LRMF,  e.g., CUR, RPCA, ... \\
    $[\bA_n, \bB_n] = \mbox{LRMF} (\bM_n,\varepsilon)$, i.e.,  $\bM_n \cong \bA_n \bB_n^{\text{T}}$
    \STATE Estimate $n$th TT rank, $R_n= \mbox{size}(\bA_n,2)$
    \STATE Reshape matrix $\bA_n$ into a 3rd-order core, as
    $\underline  {\widehat \bX}^{(n)} =\mbox{reshape}\left( \bA_n,[R_{n-1},I_n,R_n] \right)$
    \STATE Reshape the matrix $\bB_n$ into the $(n+1)$th unfolding matrix $\bM_{n+1} =\mbox{reshape} \left(\bB^{\text{T}}_n,[R_{n} I_{n+1}, \prod_{p=n+2}^N I_p] \right)$
\ENDFOR
\STATE Construct the last core as $\underline  {\widehat \bX}^{(N)} =\mbox{reshape}(\bM_N,[R_{N-1},I_N,1])$
\RETURN TT-cores: $\llangle \underline  {\widehat \bX}^{(1)}, \underline  {\widehat \bX}^{(2)}, \ldots , \underline  {\widehat \bX}^{(N)}\rrangle $.
\end{algorithmic}
}
\end{algorithm}

\section{Tucker-2/PVD Algorithms for Large-scale TT Decompositions}

\begin{figure}
(a)
\vspace{-0.1cm}
\begin{center}
\includegraphics[width=8.5cm]{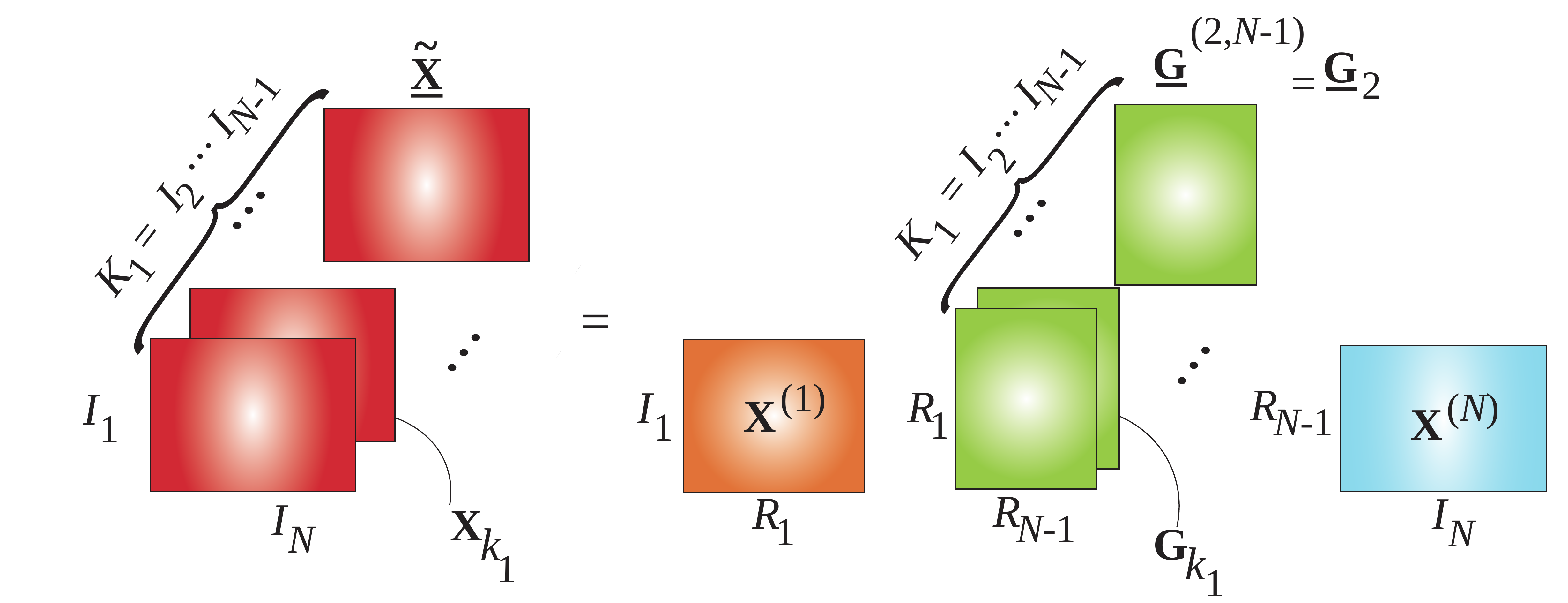}
\end{center}
(b)
\vspace{-0.5cm}
\begin{center}
\includegraphics[width=9.18cm]{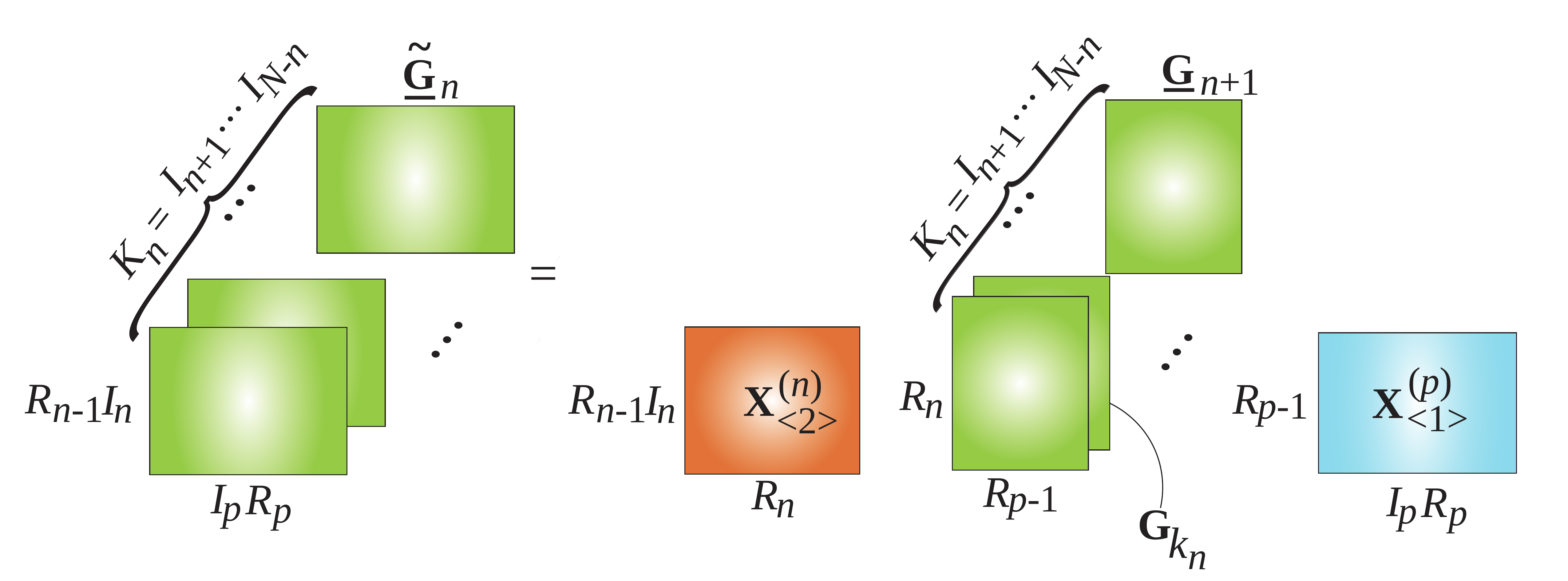}
\end{center}
(c)
\vspace{-0.6cm}
\begin{center}
\includegraphics[width=9.35cm]{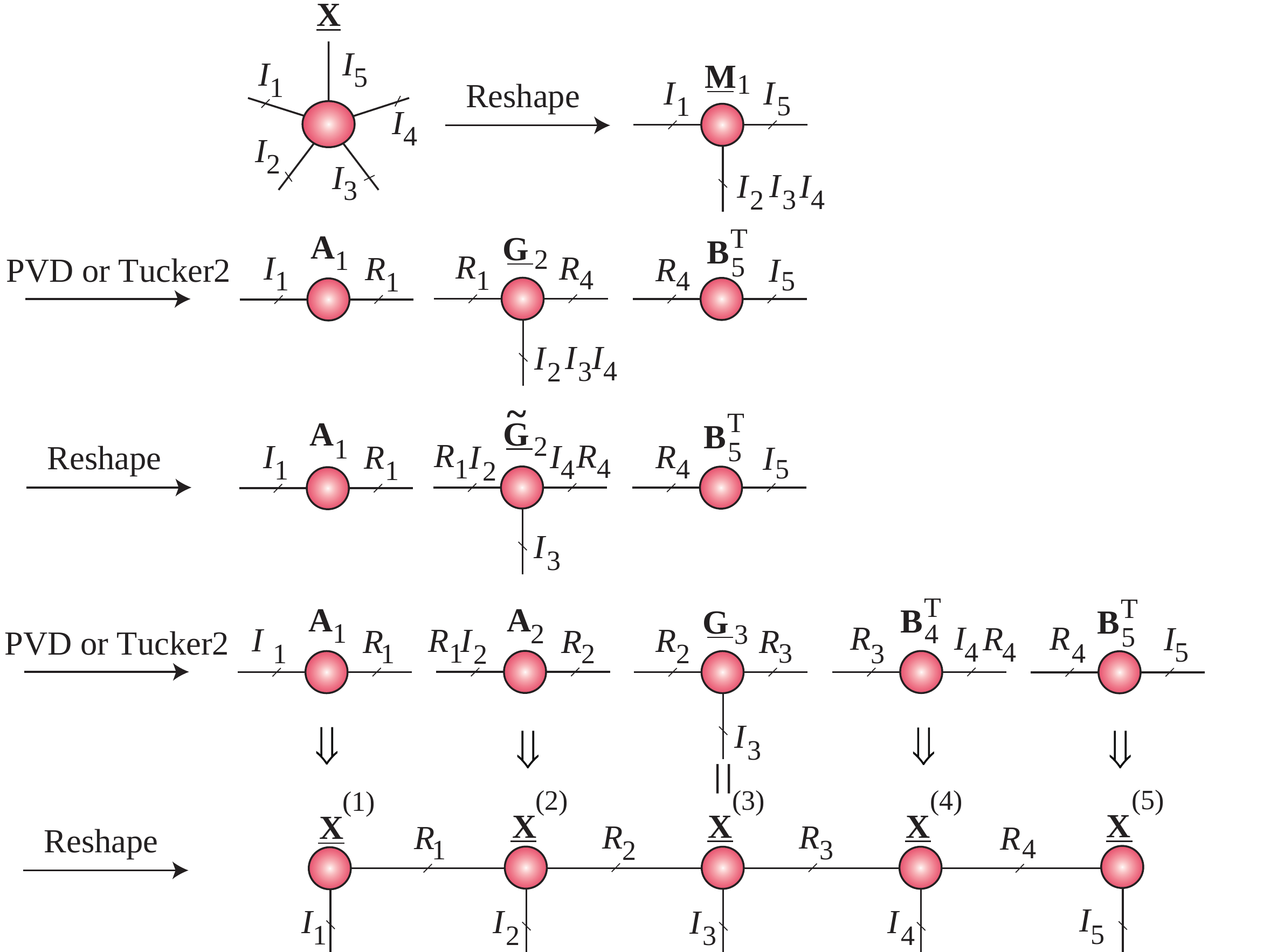}
\end{center}
\caption{{\small TT decomposition based on the Tucker-2/PVD model. (a) Extraction of the first and the last core.
 (b) The procedure can be repeated sequentially for reshaped 3rd-order tensors $\underline \bG_{n}$ (for $n=2,3,\ldots$ and $p=N-1,N-2,\ldots$). (c) Illustration of a TT  decomposition for a 5th-order data tensor, using an algorithm based on  sequential  Tucker-2/PVD decompositions.}}
\label{Fig:PVD-TT}
\end{figure}
The key idea in this approach is to reshape any
$N$th-order data tensor, $\underline \bX \in \Real^{I_1 \times I_2 \times \cdots \times I_N}$ with $N >3$, into a suitable 3rd-order tensor, e.g., $\underline {\widetilde \bX} \in \Real^{I_1 \; \times I_N \; \times \; I_2 \cdots I_{N-1}}$, in order to apply the Tucker-2 decomposition  as follows (see Algorithm \ref{alg:Tucker2} and Figure \ref{Fig:PVD-TT}(a))
\begin{equation}
\underline {\widetilde \bX} = \underline \bG^{(2,N-1)}  \times_1 \bX^{(1)} \times_2 \bX^{(N)} =   \bX^{(1)}  \times^1 \underline \bG^{(2,N-1)}  \times^1 \bX^{(N)},
\end{equation}
which, by  using frontal slices of the involved tensors, can also be expressed   in the matrix form
\be
\bX_{k_1} = \bX^{(1)} \bG_{k_1} \bX^{(N)}, \quad k_1 =1,2,\ldots, I_2 \cdots I_{N-1}.
\ee
Such representations allow us to compute the  tensor, $\underline \bG^{(2,N-1)}$, the first TT-core, $\bX^{(1)}$,  and the last TT-core, $\bX^{(N)}$.
The procedure can be repeated  sequentially  for  reshaped  tensors $\underline
{\widetilde{\bG}}_{n}=\underline \bG^{(n+1,N-n)}$ for $n=1,2,\ldots$,
 in order to extract  subsequent TT-cores in their matricized forms,
as  illustrated in Figure  \ref{Fig:PVD-TT}(b). See also the detailed step-by-step procedure  shown in Figure \ref{Fig:PVD-TT}(c).

Such a simple recursive procedure for TT decomposition can be used in conjunction with  any efficient  algorithm for
Tucker-2/PVD decompositions or the  nonnegative Tucker-2 decomposition (NTD-2) (see also Section~\ref{chap:TDs}).

\section{Tensor Train Rounding -- TT Recompression}
\label{sect:rounding}


\begin{algorithm}[t]
{
\caption{\textbf{TT Rounding (Recompression)  \cite{OseledetsTT11}}}
\label{alg:TT-rounding}
\begin{algorithmic}[1]
\REQUIRE $N$th-order tensor $\underline \bX = \llangle \underline \bX^{(1)}, \underline \bX^{(2)}, \ldots , \underline \bX^{(N)}\rrangle \in \Real^{I_1 \times I_2 \times \cdots \times I_N}$,\\ in a TT format with an overestimated  TT rank, $\brr_{TT} = \{R_1,R_2,\ldots,R_{N-1}\}$, and  TT-cores $\underline \bX \in \Real^{R_{n-1} \times I_n \times R_{n}}$,\\ absolute tolerance  $\varepsilon$, and maximum rank $R_{\max}$
\ENSURE $N$th-order tensor $\widehat{\underline \bX}$  with a reduced TT rank; the cores are rounded (reduced) according to the input tolerance
$\varepsilon$ and/or ranks bounded by $R_{\max}$, such that $\|\underline \bX - \widehat{\underline \bX}\|_F \leq \varepsilon \; \|\underline \bX \|_F $
\STATE Initialization $\underline {\widehat \bX}=\underline \bX$ and $\delta = \varepsilon/\sqrt{N-1}$
\FOR{$n=1$ to $N-1$}
    \STATE  QR decomposition $ \bX^{(n)}_{<2>} =\bQ_{n} \bR$, with $  \bX^{(n)}_{<2>} \in \Real^{R_{n-1} I_n \times R_n}$
    \STATE  Replace cores $\bX^{(n)}_{<2>}  = \bQ_n$ and  $\bX^{(n+1)}_{<1>} \leftarrow \bR \bX^{(n+1)}_{<1>}$, with $ \bX^{(n+1)}_{<1>} \in \Real^{R_{n}
    \times I_{n+1} R_{n+1}}$
\ENDFOR
\FOR{$n=N$ to $2$}
    \STATE Perform $\delta$-truncated SVD  $\bX^{(n)}_{<1>} = \bU \diag\{\mbi \sigma\} \bV^{\text{T}}$ 
    \STATE Determine minimum rank $\widehat R_{n-1}$ such that $\sum_{r>R_{n-1}} \sigma_r^2 \leq \delta^2 \|\mbi \sigma\|^2$
     \STATE  Replace cores ${\widehat \bX}^{(n-1)}_{<2>} \leftarrow {\widehat \bX}^{(n-1)}_{<2>} \widehat \bU \diag\{\widehat{\mbi \sigma}\}$ and
     ${\widehat \bX}^{(n)}_{<1>} = \widehat\bV^{\text{T}}$
\ENDFOR
\RETURN $N$th-order tensor\\ $\widehat{\bX} = \llangle \widehat{\underline \bX}^{(1)}, \widehat{\underline \bX}^{(2)}, \ldots , \widehat{\underline \bX}^{(N)}\rrangle \in \Real^{I_1 \times I_2 \times \cdots \times I_N}$, \\with reduced cores $\widehat{\underline \bX}^{(n)} \in \Real^{\widehat R_{n-1} \times I_n \times \widehat R_n}$
\end{algorithmic}
}
\end{algorithm}

\noindent
\emph{Mathematical operations in TT format produce core tensors with ranks which are not guaranteed to be optimal with respect to the desired approximation accuracy}.
For example, matrix-by-vector or matrix-by-matrix products considerably increase the TT ranks, which quickly become computationally prohibitive, so that a truncation  or low-rank TT approximations are necessary for mathematical tractability.  To  this end, the TT--rounding
 (also called  truncation or recompression)  may be used as a post-processing procedure to reduce the TT ranks. The TT rounding algorithms are typically implemented  via QR/SVD with the aim to approximate, with a desired prescribed accuracy, the  TT core tensors, $\underline \bG^{(n)}= \underline \bX^{(n)}$, by other core tensors with minimum  possible TT-ranks (see Algorithm \ref{alg:TT-rounding}). Note that  TT rounding is mathematically the same as the TT-SVD, but is more efficient owing to the  to use of TT format.

 The complexity of TT-rounding procedures is only ${\cal {O}} (N I R^3)$,  
 since all operations are performed in TT format which requires the SVD to be computed only for a relatively small matricized
 core tensor at each iteration.
 A similar approach has been developed for the HT format \cite{hTucker1,Grasedyck-rev,KressnerTobler14,espigtensorcalculus}.

\section{Orthogonalization of Tensor Train Network}
\label{sect:orthog}

The orthogonalization of core tensors is an essential procedure in many algorithms for the TT formats \cite{OseledetsTT11,Holtz-TT-2012,dolgovEIG2013,Dolgovth,KressnerEIG2014,steinlechner15,Steinlechner_phd2016}.

For convenience, we divide a TT network, which represents a tensor $\underline {\widehat{\bX}} = \llangle \underline {\widehat{\bX}}^{(1)}, \underline {\widehat{\bX}}^{(2)}, \ldots, \underline  {\widehat{\bX}}^{(N)} \rrangle \in \Real^{I_1 \times I_2 \times \cdots \times I_N}$, into sub-trains. In this way, a large-scale task is replaced by easier-to-handle sub-tasks,
whereby the aim is to extract a specific TT core or its  slices from the whole TT network.
For this purpose, the TT sub-trains can be defined  as follows 
\be
\label{TT-splitt1}
\underline {\widehat{\bX}}^{<n} &=&\llangle \underline {\widehat{\bX}}^{(1)}, \underline {\widehat{\bX}}^{(2)}, \ldots, \underline {\widehat{\bX}}^{(n-1)} \rrangle \in \Real^{I_1 \times I_2 \times \cdots \times I_{n-1} \times R_{n-1}} \\
\underline {\widehat{\bX}}^{>n} &=& \llangle \underline {\widehat{\bX}}^{(n+1)}, \underline {\widehat{\bX}}^{(n+2)}, \ldots, \underline {\widehat{\bX}}^{(N)}\rrangle \in \Real^{R_n \times I_{n+1}  \times \cdots \times I_{N}}
\label{TT-splitt2}
\ee
while the corresponding unfolding matrices, also called interface matrices, are defined by
\be
&&\widehat{\bX}^{\leq n} \in \Real^{I_1 I_2  \cdots  I_{n} \times R_{n}}, \qquad \widehat{\bX}^{> n}  \in \Real^{R_n \times I_{n+1}   \cdots I_{N}}.
\label{interface-matrices}
\ee
The left and right unfolding of the cores are defined as
\be
\widehat{\bX}_L^{(n)}= \widehat{\bX}_{<2>}^{(n)} \in \Real^{R_{n-1} I_n \times R_n} \;\; \mbox{and} \;\; \widehat{\bX}_R^{(n)}= \bX_{<1>}^{(n)} \in \Real^{R_{n-1} \times  I_n R_n} \notag \, .
\ee
\begin{algorithm}[t]
      {
        \caption{\textbf{Left-orthogonalization, right-orthog\-onalization and $n$-orthogonalization of a tensor in the TT format}}
\label{alg:TT-orthogL}
             \begin{algorithmic}[1]
\REQUIRE $N$th-order tensor $\underline {\widehat{\bX}} = \llangle \underline {\widehat{\bX}}^{(1)}, \underline {\widehat{\bX}}^{(2)}, \ldots , \underline {\widehat{\bX}}^{(N)}\rrangle \in \Real^{I_1 \times I_2 \times \cdots \times I_N}$, \\  with TT cores $\underline {\widehat{\bX}}^{(n)} \in \Real^{R_{n-1} \times I_n \times R_{n}}$ and $R_0=R_N=1$
\ENSURE  Cores $\underline {\widehat{\bX}}^{(1)}, \ldots, \underline {\widehat{\bX}}^{(n-1)}$ become left-orthogonal, while the\\ remaining cores
are right-orthogonal, except for the core $\underline {\widehat{\bX}}^{(n)}$
\FOR{$m=1$ to $n-1$}
    \STATE Perform the QR decomposition  $[\bQ, \bR] \leftarrow qr({\widehat{\bX}}^{(m)}_L)$
   for the \\unfolding cores $ {\widehat{\bX}}^{(m)}_{L} \in \Real^{R_{m-1} I_m \times R_m}$
    \STATE  Replace the cores ${\widehat{\bX}}^{(m)}_{L} \leftarrow \bQ$ and
    $\underline {\widehat{\bX}}^{(m+1)} \leftarrow \underline {\widehat{\bX}}^{(m+1)} \times_1 \bR$ 
\ENDFOR
\FOR{$m=N$ to $n+1$}
    \STATE Perform QR decomposition $[\bQ,\bR] \leftarrow qr(({\widehat{\bX}}^{(m)}_{R})^{\text{T}})$
     for the \\unfolding cores $({\widehat{\bX}}^{(m)}_{R}) \in \Real^{R_{m-1} \times I_m R_m} $,
    \STATE  Replace the cores: $\bG^{(m)}_{R} \leftarrow \bQ^{\text{T}}$ and  $\underline {\widehat{\bX}}^{(m-1)} \leftarrow \underline {\widehat{\bX}}^{(m-1)} \times_3 \bR^{\text{T}}$
\ENDFOR
\RETURN Left-orthogonal TT cores with $({\widehat{\bX}}^{(m)}_{L})^{\text{T}} {\widehat{\bX}}^{(m)}_{L} =\bI_{R_m}$ for $m=1,2,\ldots, n-1$
and right-orthogonal cores ${\widehat{\bX}}^{(m)}_{R} ({\widehat{\bX}}^{(m)}_{R})^{\text{T}} =\bI_{R_{m-1}}$ \\for $m=N, N-1, \ldots, n+1$.
\end{algorithmic}
}
 \end{algorithm}
\noindent {\bf The $n$-orthogonality of tensors.}  An $N$th-order tensor in a TT format $\underline {\widehat{\bX}} = \llangle \underline {\widehat{\bX}}^{(1)}, \ldots, \underline {\widehat{\bX}}^{(N)}\rrangle $,  is called $n$-orthogonal with $1 \leq n \leq N$, if
\be
&&({\widehat{\bX}}_L^{(m)})^{\text{T}} {\widehat{\bX}}_L^{(m)} =\bI_{R_m}, \;\;m=1, \ldots, n-1 \\
&&{\widehat{\bX}}_R^{(m)}  ({\widehat{\bX}}_R^{(m)})^{\text{T}} = \bI_{R_{m-1}},\;\; m=n+1, \ldots,N.
\ee
The tensor is called left-orthogonal if $n=N$ and right-orthogonal if $n=1$.

 When considering the $n$th TT core, it is usually assumed that all cores to the left are left-orthogonalized, and all cores to the right are right-orthogonalized.
  Notice that if a TT tensor{\footnote{By a TT-tensor we refer to as a  tensor represented in the TT format.}}, $\underline {\widehat{\bX}}$, is $n$-orthogonal then the ``left'' and ``right'' interface matrices  have orthonormal columns and rows, that is
\be
(\widehat{\bX}^{<n})^{\text{T}} \; \widehat{\bX}^{< n}  = \bI_{R_{n-1}},\qquad
\widehat{\bX}^{> n} \; (\widehat{\bX}^{>n})^{\text{T}}  =  \bI_{R_{n}}.
\ee
A tensor in a TT format can be orthogonalized efficiently using recursive QR and LQ decompositions (see Algorithm \ref{alg:TT-orthogL}).
From the above definition, for $n=N$ the algorithms perform left-orthogonalization and for $n=1$ right-orthogonalization of the whole TT network.

\section{Improved TT Decomposition Algorithm -- Alternating Single Core Update (ASCU)}
\label{sec:ascu_1}

Finally, we next present an efficient algorithm for TT decomposition, referred  to as the Alternating Single Core Update (ASCU), which sequentially optimizes  a single TT-core tensor while keeping the other TT-cores fixed in a manner similar to the modified ALS \cite{Phan_TT_part1}.

\begin{figure}[t]
\centering
\includegraphics[width=11.5cm]{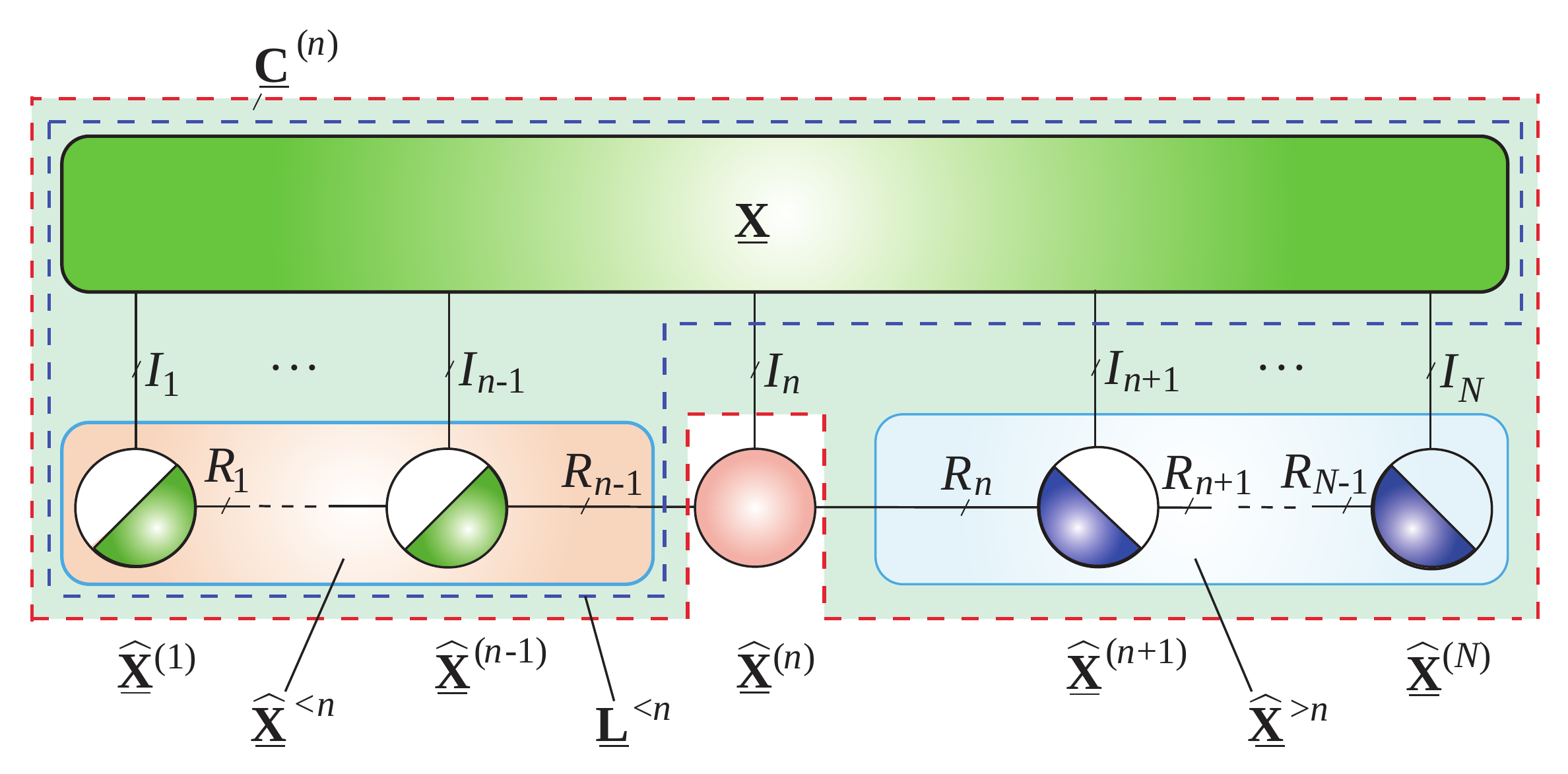}
\caption{Illustration of the contraction of tensors in the Alternating Single Core Update (ASCU) algorithm (see Algorithm \ref{alg_TT_1coreupdate}). All the cores to the left of $\underline \bX^{(n)}$ are left-orthogonal and all cores to its right are right-orthogonal.}
\label{Fig:ASCU1}
\end{figure}

Assume that the TT-tensor $\underline {\widehat {\bX}} = \llangle \underline {\widehat {\bX}}^{(1)}, \underline {\widehat {\bX}}^{(2)}, \ldots ,\underline {\widehat {\bX}}^{(N)} \rrangle $ is left- and right-orthogonalized up to $\underline {\widehat {\bX}}^{(n)}$, i.e., the unfolding matrices $\underline {\widehat {\bX}}^{(k)}_{<2>}$ for $k  = 1, \ldots, n-1$ have orthonormal columns, and $\underline {\widehat {\bX}}^{(m)}_{(1)}$ for $m = n+1, \ldots, N$ have orthonormal rows.
Then, the Frobenius norm of the TT-tensor $\underline {\widehat {\bX}}$ is  equivalent to the Frobenius norm of $\underline {\widehat {\bX}}^{(n)}$, that is, $
\|\underline {\widehat {\bX}}\|_F^2   = \|\underline {\widehat {\bX}}^{(n)}\|_F^2$, so that the Frobenius norm of the approximation error between a data tensor $\underline \bX$ and its approximate representation in the TT format $\underline {\widehat {\bX}}$ can be written as 
\be
\label{eq_cost_Gnm}
J(\underline \bX^{(n)}) &=& \|\underline \bX - \underline {\widehat {\bX}}\|_F^2 \\
&=&  \|\underline \bX\|_F^2 + \|\underline {\widehat {\bX}}\|_F^2 -2 \langle \underline \bX, \;\underline {\widehat {\bX}}\rangle \notag \\
&=& \|\underline \bX\|_F^2 + \|\underline {\widehat {\bX}}^{(n)}\|_F^2 -2 \langle \underline \bC^{(n)}, \;\underline {\widehat {\bX}}^{(n)}\rangle \notag \\
&=&  \|\underline \bX\|_F^2 - \|\underline \bC^{(n)} \|_F^2  + \| \underline \bC^{(n)} - \underline {\widehat {\bX}}^{(n)} \|_F^2, \qquad n=1,\ldots,N, \notag
\ee
where $\underline \bC^{(n)} \in \Real^{R_{n-1} \times I_n \times R_{n}}$represents a tensor contraction of $\underline \bX$ and $\underline {\widehat {\bX}}$ along all modes but the mode-$n$, as illustrated in Figure \ref{Fig:ASCU1}. The $\underline \bC^{(n)}$ can be efficiently computed through left contractions along the first $(n-1)$-modes and right contractions along the last $(N-m)$-modes, expressed as
\be
	\underline \bL^{<n} = \underline {\widehat {\bX}}^{<n} \,  \ltimesx_{n-1} \;  \underline \bX,  \qquad
	\underline \bC^{(n)}  &=& \underline \bL^{<n} \,   \rtimes_{N-n}   \; \underline {\widehat {\bX}}^{>n}.
\ee
The symbols $\ltimes_n$ and $\rtimes_m$ stand for the tensor contractions between two $N$th-order tensors along their first $n$ modes and last $m=N-n$ modes, respectively.

The optimization problem in (\ref{eq_cost_Gnm}) is usually performed subject to the following constraint
\be
\label{eq_cost_Gnm2}
	\|\underline \bX  - \underline {\widehat {\bX}}\|_F^2  \le  \varepsilon^2\,
\ee  
such that the TT-rank of $\underline {\widehat {\bX}}$ is minimum.

Observe that the constraint in (\ref{eq_cost_Gnm2}) for left- and right-orthogonalized TT-cores is equivalent to the set of sub-constraints
\be
\| \underline \bC^{(n)} - \underline {\widehat {\bX}}^{(n)} \|_F^2    \le  \varepsilon_n^2\,
\qquad n=1,\ldots,N,
\label{equ_Tn_Xn}
\ee
whereby the $n$th core $\underline \bX^{(n)} \in \Real^{R_{n-1} \times I_n \times R_n}$ should have minimum ranks $R_{n-1}$ and $R_n$.  Furthermore, $\varepsilon_n^2  = \varepsilon^2 - \|\underline \bX\|_F^2 + \|\underline \bC^{(n)} \|_F^2 $ is assumed to be non-negative.
Finally, we can formulate the following sequential optimization problem
\be
&\min  & \left(R_{n-1} \cdot R_n\right),   \notag \\
& \mbox{s.t.} \quad & \| \underline \bC^{(n)} - \underline {\widehat {\bX}}^{(n)} \|_F^2    \le  \varepsilon_n^2, \quad  n =1,2,\ldots,N.
\ee

By expressing the TT-core tensor $\underline {\widehat {\bX}}^{(n)}$ as a TT-tensor of three factors, i.e., in a Tucker-2 format given by
\be
	\underline {\widehat {\bX}}^{(n)} = \bA_n \times^1 \tilde{\underline \bX}^{(n)} \times^1 \bB_n \,,		\notag
\ee
the above optimization problem with the constraint (\ref{equ_Tn_Xn}) reduces to performing a Tucker-2 decomposition (see Algorithm \ref{alg:Tucker2}). The aim is to compute $\bA_n$, $\bB_n$  (orthogonal factor matrices) and a core tensor $\tilde{\underline \bX}^{(n)}$ which approximates tensor $\underline \bC^{(n)}$ with a minimum TT-rank-$(\tilde{R}_{n-1}, \tilde{R}_{n})$, such that 
\be
\| \underline \bC^{(n)} - \bA_n \times^1 \tilde{\underline \bX}^{(n)} \times^1 \bB_n \|_F^2    \le  \varepsilon_n^2  \, \notag \label{equ_Tn_Xn_2},
\ee
where $\bA_n \in \Real^{R_{n-1} \times \tilde{R}_{n-1}}$ and $\bB_n \in \Real^{\tilde{R}_{n} \times R_{n}}$, with $\tilde{R}_{n-1} \leftarrow R_{n-1}$ and $\tilde{R}_{n} \leftarrow R_n$.

Note that the new estimate of $\underline \bX$ is still of $N$th-order because the factor matrices $\bA_n$ and $\bB_n$ can be embedded into $\underline {\widehat {\bX}}^{(n-1)}$ and $\underline {\widehat {\bX}}^{(n+1)}$ as follows
\begin{align}
\underline {\widehat {\bX}}
&= \underline {\widehat {\bX}}^{(1)} \times^1 \cdots \times^1 (\underline {\widehat {\bX}}^{(n-1)}  \times^1 \bA_n)  \times^1 \tilde{\underline \bX}^{(n)}    \times^1 (\bB_n \times^1 \underline {\widehat {\bX}}^{(n+1)})\notag\\
&\quad \times^1 \cdots \times^1 \underline {\widehat {\bX}}^{(N)} \notag \,.
\end{align}
In this way, the three TT-cores $\underline {\widehat {\bX}}^{(n-1)}$, $\underline {\widehat {\bX}}^{(n)}$ and $\underline {\widehat {\bX}}^{(n+1)}$ are updated. Since  $\bA_n$ and $\bB_n^{\text{T}}$ have respectively  orthonormal columns and rows, the newly adjusted cores $(\underline {\widehat {\bX}}^{(n-1)} \times^1 \bA_n)$ and $(\bB_n \times^1 \underline {\widehat {\bX}}^{(n+1)})$ obey the left- and right-orthogonality conditions.
Algorithm~\ref{alg_TT_1coreupdate} outlines such a single-core update algorithm based on the Tucker-2 decomposition. In the pseudo-code, the left contracted tensor $\underline \bL^{<n}$ is computed efficiently through a progressive contraction in the form
\cite{schollwock11-DMRG,Hubig15}
\be
	\underline \bL^{<n} &=& \underline {\widehat \bX}^{(n-1)} \ltimes_2 \; \underline \bL^{<(n-1)},
\ee
where $\underline \bL^{<1}= \underline \bX$.

\begin{algorithm}[t!]
{\small
\caption{{\textbf{The Alternating Single-Core Update \hbox{Algorithm} {(two-sides rank adjustment)} \cite{Phan_TT_part1}}}
\label{alg_TT_1coreupdate}}
\begin{algorithmic}[1]
\REQUIRE{Data tensor $\underline \bX \in \Real^{I_1 \times I_2 \times \cdots \times I_N}$ and  approximation accuracy $\varepsilon$
}
\ENSURE{TT-tensor $\underline {\widehat {\bX}} = \underline {\widehat {\bX}}^{(1)} \times^1 \underline {\widehat {\bX}}^{(2)} \times^1 \cdots \times^1 \underline {\widehat {\bX}}^{(N)}$ of minimum\\ TT-rank such that $\|\underline \bX - \underline {\widehat {\bX}} \|_F^2 \le \varepsilon^2$}
\STATE Initialize $\underline {\widehat {\bX}} = \llangle \underline {\widehat {\bX}}^{(1)}, \underline {\widehat {\bX}}^{(2)}, \ldots ,\underline {\widehat {\bX}}^{(N)} \rrangle$
\REPEAT
%
\FOR{$n = 1,2 , \ldots, N-1$}  
\STATE Compute contracted tensor $\underline \bC^{(n)} =    \underline \bL^{<n}   \rtimes_{N-n} \; \underline {\widehat {\bX}}^{>n}$
\STATE  Solve a Tucker-2 decomposition $\|\underline \bC^{(n)} - \bA_n  \times^1  \underline {\widehat {\bX}}^{(n)} \times^1 \bB_n\|_F^2  \le \varepsilon^2 - \|\underline \bX\|_F^2 + \|\underline \bC^{(n)}\|_F^2$
\STATE Adjust adjacent cores\\ $\underline {\widehat \bX}^{(n-1)} \leftarrow \underline {\widehat \bX}^{(n-1)} \times^1 \bA_n$,  $\qquad  \underline {\widehat {\bX}}^{(n+1)} \leftarrow  \bB_n\times^1 \underline {\widehat {\bX}}^{(n+1)}$
\STATE  Perform left-orthogonalization of $\underline {\widehat {\bX}}^{(n)}$
\STATE Update left-side contracted tensors  \\$\underline \bL^{<n} \leftarrow \bA_n^{\text{T}} \times^1 \underline \bL^{<n}$, $\qquad  \underline \bL^{<(n+1)} \leftarrow \underline {\widehat {\bX}}^{(n)}   \ltimes_2 \; \underline \bL^{<n}$
\ENDFOR
%
\FOR{$n = N,N-1,\ldots, 2$}
\STATE Compute contracted tensor  $\underline \bC^{(n)} =    \underline \bL^{<n}   \rtimes_{N-n} \; \underline {\widehat {\bX}}^{>n}$
\STATE Solve a constrained Tucker-2 decomposition $\|\underline \bC^{(n)} - \bA_n  \times^1  \underline {\widehat {\bX}}^{(n)} \times^1 \bB_n\|_F^2  \le \varepsilon^2 - \|\underline \bX\|_F^2 + \|\underline \bC^{(n)}\|_F^2$
\STATE $\underline {\widehat \bX}^{(n-1)} \leftarrow \underline {\widehat {\bX}}^{(n-1)} \times^1 \bA_n$,  $\qquad \underline {\widehat {\bX}}^{(n+1)}\leftarrow  \bB_n\times^1 \underline {\widehat {\bX}}^{(n+1)}$
\STATE Perform right-orthogonalization of $\underline {\widehat {\bX}}^{(n)}$ 
\ENDFOR
\UNTIL{a stopping criterion is met}
\RETURN  $\llangle \underline {\widehat \bX}^{(1)}, \underline {\widehat \bX}^{(2)}, \ldots , \underline {\widehat \bX}^{(N)} \rrangle $.
\end{algorithmic}
}
\end{algorithm}

\begin{algorithm}[t!]
{\small
\caption{{\textbf{The Alternating Single-Core Update \hbox{Algorithm} {(one-side rank adjustment)} \cite{Phan_TT_part1}}}
\label{alg_TT_1coreupdate_b}}
\begin{algorithmic}[1]
\REQUIRE{Data tensor $\underline \bX \in \Real^{I_1 \times I_2 \times \cdots \times I_N}$ and  approximation accuracy $\varepsilon$
}
\ENSURE{TT-tensor $\underline {\widehat {\bX}} = \underline {\widehat {\bX}}^{(1)} \times^1 \underline {\widehat {\bX}}^{(2)} \times^1 \cdots \times^1 \underline {\widehat {\bX}}^{(N)}$ of minimum\\ TT-rank such that $\|\underline \bX - \underline {\widehat {\bX}} \|_F^2 \le \varepsilon^2$}
\STATE Initialize  TT-cores $ \underline {\widehat{\bX}}^{(n)}, \quad \forall n$
\REPEAT
%
\FOR{$n = 1, 2 , \ldots, N-1$}
\STATE Compute the contracted tensor $\underline \bC^{(n)} =    \underline \bL^{<n}   \rtimes_{N-n} \; \underline {\widehat{\bX}}^{>n}$
\STATE   Truncated SVD:\\  $\|[\underline \bC^{(n)}]_{<2>} - \bU \, \mbi \Sigma \, \bV^{\text{T}}\|_F^2  \le \varepsilon^2 - \|\underline \bX\|_F^2 + \|\underline \bC^{(n)}\|_F^2$
\STATE    Update $\underline {\widehat{\bX}}^{(n)} = \text{reshape}(\bU, R_{n-1} \times I_n \times R_{n})$
\STATE Adjust adjacent core  $\underline {\widehat{\bX}}^{(n+1)} \leftarrow  (\mbi \Sigma \, \bV^{\text{T}}) \times^1 \underline {\widehat{\bX}}^{(n+1)}$
\STATE Update left-side contracted tensors  \\$\underline \bL^{<(n+1)} \leftarrow \underline {\widehat{\bX}}^{(n)}   \ltimes_2 \; \underline \bL^{<n}$
\ENDFOR
%
\FOR{$n = N,N-1,\ldots, 2$}
\STATE Compute contracted tensor  $\underline \bC^{(n)} =    \underline \bL^{<n}   \rtimes_{N-n} \; \underline {\widehat{\bX}}^{>n}$
\STATE  Truncated SVD:\\ $\| [\underline \bC^{(n)}]_{(1)}  -  \bU \, \mbi \Sigma \, \bV^{\text{T}}  \|_F^2 \le \varepsilon^2 - \|\underline \bX\|_F^2 + \|\underline \bC^{(n)}\|_F^2$\;
\STATE $\underline {\widehat{\bX}}^{(n)}  = \text{reshape}(\bV^{\text{T}}, R_{n-1} \times I_n \times R_{n})$
\STATE $\underline {\widehat{\bX}}^{(n-1)} \leftarrow \underline {\widehat{\bX}}^{(n-1)} \times^1 (\bU \, \mbi \Sigma)$
\ENDFOR
\UNTIL{a stopping criterion is met}
\RETURN  $\llangle \underline {\widehat \bX}^{(1)}, \underline {\widehat \bX}^{(2)}, \ldots , \underline {\widehat \bX}^{(N)} \rrangle $.
\end{algorithmic}
}
\end{algorithm}

Alternatively, instead of adjusting the two TT ranks, $R_{n-1}$ and $R_{n}$, of $\underline {\widehat {\bX}}^{(n)}$, we can update only one rank, either $R_{n-1}$ or $R_{n}$, corresponding to the right-to-left or left-to-right update order procedure.
Assuming that the core tensors are updated in the left-to-right order, we need to find $\underline {\widehat {\bX}}^{(n)}$ which has a minimum rank-$R_{n}$ and satisfies  the constraints
\be
	\| \underline \bC^{(n)}  - \underline {\widehat {\bX}}^{(n)} \; \times^1 \; \bB_n \|_F^2 \le  \varepsilon_n^2, \qquad n=1,\ldots,N.  \notag \label{eq_adjust_right_rank}
\ee
This problem reduces to the truncated SVD of the mode-$\{1,2\}$ matricization of $\underline \bC^{(n)}$ with an accuracy  $\varepsilon_n^2$, that is %
\be
	[\underline \bC^{(n)}]_{<2>} \approx  \bU_n \, \mbi \Sigma \, \bV_n^{\text{T}}  \,, \notag \label{eq_factorize_Tn}
\ee 
where $\mbi \Sigma = \diag(\sigma_{n,1}, \ldots, \sigma_{n,R_{n}^{\star}})$.
Here, for the new optimized rank $R_{n}^{\star}$, the following holds%
\be
\sum_{r = 1}^{R_{n}^{\star}}  \, \sigma_{n,r}^2 \ge \|\underline \bX\|_F^2 - \varepsilon^2 > \sum_{r = 1}^{R_{n}^{\star} -1}  \, \sigma_{n,r}^2\, . \label{eq_rank_Rn}
\ee 
The core tensor $\underline {\widehat {\bX}}^{(n)}$ is then updated by reshaping $\bU_n$ to an order-3 tensor of size $R_{n-1} \times I_n \times R_{n}^{\star}$, while
the core $\underline {\widehat {\bX}}^{(n+1)}$ needs to be adjusted accordingly as 
\be
	\underline {\widehat {\bX}}^{(n+1) \star} =  \mbi \Sigma\, \bV_n^{\text{T}} \times^1 	\underline {\widehat {\bX}}^{(n+1)}\, . \label{eq_adjust_xnplus1}
\ee
When the algorithm updates the core tensors in the right-to-left order, we update $\underline {\widehat {\bX}}^{(n)}$ by using the $R_{n-1}^{\star}$ leading right singular vectors of the mode-1 matricization of $\underline \bC^{(n)}$, and adjust the core $\underline {\widehat {\bX}}^{(n-1)}$ accordingly, that is,
\begin{align}
	[\underline \bC^{(n)}]_{(1)} &\cong  \bU_n \, \mbi \Sigma \, \bV_n^{\text{T}}  \, \notag  \\
	\underline {\widehat {\bX}}^{(n) \star}  &= \text{reshape}(\bV_n^{\text{T}}, [R_{n-1}^{\star}, I_n , R_{n}])   \notag\\
	\underline {\widehat {\bX}}^{(n-1) \star}  &= \underline {\widehat {\bX}}^{(n-1)} \times^1 (\bU_n \,\mbi \Sigma)\,. \label{eq_udpate_xnminus1}
\end{align}
To summarise, the ASCU method performs a sequential update of one core and adjusts (or rotates) another core. Hence, it updates two cores at a time (for detail see  Algorithm~\ref{alg_TT_1coreupdate_b}).
%

 The ASCU algorithm can be implemented in an even more efficient way, if the data tensor $\underline \bX$ is already given in a TT format (with a non-optimal TT ranks  for the prescribed accuracy). Detailed MATLAB implementations and other variants of the TT decomposition algorithm are provided in \cite{
 Phan_TT_part1}.

\chapter{Discussion and Conclusions}

\vspace{0cm}

In Part 1 of this monograph, we have provided a systematic and example-rich guide to the basic
properties and applications of
tensor network  methodologies, and have demonstrated their promise as a tool for the
analysis of extreme-scale multidimensional data. Our main aim has been to illustrate that,
owing to the intrinsic compression ability that stems from the distributed way in which they represent
data and process information, TNs can be naturally employed for linear/multilinear dimensionality reduction.
Indeed, current applications of TNs include generalized multivariate regression, compressed
sensing, multi-way blind source separation, sparse representation and coding, feature
extraction, classification, clustering and data fusion.

With multilinear algebra as their mathematical backbone, TNs have been shown to have
intrinsic advantages over the flat two-dimensional view provided by matrices, including the
ability to model both  strong and weak couplings among multiple variables, and to cater
for multimodal, incomplete and noisy data.

In Part 2 of this monograph we introduce a scalable framework
for distributed implementation of optimization algorithms, in order to transform huge-scale
optimization problems into linked small-scale optimization sub-problems of the same type.
In that sense, TNs can be seen as a natural bridge between small-scale and very large-scale
optimization paradigms, which allows for any efficient standard numerical algorithm to be applied
to such local optimization sub-problems.

Although  research on tensor networks for dimensionality reduction and optimization problems is only emerging,
given that in many modern applications, multiway arrays (tensors) arise, either explicitly or indirectly,
through the tensorization of vectors and matrices, we foresee  this material serving  as a useful
foundation for further studies on a variety of machine learning problems for data of otherwise prohibitively
large volume, variety, or veracity. We also hope that the readers will find the approaches presented in this monograph helpful in
advancing seamlessly from numerical linear algebra to numerical multilinear algebra.



\bibliographystyle{plain}
%

\end{document}